\documentclass{amsart}

\usepackage{verbatim}
\usepackage{pst-all}

\usepackage{amsmath}
\usepackage[latin1]{inputenc}
\usepackage{amsfonts}
\usepackage{amssymb}
\usepackage{amsthm}
\usepackage{amscd}

\usepackage[latin1]{inputenc}

\usepackage[T1]{fontenc}

\usepackage{amsmath}

\usepackage{amsfonts}

\usepackage{amssymb}

\usepackage{amsthm}

\numberwithin{equation}{section}

\newcommand{\mk}{}

\newcommand{\ZZ}{\mathbb{Z}}

\newcommand{\CC}{\mathbb{C}}

\newcommand{\NN}{\mathbb{N}}

\newcommand{\QQ}{\mathbb{Q}}

\newcommand{\Glie}{\mathfrak{g}}

\newcommand{\Hlie}{\mathfrak{h}}

\newcommand{\U}{\mathcal{U}}

\newtheorem{thm}{Theorem}[section]

\newtheorem{cor}[thm]{Corollary}

\newtheorem{prop}[thm]{Proposition}

\newtheorem{lem}[thm]{Lemma}

\theoremstyle{definition}
\newtheorem{defi}[thm]{Definition}
\newtheorem{rem}[thm]{Remark}
\newtheorem*{ex}{Example}

\newtheorem*{Claim}{Claim}

\newcommand{\wt}{\mathrm{wt}}
\newcommand{\cl}{\mathrm{cl}}
\newcommand{\secref}[1]{Sect.~\ref{#1}}
\renewcommand{\epsilon}{\varepsilon}
\renewcommand{\phi}{\varphi}
\newcommand{\gr}{\mathrm{gr}}

\newenvironment{NB}{
{\bf NB}. \footnotesize
}{}
\renewenvironment{NB}{%
   \comment
}{\endcomment}

\newenvironment{NB2}{
{\bf NB}. \footnotesize
}{}
     \renewenvironment{NB2}{%
     }{}

\renewcommand{\labelenumi}{(\theenumi)}%

\usepackage[all]{xy}

\title{Level $0$ Monomial crystals}

\author{David Hernandez and Hiraku Nakajima}
\address{CNRS - UMR 8100 : Laboratoire de Math\'ematiques de Versailles, 45 avenue des Etats-Unis , Bat. Fermat, 78035 VERSAILLES, 
FRANCE}
\email{hernandez@math.cnrs.fr}
\urladdr{http://www.math.uvsq.fr/\textasciitilde hernandez}
\address{Department of Mathematics, Kyoto University, KYOTO 606-8502, JAPAN}
\email{nakajima@math.kyoto-u.ac.jp}
\urladdr{http://www.math.kyoto-u.ac.jp/\textasciitilde nakajima}

\dedicatory{Dedicated to Professor George Lusztig on his 60th birthday}

\begin{document}


\maketitle


\section*{Introduction}

In this paper we study the monomial crystal $\mathcal M$ defined by
the second author \cite{Nac}. We show that each component of $\mathcal
M$ can be embedded into a crystal $\mathcal B(\lambda)$ of an extremal
weight module $V(\lambda)$ introduced by Kashiwara \cite{kasm1}
(Theorem~\ref{emb}). 
This result was originally conjectured by Kashiwara, when 
the second author discussed the result of \cite{Nac} with him.
We prove this result by showing that the monomial crystal is equivalent
to the combinatorial crystal appeared in Kashiwara's
embedding theorem \cite{Kas-Dem}. (See Proposition~\ref{strictemb}.)
We then study the case of extremal weight modules of level $0$. We
realize the crystal $\mathcal B(\varpi_\ell)$ of a level $0$
fundamental representation via the monomial crystal
(Theorem~\ref{reallevel}). And we determine all monomials appearing in
the corresponding component of the monomial crystal for all
fundamental representations except for some fundamental
representations for $E_6^{(2)}$, $E_7^{(1)}$, $E_8^{(1)}$. Thus we
obtain explicit descriptions of the crystals in these examples. 
For classical types, we give them in terms of tableaux. For
exceptional types, we list up all monomials.
Most
of them have been calculated already in the literature
(\cite{kkmmnn,Yamane,koga,jmo,oss,Magyar,schilling,bfkl}), but we have
a few new examples in exceptional types. And our method works for
arbitrary fundamental representations in principle, though we
certainly need to use a computer with {\it huge\/} memory for the
triple node of $E_8^{(1)}$.

One of motivations of this work comes from the study of
$q$--characters of finite dimensional modules of the quantum affine
algebra, introduced by Knight~\cite{Knight},
Frenkel-Reshetikhin~\cite{Fre}, and have been intensively studied for
example in \cite{Fre2,ksoy,Naams,Naa,Nab,Nac,her02,her04,her05,cm} and
the references therein. In the combinatorial algorithm to compute
$q$--characters for arbitrary irreducible representations
\cite{Naa,Nab}, the first step was to compute ($t$--analogs of)
$q$--characters for level $0$ fundamental representations. Therefore
it would be nice if we could give their explicit forms. They can be
calculated by a computer, but we hope to see a structure by
examining their possible relations to the crystal bases.

In simply-laced type examples given in this paper, we construct {\it explicit\/}
bijections between monomials in $q$--characters, counted with
multiplicities and the crystal bases. (The existences of {\it
abstract\/} bijections are trivial as both have the same cardinality
as dimensions of modules.) In fact, the computation of the crystal
base has been done with help of explicit knowledge of
$q$--characters. This is opposite to our motivation, and we need a
further study to achieve it.

{\bf Acknowledgments :} The authors would like to thank the anonymous referee 
for comments. A part of this paper was written when the
first author visited the RIMS (Kyoto) in the summer of 2004. He
would like to thank the RIMS for his hospitality and the excellent
work conditions.

\section{Background}\label{un} In this section we give backgrounds on
quantized enveloping algebras, extremal weight modules.

\subsection{Cartan matrix} Let $C=(C_{i,j})_{1\leq i,j\leq n}$ be a
{\it generalized Cartan matrix\/}, i.e., \label{carmat} $C_{i,j}\in\ZZ$, $C_{i,i}=2$, $C_{i,j}\leq 0$ for $i\neq j$ and $C_{i,j}=0$ if and only if $C_{j,i}=0$. We set $I=\{1,\dots,n\}$ and $l=\text{rank}(C)$. In the following we suppose that $C$ is symmetrizable, that is to say that there is a matrix $D=\text{diag}(r_1,\dots,r_n)$ ($r_i\in\NN^*$)\label{ri} such that $B=DC$\label{symcar} is symmetric.

We consider a {\it realization\/} $(\Hlie, \Pi, \Pi^{\vee})$ of
$C$ (see \cite{kac}): $\Hlie$ is a $2n-l$ dimensional $\QQ$-vector
space, $\Pi=\{\alpha_1,\dots,\alpha_n\}\subset \Hlie^*$ (set of the
simple roots) and
$\Pi^{\vee}=\{\alpha_1^{\vee},\dots,\alpha_n^{\vee}\}\subset \Hlie$
(set of simple coroots) are set so that
$\alpha_j(\alpha_i^{\vee})=C_{i,j}$ for $1\leq i,j\leq n$. Let
$\Lambda_1,\dots,\Lambda_n\in\Hlie^*$ (resp.\ the
$\Lambda_1^{\vee},\dots,\Lambda_n^{\vee}\in\Hlie$) be the fundamental
weights (resp.\ coweights) :
$\Lambda_i(\alpha_j^{\vee})=\alpha_i(\Lambda_j^{\vee})=\delta_{i,j}$.


Let $P =\{\lambda \in\Hlie^* \mid \text{$\lambda(\alpha_i^{\vee})\in\ZZ$ for all $i\in I$}\}$ be the weight lattice and $P^+=\{\lambda \in P \mid \text{$\lambda(\alpha_i^{\vee})\geq 0$ for all $i\in I$}\}$ the semigroup of dominant weights.
Let $Q=\bigoplus_{i\in I} \ZZ \alpha_i\subset P$ (the root lattice) and $Q^+=\sum_{i\in I}\NN \alpha_i\subset Q$. For $\lambda,\mu\in \Hlie^*$, write $\lambda \geq \mu$ if $\lambda-\mu\in Q^+$.


\subsection{Quantized enveloping algebras}\label{qkma} In the following we suppose that $q\in\CC^*$ is not a root of unity.

Let $q_i=q^{r_i}$. For $l\in\ZZ, r\geq 0, m\geq m'\geq 0$ we introduce
the following monomials in $\ZZ[q^{\pm}]$:
$$[l]_q=\frac{q^l-q^{-l}}{q-q^{-1}}\in\ZZ[q^{\pm}],\  [r]_q!=[r]_q[r-1]_q\dots[1]_q,\ \begin{bmatrix}m\\m'\end{bmatrix}_q=\frac{[m]_q!}{[m-m']_q![m']_q!}.$$

\begin{defi} The {\it quantized enveloping algebra\/} $\U_q(\Glie)$ is the $\CC$-algebra with generators $k_h$ ($h\in \Hlie$), $x_i^{\pm}$ ($i\in I$) and relations
$$k_hk_{h'}=k_{h+h'},\ k_0=1,\ k_hx_j^{\pm}k_{-h}=q^{\pm \alpha_j(h)}x_j^{\pm},$$
$$[x_i^+,x_j^-]=\delta_{i,j}\frac{k_{r_i\alpha_i^{\vee}}-k_{-r_i\alpha_i^{\vee}}}{q_i-q_i^{-1}},$$
$$\sum_{r=0}^{1-C_{i,j}}(-1)^r\begin{bmatrix}1-C_{i,j}\\r\end{bmatrix}_{q_i}(x_i^{\pm})^{1-C_{i,j}-r}x_j^{\pm}(x_i^{\pm})^r=0 \text{ (for $i\neq j$)}.$$
\end{defi}

This algebra was introduced independently by Drinfeld and Jimbo.

We use the notation $k_i^{\pm}=k_{\pm r_i\alpha_i^{\vee}}$ and for
$l\geq 0$ we set $(x_i^{\pm})^{(l)}=(x_i^{\pm})^l/[l]_{q_i}!$.

For $J\subset I$ we denote by $\Glie_J$ the Kac-Moody algebra of Cartan matrix $(C_{i,j})_{i,j\in J}$.

Let $\U_q(\Hlie)$ the commutative subalgebra of $\U_q(\Glie)$ generated by the $k_h$ ($h\in\Hlie$).

For $V$ a $\U_q(\Hlie)$-module and $\omega\in P$ we denote by $V_{\omega}$ the weight space of weight
$\omega$ defined by
$$V_{\omega}=\{v\in V \mid \text{$k_h v=q^{\omega(h)}v$ for all $h\in \Hlie$}\}.$$ In particular for $v\in V_{\omega}$ we
have $k_i v=q_i^{\omega(\alpha_i^{\vee})}v$ and for $i\in I$ we have $x_i^{\pm}V_{\omega}\subset V_{\omega \pm
\alpha_i}$.

We say that $V$ is {\it $\U_q(\Hlie)$-diagonalizable\/} if
$V={\bigoplus_{\omega\in P}}V_{\omega}$.


\subsection{Extremal weight modules} In this section we recall the definition of extremal weight modules given by Kashiwara \cite{kasm1,kas0}.

\begin{defi} A $\U_q(\Glie)$-module $V$ is said to be {\it integrable\/} if $V$ is $\U_q(\Hlie)$-diagonalizable,
the weight subspace $V_{\omega}\subset V$ is finite dimensional for all $\omega\in P$, and for $\mu\in P$, $i\in I$ there is
$R\geq 0$ such that $V_{\mu\pm r\alpha_i}=\{0\}$ for $r\geq R$.\end{defi}

\begin{defi} For $V$ an integrable $\U_q(\Glie)$-module and $\lambda\in P$, a vector $v\in V_{\lambda}$ is called {\it extremal of weight $\lambda$\/} if there are vectors $\{v_w\}_{w\in W}$ such that $v_{\text{Id}}=v$ and
$$x_i^{\pm} v_w=0 \text{ if } \pm w(\lambda)(\alpha_i^{\vee})\geq 0
\text{ and }(x_i^{\mp})^{\pm (w(\lambda)(\alpha_i^{\vee}))} v_w=v_{s_i(w)}.$$
\end{defi}

In the same way one can define the notion of extremal elements in a crystal. Note that if $v$ is extremal of weight $\lambda$, then for $w\in W$, $v_{w}$ is extremal of weight $w(\lambda)$.

\begin{defi} For $\lambda\in P$, the {\it extremal weight module $V(\lambda)$ of extremal weight $\lambda$\/} is the $\U_q(\Glie)$-module generated by a
vector $v_{\lambda}$ with the defining relations that $v_{\lambda}$ is extremal of weight $\lambda$.\end{defi}

\begin{ex}
If $\lambda$ is dominant, $V(\lambda)$ is the
simple highest weight module of highest weight $\lambda$.
\end{ex}

\begin{thm}[\protect{\cite{kasm1}}]
For $\lambda\in P$, the module $V(\lambda)$ is
integrable and has a crystal basis $\mathcal{B}(\lambda)$.\end{thm}

Note that $u_{\lambda}\in\mathcal{B}(\lambda)$ (which represents $v_{\lambda}$) is extremal of weight $\lambda$ in the crystal
$\mathcal{B}(\lambda)$.






\section{Monomial crystal}\label{trois}

In this section we recall the definition of the monomial crystal and
show that each connected component can be embedded in the crystal of
an extremal weight module (Theorem~\ref{emb}).

In this section we suppose that $C$ is without odd
cycles, i.e., there is a function $s\colon I\rightarrow \{0,1\}$
($i\mapsto s_i$) such that $C_{i,j}\leq -1$ implies $s_i+s_j=1$. 
This situation includes all Cartan matrices of finite type and all 
Cartan matrices of affine type except $A_{2l}^{(1)}$ ($l\geq 1$).

\subsection{Construction}\label{deficrys} 
Consider formal variables $Y_{i,l}^{\pm}$, $e^{\lambda}$ ($i\in I,
l\in\ZZ, \lambda\in P$)  and let $A$ be the set of monomials of
the form $m=e^{\omega(m)}{\prod_{i\in I,l\in\ZZ}}Y_{i,l}^{u_{i,l}(m)}$ where $u_{i,l}(m)\in\ZZ$,
$\omega(m)\in P$ such that
\begin{equation}\label{eq:compat}
   {\sum_{l\in\ZZ}}u_{i,l}(m)=\omega(m)(\alpha_i^{\vee}).
\end{equation}
For $m\in A$ and $i\in I$ we set $u_i(m)={\sum_{l\in\ZZ}}u_{i,l}(m)$.

For example,
$Y_{i,l}^{\pm}e^{\pm \Lambda_i}\in A$ and
$A_{i,l}=e^{\alpha_i}Y_{i,l-1}Y_{i,l+1}{\prod_{j\neq i}}Y_{j,l}^{C_{j,i}}\in A$. 

\begin{NB}
I change $k_\lambda$ to $e^\lambda$ and make it defined for
$\lambda\in P$, rather than $\Hlie$. In your formulation, $\wt$ may
not take value in $P$ for some $r_i$. 
\end{NB}

We call $l$ the {\it grade\/} of the variable $Y_{i,l}$.

\begin{rem}\label{rem:char}
(1) If we fix a monomial $m$ and consider only monomials $m'$ which 
are products of $m$ with various $A_{i,l}^{\pm}$'s (as we shall do in this
paper), $\omega(m')$ is uniquely determined by $\omega(m)$ and
$u_{i,l}(m')$. Indeed let $z$ be a formal variable and consider the modified 
quantized Cartan matrix $C(z)=(C_{i,j}(z))_{i,j}$ defined by $C_{i,i}(z)=[2]_z$, 
and for $i\neq j$, $C_{i,j}(z)=C_{i,j}$. For $P(z)\in \ZZ[z^{\pm}]$, let
$P(z)={\sum_{l\in\ZZ}}P_lz^l$. $C(z)$ is invertible because 
$(\text{det}(C(q)))_n=1\neq 0$. Let $\tilde{C}(z)=(\tilde{C}_{i,j}(z))_{i,j}$ 
be its inverse. If $m'm^{-1}=e^{\omega(m')-\omega(m)}\prod_{i\in I, l\in\ZZ}A_{i,l}^{v_{i,l}}$ 
(with $v_{i,l}\in\ZZ$) we have 
$v_{i,l}=\sum_{j\in I, l'\in\ZZ}u_{j,l'}(m'm^{-1})(z^l\tilde{C}_{i,j}(z))_{l'}$. 
So we can safely omit $e^{\omega(m')}$.

(2) The group $A$ appears, in an equivalent form, in \cite{Nab} for
$q$--characters at roots of unity, and also in \cite{her04} to study
the $q$--characters of integrable representations of general quantum
affinizations. The additional term $e^\lambda$ (denoted by
$k_{\lambda}$ there) appears by looking at a part of a ``universal
$\mathcal{R}$-matrix''.
\end{rem}

A monomial $m$ is said to be {\it $J$-dominant\/} if for all $j\in J,
l\in\ZZ$ we have $u_{j,l}(m)\geq 0$. An $I$-dominant monomials is said
to be {\it dominant\/}. Let $B_J$ is the set of $J$-dominant monomials,
$B$ is the set of dominant monomials.


Consider the subgroup $\mathcal{M}\subset A$ defined by
$$\mathcal{M}=\{m\in A\mid u_{i,l}(m)=0\text{ if }l\equiv s_i + 1 \text{ mod } [2]\}.$$
(For the shortness of notations, we have replaced the condition $l\equiv s_i \text{ mod } [2]$ of \cite{Nac} by $l\equiv s_i+1 \text{ mod } [2]$).


Let us define $\wt\colon A\to P$ and $\epsilon_i,\phi_i, p_i,q_i\colon A\to \ZZ\cup \{\infty\}\cup\{-\infty\}$ for $i\in I$ by ($m\in A$)
\begin{gather*}
\wt(m)=\omega(m),
\\
\phi_{i,L}(m) = \sum_{l\leq L} u_{i,l}(m), 
\quad
\phi_i(m)=\max\{ \phi_{i,L}(m) \mid L\in\ZZ\}\geq 0,
\\
\epsilon_{i,L}(m)=-{\sum_{l\geq L}}u_{i,l}(m),
\quad
\epsilon_i(m) = \max\{ \epsilon_{i,L}(m) \mid L\in\ZZ\}\geq 0,
\\
p_i(m)=\max\{L\in\ZZ\mid \epsilon_{i,L}(m)=\epsilon_i(m)\}=\max\{L\in\ZZ\mid {\sum_{l<L}}u_{i,l}(m)=\phi_i(m)\},
\\
q_i(m)=\text{min}\{L\in\ZZ\mid \phi_{i,L}(m)=\phi_i(m)\}=\text{min}\{L\in\ZZ\mid -{\sum_{l> L}}u_{i,l}(m)=\epsilon_i(m)\}.
\end{gather*}
Then we define $\tilde{e}_i,\tilde{f}_i\colon A\to A\cup\{0\}$ for $i\in I$
by 
\begin{equation*}
\begin{split}
  & \tilde{e}_i(m) =
  \begin{cases}
    0 & \text{if $\epsilon_i(m) = 0$},
    \\
   mA_{i,p_i(m)-1} & \text{if $\epsilon_i(m) > 0$},
  \end{cases}
\\
  & \tilde{f}_i(m) =
  \begin{cases}
    0 & \text{if $\phi_i(m) = 0$},
    \\
   mA_{i,q_i(m)+1}^{-1} & \text{if $\phi_i(m) > 0$}.
  \end{cases}
\end{split}
\end{equation*}
By \cite{Nac,kas} $(\mathcal{M},\wt, \epsilon_i, \phi_i, \tilde{e}_i,
\tilde{f}_i)$ is a crystal (called the monomial crystal). 

\subsection{Connected components of $\mathcal{M}$ and monomial realization of highest weight crystals}\label{exa}

For $m\in \mathcal{M}$ we denote by $\mathcal{M}(m)$ the subcrystal of $\mathcal{M}$ generated by $m$.


By \cite{Nac,kas} the crystal $\mathcal{M}(m)$ is isomorphic to the
crystal $\mathcal{B}(\wt(m))$ of the highest weight module of highest
weight $\wt(m)$, if $m$ is dominant.

The aim of sections~\ref{embs} and \ref{real} is to ``generalize''
this result for general $m\in\mathcal{M}$.

\subsection{Embedding of $\mathcal{M}(m)$ into $\mathcal{B}(\lambda)$}\label{embs}

In this section we prove the following:

\begin{thm}\label{emb} For $m\in\mathcal{M}$, the crystal $\mathcal{M}(m)$ is isomorphic to a connected component of the crystal
$\mathcal{B}(\lambda)$ of an extremal weight module for some $\lambda\in P$.\end{thm}

Note that it is proved in \cite[Theorem~4.15]{bn2} that
for quantum affine algebras, all the connected components of
$\mathcal{B}(\lambda)$ are isomorphic to each other modulo shift of weight by $\delta$.

The proof is a slight modification of Kashiwara's proof of the above
mentioned result.
 
\begin{defi} A {\it shift\/} on $I$ is the data $(\leq,\phi)$ of a total ordering $\leq$ on $I$ and of a map $\phi \colon I\to \ZZ$ such that
\begin{enumerate}
\item $\phi(i)\geq \phi(j)$ for $i\leq j$,
\item if $C_{i,j}\leq -1$ and $i\leq j$, then $\phi(i)=\phi(j)+1$,
\item for $i\in I$, $s_i\equiv \phi(i) \bmod [2]$.
\end{enumerate}
\end{defi}

 For $\phi\colon I\to\ZZ$, one says that a total ordering $\leq$ on
 $I$ is {\it adapted\/} to $\phi$ if $(\leq,\phi)$ is a shift.

\begin{lem}\label{adapt} Let $\phi\colon I\rightarrow\ZZ$ such that $\phi(i)-\phi(j)\in\{\pm 1\}$ if $C_{i,j}\leq -1$ and $s_i\equiv \phi(i) \bmod [2]$ for $i\in I$. Then there is at least one total ordering on $I$ adapted to $\phi$.\end{lem}

\begin{proof} For each $r\in\ZZ$ choose a total ordering on $\{j\in I\mid \phi(j)=r\}$, and for each $(i,j)\in I^2$ such that $\phi(i)<\phi(j)$, put $i>j$.\end{proof}

 Note that in general there is at least one shift. Put $\phi(i)=s_i$, and Lemma~\ref{adapt} gives a shift $(\phi, \leq)$.

In the following we fix a shift $(\leq, \phi)$ in $I$. We put a
numbering $I=\{i_1,\dots,i_n\}$ so that $i_1 < i_2 < \dots < i_n$. 

\begin{NB}
As $\# I = n$, I omit $i_0$. (July 11)
\end{NB}

For $i\in I$, let $\mathcal B_i$ be the crystal $\mathcal B_i=\{b_i(l)\mid l\in\ZZ\}$ with $\wt(b_i(l))=l\alpha_i$ and ($j\neq i$)
$$\epsilon_i(b_i(l))=-l,\ \phi_i(b_i(l))=l,\ \tilde{e}_i(b_i(l))=b_i(l+1),\ \tilde{f}_i(b_i(l))=b_i(l-1),$$
$$\epsilon_j(b_i(l))=\phi_j(b_i(l))=-\infty,\ \tilde{e}_j(b_i(l))=\tilde{f}_j(b_i(l))=0.$$

 Let $\mathcal{B}(\infty)$ be the crystal of $\U_q^-(\Glie)$ and let
$T_{\lambda}=\{t_{\lambda}\}$ ($\lambda\in P$) be the crystal defined by $\wt(t_{\lambda})=\lambda$, $\epsilon_i(t_{\lambda})=\phi_i(t_{\lambda})=-\infty$ and $\tilde{e}_i(t_{\lambda})=\tilde{f}_i(t_{\lambda})=0$.

Let $\mathcal C$ be the crystal consisting of a single element $c$
with $\wt (c) = 0$, $\epsilon_i(c) = \varphi_i(c) = 0$, $\tilde e_i(c) =
\tilde f_i(c) = 0$.

For $m\in A$ we define the crystal $K_m=\mathcal
 C\otimes\cdots\otimes K_2 \otimes K_1\otimes K_0\otimes
 T_{\alpha}\otimes K_{-1}\otimes K_{-2}\otimes\cdots\otimes\mathcal C$
 where for $l\in\ZZ$, $K_l=\mathcal B_{i_1}\otimes \mathcal
 B_{i_2}\otimes \dots \otimes \mathcal B_{i_n}\otimes T_{\lambda(l)}$ and $\lambda(l)={\sum_{i\in I}}\lambda_i(l)\Lambda_i={\sum_{i\in I}}u_{i,2l+\phi(i)}(m)\Lambda_i$ and $\alpha=\wt(m)-{\sum_{i\in I, l\in\ZZ}}u_{i,l}(m)\Lambda_i$.

We also denote $\langle \lambda(l), \alpha_i^{\vee}\rangle$ by $\lambda_i(l)$.

\begin{defi} Let us define $\Phi_m^{\phi}\colon \mathcal{M}(m)\to K_m$
as follows: for $m'\in \mathcal{M}(m)$ with
$$m'=e^{\wt(m')}\prod_{i\in I,
k\in\ZZ}Y_{i,2k+\phi(i)}^{\lambda_i(k)}\prod_{i\in I,
k\in\ZZ}A_{i,2k+\phi(i)+1}^{z_i(k)}$$
we define $\Phi_m^{\phi}(m')=b$ by 
\begin{equation*}
\begin{split}
b&=c\otimes\cdots\otimes b_2\otimes b_1\otimes b_0\otimes
t_{\alpha}\otimes b_{-1}\otimes b_{-2}\otimes\cdots\otimes c, 
\\
& \qquad\qquad\text{where $b_l=b_{i_1}(z_{i_1}(l))\otimes \dots\otimes
b_{i_n}(z_{i_n}(l))\otimes t_{\lambda(l)}$}.
\end{split}
\end{equation*}
\end{defi}

The map $\Phi_m^{\phi}$ is well-defined as the $z_i(k)$ depend only of $m'$ (see Remark~\ref{rem:char}).

\begin{prop}\label{strictemb}
$\Phi_m^{\phi}$ is a strict embedding of the crystal.\end{prop}

When $m$ is dominant, this result appeared in \cite[8.5]{tensor} in an
equivalent form. More precisely, we parametrize
$\operatorname{Irr}\widetilde{\mathfrak Z}^\diamond$ there by
monomials as explained in \cite[\S3]{Nac}. Then the above is exactly
\cite[8.5]{tensor}.

Although the proof is exactly the same, we reproduce it here in our
current notation for the sake of the reader.

\begin{proof} The injectivity is obvious. Let $m'\in\mathcal{M}(m)$ and $b=\Phi_m^{\phi}(m')$. First we have
$$\wt(b)=\alpha+\sum_{i\in I,l\in\ZZ}u_{i,l}(m)\Lambda_i+\sum_{i\in I,l\in\ZZ}z_i(l)\alpha_i=\wt(m)+\wt(m'm^{-1}) =\wt(m').$$
 Let us prove the following formulas ($i\in I, L\in\ZZ$):

\begin{equation}\label{equaeps}\epsilon_i(b_{L-1})-\sum_{l\geq L}\wt(b_l)(\alpha_i^{\vee})=-\sum_{l\geq 2L+\phi(i)}u_{i,l}(m'),\end{equation}
\begin{equation}\label{equaphi}\phi_i(b_L)+\sum_{l<L}\wt(b_l)(\alpha_i^{\vee})=\sum_{l\leq 2L+\phi(i)}u_{i,l}(m').\end{equation}

 The equation~(\ref{equaeps}) can be checked as
\begin{equation*}
\begin{split}
& -\sum_{l\geq L}\left\{ z_i(l)+z_i(l-1)\right\}
-\sum_{l\geq L, j>i}C_{i,j}z_j(l)
-\sum_{l\geq L-1, j<i}C_{i,j}z_j(l)
-\sum_{l\geq L}\lambda_i(l)
\\
=\; &-z_i(L-1)-{\sum_{j<i}}C_{i,j}z_j(L-1)
+\sum_{l\geq L}(-{\sum_{j\in I}}C_{i,j}z_j(l)-\lambda_i(l))
\\
=\; & \epsilon_i(b_{L-1})-\sum_{l\geq L}\wt(b_l)(\alpha_i^{\vee}).
\end{split}
\end{equation*}
The equation~(\ref{equaphi}) can be checked exactly in the same way.

 The equation~(\ref{equaeps}) implies
$$\epsilon_i(b)={\max_{L\in\ZZ}}\{\epsilon_i(b_{L-1})-\sum_{l\geq L}\wt(b_l)(\alpha_i^{\vee})\}={\max_{L\in\ZZ}}\{-\sum_{l\geq 2L+\phi(i)}u_{i,l}(m')\}=\epsilon_i(m').$$
Similarly the equation (\ref{equaphi}) implies
$\phi_i(b)
=\phi_i(m').$

Let us prove the compatibility with the operators $\tilde{e}_i, \tilde{f}_i$.

If $\epsilon_i(m') = \epsilon_i(b) = 0$, then both $\tilde e_i(m')$
and $\tilde e_i(b)$ are $0$. Suppose otherwise. Then $\tilde{e}_i(b)$
is given by replacing $z_i(L_i)$ by $z_i(L_i)+1$ where $L_i=\max\{L\in\ZZ\mid \epsilon_i(b_L)-\sum_{l>L}\wt(b_l)(\alpha_i^{\vee})=\epsilon_i(b)\}$. Therefore $\tilde{e}_i(b)=\Phi_m^{\phi}(m'A_{i,2L_i+\phi(i)+1})$. But it follows from the equation~(\ref{equaeps}) that $2L_i+\phi(i)+2=p_i(m')$, and so $\tilde{e}_i(b)=\Phi_m^{\phi}(m'A_{i,p_i(m')-1})=\Phi_m^{\phi}(\tilde{e}_i(m'))$.
Similarly $\tilde f_i$ is compatible.
%
\end{proof}

Let $\mathcal B=\mathcal B_{i_1}\otimes \mathcal B_{i_2}\otimes
\dots\otimes\mathcal B_{i_n}$, and let $\mathcal{P}$ (resp.\
$\mathcal{P}^-$) be the subcrystal of  $\mathcal C\otimes \dots\otimes
\mathcal B\otimes\mathcal B$ (resp.\ of $\mathcal B\otimes \mathcal
B\otimes \dots\otimes\mathcal C$) of elements of the form $c\otimes
\dots\otimes b(0)\otimes b(0)\otimes b_l\otimes b_{l-1}\otimes
\dots\otimes b_1$ (resp.\ $b_1\otimes \dots\otimes b_{l-1}\otimes
b_l\otimes b(0)\otimes b(0)\otimes \dots\otimes c$)
where $b_{l'}\in B$ ($1\leq l'\leq l$) and $b(0)=b_{i_1}(0)\otimes
\dots\otimes b_{i_n}(0)$. 

\begin{proof}[Proof of Theorem~\ref{emb}]
By the crystal isomorphism $T_{\lambda}\otimes\mathcal B_i\simeq
\mathcal B_i\otimes T_{s_i(\lambda)}$ given by $t_{\lambda}\otimes
b_i(l)\mapsto b_i(l+\lambda(\alpha_i^{\vee}))\otimes
t_{s_i(\lambda)}$, our crystal $K_m$ is isomorphic to
\(
   \mathcal{P}\otimes T_{\lambda'}\otimes\mathcal{P}^-
\)
for some $\lambda'\in P$. 

It is known that $\mathcal P$ is isomorphic to
\(
   \bigsqcup_{\tilde e_i(b) = 0} \mathcal B(\infty)\otimes T_{\wt(b)}.  
\)
(See \cite[7.2.4]{kasb} for example.) Similarly
$\mathcal P^-$ is 
\(
   \bigsqcup_{\tilde f_i(b) = 0} T_{\wt(b)} \otimes \mathcal
   B(-\infty).
\)
Therefore
\(
   \mathcal{P}\otimes T_{\lambda'}\otimes\mathcal{P}^-
\)
is a disjoint union of various $\mathcal B(\infty)\otimes
T_\lambda\otimes B(-\infty)$. The crystal of the modified enveloping
algebra $\tilde{\U}_q(\Glie)$ is equal to
\( 
   \bigsqcup_{\lambda\in P}\mathcal{B}(\infty)\otimes
T_{\lambda}\otimes \mathcal{B}(-\infty)
\)
and its connected components can be embedded into
some $\mathcal{B}(\lambda)$ (\cite[Corollary 9.3.4]{kasm1}). Therefore
our assertion follows. \end{proof}

\section{Monomial realization of the level $0$ extremal fundamental weight crystals}\label{real}

In this section we study in more details extremal weight crystals
(Proposition~\ref{premex}) for quantum affine algebras.
We prove that the crystal of a level $0$ fundamental extremal weight module can be
realized in the monomial crystal (Theorem~\ref{reallevel}).

We omit $e^{\omega(m')}$ hereafter by Remark~\ref{rem:char}(1).

\subsection{Extremal monomials}\label{extq}

When $m$ is dominant, the component $\mathcal M(m)$ is isomorphic to
$\mathcal B(\lambda)$ where $\lambda$ is the weight of $m$.
 But the situation is different in general, as not all $m\in\mathcal{M}$ are extremal, even if the monomial is dominant or antidominant for each $i\in I$. For example in the case $D_4^{(1)}$, $m=Y_{2,0}Y_{0,3}^{-2}$ is not extremal. Indeed suppose that $m$ is extremal. Then we have
$$m_{s_2}=\tilde{f}_2(m)=Y_{2,2}^{-1}Y_{0,1}Y_{0,3}^{-2}Y_{1,1}Y_{4,1}Y_{3,1}.$$
But $(\wt(m_{s_2}))(\alpha_0^{\vee})=-1\leq 0$ and $\tilde{f}_0(m_{s_2})=Y_{0,3}^{-3}Y_{1,1}Y_{4,1}Y_{3,1}\neq 0$, and so $m_{s_2}$ is not extremal, we have a contradiction.

However we have the following consequence of Theorem~\ref{emb}.

\begin{prop}\label{premex} Let $(\phi,\leq)$ be a shift. Then for $(l_1,\dots,l_n)\in\ZZ^n$, the monomial $m=\prod_{i\in I}Y_{i,\phi(i)}^{l_i}\in\mathcal{M}$ is extremal and $\mathcal{M}(m)$ is isomorphic to the connected component of $\mathcal{B}(\wt(m))$ generated by $u_{\wt(m)}$.\end{prop}

\begin{proof} Consider the morphism $\Phi_m^{\phi}$. It follows from
Theorem~\ref{emb} that it gives an embedding $\mathcal{M}(m)\subset
\mathcal{B}(\lambda)$ where $\lambda\in P$. But for this particular
$m$ we have $\Phi_m^{\phi}(m)=c\otimes\dots\otimes b(0)\otimes
b(0)\otimes t_{\wt(m)}\otimes b(0)\otimes b(0)\otimes \dots\otimes c$
in Proposition~\ref{strictemb}.
So $m$ is sent to $u_{\wt(m)}\in\mathcal{B}(\wt(m))$ which is extremal.\end{proof}


\subsection{Monomial realization of the level $0$ extremal fundamental weight crystals}\label{numdyn}

We suppose that $C$ is of affine type. Let us number the set of simple
roots as $I=\{0,1,\dots,n\}$. We choose the extra vertices $0$ so that
$a_0=a_0^{\vee}=1$ (except $A_{2n}^{(2)}$, $a_0=2$,
$a_0^{\vee}=1$), and the index number of the vertices are the notations of
\cite{kac} (for untwisted cases $X^{(1)}$ we use the enumeration of finite type of \cite{kac} for the sub-Dynkin diagram of type $X$). This choice is unique up to an automorphism of the
Dynkin diagram. We set $I_0 = I \setminus \{ 0\}$.

We also consider a new type $A^{(2)\dagger}_{2n}$, which is the same
as $A^{(2)}_{2n}$, but we take the opposite numbering convention from
\cite{kac}, i.e., the vertex $i$ in $A^{(2)\dagger}_{2n}$ is the the
vertex $n-i$ in $A^{(2)}_{2n}$. In particular, the extra vertex $0$ is
the vertex $n$ in $A^{(2)}_{2n}$, and we have $a_0 = 1$, $a_0^\vee =
2$. We need to distinguish these as we consider the restriction of
representations to $\U_q(\Glie_{I_0})$. Note also that this convention
was taken in \cite{bn2}.

\begin{NB}
It is not good to use Bourbaki's index as there are no twisted affine
Lie algebras there.  
\end{NB}

Let $Q^{\vee}=\sum_{i\in I}\ZZ \alpha_i^{\vee}$. There is a unique
$c\in \sum_{i\in I}\NN \alpha_i^{\vee}$ such that $\{h\in Q^{\vee}\mid
\text{$\alpha_i(h)=0$ for all $i\in I$} \}=\ZZ c$. We write $c=\sum_{i\in
I}a_i^{\vee}\alpha_i^{\vee}$. In the same way one can define
$\delta=\sum_{i\in I}a_i\alpha_i\in Q$. The $a_i$ are given in \cite{kac}, 
the $a_i^{\vee}$ are the $a_i$ of the transposed Cartan matrix.

We have $\{\omega\in P\mid \text{$\omega(\alpha_i^{\vee})=0$ for all
$i\in I$} \}=\QQ\delta\cap P$. Put $P_{\cl}=P/(\QQ\delta\cap P)$. 

 Let $P^0=\{\lambda\in P\mid \lambda(c)=0\}$ be the set of level $0$ weights.

Let $\U_q(\Glie)'$ be the subalgebra of $\U_q(\Glie)$ generated by
$x_i^\pm$ and $k_h$ ($h\in \sum \QQ\alpha_i^\vee$). This has $P_\cl$
as a weight lattice. We have the corresponding definition of the
crystal. When we want to distinguish crystals of $\U_q(\Glie)$ and
$\U_q(\Glie)'$, we call the former a $P$-crystal, and the latter a
$P_\cl$-crystal.

For $i\in I_0$, let us define a {\it level $0$ fundamental weight\/}
$\varpi_i$ by $\Lambda_i-a_i^{\vee}\Lambda_0 \in P^0$ when
$\mathfrak g\neq A^{(2)\dagger}_{2n}$ and
\[
   \varpi_i = \Lambda_i - \Lambda_0 \quad (i\neq n), \qquad
   \varpi_n = 2\Lambda_n - \Lambda_0
\] 
when $\mathfrak g = A^{(2)\dagger}_{2n}$.  The corresponding extremal
weight module $V(\varpi_i)$ are called a {\it level $0$ fundamental extremal weight module}. Those representations and their crystal have been
intensively studied, see \cite{ak,beck1,bn2,kas0,kas3,ns,ns2}.

We identify these with (usual) fundamental weights of the finite
dimensional Lie algebra $\mathfrak g_{I_0}$ when $(\mathfrak g,i)\neq
(A^{(2)\dagger}_{2n},n)$. For $(\mathfrak g,i) =
(A^{(2)\dagger}_{2n},n)$, we identify $\varpi_n$ with the {\it
  twice\/} of the $n^{\mathrm{th}}$ fundamental weight.
We denote by $V_{I_0}(\varpi_i)$ the corresponding irreducible
$\mathcal U_q(\mathfrak g_{I_0})$-module, and by $\mathcal
B_{I_0}(\varpi_i)$ its crystal base, for either case.

 As $\mathcal{B}(\varpi_i)$ is connected (see \cite{kas0}), it follows from Proposition~\ref{premex} that

\begin{thm}\label{reallevel} 
Let $(\leq, \phi)$ be a shift on $I$. 
For $i\in I_0$, let $M$ be the monomial given by
$Y_{i,\phi(i)}Y_{0,\phi(0)}^{-a_i^{\vee}}$ 
for $\mathfrak g\neq A^{(2)\dagger}_{2n}$, 
$M=Y_{i,\phi(i)}Y_{0,\phi(0)}^{-1}$
for $\mathfrak g = A^{(2)\dagger}_{2n}$, $i\neq n$,
and
$M=Y_{n,\phi(n)}^2 Y_{0,\phi(0)}^{-1}$
for $\mathfrak g = A^{(2)\dagger}_{2n}$, $i = n$.
Then $M$ is extremal in $\mathcal{M}$ and
$\mathcal{M}(M)\simeq \mathcal{B}(\varpi_i)$.\end{thm}

 This result establishes a monomial realization of the level $0$ extremal fundamental weight crystals $\mathcal{B}(\varpi_i)$. We will give some examples in \secref{sec:example}.

 Not all monomials of weight $\varpi_i$ give a crystal isomorphic to
 $\mathcal{B}(\varpi_i)$ (see the example in \secref{extq}). However
 there are some other monomials which generate the same crystal as we
 will see in the next subsection.

\subsection{Other monomial realizations}

For $i\in I$, let $\theta_i\geq 0$ be the distance between $i$ and $0$, that is to say the minimum $p\geq 0$ such that there exists a sequence $\{0=j_0,j_1,\dots,j_p=i\}$ of distinct elements of $I$ satisfying $C_{j_{l},j_{l+1}}\leq -1$.

Suppose $\mathfrak g\neq A^{(2)\dagger}_{2n}$ for brevity.
\begin{cor}\label{corex} Let $i\in I_0$ and $l,l'\in\ZZ$ such that $l-l'\in\{-\theta_i,-\theta_i+2,\dots,\theta_i\}$ and $l'\equiv s_0\text{ mod }[2]$. We have $\mathcal{M}(Y_{i,l}Y_{0,l'}^{-a_i^{\vee}})\simeq\mathcal{B}(\varpi_i)$.\end{cor}

\begin{proof} It follows from Theorem~\ref{reallevel} that it suffices to show that there is a shift $(\leq, \phi)$ such that $\phi(i)=l$ and $\phi(0)=l'$. Suppose that $l-l'\leq 0$ (the proof is the same for $l-l'\geq 0$) and let $a=(\theta_i+l-l')/ 2$. Define $\phi\colon I\to \ZZ$ by
$\phi(j) = l'+\theta_j$ if $\theta_j\leq a$ and
$\phi(j)=l'+2a-\theta_j$ if $\theta_j\geq a$.
We can conclude with Lemma~\ref{adapt}.\end{proof}

 For example in all cases we have the following: 
 \begin{enumerate}
 \item if $\theta_i\in 2\ZZ$ then
$\mathcal{M}(Y_{i,0}Y_{0,0}^{-a_i^{\vee}})\simeq\mathcal{B}(\varpi_i)$,
\item  if $\theta_i\in 2\ZZ+1$ then $\mathcal{M}(Y_{i,0}Y_{0,1}^{-a_i^{\vee}})\simeq\mathcal{B}(\varpi_i)$.
 \end{enumerate}

\begin{prop}\label{otherd} 
Suppose that $C$ is of type $D_n^{(1)}$
\textup($n\geq 4$\textup) and let $i\in\{2,\dots,n-2\}$. Then
$M=Y_{i,0}Y_{0,i-1}^{-1}Y_{0,i+1}^{-1}$ is extremal and
$\mathcal{B}(M)\simeq \mathcal{B}(\varpi_i)$.
\end{prop}

\begin{proof} First suppose that $i\leq n-3$. Consider
 $$m=(\tilde{f}_2 \cdots \tilde{f}_{i-1}\tilde{f}_i)(M)=Y_{0,i+1}^{-1}Y_{1,i-1}Y_{2,i}^{-1}Y_{i+1,1}.$$
Let us define $\phi\colon I\to\ZZ$ by $\phi(0)=i+1, \phi(2)=i,
\phi(1)=\phi(3)=i-1,\phi(4)=i-2,\dots,\phi(n-2)=i-n+4,
\phi(n)=\phi(n-1)=i-n+5$. Lemma~\ref{adapt} gives a shift
$(\phi,\leq)$. So it follows from Proposition~\ref{premex} that $m$ is
extremal, and so $M=m_{s_i s_{i-1} \dots s_2}$ is extremal.

  If $i=n-2$, in the same way we consider
 $$m=(\tilde{f}_2 \cdots \tilde{f}_{i-1}\tilde{f}_i)(M)=Y_{0,n-1}^{-1}Y_{1,n-3}Y_{2,n-2}^{-1}Y_{n-1,1}Y_{n,1}.$$
 \end{proof}

In the following we will see various examples of realizations of the level $0$ extremal fundamental weight crystals.

\section{Finite dimensional crystals -- start}

\begin{NB}
I have changed the notation. The highest weight monomial should be written
$M$, $M_0$, etc. Other monomials are written by $m$, etc.   
\end{NB}

Kashiwara has shown that there are a $\U_q(\mathfrak g)'$-automorphism
$z_\ell$ of the level $0$ fundamental extremal weight module $V(\varpi_\ell)$
preserving the global crystal base, and the induced
$P_{\mathrm{cl}}$-crystal automorphism, denoted also by $z_\ell$, on
the crystal $\mathcal B(\varpi_\ell)$~\cite{kas0}. The weight of $z_\ell$ in the $P$-crystal 
is $d_\ell\delta$ where $d_\ell=\text{max}(1,a_{\ell}^{\vee}/a_{\ell})$ except $d_{\ell} = 1$ for $(\Glie,\ell) = (A^{(2)}_{2n}, n)$. The quotient
$\mathcal B(\varpi_\ell)/z_\ell$ is the crystal of the finite
dimensional irreducible $\U_q(\mathfrak g)'$-module $W(\varpi_\ell) =
V(\varpi_\ell)/(z_\ell - 1) V(\varpi_\ell)$. We denote it by $\mathcal
B(W(\varpi_\ell))$. We call $W(\varpi_\ell)$ the {\it level $0$ fundamental
representation}.

After Theorem~\ref{reallevel} it is natural to ask the followings.
\begin{enumerate}
\item Give an explicit description of monomials appearing in $\mathcal
  M(M)$.
\item Give an explicit description of the automorphism $z_\ell$.
\end{enumerate}

Note that the automorphism $z_\ell$ is defined as a composite of
operators $\tilde e_i$, $\tilde f_i$'s. But we require more explicit
description.

We do not answer these questions in general, but we give examples
in the next sections. These are motivated by known descriptions of
level $0$ crystals in terms of tableaux \cite{kkmmnn,koga,schilling}
in part, but closer to those of $q$--characters \cite{Nac}.

Before giving examples, we define $P_{\mathrm{cl}}$-crystal
automorphisms on the monomial crystal $\mathcal M$.
For $p\in\ZZ, \alpha\in \QQ\delta\cap P$ let $\tau_{2p, \alpha}$
denote the map $\tau_{2p,\alpha}\colon
\mathcal{M}\rightarrow\mathcal{M}$ defined by
$\tau_{2p,\alpha}(e^\lambda \prod Y_{i,n}^{u_{i,n}}) =
e^{\lambda+\alpha} \prod Y_{i,n+2p}^{u_{i,n}}$. This clearly preserves the
compatibility condition~(\ref{eq:compat}) and is a $P_{\mathrm{cl}}$-crystal
automorphism. In the following, we omit $\alpha$ from the notation and
denote simply by $\tau_{2p}$.

Suppose that $\mathcal M(M)$ is a monomial crystal isomorphic to
$\mathcal B(\varpi_\ell)$ such that $M$ is an extremal vector with
$\tilde{e}_i M = 0$ for all $i\in I_0$.
If we have a monomial $m\in\mathcal M(M)$ with $\wt(m) = \wt(M) +
N d_\ell\delta$ for $N\in\ZZ$, then we have $m = z_\ell^N(M)$. This follows
from \cite[\S5.2]{kas0}. In particular, if $\mathcal M(M)$ is
isomorphic to $\mathcal B(\varpi_\ell)$ and preserved under
$\tau_{2p}$, then $\tau_{2p}$ is equal to a power of $z_\ell$.

In the following examples, we answer the above questions (1),(2) in
the following manner:
\begin{enumerate}
\item First show that $\mathcal M(M)$ is invariant under $\tau_{2p}$
for some $p$. Then $\mathcal M(M)/\tau_{2p} \simeq \mathcal
B(\varpi_\ell)/z_\ell^N$ for some $N$.

\item We determine all monomials in $\mathcal M(M)/\tau_{2p}$ and
give $z_\ell$ explicitly in these monomials.
\end{enumerate}

We thus obtain explicit descriptions of crystals of some finite
dimensional representations of $\U_q(\Glie)'$ : we treat all fundamental representations except some fundamental representations for $E_6^{(2)}$,
$E_7^{(1)}$, $E_8^{(1)}$. However it is natural to hope this procedure works for any fundamental
representations with appropriate choices of the initial monomials $m$.

Note that the uniqueness of the crystal base for $W(\varpi_\ell)$ is
not known so far. But all the examples where we compare the crystal
base with those existing in the literature, we can always prove that
the crystal base is isomorphic.

\subsection{}
Let us illustrate our description in type $A^{(1)}_{2r+1}$ with $n =
2r+1$ ($r\ge 0$).

\newcommand{\yl}{15pt}
\newcommand{\yh}{8pt}
\newcommand{\ffbox}[1]{
\setbox9=\hbox{$\scriptstyle\overline{1}$}
%
\framebox[\yl][c]{\rule{0mm}{\ht9}${\scriptstyle #1}$}
}
\newcommand{\te}{\tilde e}
\newcommand{\tf}{\tilde f}
\newenvironment{aenume}{%
  \begin{enumerate}%
  \renewcommand{\theenumi}{\alph{enumi}}%
  \renewcommand{\labelenumi}{(\theenumi)}%
  }{\end{enumerate}}

\begin{NB}
When $\mathfrak g$ is of type $A^{(1)}_{n}$, the crystal $\mathcal
B(W(\varpi_i))$ can be explicitly written \cite{kkmmnn} as follows:
Let $b = (b_1,\dots b_{n+1})\in \{ 0,1 \}^{n+1}$ such that $\sum_j b_j
= i$. Then
\begin{equation*}
  \begin{split}
  & \tilde f_0(b_1,\dots, b_{n+1}) 
  = (b_1+1,b_2,\dots,b_{n},b_{n+1}-1),
\\  
  & \tilde f_j(b_1,\dots, b_{n+1}) 
  = (b_1,\dots,b_{j-1},b_j-1,b_{j+1}+1,b_{j+2},\dots,b_{n+1}), \quad j\neq 0,
  \end{split}
\end{equation*}
where we understand the right hand side is $0$ if it is not in
$\{ 0,1\}^{n+1}$.

We can easily translate the description to the monomials.
\end{NB}

Mimicking the definition in \cite{Nac,kks}, we define
\begin{equation*}
  \ffbox{k}_p = Y_{k-1,p+k}^{-1} Y_{k,p+k-1} \qquad
  \text{for $1\le k\le n+1$, $p\in\ZZ$},
\end{equation*}
where $Y_{n+1,p}$ is understood as $Y_{0,p}$.

\subsubsection{} Let us consider the first level $0$ fundamental extremal weight module 
$V(\varpi_1)$. Let $M = Y_{1,p} Y_{0,p+1}^{-1}$. We have $\mathcal M(M)
\simeq \mathcal B(\varpi_1)$ by Corollary~\ref{corex}.

 Then the crystal graph
of $\mathcal M(M)$ is given in Figure~\ref{fig:An}. 
\begin{figure}[htbp]
\centering
\psset{xunit=1mm,yunit=1mm,runit=1mm}
\psset{linewidth=0.3,dotsep=1,hatchwidth=0.3,hatchsep=1.5,shadowsize=1}
\psset{dotsize=0.7 2.5,dotscale=1 1,fillcolor=black}
\psset{arrowsize=1 2,arrowlength=1,arrowinset=0.25,tbarsize=0.7 5,bracketlength=0.15,rbracketlength=0.15}
\begin{pspicture}(0,60)(160.78,85)
\rput(10,80){$\ffbox{1}_p$}
\rput(30,80){$\ffbox{2}_p$}
\rput(90,80){$\ffbox{n}_p$}
\rput(110,80){$\ffbox{n\!+\!1}_p$}
\rput{90}(60,164.22){\parametricplot[arrows=<-]{150.26}{209.74}{ t cos 100.78 mul t sin 100.78 mul }}
\psline{->}(15,80)(25,80)
\psline{->}(35,80)(45,80)
\psline{->}(95,80)(105,80)
\psline{->}(75,80)(85,80)
\psline[linestyle=dotted,dotsep=2](70,80)(50,80)
\rput(20,82){$\scriptstyle 1$}
\rput(40,82){$\scriptstyle 2$}
\rput(80,82){$\scriptstyle n\!-\!1$}
\rput(100,82){$\scriptstyle n$}
\rput(60,65){$\scriptstyle 0[n\!+\!1]$}
\end{pspicture}
\caption{(Type $A_n^{(1)}$) the crystal $\mathcal B(\varpi_1)$ of the vector representation}
\label{fig:An}
\end{figure}
Here $0[n+1]$ means $\tilde f_0\ffbox{n\!+\!1}_p\! =
\ffbox{1}_{p+n+1}$, i.e., the suffix is shifted by $n+1$. In
particular $\mathcal M(M)$ is preserved under $\tau_{n+1}$, which has
weight $-\delta$. Therefore we have $z_1=\tau_{-n-1}$ and $\mathcal
M(M)/\tau_{n+1}\simeq \mathcal B(W(\varpi_1))$.

\begin{NB}
Note that this study includes the type $A_1^{(1)}$. Let $M=Y_{1,0}Y_{0,1}^{-1}$ and $m=\tilde{f}_1(M)=Y_{1,2}^{-1}Y_{0,1}$. We have $\mathcal B(W(\varpi_1))\simeq\mathcal M(M)/\tau_2=\mathcal M_{I_0}(M)=\{M, m\}$. The crystal graph of $\mathcal M(M)$ is given in Figure~\ref{fig:A1}.
\begin{figure}[htbp]
\centering
\psset{xunit=.55mm,yunit=.55mm,runit=.55mm}
\psset{linewidth=0.3,dotsep=1,hatchwidth=0.3,hatchsep=1.5,shadowsize=1}
\psset{dotsize=0.7 2.5,dotscale=1 1,fillcolor=black}
\psset{arrowsize=1 2,arrowlength=1,arrowinset=0.25,tbarsize=0.7 5,bracketlength=0.15,rbracketlength=0.15}
\begin{pspicture}(0,0)(220,30)
\rput(15,20){$M$}
\psline{->}(25,20)(183,20)
\rput(193,20){$m$}
\rput{0}(104,525){\parametricplot[arrows=<-]{-99.15}{-80.85}{ t cos 519.26 mul t sin 519.26 mul }}
\rput(104,25){$\scriptscriptstyle 1$}
\rput(104,0){$\scriptscriptstyle 0[2]$}
\end{pspicture}
\caption{(Type $A_1^{(1)}$) the crystal $\mathcal B(\varpi_1)$}
\label{fig:A1}
\end{figure}
\end{NB}

\subsubsection{} Next consider $\mathcal B(\varpi_\ell)$ for $\ell \le r+1$. (The
description for the remaining case $\ell > r+1$ can be obtained from
these cases by applying a diagram automorphism.) Let 
$M_0 = Y_{\ell,0}Y_{0,\ell}^{-1}$. It follows from Corollary~\ref{corex} that
$\mathcal M(M_0) \simeq \mathcal B(\varpi_\ell)$. We set
\begin{equation*}
  \begin{split}
  M_j &= Y_{\ell,2j}Y_{0,n-\ell+1+2j}^{-1}
       Y_{j,\ell+j}^{-1}Y_{j,n-\ell+1+j}
\\
   &= \left( \ffbox{1}_{n-\ell+2j}\ffbox{2}_{n-\ell+2j-2}
   \cdots \ffbox{j}_{n-\ell+2}\right)
\times
  \left(\ffbox{j\!+\!1}_{\ell-1}\ffbox{j\!+\!2}_{\ell-3}
    \cdots \ffbox{\ell}_{1-\ell+2j}\right)
\\
   &=
   \prod_{p=1}^j \ffbox{p}_{n-\ell-2p+2j+2} \times
   \prod_{p=j+1}^\ell \ffbox{p}_{\ell+1-2p+2j}
  \end{split}
\end{equation*}
with $0\le j\le \ell$. Note that $M_\ell =
Y_{\ell,n+1}Y_{0,n+1+\ell}^{-1} = \tau_{n+1}(M_0)$. Note also that $M_1
= \tau_2(M_0)$ for $\ell = r+1$.

For an increasing sequence $T = (1\le i_1 < i_2 < \dots < i_\ell\le
n+1)$ of integers (i.e., a Young tableaux of shape $(\ell)$) we assign
\begin{equation*}
   m_{T;j} = \prod_{p=1}^j \ffbox{i_p}_{n-\ell-2p+2j+2} \times
   \prod_{p=j+1}^\ell \ffbox{i_p}_{\ell+1-2p+2j}
   \qquad\text{for $0\le j\le \ell-1$}.
\end{equation*}

Then one can directly check that 
\begin{enumerate}
\item $\mathcal M_{I_0}(M_j)$ consists of $m_{T;j}$ for various
sequences $T$ (cf.\ \cite[4.6]{Nac}), 

\item $\tilde f_0(m_{T;j})$ with $T = (2,3,\dots,\ell, n+1)$ is equal
to $M_{j+1}$,

\item the automorphism $\sigma$ defined by $m_{T;j} \mapsto m_{T;j+1}$
($j$, $j+1$ are understood modulo $\ell$) is a
$P_{\mathrm{cl}}$-crystal automorphism.
\end{enumerate}
Here, for $J\subset I$ and $m\in\mathcal{M}$ we denote by
$\mathcal{M}_J(m)$ the set of monomials obtained by applying $\te_j$,
$\tf_j$ with $j\in J$ to $m$. It is a crystal for the Lie subalgebra
$\Glie_J$ associated with $J$.

From (2) all $M_j$ (and hence $m_{T;j}$ by (1)) are in $\mathcal M(M_0)$
by induction. 
Computing the weights, we find that $M_j = (z_\ell)^{-j}(M_0)$ as
explained above. In particular, $\tau_{n+1} = (z_\ell)^{-\ell}$. In
the case $\ell = r+1$, we have $\tau_2 = z_\ell^{-1}$. Therefore
$\mathcal M(M_0)/\tau_2 \simeq \mathcal B(W(\varpi_{r+1}))$.
By the same reason mentioned above, $\sigma$ is equal to $z_\ell$.
Therefore $\mathcal M(M_0)/\sigma \simeq \mathcal B(W(\varpi_\ell))$.

Let us describe Kashiwara operators $\tilde{e}_i$, $\tilde{f}_i$ in
terms of tableaux. This can be done by transfering the definition of
those operators on monomials to tableaux. For $i\neq 0$ we have
$\tilde{e}_i m_{T;j} = m_{T';j}$ or $0$. Here $T'$ is obtained from
$T$ by replacing $i$ by $i-1$. If it is not possible (say, when we
have both $i-1$ and $i$ in $T$), then it is zero. Similarly
$\tilde{f}_i = m_{T'';j}$ or $0$, where $T''$ is given by replacing
$i$ by $i+1$.
We can also describe the action of $\tilde{e}_0$, $\tilde{f}_0$ :
\begin{equation*}
\begin{split}
\tilde{e}_0(m_{T;j}) &=
        \begin{cases}
        0 &\text{ if $i_1\neq 1$ or $i_{\ell}=n+1$},
                \\  m_{(i_2,\cdots,i_{\ell},n+1);j-1}&\text{ if $i_1=1$ and $i_{\ell}\neq n+1$,}
  \end{cases}
\\
\tilde{f}_0(m_{T;j}) &=
        \begin{cases}
        0 &\text{ if $i_1 = 1$ or $i_{\ell} \neq n+1$},
                \\  m_{(1,i_1,\cdots,i_{\ell-1});j+1}&\text{ if $i_1 \neq 1$ and $i_{\ell} = n+1$.}
  \end{cases} 
\end{split}
\end{equation*}
Here we extend the definition of $m_{T;j}$ from $0\le j\le \ell-1$ to
all $j\in\mathbb Z$ so that $m_{T;j+\ell} = \tau_{n+1} m_{T;j}$.

As a corollary we get a description of $\mathcal B(W(\varpi_\ell))$ in
terms of tableaux. This coincides with one in \cite{kkmmnn}. We also get
an isomorphism of $I_0$-crystals 
$\mathcal B(W(\varpi_\ell))\simeq \mathcal B_{I_0}(\varpi_\ell)$.
This is a well-known result.

Comparing the above descriptions with the tableaux sum expressions of
$q$--characters in \cite{Nac}, we see that there is a bijection
between $\mathcal M(M_0)/\sigma \simeq \mathcal B(W(\varpi_\ell))$ and
monomials appearing the $q$--characters of $W(\varpi_\ell)$. In fact,
the bijection is simply given by putting $Y_{0,*} = 1$ in $m_{T;0}$.

\section{Finite dimensional crystals -- classical types}\label{sec:example}

In this section we treat all classical types.

\subsection{Type $D^{(1)}_n$}\label{Dn}

Let $\mathbf B = \{ 1,\dots,n,\overline{n},\dots, \overline{1}\}$. 
We give the ordering $\prec$ on the set $\mathbf B$
by
\begin{equation*}
  1 \prec 2 \prec \cdots \prec n-1 \prec
  \begin{matrix}n \\ \\ \overline{n}\end{matrix}
 \prec \overline{n-1} \prec\cdots \prec \overline{2}\prec \overline{1}.
\end{equation*}
Remark that there is no order between $n$ and $\overline{n}$.

For $p\in\ZZ$, mimicking the definition in \cite{Nac,kks}, we define 
{\allowdisplaybreaks
\begin{equation*}
\begin{aligned}[c]
& \ffbox{1}_p = Y_{0,p+2}^{-1} Y_{1,p}, \quad
  \ffbox{2}_p = Y_{0,p+2}^{-1} Y_{1,p+2}^{-1} Y_{2,p+1},
\\
& \ffbox{i}_p = Y_{i-1,p+i}^{-1} Y_{i,p+{i-1}} \qquad(3 \le i \le n-2),
\\
& \ffbox{n\!-\!1}_p
 = Y_{n-2,p+{n-1}}^{-1} Y_{n-1,p+{n-2}} Y_{n,p+{n-2}}, \qquad \ffbox{\overline{n\!-\!1}}_p = Y_{n-2,p+{n-1}} Y_{n-1,p+{n}}^{-1}
Y_{n,p+{n}}^{-1},
\\
& \ffbox{n}_p = Y_{n-1,p+{n}}^{-1} Y_{n,p+{n-2}},\qquad \ffbox{\overline{n}}_p = Y_{n-1,p+{n-2}} Y_{n,p+{n}}^{-1},
\\
& \ffbox{\overline{i}}_p = Y_{i-1,p+{2n-2-i}} Y_{i,p+{2n-1-i}}^{-1}
\qquad(3 \le i \le n-2),
\\
& \ffbox{\overline{2}}_p = Y_{0,p+2n-4} Y_{1,p+2n-4} Y_{2,p+2n-3}^{-1},
\qquad
\ffbox{\overline{1}}_p = Y_{0,p+2n-4} Y_{1,p+{2n-2}}^{-1}.
\end{aligned}
\end{equation*}

We define the {\it $i$-grade\/} $\gr_i(\ffbox{*}_p)$ as the grade of
the variable $Y_{i,*}$ appearing in $\ffbox{*}_p$. If $Y_{i,*}$ does
not appear, it is not defined. As variables appear at most once, it is
well-defined. When the suffix is clear from the context, we may omit
it and simply write $\gr_i(\ffbox{*})$.

\subsubsection{}\label{Dn1} First consider the case $\ell=1$. 
\begin{NB}
I moved this subsection here. Otherwise the reader cannot understand the notation.
\end{NB}
We take $M = Y_{1,p} Y_{0,p+2}^{-1}$. It
follows from Corollary~\ref{corex} that $\mathcal M(M) \simeq \mathcal
B(\varpi_1)$. 
The} crystal graph of $\mathcal M(M)$ is given in Figure~\ref{fig:Dn}.
\begin{figure}[htbp]
\centering
\psset{xunit=.55mm,yunit=.55mm,runit=.55mm}
\psset{linewidth=0.3,dotsep=1,hatchwidth=0.3,hatchsep=1.5,shadowsize=1}
\psset{dotsize=0.7 2.5,dotscale=1 1,fillcolor=black}
\psset{arrowsize=1 2,arrowlength=1,arrowinset=0.25,tbarsize=0.7 5,bracketlength=0.15,rbracketlength=0.15}
\begin{pspicture}(0,0)(220,50)
\rput(0,20){$\ffbox{1}_p$}
\rput(20,20){$\ffbox{2}_p$}
\psline{->}(5,20)(13,20)
\psline{->}(25,20)(33,20)
\psline{->}(75,20)(83,20)
\psline[linestyle=dashed,dash=1 1](40,20)(70,20)
\rput(90,20){$\ffbox{n\!-\!1}_p$}
\rput(110,30){$\ffbox{n}_p$}
\rput(110,10){$\ffbox{\overline{n}}_p$}
\rput(130,20){$\ffbox{\overline{n\!-\!1}}_p$}
\psline{->}(135,20)(143,20)
\psline{->}(185,20)(193,20)
\psline[linestyle=dashed,dash=1 1](150,20)(180,20)
\rput(200,20){$\ffbox{\overline{2}}_p$}
\psline{->}(205,20)(213,20)
\rput(220,20){$\ffbox{\overline{1}}_p$}
\psline{->}(95,16)(104,12)
\psline{->}(115,29)(124,25)
\psline{->}(95,25)(104,29)
\psline{->}(115,12)(124,16)
\rput{0}(0,25){\qline(-0,0)(0,0)}\rput{0}(100,-480){\parametricplot[arrows=->]{78.85}{101.15}{ t cos 517.26 mul t sin 517.26 mul }}
\rput{0}(117.26,520){\parametricplot[arrows=<-]{-101.15}{-78.85}{ t cos 517.26 mul t sin 517.26 mul }}
\rput(8,25){$\scriptscriptstyle 1$}
\rput(30,25){$\scriptscriptstyle 2$}
\rput(95,30){$\scriptscriptstyle n-1$}
\rput(138,25){$\scriptscriptstyle n-2$}
\rput(95,10){$\scriptscriptstyle n$}
\rput(122,10){$\scriptscriptstyle n-1$}
\rput(125,30){$\scriptscriptstyle n$}
\rput(78,25){$\scriptscriptstyle n-2$}
\rput(190,25){$\scriptscriptstyle 2$}
\rput(210,25){$\scriptscriptstyle 1$}
\rput(100,41){$\scriptscriptstyle 0[2n-4]$}
\rput(185,3){$\scriptscriptstyle 0[2n-4]$}
\end{pspicture}
\caption{(Type $D_n^{(1)}$) the crystal $\mathcal B(\varpi_1)$ of the vector representation}
\label{fig:Dn}
\end{figure}

We have 
\begin{equation*}
  \begin{split}
    & \tilde f_0(\ffbox{\overline{2}}_p) = 
     Y_{1,p+2n-4}Y_{0,p+2n-2}^{-1} = \ffbox{1}_{p+2n-4},
\\
    & \tilde f_0(\ffbox{\overline{1}}_p) = 
    Y_{0,p+2n-2}^{-1} Y_{1,p+2n-2}^{-1}Y_{2,p+2n-3} = \ffbox{2}_{p+2n-4}.
  \end{split}
\end{equation*}
Therefore $\mathcal M(M)$ is preserved under $\tau_{2n-4}$. Computing
weights as above, we find that $z_1=\tau_{4-2n}$ and so we have
$\mathcal M(M)/\tau_{2n-4} \simeq \mathcal B(W(\varpi_1))$.
We also get an isomorphism of $I_0$-crystals $\mathcal
B(W(\varpi_1))\simeq \mathcal B_{I_0}(\varpi_1)$.

\subsubsection{Preliminary results for crystals of finite type $D$}\label{preld}
As is illustrated in examples in type $A^{(1)}_n$, we first need to
describe the $I_0$-crystal structure on the monomials. This will be
given in this subsection. All the results on the $I_0$-crystal are
independent of the information on $Y_{0,*}$,  so we set $Y_{0,*}$ as
$1$ in this subsection. 
Note also that results can be modified in an obvious manner so that
the suffixes of $\ffbox{}_*$ can be shifted simultaneously. We will
use the results in these modified forms in later subsections.

\begin{thm}\label{find} Let $1\leq \ell\leq n-2$ and
\begin{equation*}
M = \ffbox{1}_{\ell-1}\ffbox{2}_{\ell-3}\cdots\ffbox{\ell}_{1-\ell}.
\end{equation*}
Then $\mathcal M_{I_0}(M)$ is isomorphic to $\mathcal B_{I_0}(\varpi_\ell)$ and is equal to the set of monomials 
\begin{equation*}
m_T = \ffbox{i_1}_{\ell-1}\ffbox{i_2}_{\ell-3}\cdots\ffbox{i_\ell}_{1-\ell},
\end{equation*}
indexed by the set $D_{\ell,0,0}$ of tableaux $T = (i_1,\dots,i_{\ell})$ satisfying the conditions

\textup{(1)} $i_a\in\mathbf B,  i_1 \nsucceq i_2 \nsucceq \cdots \nsucceq i_\ell$,

\textup{(2)} there is no pair $a$, $b$ such that $1\leq a < b \leq \ell$ and $i_a=k$, $i_b=\overline{k}$
                and $b-a = n -1 - k$.

\noindent
Moreover the map $T\mapsto m_T$ defines a bijection between $D_{\ell,0,0}$ and $\mathcal{M}_{I_0}(M)$.\end{thm}

\begin{NB}
  Please distinguish $\omega_\ell$ and $\varpi_\ell$.
\end{NB}

This result follows from \cite[3.4, 5.5]{Nac}. It was also proved by
Kang-Kim-Shin \cite{kks} in the present form. We briefly recall their
argument for a later purpose. They checked the following statements:
\begin{aenume}
\item The set of monomials $m_T$ with $T$ satisfying (1), but not
  necessarily (2), is preserved by $\te_i$, $\tf_i$.

\item If a monomial $m_T$ satisfies $\te_i m_T = 0$ for all $i=1,\dots, n$, then $m_T$ must be equal to $M$.

\item For a tableau $T$ satisfying (1), there exists a tableau $T'$
satisfying (1),(2) and $m_T = m_{T'}$.

\item The tableau $T$ satisfying (1),(2) is uniquely determined from
  the monomial $m_T$.
\end{aenume}

The statement (d) is not explicitly stated in \cite{kks}, but follows
from \cite[Prop.~3.2]{kks} or the argument below.

Let us give an example for the procedure (c). Suppose $n=7$ and 
\(
  T = (2,3,4,\overline{3},\overline{2}).
\)
Using the relation 
\(
  \ffbox{k}_{p}\ffbox{\overline{k}}_{p-2(n-1-k)}
  = \ffbox{k\!+\!1}_{p}\ffbox{\overline{k\!+\!1}}_{p-2(n-1-k)}
\)
several times, we get
\begin{equation*}
  \ffbox{2}_4\ffbox{3}_2 \ffbox{4}_0
  \ffbox{\overline{3}}_{-2} \ffbox{\overline{2}}_{-4}
  = 
  \ffbox{4}_4\ffbox{5}_2 \ffbox{6}_0
  \ffbox{\overline{6}}_{-2} \ffbox{\overline{5}}_{-4}.
\end{equation*}
Thus $T' = (4,5,6,\overline{6},\overline{5})$. In general, we replace
the pair
\(
  \ffbox{k}_{p}\ffbox{\overline{k}}_{p-2(n-1-k)}
\)
by
\(
  \ffbox{k\!+\!1}_{p}\ffbox{\overline{k\!+\!1}}_{p-2(n-1-k)}
\)
repeatedly from $k=1$ to $n-2$.

As we saw in examples in type $A^{(1)}_n$, we need to study a tableau
whose suffixes may {\it jump}. 
For $1\le\ell\le n-2$, $0\le r\le n-\ell-1$, $0\leq h\leq \ell$
let 
\begin{NB}
Note that $r=0$ is allowed now.
\end{NB}
\begin{equation*}
  \begin{split}
  & M_{\ell,h,r} = Y_{h,\ell-h} Y_{h,\ell-h-2r}^{-1} Y_{\ell,-2r}
\\
   =\; &
   \left(\ffbox{1}_{\ell-1} \ffbox{2}_{\ell-3} \cdots
     \ffbox{h}_{\ell-2h+1}
    \right)
  \times
   \left(\ffbox{h\!+\!1}_{\ell-2h-2r-1}
     \ffbox{h\!+\!2}_{\ell-2h-2r-3}
     \cdots \ffbox{\ell}_{1-\ell-2r}\right)
\\
   = \; &
   \prod_{p=1}^h \ffbox{p}_{\ell-2p+1} \times
   \prod_{p=h+1}^\ell \ffbox{p}_{\ell+1-2p-2r}
  \end{split}
\end{equation*}
and consider a monomial
\begin{equation*}
   m_T = \left(\ffbox{i_1}_{\ell-1}\ffbox{i_2}_{\ell-3}\cdots
     \ffbox{i_h}_{\ell-2h+1}\right)\times
    \left(
      \ffbox{i_{h\!+\!1}}_{\ell-2h-2r-1}\ffbox{i_{h\!+\!2}}_{\ell-2h-2r-3}
    \cdots\ffbox{i_\ell}_{1-\ell-2r}\right)
\end{equation*}
appearing in $\mathcal M_{I_0}(M_{\ell,h,r})$. When $h = 0$ or $\ell$,
these are obtained from $\mathcal M_{I_0}(M)$ in Theorem~\ref{find} by
the simultaneous shift of grades.

We should consider $T$ as a tableau of shape $(h,\ell-h)$ (one column with $h$ boxes and one column with $\ell - h$ boxes), where the
second column is shifted below by $h+r$ boxes. But we simply denote it
by $T = ((i_1,\dots,i_h),(i_{h+1},\dots,i_\ell))$ or by $T =
(i_1,\dots,i_\ell)$ for the sake of spaces.

From the proof of Theorem~\ref{find} in \cite{kks}, we have
\begin{equation}
\tag{D.1}
\begin{split}
   & i_1 \nsucceq i_2 \nsucceq \cdots \nsucceq i_h,
\\
   & i_{h+1}\nsucceq i_{h+2}\nsucceq \cdots \nsucceq i_\ell.    
\end{split}
\end{equation}

\begin{NB}
I have changed \verb+(\thesubsection.1)+ to D.1. (D stands for type $D$.)  For other types, use B.1, etc.
\end{NB}

Let us study the order between $i_h$ and $i_{h+1}$. The following
example shows that $i_h\nsucceq i_{h+1}$ may not be satisfied in
general: Let $n=7$, $\ell =3$, $r=n-\ell-2=2$.  Consider the starting
monomial $M=\ffbox{1}_2\ffbox{2}_{-4}\ffbox{3}_{-6}$. It gives in the
crystal $\mathcal{M}_{I_0}(M)$ the monomial
$m=\ffbox{3}_2\ffbox{4}_{-4}\ffbox{\overline{4}}_{-6}
=Y_{2,5}^{-1}Y_{3,4}Y_{3,0}^{-1}Y_{4,-1}Y_{3,2}Y_{4,3}^{-1}$.  If we
apply $\tilde{f}_3$, we get the monomial
$m'=Y_{3,6}^{-1}Y_{4,5}Y_{3,0}^{-1}Y_{4,-1}Y_{3,2}Y_{4,3}^{-1}$ and
this monomial can only be written in the form
$m'=\ffbox{4}_2\ffbox{4}_{-4}\ffbox{\overline{4}}_{-6}$.

The condition~(2) in Theorem~\ref{find} also needs to be modified
for the pair $a$, $b$ with $a\leq h$, $h+1\leq b$ as the suffix jump.
A naive guess is to replace $n-1-k$ by $n-r-k-1$, but this change does
not work as indicated by the following example:
Consider the case $n=6$, $\ell=4$, $r=1$ and the starting
monomial $m=\ffbox{1}_3\ffbox{2}_{-1}\ffbox{3}_{-3}\ffbox{4}_{-5}$.
Then consider the monomial
$m'=\ffbox{1}_3\ffbox{2}_{-1}\ffbox{3}_{-3}\ffbox{\overline{1}}_{-5}=Y_{1,3}Y_{1,5}^{-1}Y_{1,1}^{-1}Y_{3,-1}$.
For $b=4$ and $a=1$, we have $b-a=n-r-k-1=3$.
Thus this monomial violates the condition (2) of Theorem~\ref{find}.
But if we replace the pair $(i_1, i_4) = (1,\overline{1})$ to
$(2, \overline{2})$ as before, we get
$\ffbox{2}_3\ffbox{2}_{-1}\ffbox{3}_{-3}\ffbox{\overline{2}}_{-5}$,
which does not satisfies the condition~(1) of Theorem~\ref{find}.
In the original situation we can further replace the pair
$(i_2,i_4) = (2,\overline{2})$ to $(3, \overline{3})$,
and then further $(i_3,i_4) = (3,\overline{3})$ to $(4,\overline{4})$
to achieve the condition (1). But we cannot make this replacement
as $\ffbox{2}_{-1}\ffbox{\overline{2}}_{-5}\neq
\ffbox{3}_{-1}\ffbox{\overline{3}}_{-5}$.

We modify the condition~(2) as follows.
\begin{enumerate}
\def\labelenumi{(D.\theenumi)}
\def\theenumi{\arabic{enumi}}
\addtocounter{enumi}{1}
\item There is no pair $a$, $b$ such that $1\leq a < b \leq h$ and
  $i_a=k$, $i_b=\overline{k}$ and $b-a = n -1 - k$.

\item There is no pair $a$, $b$ such that $h+1\leq a < b \leq \ell$
  and $i_a=k$, $i_b=\overline{k}$ and $b-a = n -1 - k$.

\item There is no pair $a$, $b$ such that $a \leq h$ , $h+1 \leq b$ , $i_a=k$, $i_b=\overline{k}$ and $b-a = n - \max(r,1) - k$.
\end{enumerate}

The conditions (D.2,3) can be achieved without changing the
corresponding monomial by the procedure explained above. For (D.4)
(when $r\ge 1$), we replace a pair $(i_a,i_b) = (k,\overline{k})$ with
$b-a = n - \max(r,1) - k$ by $(k-1,\overline{k-1})$. If there are
several such pairs or this procedure yields a new such pair, we
replace them repeatedly starting from $k=n-1$, then $k=n-2$,..., and
finally to $k=2$. (Note that this is converse to the order of the
procedure for (D.2,3).) As $r\le n-\ell-1$, the condition (D.4) always
holds for $k=1$.

Our approach to determine all monomials appearing in $\mathcal
M_{I_0}(M_{\ell,h,r})$ is to relate them to monomials in $\mathcal
M_{I_0}(M_{\ell,h,r-1})$. Since we understand the case $r=0$, we know
a general case inductively.

In order to accomplish this approach, we first remark that the crystal
structure on the monomials can be transfered to that on the tableaux
satisfying $(\mathrm{D}.1\sim 4)$.

\begin{lem}\label{lem:tableauxcrystal}
There exists a unique crystal structure on the set of tableaux $T$
satisfying $(\mathrm{D}.1\sim 4)$ such that $T\mapsto m_T$ is a strict
morphism, i.e., it preserves $\varepsilon_k$, $\varphi_k$, $\wt$ and
commutes with $\tf_k$, $\te_k$.
\end{lem}

\begin{proof}
We transfer $\varepsilon_k$, $\varphi_k$, $\wt$ on monomials
to those on tableaux via $T\mapsto m_T$.

Let us define $\tf_k$ on tableaux. ($\te_k$ can be defined in the same
way.)
In general, $\tf_k m_T \neq 0$ can be written as $m_{T'}$ for a
tableau $T'$ which is obtained by replacing an entry $i_a$ in $T$ by a
new one according to the rule described in Figure~\ref{fig:Dn}.
To define $\tf_k$ on tableaux, we need to specify the entry $i_a$ to
be replaced. There might be ambiguity when we have a pair $(i_a,i_b) =
(k,\overline{k+1})$ with $\gr_k(\ffbox{i_a}) =
\gr_k(\ffbox{i_b})$. This happens when
\(
   b - a = n - 1 - k
\)
for $a,b\le h$ or $h+1\le a,b$ and
\(
  b - a = n - 1 - k - r
\)
for $a\le h$, $h+1\le b$.
In the first case (or the second case with $r = 0$) we replace
$k$ by $k+1$. In the second case with $r\neq 0$
we replace $\overline{k+1}$ by $\overline{k}$. Note that we are forced
to take these choices by $(\mathrm{D}.2\sim 4)$.
Now the assertion is clear.
\end{proof}

Let us prove the statement (d) after Theorem~\ref{find} as we
promised. From (a),(c) we have a surjective map $T\mapsto m_T$. Since
it commutes with $\te_i$ and $\tf_i$, the injectivity follows if we
check that $\te_i T = 0$ for all $i$ implies $T = (1,\dots,\ell)$. But
the proof of the statement (b) in \cite{kks}, in fact, gives this
statement.

Let us next define a map $\sigma_{\ell,h,r}$ from tableaux satisfying
$(\mathrm{D}.1\sim 4)$ to those where we increase $r$ by $1$, i.e.,
each $\ffbox{i_c}_{\ell-2r-2c+1}$ is replaced by
$\ffbox{i_c}_{\ell-2r-2c-1}$ for $c\ge h+1$.
Almost all the cases, $\sigma_{\ell,h,r}(T)$ is just $T$.
But the condition (D.4) is violated if there is a pair $(i_a, i_b) =
(k,\overline{k})$ such that $a\le h$, $h+1\le b$ and $b-a = n - r - k
-1$. We replace it by $(k+1,\overline{k+1})$. If there are several
such pairs or this procedure yields a new such pair, we replace them
repeatedly starting from $k=1$ to $n-r-1$.
We define $T' = \sigma_{\ell,h,r}(T)$ as the final result.
As we have
\begin{equation*}
  \ffbox{k}_{\ell-2a+1} \ffbox{\overline{k}}_{\ell-2r - 2b + 1}
  = \ffbox{k\!+\!1}_{\ell-2a+1}\ffbox{\overline{k\!+\!1}}_{\ell-2r-2b+1},
\end{equation*}
the procedure keeps the corresponding monomial unchanged if we do not
change $r$ for the map $T\mapsto m_T$.

Let us check that $\sigma_{\ell,h,r}$ intertwines $\tf_k$.
By definition, $\sigma_{\ell,h,r}\tf_k T$ is possibly different 
from $\tf_k\sigma_{\ell,h,r}T$ if there is a
pair $(i_a,i_b)$ with $a\le h$, $h+1\le b$ such that the order of $k$-grades $p =
\gr_k(\ffbox{i_a})$, $q = \gr_k(\ffbox{i_b})$ are changed by $\sigma_{\ell,h,r}$.
If both $i_a$ and $i_b$ contribute to $Y_{k,*}$ in positive or
negative powers, the rule for $\tf_k T$ is changed accordingly. (See
the proof of Lemma~\ref{lem:tableauxcrystal} how $\tf_k T$ is
defined.)
Thus it is enough to study the case when one contributes in positive,
and the other in negative.
For $k = n-1$, $n$ such a change cannot occur. As grades can only be shifted by $2$, for $k\le n-2$ we have
a possible change only when $p + 2 = q$ for $(i_a,i_b) =
(k,\overline{k})$, and $p = q$ for $(i_a,i_b) = (k+1,\overline{k+1})$.
These are equivalent to
\begin{equation*}
\begin{cases}
   b - a = n - r - k - 1 &\text{if $(i_a,i_b) = (k,\overline{k})$},
\\
   b - a = n - r - k - 2 &\text{if $(i_a,i_b) = (k+1,\overline{k+1})$}.
\end{cases}
\end{equation*}
Therefore if there is no pair $(i_a,i_b) = (k,\overline{k})$ with
$a\le h$, $h+1\le b$ and $b -a = n - r - k - 1$ for any $k$, then
$\tf_k m_T$ is unchanged when we increase $r$ by $1$.
But we have defined $\sigma_{\ell,h,r}$ exactly so that this condition
is achieved. Thus we have
\begin{equation*}
   \sigma_{\ell,h,r} \tf_k T = \tf_k \sigma_{\ell,h,r}{T}
   \qquad \text{for all $k\in I_0$}.
\end{equation*}
This equality holds even if $\tf_k T = 0$.

Similarly we define $\sigma'_{\ell,h,r}(T)$ as follows. When
$r=1$, we simply set it $T$. Assume $r > 1$ hereafter.
Suppose that there is a pair $(i_a,i_b) = (k,\overline{k})$ such that
$a\le h$, $h+1\le b$ and $b-a = n-r-k+1$. We replace it by
$(k-1,\overline{k-1})$. If there are several such pairs or this
procedure yields a new such pair, we replace them repeatedly starting
from $k=n-r$ to $3$. We define $\sigma'_{\ell,h,r}(T)$ as the final
result.
As $r\le n-\ell-1$, we have $\ell+k\ge b-a+k+1 = n - r + 2 \ge
\ell + 3$. Therefore $k \le 2$ cannot happen, so $k-1\in\mathbf B$.

These maps are somewhat similar to one defined in
\cite[Prop.~3.2]{kks}.

Now we introduce new conditions:
\begin{enumerate}
\def\labelenumi{(D.\theenumi)}
\def\theenumi{\arabic{enumi}}
\addtocounter{enumi}{4}
\item  Suppose that $i_{h+1} = k \in \{1,\dots, n-1\}$ and $i_h\succeq
i_{h+1}$. Then $i_h = k'$ is also in $\{ 1,\dots, n-1\}$,
and the successive part $(\overline{k'}, \overline{k'-1},\dots,
\overline{k})$ appears as $(i_{b'},i_{b'+1},\dots,i_{b})$ with
$n - r - k < b - h \le n- k - 1$.

\item Suppose that $i_{h+1} = \overline{k} \in \{\overline{1},\dots,
  \overline{n-1}\}$ and $i_h\succeq i_{h+1}$. Then $i_h =
  \overline{k'}$ is also in $\{ \overline{1},\dots, \overline{n-1}\}$,
  and the successive part $(k',k'+1,\dots,k)$ appears as
  $(i_{a'},i_{a'+1},\dots,i_{a})$ with
  $n-r-k \le h - a < n-k - 1$.

\item  If $i_{h+1} = n$ or $\overline{n}$, then 
$i_h \nsucceq i_{h+1}$.
\end{enumerate}

Note that (D.1) implies that the successive part in (D.5) occurs in $b'
> h+1$. This together with the second inequality (and $b+k = b'+k'$)
implies $k' < n-2$. Thus $i_h=n-1, n-2$ cannot happen in (D.5).
Similarly $i_{h+1}=\overline{n-1}, \overline{n-2}$, cannot happen in
(D.6).

\begin{defi}
Let $D_{\ell,h,r}$ be the set of tableaux $T$ satisfying
$(\mathrm{D}.1\sim 7)$.
\end{defi}

\begin{rem}
  When $r=0$, the conditions $(\mathrm{D}.1\sim 7)$ are equivalent to
  (1),(2) in Theorem~\ref{find}.
\end{rem}

\begin{prop}\label{prop:sigma}
$\sigma_{\ell,h,r}$ defines a crystal isomorphism from
$D_{\ell,h,r}$ to $D_{\ell,h,r+1}$. Its inverse is given by
$\sigma'_{\ell,h,r+1}$.
\end{prop}

As a corollary we have
\begin{thm}\label{isot}
The map $T \mapsto m_T$ induces a crystal isomorphism between
$D_{\ell,h,r}$ and $\mathcal M_{I_0}(M_{\ell,h,r})$.
\end{thm}

\begin{proof}
We first prove that the image of $D_{\ell,h,r}$ is contained in
$\mathcal M_{I_0}(M_{\ell,h,r})$ by the induction on $r$. This is true
for $r=0$ by Theorem~\ref{find}. Suppose it is true for $r$.
First note that
$\sigma_{\ell,h,r}$ maps
$T = (1,\dots,\ell)$ to $(1,\dots,\ell)$.
Take $T\in D_{\ell,h,r+1}$. By the induction hypothesis
$m_{\sigma'_{\ell,h,r+1}(T)}$ can be written as 
\[
   m_{\sigma'_{\ell,h,r+1}(T)} =
   \tf_{i_1}\tf_{i_2}\cdots \tf_{i_N} M_{\ell,h,r}
\]
for $N \ge 0$, $i_p\in I_0$. We then have
\[
   m_{T} = \tf_{i_1}\tf_{i_2}\cdots \tf_{i_N} M_{\ell,h,r+1}.
\]
This shows $m_T\in \mathcal M_{I_0}(M_{\ell,h,r+1})$.

As the crystal graph of $\mathcal M_{I_0}(M_{\ell,h,r})$ is connected
by its definition, the map is surjective.
\begin{NB}
Any $m\in \mathcal M_{I_0}(M_{\ell,h,r})$ can be written as $m =
\tf_{i_1}\tf_{i_2}\cdots \tf_{i_N} M_{\ell,h,r}$ for $N \ge 0$,
$i_p\in I_0$. We have
\(
   \tf_{i_1}\tf_{i_2}\cdots \tf_{i_N} (1,\dots,\ell)
   \mapsto m.
\)
\end{NB}

By the induction on $r$, it follows that the only tableau $T$ with
$\te_i T = 0$ for all $i\in I_0$ is the highest one $T =
(1,\dots,\ell)$. This shows that the strict crystal morphism $T\mapsto
m_T$ is injective.
\end{proof}

\begin{proof}[Proof of Proposition~\ref{prop:sigma}]
It is enough to show that $\sigma_{\ell,h,r}$ is a set theoretical
bijection, as we already observed that it is a strict crystal morphism.

When $r=0$, there is no pair $(i_a,i_b) = (k,\overline{k})$ to replace
by $(\mathrm{D}.2\sim 4)$. Thus $\sigma_{\ell,h,0}$ is just an identity.
Also $\sigma'_{\ell,h,1}$ is an identity by definition.
On the other hand, the conditions $(\mathrm{D}.1\sim 7)$ are the same
for $r= 0$ and $1$. Therefore the assertion is true for $r=0$. We assume
$r > 0$ hereafter.

Suppose $T$ satisfies $(\mathrm{D}.1\sim 7)$. We show that 
$\sigma_{\ell,h,r}(T)$ also satisfies $(\mathrm{D}.1\sim 7)$.
The condition (D.1) is clearly satisfied. The condition (D.4) with $r$
replaced by $r+1$ is satisfied by the definition of
$\sigma_{\ell,h,r}$.

We study the cases $i_{h}\succeq i_{h+1}$ and $i_h\nsucceq i_{h+1}$
separately.

Case (1): $i_h\succeq i_{h+1}$.

We assume $i_{h+1} = k \in \{ 1,\dots,n-1\}$. By (D.5) $i_h = k' \in
\{1,\dots, n-1\}$ and there exists a successive part
$(\overline{k'},\dots,\overline{k}) = (i_{b'},\dots,i_b)$ with $h+1\le
b'$, $n-r-k < b-h \le n-1-k$. The condition (D.2) automatically holds
as $i_h\in \{1,\dots,n-1\}$.

Suppose that $i_h$ is replaced during the procedure. Then in the
middle of the procedure, we find an entry $i'_B$ with $i'_B =
\overline{k}$, $B\ge h+1$, $B - h = n - r - k - 1$. As $i'_B$ is
obtained by replacing $i_B$, we have $i'_B \preceq i_B$. Therefore $B
\ge b$. But this contradicts with (D.5) as
\[
   n - r- k - 1 = B - h \ge b - h > n - r - k.
\]
Therefore $i_h$ remains unchanged during the procedure. Therefore the
procedure is performed for pairs $(K,\overline{K})$ with $K < k$, so
all $(i_h, i_{h+1},\dots,i_b)$ are also unchanged. Thus (D.5)
remains true.
Suppose (D.3) is violated, i.e., there exists $(i_A, i_B) =
(K,\overline{K})$ with $B > A \ge h+1$, $B-A = n-1-K$. As $K \ge k$,
such a pair can appear only in $(i_{h+1},\dots,i_b)$. But this
part is unchanged, so (D.3) for $r$ implies that this cannot happen.
Thus (D.3) is also satisfied.

We can similarly check the assertion when $i_{h+1}\in \{\overline{1},\dots,
\overline{n-1}\}$.

Case (2): $i_h \nsucceq i_{h+1}$.

Suppose that we apply the above procedure to a tableau $T =
((i_1,\cdots,i_h),(i_{h+1},\cdots,i_\ell))$ to get a new tableau
$T'=((j_1,\cdots,j_h),(j_{h+1},\cdots,j_\ell))$.
We separate the cases according to the order among $j_h$ and
$j_{h+1}$.

Subcase (2.1): $j_h \succeq j_{h+1}$ and $i_{h+1}\in \{1,\dots,n\}$.

As $i_{h+1}$ is unchanged, $j_h \succeq j_{h+1}$ can happen only when
$i_h$ is replaced during the procedure.  Suppose that $i_h$ is
replaced from $k'$ to $m$ with $m \ge k'+1$. Then the procedure yields a
successive part $(j_b,\dots,j_{b''}) = (\overline{m},\dots,
\overline{k'+1})$ with $b - h = n - r - m$.
We have 
\[
   m = j_h \succeq j_{h+1} = i_{h+1} \npreceq i_h = k'
\]
Thus $\overline{j_{h+1}}$ can appear only in the successive part, so
(D.5) is satisfied with $r$ replaced by $r+1$.

The condition (D.2) is automatic. Suppose that (D.3) is violated,
i.e., there exists $(j_A, j_B) = (K,\overline{K})$ with $B > A > h$,
$B-A = n-1-K$. We have $K \ge j_{h+1} = i_{h+1} \npreceq i_{h} =
k'$. Therefore $j_B$ can occur only in $B\le b''$. If $j_B$ appears
outside of the successive part, then $j_B = i_B$ and we have a
contradiction with (D.3) for the original tableau. If $j_B$ appears
in the successive part, we have
\[
   n - 1 - K = B - A < B - h = n - r - K.
\]
As $r \ge 1$, we have a contradiction.

Similarly we can check the assertion
$j_h \succeq j_{h+1}$ and $i_{h}\in \{\overline{1},\dots,\overline{n} \}$.
When $i_h\in \{1, \dots, n\}$, $i_{h+1}\in
\{\overline{1},\dots,\overline{n} \}$, the inequality $j_h \succeq
j_{h+1}$ cannot happen. Thus we checked the assertion when $j_h
\succeq j_{h+1}$.

Subcase (2.2): $j_h \nsucceq j_{h+1}$.

The conditions $(\mathrm{D}.5\sim 7)$ are satisfied by the assumption.
Let $(j_a, j_b) = (k+1,\overline{k+1})$ with $b-a = n-r-k-1$ be the
pair obtained by the last replacement in the procedure. We suppose
that (D.2) is violated, i.e., we have a pair $A < B\leq h$ such that
$j_A=K$ and $j_B=\overline{K}$ and $B-A=n-1-K$. As $i_c$ for $a < c <
b$ is unchanged by the above procedure, the condition (D.2) for $T$
implies that $j_A$ can appear only in $A\le a$. Then
$n-1-K=B-A=(a-A)+(B-b)+n-r-k-1$, so $K + a-A =b-B+r+k>k+1$. This
inequality contradicts with (D.1) as
\[
   k+1 = j_a \ge j_A + (a-A) > k+ 1.
\]
Thus (D.2) is satisfied. In the same way (D.3) is satisfied.

Next we show that $\sigma'_{\ell,h,r}(T)$ also satisfies
$(\mathrm{D}.1\sim 7)$. We may suppose $r\ge 2$. The condition (D.1)
is clearly satisfied. The condition (D.4) with $r$ replaced by $r-1$
is satisfied by the definition of $\sigma'_{\ell,h,r}$.

Suppose that we apply the above procedure to a tableau $T =
((i_1,\cdots,i_h),(i_{h+1},\cdots,i_\ell))$ to get a new tableau
$T'=((j_1,\cdots,j_h),(j_{h+1},\cdots,j_\ell))$.  Let $(i_a,i_b) =
(k,\overline{k})$ with $a\leq h$, $h+1 \leq b$, $b-a = n - r - k + 1$
be the first pair replaced in the procedure.
Suppose that (D.2) is violated, i.e., we have a pair $A < B\leq h$
such that $j_A=K$ and $j_B=\overline{K}$ and $B-A=n-1-K$. As $i_c$ for
$a< c< b$ is unchanged by the above procedure, we have $A\le a$. If
$i_A = j_A$, i.e., $i_A$ is not unchanged, we have a contradiction
with (D.2) for $T$. Therefore $i_A \ge j_A + 1$.
We have $n-1-K=B-A=(a-A)+(B-b)+n-r-k+1$, so $K + a-A=b-B+r+k-2>k-1$.
This inequality contradicts with (D.1) as
\[
   k = i_a \ge i_A + (a-A) \ge j_A + 1 + (a -A ) > k.
\]
So (D.2) is satisfied by $T'$. In the same way (D.3) is satisfied by $T'$.

In order to check the remaining conditions, we treat the cases
separately according the ordering among $i_h$, $i_{h+1}$.

Case (a): $i_h\nsucceq i_{h+1}$.

This inequality is preserved during the procedure. Therefore
we have $j_h\nsucceq j_{h+1}$, so $(\mathrm{D}.5\sim 7)$ are preserved.

Case (b): $i_h\succeq i_{h+1}$ and $i_{h+1} = k\in \{1,\dots,n-1\}$.

Take the successive part $(\overline{k'},
\overline{k'-1},\dots, \overline{k}) = (i_{b'},i_{b'+1},\dots,i_{b})$
with $i_h = k'$ as in (D.5).
Suppose that an entry in the successive part is replaced during the
procedure, i.e., we replace a pair $(i_A,i_B) = (K,\overline{K})$ with
$A \le h$, $b' \le B \le b$ with $B - A = n-r-K+1$.
The inequality in (D.5) implies
\[
  n - r - K + 1\le b - h + k - K = B - h \le B - A.
\]
So this can happen only when two inequalities are equalities, i.e., 
$n-r-k+1 = b-h$ and $A = h$. And in such case, we really replace the
pair by the definition of $\sigma'_{\ell,h,r}$.

Subcase (b.1): $i_h$ is unchanged.

As we observed above, the successive part remains unchanged. By (D.5)
we have $n-r-k < b - h \le n - k - 1$. And the case $b - h = n - r -
k+1$ is excluded as we have just observed. Therefore the left hand
side of the inequality can be improved to $n-r-k+1$. This shows that
(D.5) with $r$ replaced by $r-1$ is satisfied.

Subcase (b.2): $i_h$ is changed.

Suppose that $i_{h}$ is changed, say from ${k'}$ to $j_{h}
= {m}$ with $m \le k'-1$. Then $i_{b'} = \overline{k'}$ is replaced by
$\overline{k-1}$, $i_{b'+1}$ is replaced by $\overline{k-2}$, and so on.
This procedure continues at least until we replace $i_b$ by
$\overline{k-1}$. Thus $m < k$. This is equivalent to
$j_h < j_{h+1}$. Thus we have $(\mathrm{D}.2\sim 4)$.

If $i_h \nsucceq i_{h+1}$, we get $j_h\nsucceq j_{h+1}$ as 
the procedure preserves this inequality. Thus 
If $i_h \succeq i_{h+1}$, we have a successive part $(k',k'+1,\dots,k)
= (i_{a'},i_{a'+1},\dots,i_{a})$ with $i_{h+1} =
\overline{k}$. Therefore the procedure continues at least until
$i_{h+1}$ is replaced by $k'-1$, i.e., $m < k'$. Therefore $j_h \nsucceq
j_{h+1}$.

Case (c): $i_h\succeq i_{h+1}$ and $i_{h+1}\in
\{\overline{1},\dots,\overline{n}\}$.

This case can be proved in the same way as in case (b).

Finally it is clear that $\sigma_{\ell,h,r}$ and
$\sigma'_{\ell,h,r+1}$ are mutually inverse. All replaced pairs $(k+1,
\overline{k+1})$ are returned back to $(k,\overline{k})$. And we do
not have extra replacements by (D.4). 
\end{proof}

When $r=0$, $\mathcal M_{I_0}(M_{\ell,h,0})$ is independent of $h$.
Therefore we get a crystal isomorphism between any pair $\mathcal
M_{I_0}(M_{\ell,h,r})$ and $\mathcal M_{I_0}(M_{\ell,h',r'})$ as a composite
of various $\sigma_{\ell,h'',r''}$ and $\sigma'_{\ell,h'',r''}$.

For a later purpose we explicitly write down the crystal isomorphism 
\begin{equation*}
   \tau_{\ell,h,r} \colon 
   D_{\ell,h,r} \cong \mathcal M_{I_0}(M_{\ell,h,r}) \to
   D_{\ell,h+1,r} \cong \mathcal M_{I_0}(M_{\ell,h+1,r}).
\end{equation*}
This is the composite
\(
   \sigma_{\ell,h+1,r-1}\sigma_{\ell,h+1,r-2}\cdots
   \sigma_{\ell,h+1,0} \sigma'_{\ell,h,1}
   \sigma'_{\ell,h,2} \cdots \sigma'_{\ell,h,r}.
\)
All replaced pairs $(k-1,\overline{k-1})$ are returned back to
$(k,\overline{k})$ except for those $i_{h+1} = \overline{k-1}$. Also
we may have extra replacements for $i_h = k-1$.

Let $T = ((i_1,\dots,i_h),(i_{h+1},\dots,i_\ell))$. We describe
$\tau_{\ell,h,r}(T)$ in the following three cases separately. 
\begin{enumerate}
\def\labelenumi{(D.\theenumi)}
\def\theenumi{\alph{enumi}}
\item $i_{h+1} = k \in \{1,\dots,n-1\}$ and there is an entry $i_b =
\overline{k}$ with $n-r-k< b - h \le n-1-k$.

\item $i_{h+1} = \overline{k} \in \{\overline{1},\dots,\overline{n-1}
\}$ and there is an entry $i_a = k$ with $n-r-k \le h - a < n - 1 - k$.

\item Neither (D.a) nor (D.b) is not satisfied. 
\end{enumerate}

In the case (D.c) we simply have
\begin{equation*}
\tau_{\ell,h,r}(T)=((i_1,\cdots, i_{h+1}),(i_{h+2},\cdots, i_\ell)).
\end{equation*}
\begin{NB}
Note that we have $i_h\nsucceq i_{h+1}$ in (D.c). So
$\tau_{\ell,h,r}(T)$ also satisfies (D.1). As we have $i_{h+1} \nsucceq
i_{h+2}$, $(\mathrm{D}.5\sim 7)$ are satisfied. 
The conditions (D.3) clearly holds. Let us check the condition (D.2).
The only possibility is $(i_a, i_{h+1}) = (k,\overline{k})$ with $a\le
h$. But $h+1-a = n-1-k$ cannot happen since we are not in case (D.b).
Similarly (D.4) holds as we are not in case (D.a).
\end{NB}

Next suppose we are in the case (D.a). As was explained in the
paragraph just after (D.7), the inequalities imply $b > h+1$ and
$k < n-2$. Starting from $i_b$, we go back $i_{b-1}$, $i_{b-2}$,
...  while entries are successive. Let $i_{b''}$ be the ending entry,
so $(i_{b''},i_{b''+1},\dots,i_b)$ are successive as
$(\overline{k''},\overline{k''+1},\dots,\overline{k})$ and
$i_{b''-1}\neq \overline{k''-1}$. Also by the same reasoning as above,
we have $k'' < n-2$.
We then have
\begin{equation*}
  \begin{split}
    \tau_{\ell,h,r}(T)=((i_1,\cdots&,i_h,k''\!+\!1),
\\
  & (i_{h+2},\cdots,i_{b''-1}, \overline{k''\!+\!1},\overline{k''},
  \cdots,\overline{k\!+\!1}, i_{b+1},\cdots, i_\ell)). 
  \end{split}
\end{equation*}

Similarly in the case (D.b), we take $i_{a''}$ so that
$(i_{a''},i_{a''+1},\dots, i_a) = (k'',k''+1,\dots,k)$ and $i_{a''-1}\neq
k''-1$. We have $k < n-2$. We then have
\begin{equation*}
  \begin{split}
    \tau_{\ell,h,r}(T)=((i_1,\cdots, i_{a''-1},k''\!-\!1,\cdots,k\!-\!1,
    i_{a+1},&\cdots,i_{h}, \\
    & \overline{k''\!-\!1}), (i_{h+2},\cdots, i_\ell)).
  \end{split}
\end{equation*}

\begin{NB}
The following subsection is unnecessary, as we changed the proof. But
we keep it for a possible future use.

\subsubsection{Preliminary results for crystals of finite type $D$}\label{preld}
As is illustrated in examples in type $A^{(1)}_n$, we first need to
describe the $I_0$-crystal structure on the monomials. This will be
given in this subsection. All the results on the $I_0$-crystal are
independent of the information on $Y_{0,*}$,  so we set $Y_{0,*}$ as
$1$ in this subsection. 
Note also that results can be modified in an obvious manner so that
the suffixes of $\ffbox{}_*$ can be shifted simultaneously. We will
use the results in these modified forms in later subsections.

\begin{thm}\label{find} Let $1\leq \ell\leq n-2$ and
\begin{equation*}
M = \ffbox{1}_{\ell-1}\ffbox{2}_{\ell-3}\cdots\ffbox{\ell}_{1-\ell}.
\end{equation*}
Then $\mathcal M_{I_0}(M)$ is isomorphic to $\mathcal B_{I_0}(\omega_\ell)$ and is equal to the set of monomials 
\begin{equation*}
m_T = \ffbox{i_1}_{\ell-1}\ffbox{i_2}_{\ell-3}\cdots\ffbox{i_\ell}_{1-\ell},
\end{equation*}
indexed by the set $D_{\ell,0,0}$ of tableaux $T = (i_1,\dots,i_{\ell})$ satisfying the conditions

\textup{(1)} $i_a\in\mathbf B,  i_1 \nsucceq i_2 \nsucceq \cdots \nsucceq i_\ell$,

\textup{(2)} there is no pair $a$, $b$ such that $1\leq a < b \leq \ell$ and $i_a=k$, $i_b=\overline{k}$
                and $b-a = n -1 - k$.

\noindent
Moreover the map $T\mapsto m_T$ defines a bijection between $D_{\ell,0,0}$ and $\mathcal{M}_{I_0}(M)$.\end{thm}

This result follows from \cite[3.4, 5.5]{Nac}. It was also proved by
Kang-Kim-Shin \cite{kks} in the present form. We briefly recall their
argument for a later purpose. They checked the following statements:
\begin{aenume}
\item The set of monomials $m_T$ with $T$ satisfying (1), but not
  necessarily (2), is preserved by $\te_i$, $\tf_i$.

\item If a monomial $m_T$ satisfies $\te_i m_T = 0$ for all $i=1,\dots, n$, then $m_T$ must be equal to $M$.

\item For a tableau $T$ satisfying (1), there exists a tableau $T'$
satisfying (1),(2) and $m_T = m_{T'}$.

\item The tableau $T$ satisfying (1),(2) is uniquely determined from
  the monomial $m_T$.
\end{aenume}

\begin{NB2}
I still need to check the statement (d) in \cite{kks}
\end{NB2}

Let us give an example for the procedure (c). Suppose $n=7$ and 
\(
  T = (2,3,4,\overline{3},\overline{2}).
\)
Using the relation 
\(
  \ffbox{k}_{p}\ffbox{\overline{k}}_{p-2(n-1-k)}
  = \ffbox{k\!+\!1}_{p}\ffbox{\overline{k\!+\!1}}_{p-2(n-1-k)}
\)
several times, we get
\begin{equation*}
  \ffbox{2}_4\ffbox{3}_2 \ffbox{4}_0
  \ffbox{\overline{3}}_{-2} \ffbox{\overline{2}}_{-4}
  = 
  \ffbox{4}_4\ffbox{5}_2 \ffbox{6}_0
  \ffbox{\overline{6}}_{-2} \ffbox{\overline{5}}_{-4}.
\end{equation*}
Thus $T' = (4,5,6,\overline{6},\overline{5})$. In general, we replace
the pair
\(
  \ffbox{k}_{p}\ffbox{\overline{k}}_{p-2(n-1-k)}
\)
by
\(
  \ffbox{k\!+\!1}_{p}\ffbox{\overline{k\!+\!1}}_{p-2(n-1-k)}
\)
repeatedly from $k=1$ to $n-2$.

As we saw in examples in type $A^{(1)}_n$, we need to study a tableau
whose suffixes may {\it jump}. 
For $1\le\ell\le n-2$, $1\le r\le n-\ell-1$, $0\leq h\leq \ell$
let 
\begin{equation*}
  \begin{split}
  & M_{\ell,h,r} = Y_{h,\ell-h} Y_{h,\ell-h-2r}^{-1} Y_{\ell,-2r}
\\
   =\; &
   \left(\ffbox{1}_{\ell-1} \ffbox{2}_{\ell-3} \cdots
     \ffbox{h}_{\ell-2h+1}
    \right)
  \times
   \left(\ffbox{h\!+\!1}_{\ell-2h-2r-1}
     \ffbox{h\!+\!2}_{\ell-2h-2r-3}
     \cdots \ffbox{\ell}_{1-\ell-2r}\right)
\\
   = \; &
   \prod_{p=1}^h \ffbox{p}_{\ell-2p+1} \times
   \prod_{p=h+1}^\ell \ffbox{p}_{\ell+1-2p-2r}
  \end{split}
\end{equation*}
and consider a monomial
\begin{equation*}
   m_T = \left(\ffbox{i_1}_{\ell-1}\ffbox{i_2}_{\ell-3}\cdots
     \ffbox{i_h}_{\ell-2h+1}\right)\times
    \left(
      \ffbox{i_{h\!+\!1}}_{\ell-2h-2r-1}\ffbox{i_{h\!+\!2}}_{\ell-2h-2r-3}
    \cdots\ffbox{i_\ell}_{1-\ell-2r}\right)
\end{equation*}
appearing in $\mathcal M_{I_0}(M_{\ell,h,r})$. When $h = 0$ or $\ell$,
these are obtained from $\mathcal M_{I_0}(M)$ in Theorem~\ref{find} by
the simultaneous shift of grades.

We should consider $T$ as a tableau of shape $(h,\ell-h)$, where the
second column is shifted below by $h+r$ boxes. But we simply denote it
by $T = ((i_1,\dots,i_h),(i_{h+1},\dots,i_\ell))$ or by $T =
(i_1,\dots,i_\ell)$ for the sake of spaces.

From the proof of Theorem~\ref{find} in \cite{kks}, we have
\begin{equation}
\tag{D.1}
\begin{split}
   & i_1 \nsucceq i_2 \nsucceq \cdots \nsucceq i_h,
\\
   & i_{h+1}\nsucceq i_{h+2}\nsucceq \cdots \nsucceq i_\ell.    
\end{split}
\end{equation}

\begin{NB2}
I have changed \verb+(\thesubsection.1)+ to D.1. (D stands for type $D$.)  For other types, use B.1, etc.
\end{NB2}

Next we study the order between $i_h$ and $i_{h+1}$. The following
example shows that $i_h\nsucceq i_{h+1}$ may not be satisfied in general:
Let $n=7$, $\ell =3$, $r=n-\ell-2=2$.
Consider the starting monomial
$M=\ffbox{1}_2\ffbox{2}_{-4}\ffbox{3}_{-6}$. It gives in the crystal $\mathcal{M}_{I_0}(M)$ the monomial $m=\ffbox{3}_2\ffbox{4}_{-4}\ffbox{\overline{4}}_{-6}
=Y_{2,5}^{-1}Y_{3,4}Y_{3,0}^{-1}Y_{4,-1}Y_{3,2}Y_{4,3}^{-1}$.
If we apply $\tilde{f}_3$, we get the monomial
$m'=Y_{3,6}^{-1}Y_{4,5}Y_{3,0}^{-1}Y_{4,-1}Y_{3,2}Y_{4,3}^{-1}$
and this monomial can only be written in the form
$m'=\ffbox{4}_2\ffbox{4}_{-4}\ffbox{\overline{4}}_{-6}$.

Nevertheless $i_h \succeq i_{h+1}$ happens only in very restricted
situations as follows:

\begin{enumerate}
\def\labelenumi{(D.\theenumi)}
\def\theenumi{\arabic{enumi}}
\addtocounter{enumi}{1}
\item  Suppose that $i_{h+1} = k \in \{1,\dots, n-1\}$ and $i_h\succeq
i_{h+1}$. Then $i_h = k'$ is also in $\{ 1,\dots, n-1\}$,
and the successive part $(\overline{k'}, \overline{k'-1},\dots,
\overline{k})$ appears as $(i_{b'},i_{b'+1},\dots,i_{b})$ with
$n - r - k < b - h \le n-k$.

\item Suppose that $i_{h+1} = \overline{k} \in \{\overline{1},\dots,
  \overline{n-1}\}$ and $i_h\succeq i_{h+1}$. Then $i_h =
  \overline{k'}$ is also in $\{ \overline{1},\dots, \overline{n-1}\}$,
  and the successive part $(k',k'+1,\dots,k)$ appears as
  $(i_{a'},i_{a'+1},\dots,i_{a})$ with
  $n-r-k \le h - a < n-k$.

\item  If $i_{h+1} = n$ or $\overline{n}$, then 
$i_h \nsucceq i_{h+1}$.
\end{enumerate}

\begin{NB2}
I have changed inequalities in (D.2,3): $b-h \le n-1-k\Longrightarrow
b-h \le n-k$, $h-a < n-1-k\Longrightarrow h-a < n-k$.
\end{NB2}

Note that (D.1) implies that the successive part in (D.2) occurs in $b'
> h+1$. This together with the second inequality (and $b+k = b'+k'$)
implies $k'\le n-2$. Thus $i_h=n-1$ cannot happen in (D.2). Similarly
$i_{h+1}=\overline{n-1}$, cannot happen in (D.3).

Let us explain the inequalities in (D.2).
(See also $(\mathrm{D}.5\sim 7)$ below.)
First remark that we have $b+k = b'+k'$, so we also have
$n - r - k' < b' - h \le n-k'$. 
We have
\begin{equation*}
\begin{split}
   & \ffbox{i_{h\!+\!1}}_{\ell\!-\!2r\!-\!2h\!-\!1}
   = \ffbox{k}_{\ell\!-\!2r\!-\!2h\!-\!1} 
   = Y^{-1}_{k\!-\!1,\ell\!-\!2r\!-\!2h\!+\!k\!-\!1}
     Y_{k,\ell\!-\!2r\!-\!2h\!+\!k\!-\!2},
\\
   & \ffbox{i_{b}}_{\ell\!-\!2r\!-\!2b+\!1}
   = \ffbox{\overline{k}}_{\ell\!-\!2r\!-\!2b\!+\!1}
   = Y_{k\!-\!1,\ell\!-\!2r\!-\!2b\!+\!2n\!-\!k\!-\!1}
     Y^{-1}_{k,\ell\!-\!2r\!-\!2b\!+\!2n\!-\!k\!}.
\end{split}
\end{equation*}
The second inequality $b-h\le n-k$ and each of the followings are
equivalent:
\begin{equation*}
\begin{split}
  \gr_{k-1}(\ffbox{i_{h\!+\!1}}_{\ell\!-\!2r\!-\!2h\!-\!1})
    & \le \gr_{k-1}(\ffbox{i_{b}}_{\ell\!-\!2r\!-\!2b+\!1}),
\\
  \gr_{k}(\ffbox{i_{h\!+\!1}}_{\ell\!-\!2r\!-\!2h\!-\!1})
    & < \gr_{k}(\ffbox{i_{b}}_{\ell\!-\!2r\!-\!2b+\!1}).
\end{split}
\end{equation*}
%
\begin{NB2}
Be careful for the change.
\end{NB2}
An entry $i_c$ with $c \le h$ and $i_c = k\pm 1$ or $k$ possibly
contributes to $Y_{k,*}^\pm$ or $Y_{k-1,*}^\pm$. It is
\begin{equation*}
    \ffbox{i_c}_{\ell\!-\!2c+\!1}
   = Y_{i_c\!-\!1,\ell\!-\!2c+\!i_c\!+\!1}^{-1}
     Y_{i_c,\ell\!-\!2c+\!i_c}.
\end{equation*}
We find that the first inequality implies
\begin{equation*}
\begin{split}
   \gr_{k-1}(\ffbox{k\!-\!1}_{\ell\!-\!2c+\!1}) & > 
   \gr_{k-1}(\ffbox{i_{b}}_{\ell\!-\!2r\!-\!2b+\!1}),
\\
   \gr_{k}(\ffbox{k}_{\ell\!-\!2c+\!1}) & > 
   \gr_{k}(\ffbox{i_{b}}_{\ell\!-\!2r\!-\!2b+\!1}).
\end{split}
\end{equation*}
In fact, the first inequality is equivalent to each of these
inequalities with $c = h$.
The only remaining possibility of a contribution is the case
$i_{h+2} = k+1$ or $i_{b+1} = \overline{k-1}$. In the former (resp.\
latter) case its contribution is canceled with that of $i_{h+1}$
(resp.\ $i_b$).
Also
\begin{equation*}
\begin{split}
   & \ffbox{i_h}_{\ell\!-\!2h\!+\!1}
   = \ffbox{k'}_{\ell\!-\!2h\!+\!1} 
   = Y^{-1}_{k'\!-\!1,\ell\!-\!2h\!+\!k'\!+\!1}
     Y_{k',\ell\!-\!2h\!+\!k'},
\\
   & \ffbox{i_{b'}}_{\ell\!-\!2r\!-\!2b'+\!1}
   = \ffbox{\overline{k'}}_{\ell\!-\!2r\!-\!2b'+\!1}
   = Y_{k'\!-\!1,\ell\!-\!2r\!-\!2b'\!+\!2n\!-\!k'\!-\!1}
     Y^{-1}_{k',\ell\!-\!2r\!-\!2b'\!+\!2n\!-\!k'\!},
\\
   & \ffbox{i_{c}}_{\ell\!-\!2r\!-\!2c\!+\!1}
   = Y^{-1}_{i_c\!-\!1,\ell\!-\!2r\!-\!2c\!+\!i_c\!+\!1}
     Y_{i_c,\ell\!-\!2r\!-\!2c\!+\!i_c} \qquad
     (c\ge h+1, i_c = k', k'\pm 1).
\end{split}
\end{equation*}
We have
\begin{gather*}
  \gr_{k'}(\ffbox{k\!+\!1}_{\ell\!-\!2r\!-\!2c\!+\!1}) \le
  \gr_{k'}(\ffbox{i_{b'}}_{\ell\!-\!2r\!-\!2b'+\!1}) <
  \gr_{k'}(\ffbox{i_h}_{\ell\!-\!2h\!+\!1}),
\\
  \gr_{k'-1}(\ffbox{k}_{\ell\!-\!2r\!-\!2c\!+\!1}) \le
  \gr_{k'-1}(\ffbox{i_{b'}}_{\ell\!-\!2r\!-\!2b'+\!1}) <
  \gr_{k'-1}(\ffbox{i_h}_{\ell\!-\!2h\!+\!1}).
\end{gather*}


\begin{prop}\label{enumcond}
The set $\mathcal M_{I_0}(M_{\ell,h,r})$ consists of monomials 
$m_T$ associated with tableaux $T$ satisfying 
$(\mathrm{D}.1\sim 4)$.
\end{prop}

\begin{proof} 
For a monomial $m=\prod_{i\in I,l\in\ZZ} Y_{i,l}^{u_{i,l}(m)}$ and
$J\subset I$ we set $m^{(J)}=\prod_{i\in
J,l\in\ZZ}Y_{i,l}^{u_{i,l}(m)}$.

Let us take a monomial $m_T$ with $T$ satisfying $(\mathrm{D}.1
\sim 4)$. We assume $\tf_i m_T \neq 0$ and want to show that $\tf_i
m_T = m_{T'}$ for some $T'$ which also satisfies $(\mathrm{D}.1
\sim 4)$. To get such a $T'$, we replace an entry $i_a$ by a new one
according to the rule described in Figure~\ref{fig:Dn}. The condition
(D.1) is preserved by the proof of Theorem~\ref{find} in \cite{kks}.

Case (a): $i_{h+1} \in \{1,\dots, n-1\}$ and $i_h\succeq i_{h+1}$,
i.e., the situation described in (D.2).

We cannot replace $i_{b'}$, $i_{b'+1}$, \dots, $i_{b-1}$ as they are
successive.
If we do not replace $i_h$, $i_{h+1}$ or $i_{b}$,
the resulting monomial $\tf_i m_T$ clearly satisfies
$(\mathrm{D}.2)$. 

Subcase (a.1): we replace $i_{h+1} = k$.

As we have remarked just after (D.4), we cannot have $k
= n-1$. Thus $k \neq n-1$ and $i_{h+1}$ is replaced by $k+1$.
%
%
If the condition $i_h\nsucceq i_{h+1}$ is satisfied in the new
tableau, we are done. Otherwise we still have (D.2).

Subcase (a.2): we replace $i_h = k'$.

As we have remarked just after (D.4), we cannot have $k' =
n-1$. Thus $k' \neq n-1$ and $i_{h}$ is replaced by $k'+1$.
If $i_{b'-1} = \overline{k'+1}$, (D.2) is preserved. So we
may assume $i_{b'-1}\neq\overline{k'+1}$. Let us consider
$(m_{T})^{(k')}$. We have contributions from
$\ffbox{i_h}_{\ell\!-\!2h\!+\!1}$, 
$\ffbox{i_{b'}}_{\ell\!-\!2r\!-\!2b'+\!1}$ and possibly
from $\ffbox{i_c}_{\ell\!-\!2r\!-\!2c\!+\!1}$ with $c\ge h+1$
when $i_c = k'$ or $k'+1$. As we explained above, we have
\(
   \gr_{k'}(\ffbox{i_c}) \le \gr_{k'}(\ffbox{i_{b'}}) <
   \gr_{k'}(\ffbox{i_h}), 
\)
and the first inequality is strict if $i_c = k'$.
\begin{NB2}
 \begin{equation*}
   \begin{split}
     (m_{T})^{(k')} &= Y_{k',\ell-2h+k'}Y_{k',\ell-2r-2b'+2n-k'}^{-1}
 \\
    & \qquad\qquad
    \times \text{possibly }Y_{k',\ell-2r-2c+k'}
    \text{ or } Y_{k',\ell-2r-2c+k'+2}^{-1}. 
   \end{split}
 \end{equation*}
\end{NB2}
This implies that either $\tf_k m_T = 0$ or $i_c$ is replaced by
$\tf_k$ instead of $i_h$, and we have a contradiction.

Subcase (a.3): we replace $i_{b} = \overline{k}$ by
$\overline{k-1}$. ($\tf_i = \tf_{k-1}$)

We consider the contribution to $Y_{k-1,*}$ and use
(D.2) as in subcase~(a.2), we conclude that (D.2) is preserved.

\begin{NB2}
Detail:
\begin{equation*}
   \gr_{k-1}(\ffbox{i_{h\!+\!1}}) \le
   \gr_{k-1}(\ffbox{i_b}) < \gr_{k-1}(\ffbox{i_c}).
\end{equation*}
and
\begin{equation*}
   (m_{T})^{(k-1)}\in
   Y_{k-1,\ell-2h-2r+k-1}^{-1} 
   Y_{i,\ell-2r-2b+2n-k-1}
   \times\ZZ[Y_{i,m}^{\pm}]_{m\geq \ell-2h+k-1}.
\end{equation*}
This implies that either $\tf_k m_T = 0$ or $i_{c}$ is replaced by
$\tf_{k-1}$ instead of $i_b$, and we have a contradiction.
\end{NB2}

Similarly we have the assertion when we are in the situation described
in (D.3).

Case (b): $i_{h+1} \in \{1,\dots, n-1\}$ and $i_h\nsucceq i_{h+1}$.

\begin{NB2}
There was a gap in this case. Check carefully.
\end{NB2}

It suffices to study the case when we replace $i_h=k'$ by $k'+1$ and
$i_{h+1}$ is equal to $k'+1$. Otherwise the condition $i_h\nsucceq
i_{h+1}$ is preserved.
If there is an entry $i_b = \overline{k'+1}$ with
\(
   n-r-(k'+1) < b - h \le n -(k'+1),
\)
then (D.2) is satisfied by a new tableau $\tf_k m_T$ with
$b = b'$ and $i_b = \overline{k'+1}$.
We assume that there is no such entry. Let us consider
$(m_T)^{(k')}$. We have contributions from
\(
   \ffbox{i_h}_{\ell\!-\!2h\!+\!1},   
\)
\(
   \ffbox{i_{h\!+\!1}}_{\ell\!-\!2h\!-\!2r\!-\!1},   
\)
and possibly from
\( 
   \ffbox{i_c}_{\ell\!-\!2r\!-\!2c\!+\!1}
\)
with $c\ge h+1$
when $i_c = \overline{k'}$, $\overline{k'+1}$.
We have
\begin{equation*}
\begin{split}
   & (m_{T})^{(k')} = 
   Y_{k',\ell-2h+k'}Y_{k',\ell-2h-2r+k'}^{-1}
\\
   &\qquad\qquad\qquad \times
   \text{possibly }
   Y_{k',\ell-2r-2c+2n-k'}^{-1} \text{ or }
   Y_{k',\ell-2r-2c+2n-k'-2}.
\end{split}
\end{equation*}
\begin{NB2}
Be careful for the place of $\&$ in \verb+split+ environment. A misuse
cause a wrong spacing.
\end{NB2}
The last factor $Y_{k',\ell-2r-2c+2n-k'-2}$ appears when
$i_c=\overline{k'+1}$, so we have
\begin{equation*}
  \begin{split}
    & \gr_{k'}(\ffbox{i_c}_{\ell\!-\!2r\!-\!2c\!+\!1})
    \ge 
    \gr_{k'}(\ffbox{i_h}_{\ell-2h+1})
\\
  \text{or } &
  \gr_{k'}(\ffbox{i_c}_{\ell\!-\!2r\!-\!2c\!+\!1})
  < \gr_{k'}(\ffbox{i_{h\!+\!1}}_{\ell\!-\!2h\!-\!2r\!-\!1})
  \end{split}
\end{equation*}
by the assumption.
When we have the equality in the first inequality, $Y_{k',\ell-2h+k'}$
appears in two places, we may have two choices to replace, either
$i_h$ or $i_c$ with $c=n-r-k'+h-1$. We choose the second one, then
$i_h\nsucceq i_{h+1}$ holds. In the other cases, either $\tf_{k'} m_T
= 0$ or $i_c$ (instead of $i_h$) is replaced, so we get a
contradiction.

Similarly we have the assertion when 
$i_{h+1}\in \{\overline{1},\dots, \overline{n-1}\}$
and $i_h\nsucceq i_{h+1}$.

Case (c): $i_{h+1} = n$.

We have $i_h\nsucceq i_{h+1}$ by (D.4). We assume that
$i_h$ is replaced and $i_h\nsucceq i_{h+1}$ is violated, and lead to a
contradiction. Thus we may assume either $i_h = n-1$ is replaced
by $n$, or $i_h = \overline{n}$ is replaced by $\overline{n-1}$.

Subcase (c.1): $i_h = n-1$ is replaced by $n$.

Starting from $i_{h+1} = n$, several $\overline{n}$, $n$ may follow in
turn. We have two possibilities whether it ends with $n$ or $\overline{n}$:
\begin{equation*}
  \begin{split}
    & \ffbox{n}_{\ell\!-\!2h\!-\!2r\!-\!1}
    \ffbox{\overline{n}}_{\ell\!-\!2h\!-\!2r\!-\!3}
    \cdots
    \ffbox{{n}}_{\ell\!-\!2r\!-\!2c\!+\!1}
    = Y_{n\!-\!1,\ell\!-\!2h\!-\!2r\!+\!n\!-\!1}^{-1}
    Y_{n,\ell\!-\!2r\!-\!2c\!+\!n\!-\!1},
\\
    & \ffbox{n}_{\ell\!-\!2h\!-\!2r\!-\!1}
    \ffbox{\overline{n}}_{\ell\!-\!2h\!-\!2r\!-\!3}
    \cdots
    \ffbox{{n}}_{\ell\!-\!2r\!-\!2c\!+\!1}
    \ffbox{\overline{n}}_{\ell\!-\!2r\!-\!2c\!-\!1}
    = Y_{n\!-\!1,\ell\!-\!2h\!-\!2r\!+\!n\!-\!1}^{-1}
    Y_{n\!-\!1,\ell\!-\!2r\!-\!2c\!+\!n\!-\!3}.
  \end{split}
\end{equation*}
The grades of $Y_{n\!-\!1,*}$ are smaller than $\ell-2h+n-1$ which
is the grade of $(\ffbox{i_h}_{\ell\!-\!2h\!+\!1})^{(n-1)}$. So
either $\tf_{n-1} m_T = 0$ or $i_h$ is not replaced. We have a
contradiction, so this cannot happen.

Subcase (c.2): $i_h = \overline{n}$ is replaced by $\overline{n-1}$.

The grade of 
\(
  (\ffbox{i_h}_{\ell\!-\!2h\!+\!1})^{(n-1)}
\)
is again $\ell\!-\!2h\!+\!n\!-1$. We may have an entry $i_c = n-1$ or
$n$ with $c < h$, but the grade is always larger than
$\ell\!-\!2h\!+\!n\!-\!1$. As in case (c1), either $\tf_{n-1} m_T = 0$
or $i_h$ is not replaced. We have a contradiction, so this cannot
happen.

Now we consider a monomial $m_T$ such that $T$ satisfies
$(\mathrm{D}.1\sim 4)$ and $\te_i m_T = 0$ for all $i\in I_0$. We
want to prove $T = ((1,\dots,h),(h+1,\dots,\ell))$.

First we consider the case $i_h\nsucceq i_{h+1}$. Following the proof
of \cite[Prop.~2.1]{kks}, we suppose that there is an entry $i_c =
\overline{k}$ with $k\in \{1,\dots,n\}$ and lead to a
contradiction. We take the smallest $c$ among such entries.  We first
assume that $c\ge h+1$ and $k\neq n$. Then
\(
   (\ffbox{\overline{k}}_{\ell\!-\!2r\!-\!2c\!+\!1})^{(k)} =
   Y_{k,\ell\!-\!2r\!-\!2c\!+\!2n\!-\!k}^{-1}.
\)
As we have $\te_k m_T = 0$, we must have a contribution
$Y_{k,m}$ with $m \ge \ell\!-\!2r\!-\!2c\!+\!2n\!-\!k$. This is
possible only if we have an entry $i_{c'}$ with
$i_{c'} = k$ with
\(
   c' \le r+c-n+k.
\)
Combining this with the assumption $r\le n-\ell-1$ and $c\le \ell$, we get
\(
   c' \le k - 1.
\)
This inequality holds even if $k=n$, as $\ell\le n-2$. We also have
this even if $c'\le h$, as $c'\le c-n+k\le\ell-n+k\le k-1$.

Now we study $(m_T)^{(k-1)}$. The contribution of $i_{c'}$ cannot be
canceled by $\overline{k}$ by the comparison of grades, so we must
have $i_{c'-1} = k-1$. We continue this to get $i_1 = k - c' +
1$. Then the above inequality implies $i_1 \ge 2$, so we have
$\te_{i_1-1} m_T \neq 0$, which is a contradiction. Thus all entries are in
$\{ 1,\dots,n\}$. In this case it is easy to show that $m_T$ must be equal
to the highest weight monomial $M_{\ell,h,r}$.

In the situation in (D.3), we cannot have $\te_i m_T = 0$
for all $i$ by the same reason as above.

Next consider the situation in (D.2). We have
\(
   (\ffbox{i_h}_{\ell\!-\!2h\!+\!1})^{(k'-1)} = 
   Y^{-1}_{k'\!-\!1,\ell\!-\!2h\!+\!k'\!+\!1}.
\)
As $\te_{k'-1} m_T = 0$, we must have a contribution
$Y_{k'-1,m}$ with 
$m \ge \gr_{k'-1}(\ffbox{i_h})$.
But by the remark preceding this proposition, this is possible only when
$i_{h-1} = k'-1$. We continue this argument to get
$i_1 = k'-h + 1$. If $i_1 \neq 1$, we have $\te_{i_1-1} m_T \neq 0$. So
we must have $k' = h$. But this together with $r\le n-\ell-1$
violates (D.2) as
$\ell - h \ge b' - h > n - r - k' = n - r - h
 \ge \ell+1 - h$. This case cannot happen.
\end{proof}

\begin{NB2}
You cannot use `$\Rightarrow$' except when you use it as a logical symbol.
\end{NB2}

We now modify the procedure (c) in the proof of Theorem~\ref{find}
when the suffix may jump.
By the procedure explained before we attain
\begin{enumerate}
\def\labelenumi{(D.\theenumi)}
\def\theenumi{\arabic{enumi}}
\addtocounter{enumi}{4}
\item there is no pair $a$, $b$ such that $1\leq a < b \leq h$ and
  $i_a=k$, $i_b=\overline{k}$ and $b-a = n -1 - k$,
\item there is no pair $a$, $b$ such that $h+1\leq a < b \leq \ell$
  and $i_a=k$, $i_b=\overline{k}$ and $b-a = n -1 - k$.
\end{enumerate}

\begin{lem}
For a tableau $T$ satisfying $(\mathrm{D}.1\sim 4)$ there exists a
tableau $T'$ satisfying $(\mathrm{D}.1\sim 6)$ and $m_T = m_{T'}$.
\end{lem}

\begin{proof}
As is mentioned after Theorem~\ref{find}, we replace
the pair
\(
  \ffbox{k}_{p}\ffbox{\overline{k}}_{p-2(n-1-k)}
\)
by
\(
  \ffbox{k\!+\!1}_{p}\ffbox{\overline{k\!+\!1}}_{p-2(n-1-k)}
\)
repeatedly from $k=1$ to $n-2$.
We prove that the procedure does not violate $(\mathrm{D}.1\sim
4)$. For example, let us consider the procedure for (D.5). It is
enough to consider a possibility of violating (D.3) when $i_h$
is replaced, say from $\overline{k'}$ to $\overline{m}$ with $m >
k'$. Then $i_{a'} = k'$ with $h-{a'}=n-1-k'$ is replaced by $k'+1$,
$i_{a'+1} = k'+1$ by $k'+2$, and so on.

Suppose that $i_h \nsucceq i_{h+1}$. As $\overline{k'} \succeq
\overline{m}$, we still have $i_h \nsucceq i_{h+1}$ after the
replacing procedure. Thus $(\mathrm{D}.2\sim 4)$ hold in this case.

Next suppose $i_h \succeq i_{h+1}$. By (D.3) there exists a successive
part $(k',k'+1,\dots,k) = (i_{a'},i_{a'+1},\dots,i_a)$ with $i_{h+1} =
\overline{k}$. Then the above procedure continues at least until $i_h$
is replaced by $\overline{k+1}$, i.e., $\overline{m} \preceq
\overline{k+1}$. Therefore $i_h\nsucceq i_{h+1}$ holds after the
procedure was finished.
\end{proof}

\begin{rem}
Under (D.5,6), `$n-k$' appearing in inequalities in (D.2,3) can be
replaced by `$n-1-k$'.
\end{rem}

Next we consider the corresponding condition for the pair $a$, $b$
with $a\leq h$, $h+1\leq b$. In this case, we need to change the rule
as indicated by the following example:
Consider the case $n=6$, $\ell=4$, $r=1$ and the starting
monomial $m=\ffbox{1}_3\ffbox{2}_{-1}\ffbox{3}_{-3}\ffbox{4}_{-5}$.
Then consider the monomial
$m'=\ffbox{1}_3\ffbox{2}_{-1}\ffbox{3}_{-3}\ffbox{\overline{1}}_{-5}=Y_{1,3}Y_{1,5}^{-1}Y_{1,1}^{-1}Y_{3,-1}$.
For $b=4$ and $a=1$, we have $b-a=n-r-k-1=3$.
Thus this monomial violates the condition (2) of Theorem~\ref{find}.
But if we replace the pair $i_1 = 1$, $i_4 = \overline{1}$ to
$2$, $\overline{2}$ as before, we get
$\ffbox{2}_3\ffbox{2}_{-1}\ffbox{3}_{-3}\ffbox{\overline{2}}_{-5}$,
which does not satisfy the condition~(1) of Theorem~\ref{find}.
In the original situation we can further replace the pair
$i_2 = 2$, $i_4 = \overline{2}$ to $3$, $\overline{3}$,
and then further $(i_3,i_4) = (3,\overline{3})$ to $(4,\overline{4})$
to achieve the condition (1). But we cannot make this replacement
as $\ffbox{2}_{-1}\ffbox{\overline{2}}_{-5}\neq
\ffbox{3}_{-1}\ffbox{\overline{3}}_{-5}$.

Thus we need to modify either of the conditions (1), (2).
We keep the condition (1) and the change (2) to
\begin{enumerate}
\def\labelenumi{(D.\theenumi)}
\def\theenumi{\arabic{enumi}}
\addtocounter{enumi}{6}
\item there is no pair $a$, $b$ such that $a \leq h$ , $h+1 \leq b$ , $i_a=k$, $i_b=\overline{k}$ and $b-a = n - r - k$.
\end{enumerate}

Note that this condition coincides with one in \cite{kks} when
we write $T$ as a tableau of shape $(h,\ell-h)$ where the second
column is shifted below by $h+r$ boxes.

\begin{lem}\label{subcrys}
For a tableau $T$ satisfying $(\mathrm{D}.1\sim 6)$ there exists a
tableau $T'$ satisfying $(\mathrm{D}.1\sim 7)$ and $m_T = m_{T'}$.
\end{lem}

\begin{defi}
Let $D_{\ell,h,r}$ be the set of tableaux $T$ satisfying
$(\mathrm{D}.1\sim 7)$.
\end{defi}

\begin{proof}
Let $T$ be a tableau satisfying $(\mathrm{D}.1\sim 6)$.
If we have a pair $a$, $b$ violating (D.7),
we replace the tableau by using the relation
\(
   \ffbox{k}_{\ell-2a+1} \ffbox{\overline{k}}_{\ell-2b+1-2r}
   = \ffbox{k\!-\!1}_{\ell-2a+1} \ffbox{\overline{k\!-\!1}}_{\ell-2b+1-2r}.
\)
If there are several such pairs or this procedure yields a new such
pair, we replace them repeatedly starting from $k=n-1$, then
$k=n-2$,..., and finally to $k=2$. (Note that this is converse to the
order of the procedure among entries in the same column.) As $r\le
n-\ell-1$, the condition (D.7) always holds for $k=1$.

It is clear that (D.1) is preserved when we finish the
above procedure. Let us check that $(\mathrm{D}.2\sim 6)$ are
preserved.
  
Suppose that we apply the above procedure to a tableau $T =
((i_1,\cdots,i_h),(i_{h+1},\cdots,i_\ell))$ starting from changing a
pair $a\leq h$, $h+1 \leq b$ such that $i_a=k$, $i_b=\overline{k}$ and
$b-a = n - r - k$, and finally get a new tableau
$T'=((j_1,\cdots,j_h),(j_{h+1},\cdots,j_\ell))$.
Suppose that (D.5) is violated, i.e., we have a pair $A < B\leq h$
such that $j_A=K$ and $j_B=\overline{K}$ and $B-A=n-1-K$. As $i_m$ for
$a+1\leq m\leq b-1$ is unchanged by the above procedure, we have $A\le
a$. We have $n-1-K=B-A=(a-A)+(B-b)+n-r-k$, so
$a-A=b-B+r+k-K-1>k-K$. This inequality contradicts with (D.1) as $k =
j_a \ge j_A + (a-A) > k$. So (D.5) is satisfied by $T'$. In the same
way (D.6) is satisfied by $T'$.

We treat the cases separately according whether $i_h$, $i_{h+1}$ are
changed or not.

Case (1): Both $i_h$ and $i_{h+1}$ are changed during the procedure.

In this case we have $i_h\in \{2,\dots,n-1\}$, $i_{h+1}\in
\{\overline{2},\dots,\overline{n-1}\}$. Then $(\mathrm{D}.2\sim 4)$
imply $i_h\nsucceq i_{h+1}$. This inequality is preserved during the
procedure, so $(\mathrm{D}.2\sim 4)$ are preserved.

Case (2): Both $i_h$ and $i_{h+1}$ are unchanged during the procedure.

If $i_h\nsucceq i_{h+1}$, we clearly have $(\mathrm{D}.2\sim 4)$. Thus
we assume $i_h\succeq i_{h+1}$. Suppose $i_{h+1} = k\in
\{1,\dots,n-1\}$ and take the successive part $(\overline{k'},
\overline{k'-1},\dots, \overline{k}) = (i_{b'},i_{b'+1},\dots,i_{b})$
with $i_h = k'$ as in (D.2).
If (D.2) is violated during the procedure, we replace a pair
$(i_A,i_B) = (K,\overline{K})$ with $A < h$, $b' < B < b$ with
$B - A = n-r-K$. 
But this contradicts with the inequality in (D.2) as 
\(
  n-r-K < b- h + k - K = B-h < B - A.
\)
So this cannot happen. The case $i_{h+1}\in
\{\overline{1},\dots,\overline{n-1}\}$ can be proved in the same way.

Case (3): $i_h$ is unchanged, $i_{h+1}$ is changed.

Suppose that $i_{h+1}$ is changed, say from $\overline{k}$ to $j_{h+1}
= \overline{m}$ with $m < k$. Then $i_a = k$ with $h-a=n-r-k$ is
replaced by $k-1$, $i_{a-1}$ is replaced by $k-2$, and so on.
If $i_h \nsucceq i_{h+1}$, we get $j_h\nsucceq j_{h+1}$ as 
the procedure preserves this inequality. Thus we have
$(\mathrm{D}.2\sim 4)$.
If $i_h \succeq i_{h+1}$, we have a successive part $(k',k'+1,\dots,k)
= (i_{a'},i_{a'+1},\dots,i_{a})$ with $i_{h+1} =
\overline{k'}$. Therefore the procedure continues at least until
$i_{h}$ is replaced by $k'-1$, i.e., $m < k'$. Therefore $j_h \nsucceq
j_{h+1}$.

Case (4): $i_h$ is changed, $i_{h+1}$ is unchanged.

If $i_h = k'$ is replaced, there exists $i_{b'} = \overline{k'}$ with
$b'-h = n-r-k'$. But it contradicts with an inequality in (D.2), we
have $i_h\nsucceq i_{h+1}$. Then we have $(\mathrm{D}.2\sim 4)$ as in (3).
\end{proof}

\begin{NB2}
The following paragraph explains the motivation of the definition of
$\tau_{\ell,h,r}$. But I do not quite understand. 
\end{NB2}

+++++++++++++++++++++++++++++++++++++++++++++++++++++++++++++++++++++

Now we describe explicitly the crystal isomorphism between $\mathcal{M}_{I_0}(M_{\ell,h,r})$ and $\mathcal{M}_{I_0}(M_{\ell,h+1,r})$ for $0\leq h\leq \ell -1$. Let us first explain the idea of the construction. When $h$ vary without changing the $i_a$, one property (\thesubsection.1,2,3,4,5,6,7) may not be satisfied. For example if $i_{H+1}\preceq i_H$ at $h=H$ we have to move $i_{H+1}$ so that $i_H\nsucceq i_{H+1}$ at $h=H+1$. To do it, the idea is to replace a product of boxes $\ffbox{i}_{L_1}\ffbox{\overline{i}}_{L_2}$ by a product $\ffbox{i+1}_{L_1}\ffbox{\overline{i+1}}_{L_2'}$ in certain situations. However for $a\leq b$, $d=b-a$, $L_1\in\ZZ$, $L_2=L_1-2d$ and $L_2'=L_2-2r$, the map $\{\ffbox{i}_{L_1}\ffbox{\overline{i}}_{L_2}/1\leq i\leq n-1\}\rightarrow \{\ffbox{i}_{L_1}\ffbox{\overline{i}}_{L_2'}/2\leq i\leq n\}$ such that $\ffbox{i}_{L_1}\ffbox{\overline{i}}_{L_2}\mapsto\ffbox{i+1}_{L_1}\ffbox{\overline{i+1}}_{L_2'}$ is not well-defined because of the identities involving products of boxes. We can define a map in the following way. Let $j_1=n-1-d-r$ and $j_2=n-1-d$. We have the relations
\begin{equation*}\ffbox{j_1}_{L_1}\ffbox{\overline{j_1}}_{L_2'}=\ffbox{j_1+1}_{L_1}\ffbox{\overline{j_1+1}}_{L_2'}\text{ , }\ffbox{j_2}_{L_1}\ffbox{\overline{j_2}}_{L_2}=\ffbox{j_2+1}_{L_1}\ffbox{\overline{j_2+1}}_{L_2}.
\end{equation*}
So there is a bijection $\{\ffbox{i}_{L_1}\ffbox{\overline{i}}_{L_2}/1\leq i\leq n\}\rightarrow \{\ffbox{i}_{L_1}\ffbox{\overline{i}}_{L_2'}/1\leq i\leq n\}$ defined by
\begin{enumerate}
    \item if $i\leq j_1$, $\ffbox{i}_{L_1}\ffbox{\overline{i}}_{L_2}\mapsto \ffbox{i}_{L_1}\ffbox{\overline{i}}_{L_2'}$,
    
    \item if $j_1+1\leq i\leq j_2-1$, $\ffbox{i}_{L_1}\ffbox{\overline{i}}_{L_2}\mapsto \ffbox{i+1}_{L_1}\ffbox{\overline{i+1}}_{L_2'}$,

    \item if $j_2+1\leq i$, $\ffbox{i}_{L_1}\ffbox{\overline{i}}_{L_2}\mapsto \ffbox{i}_{L_1}\ffbox{\overline{i}}_{L_2'}$.
\end{enumerate}
The idea of the following definition is to use this transformation.













\begin{defi}\label{defitau} Let $0\leq h\leq \ell$ and $T=((i_1,\cdots,i_h),(i_{h+1},\cdots,i_\ell))\in D_{\ell,h,r}$. We define $\tau_{\ell,h,r}(T)$ in the following way.

\textup{(1)} If $i_{h+1}=i\in\{1,\cdots,n\}$ and there is $h'\geq h+2$ such that $i_{h'}=\overline{i}$ and $n-h'+h+1-r\leq i\leq n-1-h'+h$. Then $h'$ is unique and there is a unique minimal $h''\in\{h+2,\cdots ,h'\}$ such that $(i_{h''},i_{h''+1},\cdots,i_{h'})=(\overline{i+(h'-h'')},\overline{i+(h'-h'')-1},\cdots,\overline{i})$. We set 
\begin{equation*}
\tau_{\ell,h,r}(T)=((i_1,\cdots, i_h,i+(h'-h'')+1),
\end{equation*}
\begin{equation*}
(i_{h+2},\cdots,i_{h''-1},
\overline{i+(h'-h'')+1},\overline{i+(h'-h'')},\cdots,\overline{i+1}, i_{h'+1},\cdots, i_\ell)).
\end{equation*}
\begin{NB2}
  You should not separate a single equation into two
  \verb+\begin{equation*}..\end{equation*}+. Read the manual of TeX.
  You must use split or align environment.
\end{NB2}

\textup{(2)} If $i_{h+1}=\overline{i}\in\{\overline{1},\cdots,\overline{n}\}$ and there is $h' < h$ such that $i_{h'}=i$ and $n+h'-h-r\leq i\leq n-2-h+h'$. Then $h'$ is unique and there is a unique minimal $h''\leq h'$ such that $(i_{h''},i_{h''+1},\cdots,i_{h'})=(i-(h'-h''),1+i-(h'-h''),\cdots,i)$. We set 
\begin{equation*}
\tau_{\ell,h,r}(T)=((i_1,\cdots, i_{h''-1},i-(h'-h'')-1,\cdots,i-1,i_{h'+1},\cdots,i_h),
\end{equation*}
\begin{equation*}
(\overline{i-(h'-h'')-1},i_{h+2},\cdots, i_\ell)).
\end{equation*}

\textup{(3)} If the conditions (1) or (2) are not satisfied, we set 
\begin{equation*}
\tau_{\ell,h,r}(T)=((i_1,\cdots, i_{h+1}),(i_{h+2},\cdots, i_\ell)).
\end{equation*}
\end{defi}
The map $\tau_{\ell,h,r}$ is well-defined : for example for $T$ satisfying the condition (1) of the definition \ref{defitau}, we have $i+(h'-h'')+1\preceq n-1$ and $\overline{n-2}\preceq i_{h''}$ because $i+h'-(h+2)+1\leq n-2$. It is also clear that if $T\in D_{\ell,h,r}$, then $\tau_{\ell,h,r}(T)\in D_{\ell,h+1,r}$.



+++++++++++++++++++++++++++++++++++++++++++++++++++++++++++++++++++++

Let $T = ((i_1,\dots,i_h), (i_{h+1},\dots,i_\ell))\in
D_{\ell,h,r}$. We divide the situation into three cases:
\begin{enumerate}
\def\labelenumi{(D.\theenumi)}
\def\theenumi{\alph{enumi}}
\item $i_{h+1} = k \in \{1,\dots,n-1\}$ and there is an entry $i_b =
\overline{k}$ with $n-r-k< b - h \le n-1-k$.

\item $i_{h+1} = \overline{k} \in \{\overline{1},\dots,\overline{n-1}
\}$ and there is an entry $i_a = k$ with $n-r-k \le h - a < n - 1 - k$.

\item Neither (D.a) nor (D.b) is not satisfied. 
\end{enumerate}

In the case (D.c) we simply define
\begin{equation*}
\tau_{\ell,h,r}(T)=((i_1,\cdots, i_{h+1}),(i_{h+2},\cdots, i_\ell)).
\end{equation*}
Note that we have $i_h\nsucceq i_{h+1}$ in (D.c). So
$\tau_{\ell,h,r}(T)$ also satisfies (D.1). As we have $i_{h+1} \nsucceq
i_{h+2}$, $(\mathrm{D}.2\sim 4)$ are satisfied. 
The conditions (D.6) clearly holds. Let us check the condition (D.5).
The only possibility is $(i_a, i_{h+1}) = (k,\overline{k})$ with $a\le
h$. But $h+1-a = n-1-k$ cannot happen since we are not in case (D.b).
Similarly (D.7) holds as we are not in case (D.a).

Next suppose we are in the case (D.a). As was explained in the
paragraph just after (D.4), the inequalities imply $b > h+1$ and
$k\neq n-1$. Starting from $i_b$, we go back $i_{b-1}$, $i_{b-2}$,
...  while entries are successive. Let $i_{b''}$ be the ending entry,
so $(i_{b''},i_{b''+1},\dots,i_b)$ are successive as
$(\overline{k''},\overline{k''+1},\dots,\overline{k})$ and
$i_{b''-1}\neq \overline{k''-1}$. Also by the same reasoning as above,
we have $k''\le n-2$.
We define $\tau_{\ell,h,r}(T)$ by
\begin{equation*}
  \begin{split}
    \tau_{\ell,h,r}(T)=((i_1,\cdots&,i_h,k''\!+\!1),
\\
  & (i_{h+2},\cdots,i_{b''-1}, \overline{k''\!+\!1},\overline{k''},
  \cdots,\overline{k\!+\!1}, i_{b+1},\cdots, i_\ell)). 
  \end{split}
\end{equation*}
Let us check that the new tableau satisfies $(\mathrm{D}.1\sim 7)$.
When $i_h\succeq i_{h+1}$, i.e., in the situation described in (D.2), we 
have $b'' \ge b'$ and $k'' \ge k' = i_h$. If $i_h \nsucceq i_{h+1}$,
we have $k'' \ge k = i_{h+1} \ge i_h$. Thus (D.1) is satisfied by the
new tableau. 
Note that we have $k = i_{h+1} \nsucceq i_{h+2}$ by (D.1). Therefore
even if we have $k''+1 \succeq i_{h+2}$, we find a successive part as
in (D.2).
It is clear that (D.5) is satisfied.
The property (D.6) could be violated by a pair $j_A=K$,
$j_B=\overline{K}$ with $h+2\leq A < B\leq b$. But we have
$K=k+1+b-B$, so $B-A+K=b+1+k-A\leq n+h-A\leq n-2$. Thus
$B-A = n-K-1$ cannot happen.
The property (D.7) could be violated by $j_A=K$, $j_B=\overline{K}$ with
$A\leq h+1$ and $h+2\leq B\leq b$. Since we have $K=k+1+b-B$,
$B-A+K+r=b+1+k-A+r >  n+h-A+1\geq n$, and so $B-A = n-r-K$ cannot happen.

Similarly in the case (D.b), we take $i_{a''}$ so that
$(i_{a''},i_{a''+1},\dots, i_a) = (k'',k''+1,\dots,k)$ and $i_{a''-1}\neq
k''-1$. We have $k\le n-3$. We then define
\begin{equation*}
  \begin{split}
    \tau_{\ell,h,r}(T)=((i_1,\cdots, i_{a''-1},k''\!-\!1,\cdots,k\!-\!1,
    i_{a+1},&\cdots,i_{h}, \\
& \overline{k''\!-\!1}), (i_{h+2},\cdots, i_\ell)).
  \end{split}
\end{equation*}

\begin{lem}\label{lem:tau}
Let $0\le h\le \ell-1$, $T_1\in D_{\ell,h,r}$ and $j\in I_0$.
Assume that $\tf_j m_{T_1} \neq 0$ and it is equal to $m_{T_2}$
for some $T_2\in D_{\ell,h,r}$. Then
$\tf_j m_{\tau_{\ell,h,r}(T_1)} = m_{\tau_{\ell,h,r}(T_2)}$.
\end{lem}

\begin{proof}
As we have
\(
  \tau_{\ell,h,r}((1,\cdots,h),(h+1,\cdots,\ell))
   =((1,\cdots,h+1),(h+2,\cdots,\ell)),
\)
it is enough to show that $\tf_j$ commutes with $\tau_{\ell,h,r}$.
Let us denote $T_1$, $T_2$, $T_3 = \tau_{\ell,h,r}(T_1)$ 
by 
\(
   ((i_1,\cdots,i_h),(i_{h+1},\cdots,i_\ell)),
\)
\(
   ((i_1',\cdots,i_h'),(i_{h+1}',\cdots,i_\ell')),
\)
\(
   ((i_1'',\cdots,i_{h+1}''),(i_{h+2}'',\cdots,i_\ell''))
\)
respectively.
We suppose that $T_2$ is obtained from $T_1$ by replacing
an entry $i_c$ by $i_c'$ by the rule described in Figure~\ref{fig:Dn}.
We have either $i_c = j$ or $\overline{j+1}$.

First we study the case (D.a). We take $i_{b''}$ and $k''$ as before.
We separate cases according to the value of $j$.

Case (a.1): $j\neq k-1,k,k'',k''+1$.

We clearly have the assertion.

Case (a.2): $j=k-1$.

Let us first consider $\tf_{k-1} m_{T_1}$.
Suppose that we replace $\overline{k}$. The only possibility is
$i_b$. But this is not possible by a comparison of grades of
$Y_{k-1,*}$, as explain in the paragraph preceding
Proposition~\ref{enumcond}.
Next suppose we replace $i_c = k-1$. We have $c\le h$ by (D.1). If
$i_{b+1} = \overline{k-1}$, we have $\tf_{k-1} m_{T_1} = 0$ by a
comparison of grades as before. So this case is excluded by the
assumption. We have $i_{b+1} \neq \overline{k-1}$, and then
$\tf_{k-1} m_{T_1}$ is obtained by replacing $i_c = k-1$ by $k$.

Now we consider $\tf_{k-1} m_{T_2}$. We have $i_{b+1}\neq
\overline{k}$ by (D.1) and $i_{b+1}\neq \overline{k-1}$ by above.
Therefore $\tf_{k-1} m_{T_2}$ is obtained by replacing $i_c = k-1$ by
$k$. This checks the assertion.

Case (a.3):  $j=k$.

An entry $i_c = \overline{k+1}$ could appear only at $i_{b+1}$, but we
cannot replace it by $\overline{k}$ as $i_b = \overline{k}$. Thus the
only possibility is $i_c = k$. It can appear at $c\le h+1$.

Subcase (a.3.1): $c\le h$.

We have $\gr_{k}(\ffbox{i_{h\!+\!1}}) < \gr_k(\ffbox{i_b}) <
\gr_k(\ffbox{i_c})$. In order to replace $i_c$, the contribution from
${i_b}$ must be canceled with ${i_{b-1}} = {\overline{k\!+\!1}}$,
i.e., $b''\le b'-1$. In the new tableau $T_3$, we have contributions
from ${i''_b} = {i_b}$ and  ${i_c}$, but the first one is still
canceled with ${i_{b-1}}$. Therefore $\tf_j m_{T_3}$ is obtained by
replacing $i''_c = i_c$ by $k+1$.

Subcase (a.3.2): $c=h+1$, i.e., we replace $i_{h+1} = k$ by $k+1$.

This can happen only when $i_{h+2} \neq k+1$ and
\(
  \varphi_{k,L}(m_{T_1}) \le 1,
\)
while $\varphi_{k,L}(m_{T_1}) = 0$ if $L < \gr_k(\ffbox{i_{h\!+\!1}})$
and $= 1$ if $L = \gr_k(\ffbox{i_{h\!+\!1}})$ are satisfied by the
assumption.   
Then
\begin{equation*}
\begin{split}
   & T_2=((i_1,\cdots,i_h),\\
   & \qquad\qquad\qquad (k+1,i_{h+2},\cdots,
   i_{b''-1}, \overline{k''},\cdots,\overline{k+1},
     \overline{k}, i_{b+1},\cdots, i_{\ell})), \\ 
   & \tau_{\ell,h,r}(T_2) = ((i_1,\cdots,i_h,k''+1), \\
   &\qquad\qquad\qquad
   (i_{h+2},\cdots,i_{b''-1},\overline{k''+1},\cdots,\overline{k+2},
   \overline{k},i_{b+1},\cdots,i_\ell)),\\
   & T_3 = ((i_1,\cdots,i_h,k''+1),
\\
   &\qquad\qquad\qquad
   (i_{h+2},\cdots,i_{b''-1},\overline{k''+1},\cdots,\overline{k+1},
   i_{b+1},\cdots,i_\ell)).
\end{split}
\end{equation*}
Now consider $\tf_k m_{T_3}$. As $i_{h+2}\neq k+1$, the only 
contributions to $(m_{T_3})^{(k)}$ are $i''_b = \overline{k+1}$ and
possibly $i_c$ with $c\le h$ or $i''_{h+1} = k''+1$ when $k'' = k$.
But $\gr_k(\ffbox{i_c}) \ge \gr_k(\ffbox{i''_{h\!+\!1}}) >
\gr_k(\ffbox{i''_b})$, thus the condition
\(
  \varphi_{k,L}(m_{T_3}) \le 1
\)
is preserved. Therefore $\tf_k m_{T_3}$ is obtained by replacing
$i''_b = \overline{k+1}$ by $\overline{k}$. This is equal to
$\tau_{\ell,h,r}(T_2)$.

Case (a.4): $j=k''$.

We may assume $b''\neq b$, as $b'' = b$ case was already treated in
(a.2).

Consider $(m_{T_1})^{k''}$. An entry $\overline{k''+1}$ does not
appear by our choice of $k''$. Therefore we replace $i_c = k''$ by
$k''+1$. It appears either $c\le h$ or $h+2\le c\le b''-1$, or
both. But in the first case we have $\gr_{k''}(\ffbox{i_c}) >
\gr_{k''}(\ffbox{i_{b''}})$, so either $\tf_{k''}(m_{T_1}) = 0$ or we
replace another $i_{c'} = k''$ with $h+2\le c'\le b''-1$. We have a
contradiction.
In the second case we consider $\tf_{k''} m_{T_3}$. We have
\begin{equation*}
\begin{split}
  \gr_{k''}(\ffbox{i''_c}) = \ell-2r-2c+k''
   & < \gr_{k''}(\ffbox{i''_{b''\!\!\!+\!1}}) = \ell-2r-2b''+2n-k''-2
\\
   & < \gr_{k''}(\ffbox{i''_{h\!+\!1}}) = \ell-2h+k'',
\end{split}
\end{equation*}
where the middle part contribute only when $b''+1=b$. The
contribution of $i''_{b''}$ is canceled by $i''_{b''+1}$ by the
assumption $b''\neq b$.
Therefore we replace $i''_c$ by $k''+1$ and get the assertion.

Case (a.5): $j=k''+1$.

Consider $\tf_j(m_{T_1})$. We replace either $i_c = k''+1$ or
$\overline{k''+2}$.

Subcase (a.5.1): we replace $i_c = k''+1$ by $k''+2$.

The entry $i_c = k''+1$ appears in $h+2\leq c\leq b''-1$.
We have
\begin{equation*}
  \begin{split}
    & \gr_{k''+1}(\ffbox{i_c}) < \gr_{k''+1}(\ffbox{i_{b''}}),
\\
    & \gr_{k''+1}(\ffbox{i''_c}) < \gr_{k''+1}(\ffbox{i''_{h\!+\!1}}),
  \end{split}
\end{equation*} 
and $i_c = i''_c = k''+1$ are replaced by $k''+2$ in
$\tf_{k''+1}m_{T_1}$, $\tf_{k''+1}(m_{T_3})$ respectively. Therefore
we have the assertion.

Subcase (a.5.2): we replace $i_c=\overline{k''+2}$ by 
$\overline{k''+1}$.

The entry $i_c = \overline{k''+2}$ appears at $c=b''-1$. Therefore
{\allowdisplaybreaks
\begin{equation*}
  \begin{split}
   & T_1=((i_1,\cdots,i_h),
   (k,i_{h+2},\cdots,
   i_{b''-2}, \overline{k''+2},\overline{k''}, 
   \cdots, \overline{k}, i_{b+1},\cdots, i_{\ell})), \\ 
   & T_2 = ((i_1,\cdots,i_h),(k,i_{h+2},\cdots,i_{b''-2},
   \overline{k''+1},\overline{k''},\cdots,\overline{k},
   i_{b+1},\cdots,i_\ell)),
\\
   & \tau_{\ell,h,r}(T_2) = ((i_1,\cdots,i_h,k''+2),
\\
   &\qquad\qquad\qquad
   (i_{h+2},\cdots,i_{b''-2},\overline{k''+2},\overline{k''+1},
   \cdots,\overline{k+1},
   i_{b+1},\cdots,i_\ell)),
\\
   & T_3 = ((i_1,\cdots,i_h,k''+1),
\\
   &\qquad\qquad\qquad
   (i_{h+2},\cdots,i_{b''-2},\overline{k''+2},\overline{k''+1},
   \cdots,\overline{k+1},
   i_{b+1},\cdots,i_\ell)).
  \end{split}
\end{equation*}
Therefore} $\tilde{f}_{k''+2}(m_{T_2})$ is obtained by replacing
$i_{h+1}''=k''+1$ by $k''+2$, and we get the assertion.

In the case (D.b) we can prove the assertion as in (D.a).

Finally suppose that $T_1$ is in the case (D.c). By $(\mathrm{D}.5\sim
7)$, we have $i_h\nsucceq i_{h+1}$. We assume $i_{h+1} = k\in
\{1,\cdots, n\}$, but the case $i_{h+1}\in
\{\overline{1},\cdots,\overline{n}\}$ can be studied exactly in the
same way.
We then separate the cases according to the values of $j$ and $i_c$.

Case (c.1): $j\neq k-1$ or $k$.

The assertion clearly holds. 

Case (c.2): we replace $i_c = k$ by $k+1$. ($j=k$)

The entry $i_c$ appears at $c = h+1$. We must have $i_h < k$.
Let us study other contributions $\overline{k}$, $\overline{k+1}$.
If neither appears, or both appear, then the assertion clearly holds.

Subcase (c.2.1): There exists $b$ such that $i_b = \overline{k+1}$ and
$i_{b+1}\neq \overline{k}$.

We have
\begin{equation*}
  \begin{split}
  & (m_{T_1})^{(k)}=Y_{k,\ell-2h-2r+k-2}Y_{k,\ell-2r-2b+2n-k-2},
\\
  &  (m_{T_3})^{(k)}=Y_{k,\ell-2h+k-2}Y_{k,\ell-2r-2b+2n-k-2}.
  \end{split}
\end{equation*}
As we replace $i_{h+1}$, we must have
$\ell-2h-2r+k-2 \ge \ell-2r-2b+2n-k-2$. Therefore we have
$\ell-2h+k-2 > \ell-2r-2b+2n-k-2$. Thus we replace
$i''_{h+1} = k$ by $k+1$, and we have the assertion.

Subcase (c.2.2): There exists $b$ such that $i_b = \overline{k}$ and
$i_{b-1}\neq \overline{k+1}$.

We have
\begin{equation*}
  \begin{split}
  & (m_{T_1})^{(k)}=Y_{k,\ell-2h-2r+k-2}Y_{k,\ell-2r-2b+2n-k}^{-1},
\\
  &  (m_{T_3})^{(k)}=Y_{k,\ell-2h+k-2}Y_{k,\ell-2r-2b+2n-k}^{-1}.
  \end{split}
\end{equation*}
As (D.a) does not hold, we have either $b-h \le n - r -k$ or $b-h >
n-1-k$. The second inequality contradicts with $\tf_k m_{T_1} \neq 0$.
Thus we have the first inequality. Then $\tf_k(m_{T_3})$ is
obtained by replacing $i''_{h+1} = k$ by $k+1$. We have the assertion.

Case (c.3): we replace $i_c = \overline{k+1}$ by $\overline{k}$. ($j=k$)

An entry $i_c = \overline{k+1}$ can appear only in $c\ge h+2$. We must
have $i_{c+1}\neq\overline{k}$.

Subcase (c.3.1): $i_{h+2} = k+1$.

We have $c-h \neq n-k$ thanks to $i_{c} = \overline{k+1}$ and (D.6). 
Let us study $\tf_k m_{T_3}$. We have
\begin{NB2}
\begin{equation*}
   (m_{T_2})^{(k)}=Y_{k,\ell-2h+k-2}Y_{k,\ell-2r-2h+k-2}^{-1}
   Y_{k,\ell-2r-2c+2n-k-2}
\end{equation*}
with $\ell-2r-2h+k-2\neq \ell-2r-2c+2n-k-2$.
\end{NB2}
\begin{gather*}
    \gr_k(\ffbox{i''_{h\!+\!1}}) = \ell-2h+k-2, \qquad
   \gr_k(\ffbox{i''_{h\!+\!2}}) = \ell-2r-2h+k-2,
\\
   \gr_k(\ffbox{i''_{c}}) = \ell-2r-2c+2n-k-2
\end{gather*}
with $\gr_k(\ffbox{i''_{h\!+\!2}}) \neq \gr_k(\ffbox{i''_{c}})$.

Subsubcase (c.3.1.1): 
\(
   \gr_k(\ffbox{i''_{h\!+\!2}}) < \gr_k(\ffbox{i''_{c}})
   < \gr_k(\ffbox{i''_{h\!+\!1}}).
\)

Then $\tf_k(m_{T_3})$ is obtained by replacing $i''_{h+1} = k$ by
$k+1$. On the other hand, $\tau_{\ell,h,r}(T_2)$ is obtained by
replacing the pair $(i'_{h+1},i'_b) = (k,\overline{k})$ by 
$(k+1,\overline{k+1})$. The results are the same.

Subsubcase (c.3.1.2): 
\(
   \gr_k(\ffbox{i''_{c}}) < \gr_k(\ffbox{i''_{h\!+\!2}})
\)
or
\(
   \gr_k(\ffbox{i''_{h\!+\!1}}) \le 
   \gr_k(\ffbox{i''_{c}}).
\)

Then $\tilde{f}_k(m_{T_3})$ is obtained by replacing
$i_c''=\overline{k+1}$ by $\overline{k}$. On the other hand, the
inequality in (D.a) is violated for $(i'_{h+1},i'_c) =
(k,\overline{k})$, so $T_2$ is of type (D.c). Therefore
$\tau_{\ell,h,r}(T_2)$ is the same as $T_2$. The results are the same.

Subcase (c.3.2): $i_{h+2}\neq k+1$.

We have $(m_{T_1})^{(i)}=Y_{i,-2h+i-1-2r}Y_{i,2n-i-1-2r-2t}$ and so
$-2h+i-1-2r < 2n-i-1-2r-2t$. We study the two cases $-2h+i-1-2r<
2n-i-1-2r-2t<-2h+i-1$ and $-2h+i-1 \leq 2n-i-1-2r-2t$ separately as
above and we get the assertion.

Case (c.4): $j=k-1$.

\begin{NB2}
I think that the proof was wrong in this case.   
\end{NB2}

We replace either $i_c = k-1$ or $\overline{k}$.

Subcase (c.4.1): we replace $i_c=k-1$ by $k$.

In this case, $i_c = k-1$ can appear only at $c=h$.
As
\(
   \gr_{k-1}(\ffbox{i_h}) > \gr_{k-1}(\ffbox{i_{h+1}}),
\)
we must have an entry $i_b = \overline{k}$ with
\(
   \gr_{k-1}(\ffbox{i_h}) > 
   \gr_{k-1}(\ffbox{i_b}) \ge \gr_{k-1}(\ffbox{i_{h+1}}).
\)
Also we have $b-h\neq n-k$ by (D.7). This means that the equality
cannot hold in the second inequality. Then these contradict with the
assumption that $T_1$ is of type (D.c).

Subcase (c.4.2): we replace $i_c=\overline{k}$ by $\overline{k-1}$.
The entry $i_c=\overline{k}$ appears in $c\geq h+2$ and we have
$i_{c+1}\neq \overline{k-1}$.

Let us consider $(m_{T_3})^{(k-1)}$. The contributions are
$i_c = \overline{k}$, $i_{h+1} = k$ and possibly $i_h = k-1$. Even if
we have $i_h = k-1$, its contribution is killed by $i_{h+1} = k$. So
$\tf_{k-1}(m_{T_3})$ is obtained by replacing $i_c$ by $\overline{k-1}$.
On the other hand, $T_2$ is of type (D.c), so is not changed under
$\tau_{\ell,h,r}$. We therefore the assertion.
\end{proof}

\begin{thm}\label{gentau}
  Let $0\leq h\leq \ell-1$.  The map $\tau_{\ell,h,r}$ induces a
  crystal isomorphism from $\mathcal M_{I_0}(M_{\ell,h,r})$ to
  $\mathcal M_{I_0}(M_{\ell,h+1,r})$.
\end{thm}

\begin{proof}
First we have $\tau_{\ell,h,r}((1,\dots,h),(h+1,\dots,\ell))
= ((1,\dots,h+1),(h+2,\dots,\ell))$ and
$m_{((1,\dots,h),(h+1,\dots,\ell))} = M_{\ell,h,r}$,
$m_{((1,\dots,h+1),(h+2,\dots,\ell))} = M_{\ell,h+1,r}$.

Suppose that a monomial $m$ can be expressed in two ways as $m = m_T =
m_{T'}$. We write $\tf_{i_1}\tf_{i_2}\dots \tf_{i_N} M_{\ell,h,r} = m$
for $N \ge 0$, $i_p\in I_0$. Using Lemma~\ref{lem:tau} inductively, we
get 
\[
  m_{\tau_{\ell,h,r}(T)} = 
  \tf_{i_1}\tf_{i_2}\dots \tf_{i_N}
  m_{\tau_{\ell,h,r}((1,\dots,h),(h+1,\dots,\ell))}
  = m_{\tau_{\ell,h,r}(T')}.
\]
This shows that $\tau_{\ell,h,r}$ induces a map from $\mathcal
  M_{I_0}(M_{\ell,h,r})$ to $\mathcal M_{I_0}(M_{\ell,h+1,r})$.
We already know that both $\mathcal M_{I_0}(M_{\ell,h,r})$ and
$\mathcal M_{I_0}(M_{\ell,h+1,r})$ are isomorphic to the crystal of
the finite dimensional fundamental representation corresponding to
the vertex $\ell$. Therefore there exists a (unique) crystal isomorphism
$\tau\colon \mathcal M_{I_0}(M_{\ell,h,r})\to
\mathcal M_{I_0}(M_{\ell,h+1,r})$. Both $\tau$ and $\tau_{\ell,h,r}$ maps 
the highest weight vector $M_{\ell,h,r}$ to $M_{\ell,h+1,r}$ and
satisfies $\tf_i \tau(m) = \tau \tf_i(m)$ if $\tf_i(m)\neq 0$.
As such a map is unique, we have $\tau = \tau_{\ell,h,r}$.
\end{proof}

\begin{NB2}
The following was original.
\end{NB2}

**********************************************************************

\begin{thm}\label{gentau} For $0\leq h\leq \ell-1$ the map $m_{T}\mapsto m_{\tau_{\ell,h,r}(T)}$ is a crystal isomorphism from $\mathcal M_{I_0}(M_{\ell,h,r})$ to $\mathcal M_{I_0}(M_{\ell,h+1,r})$.\end{thm}

\begin{proof} For the shortness of formulas we suppose in this proof that the suffixes of $m_t$ start with $0$. First we have $\tau_{\ell,h,r}((1,\cdots,h),(h+1,\cdots,\ell))=((1,\cdots,h+1),(h+2,\cdots,\ell))$. Let $T_1=((i_1,\cdots,i_h),(i_{h+1},\cdots,i_\ell))\in B_{\ell,h,r}$ and $j\in I$ such that $\tilde{f}_j (m_{T_1})=m_{T_3}\neq 0$ where $T_3=((i_1'',\cdots,i_h''),(i_{h+1}'',\cdots,i_\ell''))\in D_{\ell,h,r}$ is obtained from $T_1$ by changing one $j$ to $j+1$ or one $\overline{j+1}$ to $\overline{j}$. Let $t$ such that $i_t\in\{j,\overline{j+1}\}$ and $i_t''\in\{j+1,\overline{j}\}$. Let $T_2=((i_1',\cdots,i_{h+1}'),(i_{h+2}',\cdots,i_\ell'))=\tau_{\ell,h,r}(T_1)$. We suppose that $m_{T_1}$ and $m_{T_2}$ correspond to the same vector in $\mathcal{B}_{I_0}(\omega_l)$ and we prove that $m_{\tau_{\ell,h,r}(T_3)}=\tilde{f}_j(m_{T_2})$. 

First suppose that $T_1$ satisfies the condition (1) of the definition \ref{defitau} (we use the notations $i,h',h''$ of this definition). We have different cases depending of the value of $j$.

Case $(\alpha)$ : $j\notin\{i-1,i,i+(h'-h''),i+(h'-h'')+1\}$. We have clearly the result.

Case $(\beta)$ : $j=i-1$. If $i_t=\overline{i}$, we have $t=h'$ and $i_{t+1}\neq \overline{i-1}$. As $n-h'+h+1-r\leq i\leq n-1-h'+h$, we have $-2h-2r+i<-2h'-2r+2n-i< -2H+i$ for all $H\leq h$. So 
\begin{equation*}
(m_{T_1})^{(i-1)}\in Y_{i-1,-2h-2r+i}^{-1}Y_{i-1,-2h'-2r+2n-i}\times \ZZ[Y_{i-1,L}^{\pm}]_{L>-2h'-2r+2n-i},
\end{equation*}
and $i_t=i_t''$, contradiction. So $i_t=i-1$ and $t\leq h$. If $i_{h'+1}=\overline{i-1}$, we have 
\begin{equation*}
(m_{T_1})^{(i-1)}= Y_{i-1,-2h-2r+i}^{-1}Y_{i-1,-2t+i},\end{equation*}
with $-2h-2r+i<-2t+i$, and $\tilde{f}_j(m_{T_1})= 0$, contradiction. So $i_{h'+1}\neq\overline{i-1}$ and
\begin{equation*}
(m_{T_1})^{(i-1)}= Y_{i-1,-2h-2r+i}^{-1}Y_{i-1,-2h'-2r+2n-i}Y_{i-1,-2t+i},\end{equation*}
with $-2h-2r+i<-2h'-2r+2n-i<-2t+i$. So $(m_{T_2})^{(i-1)}=Y_{i-1,-2t+i}$ and $\tilde{f}_j(m_{T_2})$ is obtained from $m_{T_2}$ by replacing $i_t'=i-1$ by $i$, and we get the result. 

Case $(\gamma)$ :  $j=i$. If $i_t=\overline{i+1}$, we have $t=h'-1$, contradiction because $i_{h'}=\overline{i}$. So $i_t=i$ and $t\leq h+1$. 
Suppose that $t\leq h$. If $i_{h'-1}\neq \overline{i+1}$, we have 
\begin{equation*}
(m_{T_1})^{(i)}=Y_{i,-2h+i-1-2r}^{\epsilon}Y_{i,2n+1-i-2r-2h'}^{-1}Y_{i,-2t+1+i},
\end{equation*}
where $\epsilon\in \{0,1\}$, so $i_t''=i_t$, contradiction. So $i_{h'-1}=\overline{i+1}$ and $h''\leq h'-1$. So
\begin{equation*}
(m_{T_1})^{(i)}=Y_{i,-2h+i-1-2r}^{\epsilon}Y_{i,-2t+1+i},
\end{equation*}
where $\epsilon\in\{0,1\}$, and
\begin{equation*}
(m_{T_2})^{(i)}=Y_{i,-2t+1+i}Y_{i,2n-1-i-2h'-2r},
\end{equation*}
with $-2t+1+i > 2n-1-i-2h'-2r$, so $\tilde{f}_j(m_{T_2})$ is obtained from $m_{T_2}$ by replacing $i_t'=i$ by $i+1$ and we get the result. 
Suppose that $t=h+1$. So $i_{h+2}\neq i+1$. As 
\begin{equation*}
m_{T_1}^{(i)}\in Y_{i,-2h+i-1-2r}\ZZ[Y_{i,L}^{\pm}]_{L\geq 2n+1-i-2h'-2r},
\end{equation*}
we have $\phi_i(m_{T_1})=1$. We have $T_3=((i_1,\cdots,i_h),(i+1,i_{h+2},\cdots,i_{\ell}))$, and
\begin{equation*}
\begin{split}
\tau_{\ell,h,r}(T_3) &= ((i_1,\cdots,i_h,i+(h'-h'')+1),
\\
&\qquad\qquad(i_{h+2},\cdots,i_{h''-1},\overline{i+(h'-h'')+1},\cdots,\overline{i+2},\overline{i},i_{h'+1},\cdots,i_\ell)),
\end{split}
\end{equation*}
and
\begin{equation*}
\begin{split}
T_2 = &((i_1,\cdots,i_h,i+(h'-h'')+1),
\\
&\qquad\qquad(i_{h+2},\cdots,i_{h''-1},\overline{i+(h'-h'')+1},\cdots,\overline{i+1},i_{h'+1},\cdots,i_\ell)).
\end{split}
\end{equation*}
If $h''<h'$ we have
\begin{equation*}
m_{T_2}^{(i)}\in Y_{i,2n-i-1-2h'-2r}\times\ZZ[Y_{i,L}^{\pm}]_{L\geq 2n+1-i-2h'-2r}.
\end{equation*}
As $\phi_i(m_{T_2})=\phi_i(m_{T_1})=1$, $\tilde{f}_i(m_{T_2})$ is obtained from $m_{T_2}$ by replacing $i_{h'}'=\overline{i+1}$ by $\overline{i}$ and we get the result. If $h''=h'$, we have 
\begin{equation*}
m_{T_2}^{(i)}\in Y_{i,2n-i-1-2h'-2r}Y_{i,-2h+i+1}^{-1}\times\ZZ[Y_{i,L}^{\pm}]_{L\geq 2n+1-i-2h'-2r},
\end{equation*}
with $-2h+i+1 > 2n+1-i-2r-2h'$, and we get the result with the same argument.


Case $(\delta)$ : $j=i+h'-h''$ with $h'>h''$. If $i_t=\overline{i+h'-h''+1}$ we have $t=h''-1$ and $i_{h''}=\overline{i+h'-h''}$, contradiction. So $i_t=i+h'-h''$, and $t\leq h$ or $h+2\leq t\leq h''-1$. Suppose that $t\leq h$. As for all $H\geq h+2$, $-2r-2H+2+i+h'-h''+1 < 2n+1-i-h'-h''-2r$, we have
\begin{equation*}
(m_{T_1})^{(j)}\in Y_{j,-2t+1+i+h'-h''}Y_{j,2n+1-i-h'-h''-2r}^{-1}\times\ZZ[Y_{j,L}^{\pm}]_{L<2n+1-i-h'-h''-2r},
\end{equation*}
with $-2t+1+i+h'-h'' > 2n+1-i-h'-h''-2r$, so $i_t=i_t''$, contradiction. So $h+2\leq t\leq h'' - 1$ and we have
\begin{equation*}
(m_{T_1})^{(j)}\in Y_{j,-2t+1+i+h'-h''-2r}Y_{j,2n+1-i-h'-h''-2r}^{-1}\times\ZZ[Y_{j,L}^{\pm}]_{L>2n+1-i-h'-h''-2r},
\end{equation*}
with $-2t+1+i+h'-h''-2r < 2n+1-i-h'-h''-2r$. We have
\begin{equation*}
(m_{T_2})^{(j)}\in Y_{j,-2t+1+i+h'-h''-2r}Y_{j,i+h'-h''+1-2h}^{-1}\times\ZZ[Y_{j,L}^{\pm}]_{L>2n+1-i-h'-h''-2r},
\end{equation*}
with $-2t+1+i+h'-h''-2r < i+h'-h''-2h+1$. So $\tilde{f}_j(m_{T_2})$ is obtained from $m_{T_2}$ by replacing $i_t'=j$ by $j+1$, and we get the result.


Case $(\epsilon)$ : $j=i+(h'-h'')+1$. Suppose that $i_t=i+(h'-h'')+1$. We have $h+2\leq t\leq h''-1$. The box $\ffbox{j}_{-2t+2-2r}$ gives the term $Y_{j,-2t+2-2r+i+h'-h''}$ in $m_{T_1}$, and we have
\begin{equation*}
(m_{T_2})^{(j)} = (m_{T_1})^{(j)}Y_{j,i+h'-h''-2h}Y_{j,-2r-h''+2n-i-h'}^{-1},
\end{equation*} 
where $i+h'-h''-2h>-2r-h''+2n-i-h'>-2t+2-2r+i+h'-h''$. So in both $m_{T_1}$ and $m_{T_2}$, we replace the box $\ffbox{j}_{-2t+2-2r}$ by $\ffbox{j+1}_{-2t+2-2r}$ when we apply $\tilde{f}_j$, and we get the result. Now suppose that $i_t=\overline{j+1}$. We have $t=h''-1$. So
\begin{equation*}
T_3 = ((i_1,\cdots,i_h),(i,i_{h+2},\cdots,i_{h''-2},\overline{j},\overline{j-1},\cdots,\overline{i},i_{h'+1},\cdots,i_\ell)),
\end{equation*}
\begin{equation*}
\begin{split}
\tau_{\ell,h,r}(T_3) = &((i_1,\cdots,i_h,j+1),
\\
&\qquad\qquad(i_{h+2},\cdots,i_{h''-2},\overline{j+1},\overline{j},\cdots,\overline{i+1},i_{h'+1},\cdots,i_\ell)).
\end{split}
\end{equation*}
But we have
\begin{equation*}
T_2 = ((i_1,\cdots,i_h,j),(i_{h+2},\cdots,i_{h''-2},\overline{j+1},\cdots,\overline{i+1},i_{h'+1},\cdots,i_\ell)),
\end{equation*}
so $\tilde{f}_j(m_{T_2})$ is obtained from $m_{T_2}$ by replacing $i_{h+1}'=j$ by $j+1$, and we get the result.


Secondly suppose that $T_1$ satisfies the condition (2) of the definition \ref{defitau}. We get the result with arguments analogous to the above arguments for the condition (1).

Finally suppose that $T_1$ satisfies the condition (3) of the definition \ref{defitau}. We have $i_h\nsucceq i_{h+1}$. Let us study the different cases.

Case $(\alpha)$ : $i_{h+1}=i\in\{1,\cdots,n\}$ and $j\notin \{i-1,i\}$. The result is clear. 

Case $(\beta)$ : $i_{h+1}=i\in\{1,\cdots,n\}$ and $j = i$. Suppose that $i_t=i$. So we have $t = h+1$ and $i_h<i$. If there is no $h'\geq h+2$ such that $i_{h'}=\overline{i+1}$ or $i_{h'}=\overline{i}$ the result is clear. If there is $h'\geq h+2$ such that $i_{h'}=\overline{i+1}$ and $i_{h'+1}=\overline{i}$ the result is clear. If there is $h'\geq h+2$ such that $i_{h'}=\overline{i+1}$ and $i_{h'+1}\neq \overline{i}$, we have $(m_{T_1})^{(i)}=Y_{i,-2h+i-1-2r}Y_{i,-2h'-2r+2n-i-1}$. So $-2h+i-1-2r\geq -2h'-2r+2n-i-1$. So we have $(m_{T_2})^{(i)}=Y_{i,-2h+i-1}Y_{i,-2h'-2r+2n-i-1}$ with $-2h+i-1>-2h'-2r+2n-i-1$, and $\tilde{f}_i(m_{T_2})$ is obtained from $m_{T_2}$ by replacing $i_t'=i$ by $i+1$, and the result is clear. If there is $h'\geq h+2$ such that $i_{h'}=\overline{i}$ and $i_{h'-1}\neq \overline{i+1}$, we have $(m_{T_2})^{(i)}=Y_{i,-2h+i-1}Y_{i,-2h'-2r+2n-i+1}^{-1}$ so $\tilde{f}_j(m_{T_2})$ is obtained from $m_{T_2}$ by replacing $i_t'=i$ by $i+1$ and the result is cl!
 ear. 
Suppose that $i_t=\overline{i+1}$. So we have $t\geq h+2$ and $i_{t+1}\neq\overline{i}$. Suppose moreover that $i_{h+2}=i+1$. We have $t-h\neq n-i$ because $i_t=\overline{i+1}$. We have $(m_{T_2})^{(i)}=Y_{i,-2h+i-1}Y_{i,-2r-2h+i-1}^{-1}Y_{i,2n-i-1-2r-2t}$ with $2n-i-1-2r-2t\neq -2r-2h+i-1 $. If $-2r-2h+i-1 < 2n-i-1-2r-2t<-2h+i-1$, $\tilde{f}_j(m_{T_2})$ is obtained from $m_{T_2}$ by replacing $i_{h+1}'=i$ by $i+1$, and $\tau_{\ell,h,r}(T_3)$ is obtained from $T_3$ by replacing $(i_{h+1}'',i_t'')=(i,\overline{i})$ by $(i+1,\overline{i+1})$ and so we get the result. If $-2r-2h+i-1 > 2n-i-1-2r-2t$ or $2n-i-1-2r-2t \geq -2h+i-1$, then $\tilde{f}_j(m_{T_2})$ is obtained from $m_{T_2}$ by replacing $i_t'=\overline{i+1}$ by $i_t'=\overline{i}$, and $T_3$ satisfies the condition (3) of the definition \ref{defitau} and so we get the result. Suppose that $i_{h+2}\neq i+1$, we have $(m_{T_1})^{(i)}=Y_{i,-2h+i-1-2r}Y_{i,2n-i-1-2r-2t}$ and so $-2h+i-1-2r < 2n-i-1-2r-2t$. So we have two !
 cases $-2h+i-1-2r< 2n-i-1-2r-2t<-2h+i-1$ and $-2h+i-1 \leq 2n-i-1-2r-2t$ as above and we get the result.

Case $(\gamma)$ : $i_{h+1}=i\in\{1,\cdots,n\}$ and $j = i-1$. Then $T_3$ satisfies the property (3) of the definition \ref{defitau}. Suppose that $i_t=i-1$. We have $t=h$. If moreover there is no $h'$ such that $i_{h'}=\overline{i}$, or there is $h'$ such that $i_{h'}=\overline{i}$ and $i_{h'+1}=\overline{i-1}$, $\tilde{f}_j(m_{T_2})$ is obtained from $m_{T_2}$ by replacing $i_t'=i-1$ by $\overline{i-1}$, and so the result is clear. If there is $h'$ such that $i_{h'}=\overline{i}$ and $i_{h'+1}\neq\overline{i-1}$, we have 
\begin{equation*}
(m_{T_1})^{(i-1)} = Y_{i-1,i-2h}Y_{i-1,-2h-2r+i}^{-1}Y_{i-1,-2r-2h'+2n-i},
\end{equation*}
 with $i-2h\leq 2n-2h'-2r-i$ or $i-2h-2r > 2n-2h'-2r-i$, so $i_t''=i_t$, contradiction. Suppose that $i_t=\overline{i}$. We have $t\geq h+2$ and $i_{t+1}\neq \overline{i-1}$. If $i_{h}\prec i-1$, $\tilde{f}_j(m_{T_2})$ is obtained from $m_{T_2}$ by replacing $i_t'=\overline{i}$ by $\overline{i-1}$ and we get the result. If $i_h=i-1$, we have $(m_{T_2})^{(i-1)}=Y_{i-1,2n-2t-2r-i}$ and $\tilde{f}_j(m_{T_2})$ is obtained from $m_{T_2}$ by replacing $i_t'=\overline{i}$ by $\overline{i-1}$, and we get the result.

The cases $i_{h+1}\in\{\overline{n},\cdots,\overline{1}\}$ are treated in an analog way.\end{proof}

 }}_{L+2(n-i-1)}$ appears in one of the decomposition of $M$, and one product $\ffbox{i+1}_L\ffbox{\overline{i+1}}_{L+2(n-i-1)}$ in the other (with the same $p$ and $q$). Contradiction because it follows from the conditions (2), (3) and (4) that only one of this product is allowed.\end{proof}






\end{NB}

\subsubsection{}\label{Dn2}
Now we study $\mathcal B(\varpi_\ell)$ for $2\le \ell\le n-2$. Let
\(
  M_{0,0} = Y_{\ell,0} Y_{0,\ell-1}^{-1} Y_{0,\ell+1}^{-1}
    = \ffbox{1}_{\ell-1}\ffbox{2}_{\ell-3} \cdots
      \ffbox{\ell}_{1-\ell}.
\)
Then $\mathcal M(M_{0,0}) \simeq \mathcal B(\varpi_\ell)$ by
Proposition~\ref{otherd}. 

For $0\le j\le\ell$, we set
\begin{equation*}
\begin{split}
  M_{j,0} &= \left( \ffbox{1}_{2n-\ell+2j-5}\ffbox{2}_{2n-\ell+2j-7}
   \cdots \ffbox{j}_{2n-\ell-3}\right)
\times
  \left(\ffbox{j\!+\!1}_{\ell-1}\ffbox{j\!+\!2}_{\ell-3}
    \cdots \ffbox{\ell}_{1-\ell+2j}\right)
\\
   &=
   \prod_{a=1}^j \ffbox{a}_{2n-\ell-2a+2j-3} \times
   \prod_{a=j+1}^\ell \ffbox{a}_{\ell-2a+2j+1}
\\
 &=
 \begin{cases}
   Y_{\ell,0} Y_{0,\ell-1}^{-1} Y_{0,\ell+1}^{-1} & \text{if $j=0$},
\\
   Y_{\ell,2} Y_{0,\ell+1}^{-1} Y_{0,2n-\ell-1}^{-1}
   Y_{1,\ell+1}^{-1} Y_{1,2n-\ell-3} &\text{if $j=1$},
\\
   Y_{\ell,2j} Y_{0,2n-\ell+2j-5}^{-1}Y_{0,2n-\ell+2j-3}^{-1}
       Y_{j,\ell+j}^{-1}Y_{j,2n-\ell+j-4}
   &\text{otherwise.}
 \end{cases}
\end{split}
\end{equation*}
Note that $M_{\ell,0} = \tau_{2n-4}(M_{0,0})$. 

For a tableau $T = ((i_1,\dots,i_j),(i_{j+1},\dots,i_\ell))$ we define
$m_{T;j,0}$ by replacing the $a^{\mathrm{th}}$-entry by $i_a$.

\begin{Claim}
We have $M_{j,0}\in \mathcal M(M_{0,0})$ for $0\leq j \le \ell$. 
\end{Claim}

In fact, by Theorem~\ref{isot} we have $m_{T;j,0}$
with $T = (3,\dots,\ell+1,\overline{2})$ is contained in $\mathcal
M_{I_0}(M_{j,0})$ as $M_{j,0}=m_{(1,\cdots,\ell);j,0}$. Then we get $\tf_0 m_{T;j,0} = m_{T';j+1,0}$ with
$T' = (1,3,4,\dots,\ell+1)$. Again by Theorem~\ref{isot} this is contained in $\mathcal M_{I_0}(M_{j+1,0})$ as $M_{j+1,0}=m_{(1,\cdots,\ell);j+1,0}$. By
induction we obtain the claim.

\begin{NB}
we have
\begin{equation*}M_{j+1,0}=\tilde{e}_\ell\tilde{e}_{\ell-1}\cdots\tilde{e}_2\tilde{f}_0
\tilde{f}_2^2\tilde{f}_3^2\cdots\tilde{f}_\ell^2\tilde{f}_{\ell+1}\tilde{f}_{\ell+2}\cdots\tilde{f}_{n-1}\tilde{f}_{n+1}\tilde{f}_n\tilde{f}_{n-1}\cdots\tilde{f}_{\ell+1}\tilde{f}_1\tilde{f}_2\cdots\tilde{f}_\ell M_{j,0}.\end{equation*}
We have $M_{j,0}\in \mathcal M(M_{0,0})$ for all $j$ by induction.
\end{NB}

We have $\wt(M_{j,0}) = \varpi_\ell - j\delta$. Thus $M_{1,0} =
z_\ell^{-1}(M_{0,0})$. As $M_{\ell,0} = \tau_{2n-4}(M_{0,0})$, $\mathcal B(M_{0,0})$ is preserved under $\tau_{2n-4}$ and we have
$(z_\ell)^{-\ell} = \tau_{2n-4}$.
As in type $A^{(1)}_n$, it is enough to study $\mathcal
M(M_{0,0})/\tau_{2n-4}$.
We extend the definition of $M_{j,0}$ from $0\le j \le \ell$ to all
$j\in \ZZ$ so that $M_{j+\ell,0} = \tau_{2n-4} M_{j,0}$.
The same applies to other various other monomials introduced below
though we do not mention it hereafter.

\begin{NB}
Let $0\le j\le\ell$.    
\end{NB}
If we apply $\tilde e_0$ to $M_{j,0}$, we get
the monomial given by replacing $\ffbox{2}_*$ by
$\ffbox{\overline{1}}_{4-2n+*}$, that is
\begin{equation*}
\begin{split}
  \tilde e_0(M_{j,0})
  = 
  & \prod_{a=3}^{j} \ffbox{a}_{2n-\ell-2a+2j-3}
   \times
   \prod_{a=\max(j+1,3)}^\ell \ffbox{a}_{\ell-2a+2j+1}
\\
  & \qquad \times
  \begin{cases}
   \ffbox{1}_{\ell-1} \ffbox{\overline{1}}_{\ell-2n+1} & \text{if $j=0$},
\\
   \ffbox{1}_{2n-\ell-3}
   \ffbox{\overline{1}}_{\ell-2n+3} &\text{if $j=1$},
\\
   \ffbox{1}_{2n-\ell+2j-5}
   \ffbox{\overline{1}}_{-\ell+2j-3}
   &\text{if $2\le j\le\ell$}.
  \end{cases}
\end{split}
\end{equation*}
Let $N_{j,1}$ denote $\ffbox{1}_* \ffbox{\overline{1}}_*$ in the right
hand side. We have
\begin{equation*}
   N_{j,1} =
\begin{cases}
    Y_{0,\ell-3} Y_{0,\ell+1}^{-1} & \text{if $j=0$},
\\
    Y_{0,\ell-1} Y_{0,2n-\ell-1}^{-1}
      Y_{1,\ell+1}^{-1} Y_{1,2n-\ell-3} & \text{if $j=1$},
\\
    Y_{0,2n-\ell+2j-7} Y_{0,2n-\ell+2j-3}^{-1} & \text{if
    $2\le j\le \ell$}.
\end{cases}
\end{equation*}

We define $M_{j,1}$ by replacing $\ffbox{a}_*$ by $\ffbox{a\!-\!2}_*$ for
$a\ge 3$ and $\ffbox{2}_*$ by $\ffbox{\overline{1}}_{4-2n+*}$ 
in $M_{j,0}$, that is
\begin{equation*}
\begin{split}
  M_{j,1}
  = 
  & \prod_{a=1}^{j-2} \ffbox{a}_{2n-\ell-2a+2j-7}
   \times
   \prod_{a=\max(j-1,1)}^{\ell-2} \ffbox{a}_{\ell-2a+2j-3}
   \times N_{j,1}.
\end{split}
\end{equation*}
\begin{NB}
We extend the definition of $M_{j,1}$ so that $M_{j+\ell,1} =
\tau_{2n-4}(M_{j,1})$ for all $j\in\ZZ$.
\end{NB}

We have
\begin{equation*}
\begin{split}
   \wt(M_{j,1}) &= \varpi_\ell - j\delta + \alpha_0
   + \sum_{a=1}^{\ell-2} \left(\alpha_a + \alpha_{a+1}\right)
\\
   &= \varpi_\ell - (j-1)\delta - \alpha_{\ell-1} 
   - 2\alpha_\ell - 2\alpha_{\ell+1} - \cdots - 2\alpha_{n-2}
   - \alpha_{n-1} - \alpha_n
\\ 
   &= \varpi_{\ell-2} - (j-1)\delta .
\end{split}
\end{equation*}

We recursively define $M_{j,k}$ by replacing $\ffbox{a}_*$ by
$\ffbox{a\!-\!2}_*$ for $a\ge 3$ and $\ffbox{2}_*$ by
$\ffbox{\overline{1}}_{4-2n+*}$  in $M_{j,k-1}$ until all boxes are
either $\ffbox{1}_*$ or $\ffbox{\overline{1}}_*$. We have
$k = 0, \dots, \lfloor \ell/2\rfloor$ where $\lfloor \ell/2\rfloor$ is
the largest integer which does not exceed $\ell/2$ (the integer part
of $\ell/2$). We define $N_{j,k}$ in the same way. We have
\(
   \wt(M_{j,k}) = \varpi_{\ell-2k} - (j-k)\delta.
\)

Let us give $M_{j,k}$, $N_{j,k}$ 
\begin{NB}
for $0\le j < \ell$, $0\le k\le \lfloor \ell/2\rfloor$   
\end{NB}
 explicitly.
\begin{enumerate}
\item $k < \lfloor j/2\rfloor$: 
\begin{equation*}
\begin{split}
  N_{j,k}
  & = \prod_{a=1}^k \left( \ffbox{1}_{2n-\ell-4a+2j-1}
          \ffbox{\overline{1}}_{-\ell-4a+2j+1} \right)
\\
  &= Y_{0,2n-\ell-4k+2j-3} Y_{0,2n-\ell+2j-3}^{-1}
  ,
\\
  M_{j,k}
  & = N_{j,k} \times 
    \prod_{a=1}^{j-2k} \ffbox{a}_{2n-\ell-4k-2a+2j-3}
    \times
    \prod_{a=j-2k+1}^{\ell-2k} \ffbox{a}_{\ell-2(a-j+2k)+1}
  ,
\\
  &=
  \begin{aligned}[t]
  & Y_{\ell-2k,2j-2k}
  Y_{0,2n-\ell-4k+2j-5}^{-1} Y_{0,2n-\ell+2j-3}^{-1}
\\
  & \qquad\times
  Y_{j-2k,j+\ell-2k}^{-1} Y_{j-2k,2n-\ell+j-2k-4}
  \end{aligned}
\end{split}
\end{equation*}

\begin{NB} The following case can be absorbed in the case
$j$ even and $k\ge j/2$.

\item $j$ is even and $k = j/2$:
\begin{equation*}
\begin{split}
   N_{j,j/2} &= Y_{0,2n-\ell-3} Y_{0,2n-\ell+2j-3}^{-1},
\\
   M_{j,j/2} &= N_{j,j/2} \times
   \prod_{a=1}^{\ell-j} \ffbox{a}_{\ell-2a+1}
\\
   &=
   \begin{cases}
   Y_{\ell-j,j}Y_{0,\ell-1}^{-1}Y_{0,\ell+1}^{-1}
   Y_{0,2n-\ell-3} Y_{0,2n-\ell+2j-3}^{-1}
   & \text{if $0\le j\le \ell-2$}
   \\
   Y_{1,\ell-1} Y_{0,\ell+1}^{-1}
   Y_{0,2n-\ell-3} Y_{0,2n-\ell+2j-3}^{-1}
    & \text{if $j = \ell - 1$}
   \end{cases}
\end{split}
\end{equation*}
\end{NB}

\item $j$ is odd and $k = (j-1)/2$:
\begin{equation*}
\begin{split}
   N_{j,(j-1)/2} &= Y_{0,2n-\ell-5} Y_{0,2n-\ell+2j-3}^{-1},
\\
   M_{j,(j-1)/2} &= N_{j,(j-1)/2} \times
   \ffbox{1}_{2n-\ell-3}
   \times
   \prod_{a=2}^{\ell-j+1} \ffbox{a}_{\ell-2a+3}
\\
   &=
   Y_{0,\ell+1}^{-1} Y_{0,2n-\ell+2j-3}^{-1} Y_{1,\ell+1}^{-1}
   Y_{1,2n-\ell-3} Y_{\ell-j+1,j+1},
\end{split}
\end{equation*}

\begin{NB}
This case can be absorbed into $j$ odd, $k\ge (j-1)/2$ below, if we
consider only $M_{j,k}$, not $N_{j,k}$. 
\end{NB}

\item $j$ is even and $k \ge j/2$:
\begin{equation*}
\begin{split}
   N_{j,k} & = N_{j,j/2}
   \times
    \prod_{a=1}^{k-j/2} \left( \ffbox{1}_{\ell-4a+3}
      \ffbox{\overline{1}}_{\ell- 4a - 2n + 5}\right)
\\
   &= Y_{0,\ell-4k+2j+1} Y_{0,\ell+1}^{-1}
      Y_{0,2n-\ell-3}Y_{0,2n-\ell+2j-3}^{-1}
   ,
\\
  M_{j,k} &=
     N_{j,k} \times \prod_{a=1}^{\ell-2k} \ffbox{a}_{\ell-2a-4k+2j+1}
\\
  &= 
  \begin{cases}
    Y_{0,-\ell+2j+1} Y_{0,\ell+1}^{-1}
      Y_{0,2n-\ell-3}Y_{0,2n-\ell+2j-3}^{-1} & \text{if $k = \ell/2$},
  \\
    Y_{0,\ell+1}^{-1} Y_{0,2n-\ell-3}
    Y_{0,2n-\ell+2j-3}^{-1} Y_{1,-\ell+2j+1} & \text{if $k = (\ell-1)/2$},
  \\
  \begin{aligned}[c]
    & Y_{0,\ell-4k+2j-1}^{-1}Y_{0,\ell+1}^{-1}
      Y_{0,2n-\ell-3}Y_{0,2n-\ell+2j-3}^{-1}
    Y_{\ell-2k,2j-2k}
  \end{aligned}
   & \text{otherwise},
  \end{cases}
\end{split}
\end{equation*}

\item $j$ is odd and $k\ge (j+1)/2$:
\begin{equation*}
\begin{split}
   N_{j,k} &= 
   \begin{aligned}[t]
   & N_{j,(j-1)/2} \times \ffbox{1}_{2n-\ell-3}
   \ffbox{\overline{1}}_{\ell-2n+3}
   \\ & \qquad
   \times
   \prod_{a=1}^{k-(j+1)/2}
    \left( \ffbox{1}_{\ell-4a+1}
      \ffbox{\overline{1}}_{\ell-4a-2n+3} \right)
   \end{aligned}
\\
   &= Y_{0,\ell-4k+2j+1} Y_{0,2n-\ell+2j-3}^{-1}
   Y_{1,\ell+1}^{-1} Y_{1,2n-\ell-3}
   ,
\\
   M_{j,k} &= N_{j,k} \times
    \prod_{a=1}^{\ell-2k} \ffbox{a}_{\ell-2a-4k+2j+1}
\\
  &= 
  \begin{cases}
  Y_{0,-\ell+2j+1} Y_{0,2n-\ell+2j-3}^{-1}
   Y_{1,\ell+1}^{-1} Y_{1,2n-\ell-3} & \text{if $k=\ell/2$},
  \\
  Y_{0,2n-\ell+2j-3}^{-1}
   Y_{1,-\ell+2j+1} Y_{1,\ell+1}^{-1} Y_{1,2n-\ell-3}
  & \text{if $k=(\ell-1)/2$},
  \\
  Y_{0,\ell-4k+2j-1}^{-1} Y_{0,2n-\ell+2j-3}^{-1}
     Y_{1,\ell+1}^{-1} Y_{1,2n-\ell-3}
     Y_{\ell-2k,2j-2k} & \text{otherwise}.
  \end{cases}
\end{split}
\end{equation*}
\end{enumerate}
All $M_{j,k}$ satisfy $\te_i M_{j,k} = 0$ for all $i\in I_0$. The
monomials appearing in $\mathcal M_{I_0}(M_{j,k})\cong \mathcal
B_{I_0}(\varpi_{\ell-2k})$ can be described as in the previous
subsection. Indeed for 
\begin{NB}
$0\le j < \ell$, $0\le k < \ell/2$,
\end{NB}
$k\neq \ell/2$ let us define a monomial $m_{T;j,k}$ associated with a
tableau $T = ((i_1,\dots,i_{j-2k}),(i_{j-2k+1},\cdots,i_{\ell-2k})\in
D_{\ell - 2k,j - 2k,n - \ell -2}$ by
\begin{NB}
For all $0\le k\leq\lfloor \ell/2\rfloor$, we have $M_{\ell,k} =
\tau_{2n-4}(M_{0,k})$ . For $0\le j \leq \ell$, $0\le k\leq\lfloor \ell/2\rfloor-1$, we have 
\begin{equation*}M_{j,k+1}=\tilde{e}_{\ell-2k-2}\tilde{e}_{\ell-2k-3}\cdots\tilde{e}_{1}\tilde{e}_{\ell-2k-1}\tilde{e}_{\ell-2k-2}\cdots\tilde{e}_2\tilde{e}_0 M_{j,k}.\end{equation*}
And so for $0\le j \leq \ell$, $0\le k\leq\lfloor \ell/2\rfloor$ we have $M_{j,k}\in\mathcal M(m)$ and for $0\le j < \ell$, $0\le k\leq\lfloor \ell/2\rfloor$ we have $z_\ell^{-1}(M_{j,k})=M_{j+1,k}$. For $0\le j < \ell$, $0\le k < \ell/2$, let us define a monomial $m_{T;j,k}$ associated with a tableau $T = ((i_1,\dots,i_{j-2k}),(i_{j-2k+1},\cdots,i_{\ell-2k})\in D_{\ell - 2k,j - 2k,n - \ell -2}$ by
\end{NB}
\begin{enumerate}
    \item $k < \lfloor j/2\rfloor$: 
                \begin{equation*}
                    m_{T;j,k} =  N_{j,k} \prod_{a=1}^{j-2k} \ffbox{i_a}_{2n-\ell-4k-2a+2j-3} \times \prod_{a=j-2k+1}^{\ell-2k}                          \ffbox{i_a}_{\ell-2(a-j+2k)+1},
                \end{equation*}

    \item $j$ is odd and $k = (j-1)/2$:
                \begin{equation*}
                m_{T;j,(j-1)/2} = N_{j,k} \times \ffbox{i_1}_{2n-\ell-3}
                \prod_{a=2}^{\ell-j+1} \ffbox{i_a}_{\ell-2a+3},
                \end{equation*}

    \item $j$ is even and $k \ge j/2$:
                \begin{equation*}
                m_{T;j,k} = N_{j,k} \times \prod_{a=1}^{\ell-2k} \ffbox{i_a}_{\ell-2a-4k+2j+1},
                \end{equation*}

    \item $j$ is odd and $k\ge (j+1)/2$:
                \begin{equation*}
                m_{T;j,k} = N_{j,k}\times
                \prod_{a=1}^{\ell-2k} \ffbox{i_a}_{\ell-2a-4k+2j+1}.
                \end{equation*}
\end{enumerate}
(For the case (4) the $Y_{1,\ell+1}^{-1} Y_{1,2n-\ell-3}$ does not
change anything because all other $Y_{1,r}^{\pm}$ satisfy $r < \ell+1$.) 
For $k = \ell/2$ we set $D_{0,j-\ell,n-\ell-2} = \{ \emptyset\}$ and
define $m_{\emptyset;j,k}$ by the same formula as in (3),(4) where the
last product is understood as $1$.
If $k > j/2$, we set $D_{\ell-2k,j-2k,n-\ell-2} = D_{\ell-2k,0,0}$,
i.e., the set of tableaux whose suffixes do not jump.

As $M_{j,0}\in \mathcal M(M_{0,0})$, it becomes clear by induction on $k$ that all $M_{j,k}$ are in
$\mathcal M(M_{0,0})$, and so by using Theorem~\ref{isot} all $m_{T;j,k}$ are in
$\mathcal M(M_{0,0})$ (the argument is similar to one for the type $A^{(1)}_n$). 

As $\wt(M_{j,0}) = \varpi_\ell - j\delta$, we have $z_\ell^{-1}
(M_{j,0}) = M_{j+1,0}$ by the reason explained in the beginning of
this section. Then from the definition of $M_{j,k}$ we have
$z_\ell^{-1}(M_{j,k}) = M_{j+1,k}$. Let us consider
$\tau_{\ell-2k,j-2k,n-\ell-2}$, where we understand it as the identity
map when $k > j/2$. It maps $M_{j,k}$ to $M_{j+1,k}$ and respects the
$I_0$-crystal structure. Since such a map is unique, we have
$z_\ell^{-1} = \tau_{\ell-2k,j-2k,n-\ell-2}\colon \mathcal M_{I_0}(M_{j,k})\to 
\mathcal M_{I_0}(M_{j+1,k})$.

\begin{NB}
We extend the definition of $m_{T;j,k}$ so that $m_{T;\ell,k} =
\tau_{2n-4}(m_{T;0,k})$ for all $0\le k\leq\lfloor \ell/2\rfloor$. Note that $m_{(1,2,\cdots,\ell);j,k}=M_{j,k}$. So it follows from Section \ref{preld} that the set of monomials $m_{T;j,k}$ defined above is $\mathcal M_{I_0}(M_{j,k})\simeq\mathcal B_{I_0}(\omega_{\ell-2k})$ (for the case (4) the $Y_{1,\ell+1}^{-1} Y_{1,2n-\ell-3}$ does not change anything because all other $Y_{1,r}^{\pm}$ appearing in the boxes satisfy $r < \ell+1$). Besides if $k=\ell/2$, $\mathcal M_{I_0}(M_{j,k})=\{M_{j,k}\}$.

So all monomials of $\mathcal M(M_{0,0})$ are connected to either $M_{j,k}$ ($0\le j < \ell$, $0\le k\leq\lfloor \ell/2\rfloor$) or their $\tau_{2n-4}$ images in the $I_0$-crystal. Thus
\begin{equation*}
    \mathcal M(M_{0,0})/\tau_{2n-4} = \bigsqcup_{0\le j < \ell, 0\le k\leq\lfloor \ell/2\rfloor} \mathcal
    M_{I_0}(M_{j,k}).
\end{equation*}
Moreover for $0\leq j\leq \ell-1$, $0\leq k<\lfloor \ell/2\rfloor$ we have
\begin{equation*}
    (z_\ell)^{-1}(m_{T;j,k})=m_{\tau_{\ell-2k,j-2k,n-\ell-2}(T);j+1,k}.
\end{equation*} 
So we get
\end{NB}

We can describe Kashiwara operators $\tilde{e}_i$, $\tilde{f}_i$ in
terms of tableaux as in type $A^{(1)}_n$. For $i\neq 0$, it is
basically explained in the proof of
Lemma~\ref{lem:tableauxcrystal}. So let us consider the case
$\tilde{e}_0$, $\tilde{f}_0$. 
We get that $\tilde{e}_0(m_{T;j,k})$ is equal to
\begin{equation*}
\begin{split} 
        \begin{cases}
        m_{(i_3,\cdots,i_{\ell-2k});j,k+1}&\text{ if $i_2=2$ and $i_{\ell-2k-1}\nsucceq \overline{2}$,}
                \\
                m_{(i_1,\cdots,i_{\ell-2k},\overline{2},\overline{1});j-2,k-1}&\text{
                if $i_2\npreceq 2$, $i_{\ell-2k}\nsucceq \overline{2}$
                and $k > 0$},
                \\ m_{(i_2,\cdots,i_{\ell},\overline{3-i_1});j-1,0}&\text{ if $i_1\preceq 2$, $i_2\npreceq 2$, $i_{\ell}\nsucceq \overline{2}$ and $k=0$,}
                \\0 &\text{ otherwise,}
  \end{cases} 
\end{split}
\end{equation*}
and that $\tilde{f}_0(m_{T;j,k})$ is equal to
\begin{equation*}
\begin{split} 
        \begin{cases}
        m_{(1,2,i_1,\cdots,i_{\ell-2k});j,k - 1}&\text{ if $i_1 \npreceq 2$, $i_{\ell-2k-1}\nsucceq \overline{2}$ and $k > 0$},
                \\ m_{(i_1,\cdots,i_{\ell-2k-2});j+2,k+1}&\text{ if $i_{\ell-2k-1} = \overline{2}$ and $i_2\npreceq 2$},
                \\ m_{(3-\overline{i_{\ell}},i_1,\cdots,i_{\ell-1});j+1,0}&\text{ if $i_1 \npreceq 2$, $i_{\ell-1}\nsucceq \overline{2}$, $i_{\ell}\succeq \overline{2}$ and $k = 0$,}
                \\0 &\text{ otherwise,}
  \end{cases} 
\end{split}
\end{equation*}
where we denote $\overline{\overline{1}}=1$ and
$\overline{\overline{2}}=2$, 
and we extend the definition of $m_{T;j,k}$ from $0\le j\le \ell-1$ to
all $j\in\mathbb Z$ so that $m_{T;j+\ell} = \tau_{2n-4} m_{T;j}$.
We understood that the condition
`$i_2 = 2$' is not satisfied when $k = \lfloor \ell/2\rfloor$ (and
hence there is no entry $i_2$). Similarly `$i_2\npreceq 2$' is
satisfied when $k = \lfloor \ell/2\rfloor$. The same rules apply to
other conditions. And we will use the same conventions for other
classical types.


As we have checked the stability for operators $\tilde{e}_0$, $\tilde{f}_0$, all the monomials
appearing in $\mathcal M(M_{0,0})$ are the $\tau_{2n-4}$-images of the $m_{T;j,k}$.

In particular, we can describe the $I_0$-crystal structure of
$\mathcal B(W(\varpi_\ell))$ as
\begin{equation*}
\begin{split}
    \mathcal B(W(\varpi_\ell)) & \simeq\mathcal M(M_{0,0})/z_\ell
    = \mathcal M_{I_0}(M_{0,0})\sqcup \mathcal M_{I_0}(M_{0,1})\sqcup \cdots\sqcup 
    \mathcal M_{I_0}(M_{0,\lfloor \ell/2\rfloor})
\\
   & \simeq
   \mathcal B_{I_0}(\varpi_\ell) \sqcup
   \mathcal B_{I_0}(\varpi_{\ell-2}) \sqcup \cdots \sqcup
   \begin{cases}
   \mathcal B_{I_0}(\varpi_1) & \text{if $\ell$ is odd},
   \\
   \mathcal B_{I_0}(0) & \text{if $\ell$ is even}.
   \end{cases}
\end{split}
\end{equation*}
In fact, this last result is well-known.

As an application of the description of what we just obtained, we construct
an explicit bijection between two sets of monomials, one is $\mathcal
M(M_{0,0})/z_\ell$, the other is those appearing in the
$q$--characters of $W(\varpi_\ell)$ counted with multiplicities.
Recall the conditions (1),(2) in Theorem~\ref{find}.
In \cite{Nac} we proved that the $q$--character of $W(\varpi_\ell)$ is
given by the sum of monomials corresponding to $T =
(i_1,\dots,i_\ell)$ satisfying (1) alone. We then defined $l(T)$ as
the number of pairs as in (2). Now we define the bijection
\begin{multline*}
  \{ T = (i_1,\dots,i_\ell) \mid \text{$T$ satisfies (1)}, l(T) = d \}
\\
  \longleftrightarrow
  \{ T' = (i_1,\dots,i_{\ell-2d}) \mid \text{$T'$ satisfies (1),(2)} \}
\end{multline*}
by letting $T'$ be the tableaux obtained by removing all the pairs
violating (2) in $T$.

This bijection cannot be expressed in terms of monomials in a simple
way unlike type $A$ case.

As another application, we get a description of the crystal $\mathcal
B(W(\varpi_\ell))$ in terms of tableaux. Namely we identify it with
$\{ m_{T;0,k} \mid 0\le k\le \lfloor \ell/2\rfloor\}$. Then we express
$\te_0 m_{T;0,k}$, $\tf_0 m_{T;0,k}$ as $m_{T';0,k'}$, $m_{T'';0,k''}$
by the above formula composed with the crystal automorphism
$\tau_{\ell,h,r}$ for suitable $h$, $r$.  
This description is similar to one in
\cite{koga,schilling}, probably the same if we use the isomorphism
between our $D_{\ell,0,0}$ and Kashiwara-Nakashima's tableaux
\cite{KN} in \cite{kks}. Note that the uniqueness of the crystal base
of $W(\varpi_\ell)$ was proved in \cite{koga}.

\subsubsection{Spin representations}\label{Dn3} Finally we consider the case
$\ell = n-1$ or $n$. Following \cite{Nac,kks} we define the half size
numbered box as
\newcommand{\fhbox}[1]{
\setbox9=\hbox{$\scriptstyle\overline{n\!\!-\!\!1}$}
\framebox[\yh][c]{\rule{0mm}{\ht9}${\scriptstyle #1}$}
}
\begin{equation*}
\begin{split}
   \fhbox{i}_p & = 
   \begin{cases}
      Y_{1,p-1} & \text{if $i=1$},
   \\
      Y_{1,p+1}^{-1} Y_{2,p} Y_{0,p+1}^{-1} & \text{if $i=2$},
   \\
     Y_{i-1,p+{i-1}}^{-1} Y_{i,p+{i-2}} 
       & \text{if $3\le i\le n-2$},
   \\
     Y_{n-2,p+{n-2}}^{-1} & \text{if $i=n-1$},
   \\
     Y_{n,p+{n-1}} & \text{if $i=n$},
   \end{cases}
\\
   \fhbox{\overline{i}}_p &= 
   \begin{cases}
   Y_{0,p+2n-1} & \text{if $i=1$},
   \\
     1 & \text{if $2\le i\le n-2$},
   \\
   Y_{n-1,p+{n+1}}^{-1}Y_{n,p+{n+1}}^{-1} 
     & \text{if $i=n-1$},
   \\
   Y_{n-1,p+{n-1}} & \text{if $i=n$}.
   \end{cases}
\end{split}
\end{equation*}
Let $M = Y_{\ell,0} Y_{0,n-2}^{-1}
=
\prod_{a=1}^{n-1} \fhbox{a}_{n+1-2a} \times
\fhbox{n}_{1-n}$ ($\ell=n$) or
$ 
\prod_{a=1}^{n-1} \fhbox{a}_{n+1-2a} \times
\fhbox{\overline{n}}_{1-n}$ ($\ell = n-1$).
We have $\mathcal M(M) \simeq \mathcal B(\varpi_\ell)$ by
Corollary~\ref{corex}.

Let $\mathcal B_{\operatorname{sp}}^+$ (resp.\ $\mathcal
B_{\operatorname{sp}}^-)$ be the set of tableaux $T = (i_1,\dots,i_n)$
satisfying
\begin{enumerate}
    \item $i_a\in \mathbf B, \; i_1 \prec i_2 \prec \dots \prec i_n$,
    \item $i$ and $\overline{i}$ do not appear simultaneously,
    \item if $i_a = n$, $n-a$ is even (resp.\ odd),
    \item if $i_a = \overline{n}$, $n-a$ is odd (resp.\ even).
\end{enumerate}
We define $m_T$ by
\begin{equation*}
  m_T = \prod_{a=1}^n \fhbox{i_a}_{n+1-2a}.
\end{equation*}
Then $\mathcal B_{I_0}(M)$ is $\{ m_T \mid T\in 
\mathcal B_{\operatorname{sp}}^\pm\}$, where $\pm$ is $-$ if $\ell =
n-1$ and $+$ if $\ell = n$. Let
$T=(3,4,5,\dots,n-1,n,\overline{2},\overline{1})$ for $\ell = n$ or
$T=(3,4,5,\dots,n-1,\overline{n},\overline{2},\overline{1})$ for $\ell
= n-1$. Then
\(
   m_T = Y_{2,n+1}^{-1} Y_{\ell,4} Y_{0,n}.
\)
Applying $\tilde f_0$ to $m_T$, we get $Y_{\ell,4} Y_{0,n+2}^{-1} =
\tau_4(M)$. As this has weight $\wt(M) - \delta$, it follows that
$\tau_4(M) = z_\ell^{-1}(M)$ as before. As a consequence, we have $z_\ell =
\tau_{-4}$ and
\(
  \mathcal B(W(\varpi_\ell)) \simeq
   \mathcal M(M)/\tau_4.
\)

We describe the action of $\tilde{e}_0$, $\tilde{f}_0$. We have
\begin{equation*}
\begin{split}
\tilde{e}_0(m_{T}) &=
        \begin{cases}
         \tau_{-4}(m_{(i_3,\cdots,i_n,\overline{2},\overline{1})}) &\text{ if $i_2 = 2$},
                \\  0&\text{ otherwise,}
  \end{cases}
\\
\tilde{f}_0(m_{T}) &=
        \begin{cases}
        \tau_4(m_{(1,2,i_1,\cdots,i_{n-2})}) &\text{ if $i_{n-1}=\overline{2}$,}
                \\  0&\text{ otherwise.}
  \end{cases} 
\end{split}
\end{equation*}
So the above are all the monomials in
$\mathcal M(M)/\tau_4$. So we recover a well-known result
\(
  \mathcal B(W(\varpi_\ell)) \simeq
   \mathcal B_{I_0}(\varpi_\ell).
\)
The map $Y_{0,*} \mapsto 1$ gives a bijection between $\{ m_T \mid
T\in \mathcal B_{\operatorname{sp}}^\pm\}$ and the monomials appearing
in $q$--characters of $W(\varpi_\ell)$, where all the multiplicities
are $1$ in these cases.

\subsection{Type $B_n^{(1)}$}\label{bn}
We can describe monomial crystals of other classical types by a
similar method. We just state the result without proofs.

Let $\mathbf B = \{ 1,\dots,n,0,\overline{n},\dots, \overline{1}\}$. 
We give the ordering $\prec$ on the set $\mathbf B$
by
\begin{equation*}
  1 \prec 2 \prec \cdots \prec n \prec 0 \prec \overline{n} \prec\cdots \prec \overline{2}\prec \overline{1}.
\end{equation*}

For $p\in \ZZ$, we define
{\allowdisplaybreaks
\begin{equation*}
\begin{aligned}[c]
    & \ffbox{1}_p = Y_{0,p+2}^{-1} Y_{1,p}, \quad
     \ffbox{2}_p = Y_{0,p+2}^{-1} Y_{1,p+2}^{-1} Y_{2,p+1},
    \\
    & \ffbox{i}_p = Y_{i-1,p+i}^{-1} Y_{i,p+{i-1}} \qquad(3 \le i \le n-1),
    \\
    & \ffbox{n}_p
     = Y_{n-1,p+n}^{-1} Y_{n,p+{n-1}}^2,
    \\
    & \ffbox{0}_p = Y_{n,p+{n}+1}^{-1} Y_{n,p+{n-1}},
    \\
    & \ffbox{\overline{n}}_p = Y_{n-1,p+n} Y_{n,p+{n}+1}^{-2},
    \\
    & \ffbox{\overline{i}}_p = Y_{i-1,p+{2n-i}} Y_{i,p+{2n+1-i}}^{-1}
    \qquad(3 \le i \le n-1),
    \\
    & \ffbox{\overline{2}}_p = Y_{0,p+2n-2} Y_{1,p+2n-2} Y_{2,p+2n-1}^{-1},
    \qquad
    \ffbox{\overline{1}}_p = Y_{0,p+2n-2} Y_{1,p+{2n}}^{-1}.
\end{aligned}
\end{equation*}}

\subsubsection{}\label{Bn1} 
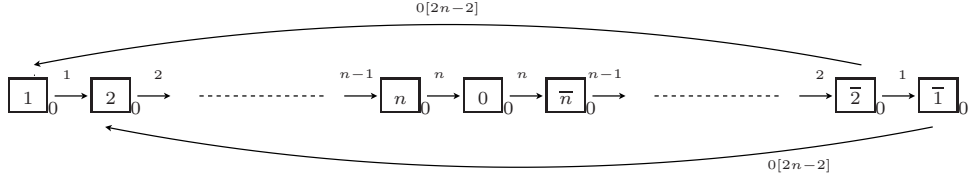
\begin{figure}[htbp]
\centering
\psset{xunit=.55mm,yunit=.55mm,runit=.55mm}
\psset{linewidth=0.3,dotsep=1,hatchwidth=0.3,hatchsep=1.5,shadowsize=1}
\psset{dotsize=0.7 2.5,dotscale=1 1,fillcolor=black}
\psset{arrowsize=1 2,arrowlength=1,arrowinset=0.25,tbarsize=0.7 5,bracketlength=0.15,rbracketlength=0.15}
\begin{pspicture}(0,0)(220,50)
\rput(0,20){$\ffbox{1}_0$}
\rput(20,20){$\ffbox{2}_0$}
\psline{->}(5,20)(13,20)
\psline{->}(25,20)(33,20)
\psline{->}(75,20)(83,20)
\psline[linestyle=dashed,dash=1 1](40,20)(70,20)
\rput(90,20){$\ffbox{n}_0$}
\rput(110,20){$\ffbox{0}_0$}
\rput(130,20){$\ffbox{\overline{n}}_0$}
\psline{->}(135,20)(143,20)
\psline{->}(185,20)(193,20)
\psline[linestyle=dashed,dash=1 1](150,20)(180,20)
\rput(200,20){$\ffbox{\overline{2}}_0$}
\psline{->}(205,20)(213,20)
\rput(220,20){$\ffbox{\overline{1}}_0$}
\psline{->}(95,20)(103,20)
\psline{->}(115,20)(123,20)
\rput{0}(0,25){\qline(-0,0)(0,0)}\rput{0}(100,-480){\parametricplot[arrows=->]{78.85}{101.15}{ t cos 517.26 mul t sin 517.26 mul }}
\rput{0}(117.26,520){\parametricplot[arrows=<-]{-101.15}{-78.85}{ t cos 517.26 mul t sin 517.26 mul }}
\rput(8,25){$\scriptscriptstyle 1$}
\rput(30,25){$\scriptscriptstyle 2$}
\rput(98,25){$\scriptscriptstyle n$}
\rput(138,25){$\scriptscriptstyle n-1$}
\rput(118,25){$\scriptscriptstyle n$}
\rput(78,25){$\scriptscriptstyle n-1$}
\rput(190,25){$\scriptscriptstyle 2$}
\rput(210,25){$\scriptscriptstyle 1$}
\rput(100,41){$\scriptscriptstyle 0[2n-2]$}
\rput(185,3){$\scriptscriptstyle 0[2n-2]$}
\end{pspicture}
\caption{(Type $B_n^{(1)}$) the crystal $\mathcal B(\varpi_1)$}
\label{fig:Bn}
\end{figure}

First consider the case $\ell=1$.
Let $M = Y_{1,0}Y_{0,2}^{-1}$. It follows from Corollary~\ref{corex}
that $\mathcal M(M) \simeq \mathcal B(\varpi_\ell)$. The crystal
graph of $\mathcal M(M)$ is given in Figure~\ref{fig:Bn}.
We find $\tau_{2n-2} = z_{\ell}^{-1}$ and $\mathcal
M(M)/\tau_{2n-2}=\mathcal M_{I_0}(M)$.  

\subsubsection{Preliminary results for crystals of finite type $B$}\label{secfinb}

\begin{NB}
As for type $D^{(1)}_n$, we first need to
describe the $I_0$-crystal structure on the monomials. This will be
given in this subsection. All the results on the $I_0$-crystal are
independent of the information on $Y_{0,*}$,  so we set $Y_{0,*}$ as
$1$ in this subsection. 
Note also that results can be modified in an obvious manner so that
the suffixes of $\ffbox{}_*$ can be shifted simultaneously. We will
use the results in these modified forms in later subsections.

The following result was proved in \cite{kks}.

\begin{thm}\label{debb} Let $1\leq \ell\leq n-1$. Consider the monomial 
\begin{equation*}M=\ffbox{1}_{\ell-1}\ffbox{2}_{\ell-3}\cdots\ffbox{\ell}_{-(\ell-1)}.\end{equation*}
Then $\mathcal M_{I_0}(M)\simeq\mathcal B_{I_0}(\varpi_\ell)$ and is equal to the set of monomials
\begin{equation*}m_T=\ffbox{i_1}_{\ell-1}\ffbox{i_2}_{\ell-3}\cdots\ffbox{i_\ell}_{-(\ell-1)},\end{equation*}
indexed by the set of tableaux $T = (i_1,\dots,i_{\ell})$ satisfying the conditions :

    \textup{(1)} $i_a\in\mathbf B,  i_1 \prec i_2 \prec \cdots \prec i_\ell$, but $0$ can be repeated,

    \textup{(2)} there is no pair $a$, $b$ such that $1\leq a < b \leq \ell$ and $i_a=k$, $i_b=\overline{k}$
                and $b-a = n - k$.
\end{thm}


\end{NB}

Let $1\leq \ell\leq n-1$, $0\leq r\leq n-\ell$ and $0\leq h\leq \ell$. Consider the monomial
\begin{equation*}
  \begin{split}
  & M_{\ell,h,r} = Y_{h,\ell-h} Y_{h,\ell-h-2r}^{-1} Y_{\ell,-2r}
\\
   =\; &
   \left(\ffbox{1}_{\ell-1} \ffbox{2}_{\ell-3} \cdots
     \ffbox{h}_{\ell-2h+1}
    \right)
  \times
   \left(\ffbox{h\!+\!1}_{\ell-2h-2r-1}
     \ffbox{h\!+\!2}_{\ell-2h-2r-3}
     \cdots \ffbox{\ell}_{1-\ell-2r}\right)
\\
   = \; &
   \prod_{p=1}^h \ffbox{p}_{\ell-2p+1} \times
   \prod_{p=h+1}^\ell \ffbox{p}_{\ell+1-2p-2r}.
  \end{split}
\end{equation*}
For $T=((i_1,\cdots,i_h),(i_{h+1},\cdots,i_\ell))$ such that $i_p\in\mathbf B$, we define the monomial
\begin{equation*}
m_{T} = \ffbox{i_1}_{\ell-1}\ffbox{i_2}_{\ell-3}\cdots\ffbox{i_h}_{\ell-2h+1}
      \ffbox{i_{h\!+\!1}}_{\ell-1-2h-2r}\ffbox{i_{h\!+\!2}}_{\ell-2h-3-2r}\cdots\ffbox{i_\ell}_{-\ell+1-2r}.
\end{equation*}

Let $B_{\ell,h,r}$ be the set of tableaux $T$ satisfying the following
conditions 
\begin{enumerate}
\def\labelenumi{(B.\theenumi)}
\def\theenumi{\arabic{enumi}}
%
\item  $i_a\in\mathbf B$, $i_1 \prec i_2 \prec \cdots \prec i_h$ but $0$ can be repeated, and $i_{h+1}\prec i_{h+2}\prec \cdots \prec i_{\ell}$ but $0$ can be repeated.

\item There is no pair $a$, $b$ such that $1\leq a < b \leq h$ and
  $i_a=k$, $i_b=\overline{k}$ and $b-a = n - k$.

\item There is no pair $a$, $b$ such that $h+1\leq a < b \leq \ell$
  and $i_a=k$, $i_b=\overline{k}$ and $b-a = n - k$.

\item There is no pair $a$, $b$ such that $a \leq h$ , $h+1 \leq b$ , $i_a=k$, $i_b=\overline{k}$ and $b-a = n + 1 - \max(r,1) - k$.

\item  Suppose that $i_{h+1} = k \in \{1,\dots, n\}$ and $i_h\succeq
i_{h+1}$. Then $i_h = k'$ is also in $\{ 1,\dots, n\}$,
and the successive part $(\overline{k'}, \overline{k'-1},\dots,
\overline{k})$ appears as $(i_{b'},i_{b'+1},\dots,i_{b})$ with
$n - r - k + 1 < b - h \le n- k$.

\item Suppose that $i_{h+1} = \overline{k} \in \{\overline{1},\dots,
  \overline{n}\}$ and $i_h\succeq i_{h+1}$. Then $i_h = \overline{k'}$ is also in $\{ \overline{1},\dots, \overline{n}\}$,
  and the successive part $(k',k'+1,\dots,k)$ appears as
  $(i_{a'},i_{a'+1},\dots,i_{a})$ with
  $n-r-k+1 \le h - a < n-k$.

\item  If $i_{h+1} = 0$, then $i_h \preceq 0$.
\end{enumerate}

Note that the conditions above are the same as the ones in \cite{kks}
when $r=0$.

For $T = ((i_1,\dots,i_h),(i_{h+1},\dots,i_\ell))\in B_{\ell,h,r}$ we define the tableau $\tau_{\ell,h,r}(T)$ in the following three cases separately. 
\begin{enumerate}
\def\labelenumi{(B.\theenumi)}
\def\theenumi{\alph{enumi}}
\item $i_{h+1} = k \in \{1,\dots,n\}$ and there is an entry $i_b = \overline{k}$ with $n-r-k+1< b - h \le n-k$.

\item $i_{h+1} = \overline{k} \in \{\overline{1},\dots,\overline{n}
\}$ and there is an entry $i_a = k$ with $n-r-k+1 \le h - a < n - k$.

\item Neither (B.a) nor (B.b) is not satisfied. 
\end{enumerate}

In the case (B.a), let $b''$ such that 
$(i_{b''},i_{b''+1},\dots,i_b)$ are successive as
$(\overline{k''},\overline{k''+1},\dots,\overline{k})$ and
$i_{b''-1}\neq \overline{k''-1}$. We have $k'' < n-1$. We set
\begin{equation*}
  \begin{split}
    \tau_{\ell,h,r}(T)=((i_1,\cdots&,i_h,k''\!+\!1),
\\
  & (i_{h+2},\cdots,i_{b''-1}, \overline{k''\!+\!1},\overline{k''},
  \cdots,\overline{k\!+\!1}, i_{b+1},\cdots, i_\ell)). 
  \end{split}
\end{equation*}

Similarly in the case (B.b), we take $i_{a''}$ so that
$(i_{a''},i_{a''+1},\dots, i_a) = (k'',k''+1,\dots,k)$ and $i_{a''-1}\neq
k''-1$. We have $k < n-1$. We then set
\begin{equation*}
  \begin{split}
    \tau_{\ell,h,r}(T)=((i_1,\cdots, i_{a''-1},k''\!-\!1,\cdots,k\!-\!1,
    i_{a+1},&\cdots,i_{h}, \\
    & \overline{k''\!-\!1}), (i_{h+2},\cdots, i_\ell)).
  \end{split}
\end{equation*}

In the case (B.c) we set
\begin{equation*}
\tau_{\ell,h,r}(T)=((i_1,\cdots, i_{h+1}),(i_{h+2},\cdots, i_\ell)).
\end{equation*}

\begin{thm}
\textup{(1)} The map $T \mapsto m_T$ induces a crystal isomorphism between
$B_{\ell,h,r}$ and $\mathcal M_{I_0}(M_{\ell,h,r})$.

\textup{(2)} $\tau_{\ell,h,r}$ induces a crystal isomorphism
$\mathcal M_{I_0}(M_{\ell,h,r})$ to $\mathcal M_{I_0}(M_{\ell,h+1,r})$.
\end{thm}

\subsubsection{} Now we study $\mathcal{B}(\varpi_\ell)$ for $2\leq \ell\leq n-1$. Let
\(
  M_{0,0} = Y_{\ell,0} Y_{0,\ell-1}^{-1} Y_{0,\ell+1}^{-1}
    = \ffbox{1}_{\ell-1}\ffbox{2}_{\ell-3} \cdots
      \ffbox{\ell}_{1-\ell}.
\)
For $\ell\neq n-1$ we have $\tilde{f}_2\tilde{f}_3\cdots\tilde{f}_\ell
M_{0,0}=Y_{0,\ell+1}^{-1}Y_{1,\ell-1}Y_{2,\ell}^{-1}Y_{\ell+1,1}$ and
for $\ell=n-1$ we have $\tilde{f}_2\tilde{f}_3\cdots\tilde{f}_\ell
M_{0,0}=Y_{0,\ell+1}^{-1}Y_{1,\ell-1}Y_{2,\ell}^{-1}Y_{n,1}^2$. By
a method similar to the proof of Proposition~\ref{otherd} we have
  $\mathcal M(M_{0,0}) \simeq \mathcal B(\varpi_\ell)$.

\begin{NB}
For $0\le j\le\ell$, we set
\begin{equation*}
\begin{split}
  M_{j,0} &= \left( \ffbox{1}_{2n-\ell+2j-3}\ffbox{2}_{2n-\ell+2j-5}
   \cdots \ffbox{j}_{2n-\ell-1}\right)
\times
  \left(\ffbox{j\!+\!1}_{\ell-1}\ffbox{j\!+\!2}_{\ell-3}
    \cdots \ffbox{\ell}_{1-\ell+2j}\right)
\\
 &= \begin{cases}
   Y_{\ell,0} Y_{0,\ell-1}^{-1} Y_{0,\ell+1}^{-1} & \text{if $j=0$},
\\
   Y_{\ell,2} Y_{0,\ell+1}^{-1} Y_{0,2n-\ell+1}^{-1}
   Y_{1,\ell+1}^{-1} Y_{1,2n-\ell-1} &\text{if $j=1$},
\\
   Y_{\ell,2j} Y_{0,2n-\ell+2j-3}^{-1}Y_{0,2n-\ell+2j-1}^{-1}
       Y_{j,\ell+j}^{-1}Y_{j,2n-\ell+j-2}
   &\text{otherwise.}
 \end{cases}
 \end{split}
\end{equation*}

\begin{Claim}
We have $M_{j,0}\in \mathcal M(M_{0,0})$ for $0\leq j \le \ell$. 
\end{Claim}
It is clear by induction on $j$ because for $0\leq j < \ell$, we have
\begin{equation*}
M_{j+1,0}=\tilde{e}_\ell\tilde{e}_{\ell-1}\cdots\tilde{e}_2\tilde{f}_0\tilde{f}_2^2\tilde{f}_3^2\cdots\tilde{f}_\ell^2\tilde{f}_{\ell+1}\tilde{f}_{\ell+2}\cdots\tilde{f}_{n}\tilde{f}_n\tilde{f}_{n-1}\cdots\tilde{f}_{\ell+1}\tilde{f}_1\tilde{f}_2\cdots\tilde{f}_\ell M_{j,0}.
\end{equation*}
For weight reason, for $0\leq j\leq \ell-1$ we have $z_\ell^{-1}(M_{j,0})=M_{j+1,0}$. As $M_{\ell,0} =\tau_{2n-2}(M_{0,0})$, $\mathcal M(M_{0,0})$ is preserved under $\tau_{2n-2}$ and $(z_\ell)^{-\ell}=\tau_{2n-2}$.

For $0\le j < \ell$, $0\le k\le
\lfloor \ell/2\rfloor$ let us define the monomial $M_{j,k}$ by
\begin{enumerate}
    \item $k < \lfloor j/2\rfloor$: 
                \begin{equation*}
                M_{j,k}=
                                Y_{\ell-2k,2j-2k}
                Y_{0,2n-\ell-4k+2j-3}^{-1} Y_{0,2n-\ell+2j-1}^{-1}
                Y_{j-2k,j+\ell-2k}^{-1} Y_{j-2k,2n-\ell+j-2k-2},
                \end{equation*}

    \item $j$ is odd and $k = (j-1)/2$:
                \begin{equation*}
                M_{j,(j-1)/2}= 
                                Y_{0,\ell+1}^{-1} Y_{0,2n-\ell+2j-1}^{-1} Y_{1,\ell+1}^{-1}
                Y_{1,2n-\ell-1} Y_{\ell-j+1,j+1},
                \end{equation*}

    \item $j$ is even and $k \ge j/2$:
                \begin{equation*}
                \begin{split}
                M_{j,k} 
                &= 
                \begin{cases}
                Y_{0,-\ell+2j+1} Y_{0,\ell+1}^{-1}
                    Y_{0,2n-\ell-1}Y_{0,2n-\ell+2j-1}^{-1} & \text{if $k = \ell/2$},
                \\
                Y_{0,\ell+1}^{-1} Y_{0,2n-\ell-1}
                Y_{0,2n-\ell+2j-1}^{-1} Y_{1,-\ell+2j+1} & \text{if $k = (\ell-1)/2$},
                \\
                \begin{aligned}[c]
                & Y_{0,\ell-4k+2j-1}^{-1}Y_{0,\ell+1}^{-1}
            Y_{0,2n-\ell-1}Y_{0,2n-\ell+2j-1}^{-1}
                Y_{\ell-2k,2j-2k}
                \end{aligned}
                & \text{otherwise},
                \end{cases}
                \end{split}
                \end{equation*}

    \item $j$ is odd and $k\ge (j+1)/2$:
                \begin{equation*}
                \begin{split}
                M_{j,k} 
                &= 
                \begin{cases}
                Y_{0,-\ell+2j+1} Y_{0,2n-\ell+2j-1}^{-1}
                Y_{1,\ell+1}^{-1} Y_{1,2n-\ell-1} & \text{if $k=\ell/2$},
                \\
                Y_{0,2n-\ell+2j-1}^{-1}
                Y_{1,-\ell+2j+1} Y_{1,\ell+1}^{-1} Y_{1,2n-\ell-1}
                & \text{if $k=(\ell-1)/2$},
                \\
                Y_{0,\ell-4k+2j-1}^{-1} Y_{0,2n-\ell+2j-1}^{-1}
                Y_{1,\ell+1}^{-1} Y_{1,2n-\ell-1}
                Y_{\ell-2k,2j-2k} & \text{otherwise}.
                \end{cases}
                \end{split}
                \end{equation*}
\end{enumerate}
We extend the definition of $M_{j,k}$ so that $M_{\ell,k} = \tau_{2n-2}(M_{0,k})$ for all $0\le k\leq\lfloor \ell/2\rfloor$. For $0\le j \leq \ell$, $0\le k\leq\lfloor \ell/2\rfloor-1$, we have
\begin{equation*}
M_{j,k+1}=\tilde{e}_{\ell-2k-2}\tilde{e}_{\ell-2k-3}\cdots\tilde{e}_{1}\tilde{e}_{\ell-2k-1}\tilde{e}_{\ell-2k-2}\cdots\tilde{e}_2\tilde{e}_0 M_{j,k}.
\end{equation*}
So for $0\le j \leq \ell$, $0\le k\leq\lfloor \ell/2\rfloor$ we have
$M_{j,k}\in\mathcal M(M_{0,0})$ and for $0\le j < \ell$, $0\le
k\leq\lfloor \ell/2\rfloor$ we have $z_\ell^{-1}(M_{j,k})=M_{j+1,k}$.
\end{NB}

For $0\le j < \ell$, $0\le k < \ell/2$, let us define the monomial
$m_{T;j,k}$ associated with $T =
((i_1,\dots,i_{j-2k}),(i_{j-2k+1},\cdots,i_{\ell-2k}))\in
B_{\ell-2k,j-2k,n-\ell-1}$ by
\begin{enumerate}
    \item $k < \lfloor j/2\rfloor$: 
                \begin{equation*}
                \begin{split}
                    m_{T;j,k}
                = &Y_{0,2n-\ell-4k+2j-1} Y_{0,2n-\ell+2j-1}^{-1} 
                \prod_{a=1}^{j-2k} \ffbox{i_a}_{2n-\ell-4k-2a+2j-1}
                \\ & \qquad\qquad\times
                \prod_{a=j-2k+1}^{\ell-2k} \ffbox{i_a}_{\ell-2(a-j+2k)+1},
                \end{split}
                \end{equation*}

    \item $j$ is odd and $k = (j-1)/2$:
                \begin{equation*}
                m_{T;j,(j-1)/2} = Y_{0,2n-\ell-3} Y_{0,2n-\ell+2j-1}^{-1} \ffbox{i_1}_{2n-\ell-1}
                \prod_{a=2}^{\ell-j+1} \ffbox{i_a}_{\ell-2a+3},
                \end{equation*}

    \item $j$ is even and $k \ge j/2$:
                \begin{equation*}
                m_{T;j,k} = Y_{0,\ell-4k+2j+1} Y_{0,\ell+1}^{-1}
            Y_{0,2n-\ell-1}Y_{0,2n-\ell+2j-1}^{-1}\prod_{a=1}^{\ell-2k} \ffbox{i_a}_{\ell-2a-4k+2j+1},
                \end{equation*}

    \item $j$ is odd and $k\ge (j+1)/2$:
                \begin{equation*}
                m_{T;j,k} = Y_{0,\ell-4k+2j+1} Y_{0,2n-\ell+2j-1}^{-1}
                Y_{1,\ell+1}^{-1} Y_{1,2n-\ell-1}
                \prod_{a=1}^{\ell-2k} \ffbox{i_a}_{\ell-2a-4k+2j+1}.
                \end{equation*}
\end{enumerate}
For $\ell = k/2$ we set $B_{0,j-\ell,n-\ell-1} = \{ \emptyset\}$ and
define $m_{\emptyset;j,k}$ by the same formula as in (3),(4) where the
last product is understood as $1$.
We extend the definition of $m_{T;j,k}$ for all $j\in\mathbb Z$ so
that $m_{T;j+\ell,k} = \tau_{2n-2} m_{T;j,k}$. 

\begin{NB}
We extend the definition of $m_{T;j,k}$ so that $m_{T;\ell,k} = \tau_{2n-2}(m_{T;0,k})$ for all $0\le k\leq\lfloor \ell/2\rfloor$. Note that $m_{(1,2,\cdots,\ell);j,k}=M_{j,k}$. So it follows from Section~\ref{secfinb} that the set of monomials $m_{T;j,k}$ defined above is $\mathcal M_{I_0}(M_{j,k})\simeq \mathcal B_{I_0}(\omega_{\ell-2k})$ (for the case (4) the $Y_{1,\ell+1}^{-1} Y_{1,2n-\ell-3}$ does not change anything because all other $Y_{1,r}^{\pm}$ appearing in the boxes satisfy $r < \ell+1$). Besides if $k=\ell/2$, $\mathcal M_{I_0}(M_{j,k})=\{M_{j,k}\}$.
We extend the definition of $m_{T;j,k}$ for all $j\in\mathbb Z$ so
that $m_{T;j+\ell,k} = \tau_{2n} m_{T;j,k}$. 
\end{NB}

We describe the action of $\tilde{e}_0$, $\tilde{f}_0$. We get that $\tilde{e}_0(m_{T;j,k})$ is equal to
\begin{equation*}
\begin{split} 
        \begin{cases}
        m_{(i_3,\cdots,i_{\ell-2k});j,k+1}&\text{ if $i_2=2$ and $i_{\ell-2k-1}\nsucceq \overline{2}$,}
                \\ m_{(i_1,\cdots,i_{\ell-2k},\overline{2},\overline{1});j-2,k-1}&\text{ if $i_2\npreceq 2$, $i_{\ell-2k}\nsucceq \overline{2}$ and $k > 0$},
                \\ m_{(i_2,\cdots,i_{\ell},\overline{3-i_1});j-1,0}&\text{ if $i_1\preceq 2$, $i_2\npreceq 2$, $i_{\ell}\nsucceq \overline{2}$ and $k=0$,}
                \\0 &\text{ otherwise,}
  \end{cases} 
\end{split}
\end{equation*}
and that $\tilde{f}_0(m_{T;j,k})$ is equal to
\begin{equation*}
\begin{split} 
        \begin{cases}
        m_{(1,2,i_1,\cdots,i_{\ell-2k});j,k - 1}&\text{ if $i_1 \npreceq 2$, $i_{\ell-2k-1}\nsucceq \overline{2}$ and $k > 0$},
                \\ m_{(i_1,\cdots,i_{\ell-2k-2});j+2,k+1}&\text{ if $i_{\ell-2k-1} = \overline{2}$ and $i_2\npreceq 2$},
                \\ m_{(3-\overline{i_{\ell}},i_1,\cdots,i_{\ell-1});j+1,0}&\text{ if $i_1 \npreceq 2$, $i_{\ell-1}\nsucceq \overline{2}$, $i_{\ell}\succeq \overline{2}$ and $k = 0$,}
                \\0 &\text{ otherwise.}
  \end{cases} 
\end{split}
\end{equation*}
So all monomials of $\mathcal M(M_{0,0})$ are connected to either $M_{j,k}$ ($0\le j < \ell$, $0\le k\leq\lfloor \ell/2\rfloor$) or their $\tau_{2n-2}$ images in the $I_0$-crystal, thus
\begin{equation*}\mathcal M(M_{0,0})/\tau_{2n-2}=\bigsqcup_{0\le j < \ell, 0\le k\leq\lfloor \ell/2\rfloor} \mathcal M_{I_0}(M_{j,k}).\end{equation*}
Moreover for $0\leq j\leq \ell-1$, $0\leq k<\lfloor \ell/2\rfloor$ we have
\begin{equation*}
    (z_\ell)^{-1}(m_{T;j,k})=m_{\tau_{\ell-2k,j-2k,n-\ell-1}(T);j+1,k}.
\end{equation*}
We have $\tau_{2n-2} = (z_\ell)^{-\ell}$, and all
monomials in $\mathcal M(M_{0,0})/\tau_{2n-2}$ are written as $m_{T;j,k}$. The crystal
automorphism $z_\ell$ is given by $\tau_{\ell-2k,j-2k,n-\ell-1}^{-1}$.

So we get
\begin{equation*}
    \mathcal B(W(\varpi_\ell))
   \simeq
   \mathcal B_{I_0}(\varpi_\ell) \sqcup
   \mathcal B_{I_0}(\varpi_{\ell-2}) \sqcup \cdots \sqcup
   \begin{cases}
   \mathcal B_{I_0}(\varpi_1) & \text{if $\ell$ is odd},
   \\
   \mathcal B_{I_0}(0) & \text{if $\ell$ is even}.
   \end{cases}
\end{equation*}

Our crystal structure described here is probably the same as one in
\cite{koga} if we use the isomorphism between our $B_{\ell,0,0}$ and
Kashiwara-Nakashima's tableaux \cite{KN} in \cite{kks}.
Note that the uniqueness of the crystal base of $W(\varpi_\ell)$ was
proved in \cite{koga}.

\subsubsection{} Finally we consider the case $\ell=n$. Let $M = Y_{n,0}Y_{0,n-1}^{-1}$. It follows from Corollary~\ref{corex} that $\mathcal M(M) \simeq \mathcal B(\varpi_\ell)$. 

\begin{NB}
As $\tilde{e}_n\tilde{e}_{n-1}\tilde{e}_n\tilde{e}_{n-2}\tilde{e}_{n-3}\cdots\tilde{e}_1\tilde{e}_{n-1}\tilde{e}_{n-2}\cdots\tilde{e}_2\tilde{e}_0 M = \tau_{-4}(M)$, $\mathcal M (M)$ is preserved under $\tau_4$. Moreover for weight reason $z_\ell=\tau_{-4}$ and so $\mathcal{M}(M)/\tau_4\simeq \mathcal{B}(W(\varpi_\ell))$. 

Mimicking the definition in \cite{Nac}, we define
\end{NB}
Let
\begin{equation*}
\begin{split}
   \fhbox{i}_p & = 
   \begin{cases}
      Y_{1,p-1} & \text{if $i=1$},
   \\
      Y_{1,p+1}^{-1} Y_{2,p} Y_{0,p+1}^{-1} & \text{if $i=2$},
   \\
     Y_{i-1,p+{i-1}}^{-1} Y_{i,p+{i-2}} 
       & \text{if $3\le i\le n-1$},
   \\
     Y_{n-1,p+{n-1}}^{-1} & \text{if $i=n$},
   \\
     Y_{n,p+n} & \text{if $i=0$},
   \end{cases}
\\
   \fhbox{\overline{i}}_p &= 
   \begin{cases}
   Y_{0,p+2n+1} & \text{if $i=1$},
   \\
     1 & \text{if $2\le i\le n-1$},
   \\
   Y_{n,p+{n+2}}^{-2} 
     & \text{if $i=n$}.\end{cases}
\end{split}
\end{equation*}
Then the monomials appearing in $\mathcal M_{I_0}(M)$ are 
$m_T=\prod_{a=1}^{n+1} \fhbox{i_a}_{n+2-2a}$ associated with a tableau
$T = (i_1,\dots,i_{n+1})$ satisfying the conditions
\begin{enumerate}
    \item $i_a\in \mathbf B, \; i_1 \prec i_2 \prec \dots \prec i_{n+1}$,
    \item $i$ and $\overline{i}$ do not appear simultaneously.
\end{enumerate}
We have $z_\ell = \tau_{-4}$. 
We describe the action of $\tilde{e}_0$, $\tilde{f}_0$ : we have
\begin{equation*}
\begin{split}
\tilde{e}_0(m_{T}) &=
        \begin{cases}
         \tau_{-4}(m_{(i_3,\cdots,i_{n+1},\overline{2},\overline{1})}) &\text{ if $i_2 = 2$},
                \\  0&\text{ otherwise,}
  \end{cases}
\\
\tilde{f}_0(m_{T}) &=
        \begin{cases}
        \tau_4(m_{(1,2,i_1,\cdots,i_{n-1})}) &\text{ if $i_n=\overline{2}$,}
                \\  0&\text{ otherwise.}
  \end{cases} 
\end{split}
\end{equation*}
So all monomials in
$\mathcal M(M)/\tau_4$ are written as $m_T$.
As an application, we recover a known result
$\mathcal B(W(\varpi_\ell))=\mathcal B_{I_0}(\varpi_\ell)$. 

By the condition (2) there is always an entry $i_a = 0$. If we remove
this entry, we get the tableaux description in \cite{KN}.




\subsection{Type $C_n^{(1)}$}

Let $\mathbf B = \{ 1,\dots,n,\overline{n},\dots, \overline{1}\}$. 
We give the ordering $\prec$ on the set $\mathbf B$
by
\begin{equation*}
  1 \prec 2 \prec \cdots \prec n \prec \overline{n} \prec\cdots \prec \overline{2}\prec \overline{1}.
\end{equation*}

For $p\in \ZZ$, we define
{\allowdisplaybreaks
\begin{equation*}
\begin{aligned}[c]
    & \ffbox{i}_p = Y_{i-1,p+i}^{-1} Y_{i,p+{i-1}} \qquad(1 \le i \le n),
    \\
    & \ffbox{\overline{i}}_p = Y_{i-1,p+{2n-i}} Y_{i,p+{2n+1-i}}^{-1}
    \qquad(1 \le i \le n).
\end{aligned}
\end{equation*}}

\subsubsection{} First consider the case $\ell=1$. Let $M = Y_{0,1}^{-1}Y_{1,0}$. It follows from Corollary~\ref{corex} that $\mathcal M(M) \simeq \mathcal B(\varpi_\ell)$.
The crystal graph of $\mathcal M(M)$ is given in Figure~\ref{fig:Cn}.
\begin{figure}[htbp]
\centering
\psset{xunit=.55mm,yunit=.55mm,runit=.55mm}
\psset{linewidth=0.3,dotsep=1,hatchwidth=0.3,hatchsep=1.5,shadowsize=1}
\psset{dotsize=0.7 2.5,dotscale=1 1,fillcolor=black}
\psset{arrowsize=1 2,arrowlength=1,arrowinset=0.25,tbarsize=0.7 5,bracketlength=0.15,rbracketlength=0.15}
\begin{pspicture}(0,0)(220,30)
\rput(0,20){$\ffbox{1}_0$}
\rput(20,20){$\ffbox{2}_0$}
\psline{->}(5,20)(13,20)
\psline{->}(25,20)(33,20)
\psline{->}(85,20)(93,20)
\psline[linestyle=dashed,dash=1 1](40,20)(80,20)
\rput(100,20){$\ffbox{n}_0$}
\psline{->}(105,20)(113,20)
\rput(120,20){$\ffbox{\overline{n}}_0$}
\psline{->}(125,20)(133,20)
\psline{->}(185,20)(193,20)
\psline[linestyle=dashed,dash=1 1](140,20)(180,20)
\rput(200,20){$\ffbox{\overline{2}}_0$}
\psline{->}(205,20)(213,20)
\rput(220,20){$\ffbox{\overline{1}}_0$}
\rput{0}(110.26,520){\parametricplot[arrows=<-]{-102.15}{-77.85}{ t cos 520.26 mul t sin 520.26 mul }}
\rput(8,25){$\scriptscriptstyle 1$}
\rput(30,25){$\scriptscriptstyle 2$}
\rput(108,25){$\scriptscriptstyle n$}
\rput(129,25){$\scriptscriptstyle n-1$}
\rput(89,25){$\scriptscriptstyle n-1$}
\rput(190,25){$\scriptscriptstyle 2$}
\rput(210,25){$\scriptscriptstyle 1$}
\rput(110,3){$\scriptscriptstyle 0[2n]$}
\end{pspicture}
\caption{(Type $C_n^{(1)}$) the crystal $\mathcal B(\varpi_1)$}
\label{fig:Cn}
\end{figure}
We find $\tau_{2n} = z_{\ell}^{-1}$ and $\mathcal
M(M)/\tau_{2n}=\mathcal M_{I_0}(M)$.  

\subsubsection{Preliminary results for crystals of finite type $C$}\label{secfinc}

Let $1\leq \ell\leq n$, $0\leq r\leq n-\ell$ and $0\leq h\leq \ell$. Consider the monomial
\begin{NB}
$\ell = n$ is included now.  
\end{NB}
\begin{equation*}
  \begin{split}
  & M_{\ell,h,r} = Y_{h,\ell-h} Y_{h,\ell-h-2r}^{-1} Y_{\ell,-2r}
\\
   =\; &
   \left(\ffbox{1}_{\ell-1} \ffbox{2}_{\ell-3} \cdots
     \ffbox{h}_{\ell-2h+1}
    \right)
  \times
   \left(\ffbox{h\!+\!1}_{\ell-2h-2r-1}
     \ffbox{h\!+\!2}_{\ell-2h-2r-3}
     \cdots \ffbox{\ell}_{1-\ell-2r}\right)
\\
   = \; &
   \prod_{p=1}^h \ffbox{p}_{\ell-2p+1} \times
   \prod_{p=h+1}^\ell \ffbox{p}_{\ell+1-2p-2r}.
  \end{split}
\end{equation*}
For $T=((i_1,\cdots,i_h),(i_{h+1},\cdots,i_\ell))$ such that $i_p\in\mathbf B$, we define the monomial
\begin{equation*}
m_{T} = \ffbox{i_1}_{\ell-1}\ffbox{i_2}_{\ell-3}\cdots\ffbox{i_h}_{\ell-2h+1}
      \ffbox{i_{h\!+\!1}}_{\ell-1-2h-2r}\ffbox{i_{h\!+\!2}}_{\ell-2h-3-2r}\cdots\ffbox{i_\ell}_{-\ell+1-2r}.
\end{equation*}
Let $C_{\ell,h,r}$ be the set of tableaux $T$ satisfying the following
conditions 
\begin{enumerate}
\def\labelenumi{(C.\theenumi)}
\def\theenumi{\arabic{enumi}}
\item  $i_a\in\mathbf B$, $i_1 \prec i_2 \prec \cdots \prec i_h$, and $i_{h+1}\prec i_{h+2}\prec \cdots \prec i_{\ell}$.

\item There is no pair $a$, $b$ such that $1\leq a < b \leq h$ and
  $i_a=k$, $i_b=\overline{k}$ and $b-a = n - k$.

\item There is no pair $a$, $b$ such that $h+1\leq a < b \leq \ell$
  and $i_a=k$, $i_b=\overline{k}$ and $b-a = n - k$.

\item There is no pair $a$, $b$ such that $a \leq h$ , $h+1 \leq b$ , $i_a=k$, $i_b=\overline{k}$ and $b-a = n + 1 - \max(r,1) - k$.

\item  Suppose that $i_{h+1} = k \in \{1,\dots, n\}$ and $i_h\succeq
i_{h+1}$. Then $i_h = k'$ is also in $\{ 1,\dots, n\}$,
and the successive part $(\overline{k'}, \overline{k'-1},\dots,
\overline{k})$ appears as $(i_{b'},i_{b'+1},\dots,i_{b})$ with
$n - r - k + 1 < b - h \le n- k$.

\item Suppose that $i_{h+1} = \overline{k} \in \{\overline{1},\dots,
  \overline{n}\}$ and $i_h\succeq i_{h+1}$. Then $i_h = \overline{k'}$ is also in $\{ \overline{1},\dots, \overline{n}\}$,
  and the successive part $(k',k'+1,\dots,k)$ appears as
  $(i_{a'},i_{a'+1},\dots,i_{a})$ with
  $n-r-k+1 \le h - a < n-k$.\end{enumerate}

Note that the conditions above are the same as the ones in \cite{kks}
when $r=0$. Note also that we only have $r = 0$ when $\ell = n$.

For $T = ((i_1,\dots,i_h),(i_{h+1},\dots,i_\ell))\in C_{\ell,h,r}$ we define the tableau $\tau_{\ell,h,r}(T)$ in the following three cases separately. 
\begin{enumerate}
\def\labelenumi{(C.\theenumi)}
\def\theenumi{\alph{enumi}}
\item $i_{h+1} = k \in \{1,\dots,n\}$ and there is an entry $i_b = \overline{k}$ with $n-r-k+1< b - h \le n-k$.

\item $i_{h+1} = \overline{k} \in \{\overline{1},\dots,\overline{n}
\}$ and there is an entry $i_a = k$ with $n-r-k+1 \le h - a < n - k$.

\item Neither (C.a) nor (C.b) is not satisfied. 
\end{enumerate}

In the case (C.a), let $b''$ such that 
$(i_{b''},i_{b''+1},\dots,i_b)$ are successive as
$(\overline{k''},\overline{k''+1},\dots,\overline{k})$ and
$i_{b''-1}\neq \overline{k''-1}$. We have $k'' < n-1$. We set
\begin{equation*}
  \begin{split}
    \tau_{\ell,h,r}(T)=((i_1,\cdots&,i_h,k''\!+\!1),
\\
  & (i_{h+2},\cdots,i_{b''-1}, \overline{k''\!+\!1},\overline{k''},
  \cdots,\overline{k\!+\!1}, i_{b+1},\cdots, i_\ell)). 
  \end{split}
\end{equation*}

Similarly in the case (C.b), we take $i_{a''}$ so that
$(i_{a''},i_{a''+1},\dots, i_a) = (k'',k''+1,\dots,k)$ and $i_{a''-1}\neq
k''-1$. We have $k < n-1$. We then set
\begin{equation*}
  \begin{split}
    \tau_{\ell,h,r}(T)=((i_1,\cdots, i_{a''-1},k''\!-\!1,\cdots,k\!-\!1,
    i_{a+1},&\cdots,i_{h}, \\
    & \overline{k''\!-\!1}), (i_{h+2},\cdots, i_\ell)).
  \end{split}
\end{equation*}

In the case (C.c) we set
\begin{equation*}
\tau_{\ell,h,r}(T)=((i_1,\cdots, i_{h+1}),(i_{h+2},\cdots, i_\ell)).
\end{equation*}

\begin{thm}
\textup{(1)} The map $T \mapsto m_T$ induces a crystal isomorphism between
$C_{\ell,h,r}$ and $\mathcal M_{I_0}(M_{\ell,h,r})$.

\textup{(2)} $\tau_{\ell,h,r}$ induces a crystal isomorphism
$\mathcal M_{I_0}(M_{\ell,h,r})$ to $\mathcal M_{I_0}(M_{\ell,h+1,r})$.
\end{thm}

\subsubsection{} Now we study $\mathcal{B}(\varpi_\ell)$ for $1\leq \ell\leq n$. Let
\(
  M_{0} = Y_{\ell,0} Y_{0,\ell}^{-1}
    = \ffbox{1}_{\ell-1}\ffbox{2}_{\ell-3} \cdots
      \ffbox{\ell}_{1-\ell}.
\)
It follows from Corollary~\ref{corex} that $\mathcal M(M_{0}) \simeq \mathcal B(\varpi_\ell)$.

For $0\le j < \ell$, let us define the monomial
$m_{T;j}$ associated with $T =
((i_1,\dots,i_{\ell - j}),(i_{\ell-j+1},\cdots,i_{\ell}))\in
C_{\ell,\ell-j,n-\ell}$ by
\begin{equation*}
m_{T;j} = \prod_{a=1}^{\ell-j} \ffbox{i_a}_{-2j+\ell+1-2a}\times \prod_{a=\ell-j+1}^\ell \ffbox{i_a}_{3\ell+1-2n-2j-2a}.
\end{equation*}
We extend the definition of $m_{T;j}$ for all $j\in\mathbb Z$ so
that $m_{T;j+\ell} = \tau_{2n} m_{T;j}$.

We describe the action of $\tilde{e}_0$, $\tilde{f}_0$ by computation on monomials. We get that $\tilde{e}_0(m_{T;j})$ is equal to
\begin{equation*}
\begin{split} 
        \begin{cases}
        m_{(i_2,\cdots,i_{\ell},\overline{1});j+1}&\text{ if $i_1 = 1$ and $i_{\ell}\neq \overline{1}$,}
                \\0 &\text{ otherwise,}
  \end{cases} 
\end{split}
\end{equation*}
and that $\tilde{f}_0(m_{T;j})$ is equal to
\begin{equation*}
\begin{split} 
        \begin{cases}
        m_{(1,i_1,\cdots,i_{\ell-1});j-1}&\text{ if $i_1 \neq 1$ and $i_{\ell} = \overline{1}$,}
                \\0 &\text{ otherwise.}
  \end{cases} 
\end{split}
\end{equation*}

We have $\tau_{2n} = z_{\ell}^{-\ell}$ and all monomials in $\mathcal
M(M_{0})/\tau_{2n}$ are written as $m_{T;j}$. The case $\ell = n$ is exceptional. We have $\tau_2 = z_{n}^{-1}$, so $\mathcal M(M_0)/\tau_2 \simeq \mathcal B(W(\varpi_n))$.
\begin{NB}
We use the same formula to define $m_{T;\ell} = \tau_{2n}(m_{T ;
0})$. For $0\le j\le \ell$ we define $M_{j}=m_{T_0;j}$ where
$T_0=(1,\cdots,\ell)$.\end{NB}
The $P_{\mathrm{cl}}$-crystal automorphism $z_\ell$ is given by
$\tau_{\ell,\ell-j-1,n-\ell}^{-1}\begin{NB}
\colon \mathcal M_{I_0}(M_{j})\to
\mathcal M_{I_0}(M_{j+1})
\end{NB}
$ 
\begin{NB}
(for $0\le j <\ell$)
\end{NB}
.

As an application, we have
\begin{equation*}
    \mathcal B(W(\varpi_\ell))
    \simeq \mathcal B_{I_0}(\varpi_\ell).
\end{equation*}
\begin{NB}
\begin{equation*}
\simeq\mathcal M(M_{0})/z_\ell=\mathcal M_{I_0}(M_{0})  
\end{equation*}
\end{NB}


A conjectural description of the crystal of $\mathcal B(W(\varpi_\ell))$ was
proposed in \cite{oss}. As their description is given by relating the
crystal to an $A^{(1)}_{2n+1}$-crystal, it is not clear, at least to 
authors, whether their conjecture is true or not.

\begin{NB}
\subsubsection{} Finally we consider the case $\ell=n$. Let $M =
Y_{n,0}Y_{0,n}^{-1}$. It follows from Corollary~\ref{corex} that
$\mathcal M(M) \simeq \mathcal B(\varpi_\ell)$. 
As $\tilde{e}_n\tilde{e}_{n-1}^2\tilde{e}_{n-2}^2\cdots\tilde{e}_1^2\tilde{e}_0 M = \tau_{-2}(M)$, $\mathcal M(M)$ is preserved under $\tau_2$. Moreover by weight calculation $z_l=\tau_{-2}$ and so $\mathcal M(M)/\tau_2\simeq \mathcal B(W(\varpi_\ell))$.

The monomials appearing in $\mathcal M_{I_0}(M)$ are 
$m_T=\prod_{a=1}^{n} \ffbox{i_a}_{n+1-2a}$ associated with a tableau
$T = (i_1,\dots,i_n)$ satisfying $i_a\in \mathbf B$ and $i_1 \prec i_2 \prec \dots \prec i_{n+1}$. These and their $\tau_2$-images are the monomials appearing in $\mathcal M(M)$. As an application, we have
$\mathcal M(M)/\tau_2=\mathcal M_{I_0}(M)$.
\end{NB}

\subsection{Type $A_{2n}^{(2)}$ ($n\geq 1$)}

Let $\mathbf B = \{ 1,\dots,n,\overline{n},\dots, \overline{1}\}$. 
We give the ordering $\prec$ on the set $\mathbf B$
by
\begin{equation*}
  1 \prec 2 \prec \cdots \prec n \prec \overline{n} \prec\cdots \prec \overline{2}\prec \overline{1}.
\end{equation*}

For $p\in \ZZ$, we define
{\allowdisplaybreaks
\begin{equation*}
\begin{aligned}[c]
    & \ffbox{1}_p = Y_{1,p}Y_{0,p+1}^{-2} , \quad \ffbox{\overline{1}}_p=Y_{0,p+2n-1}^2Y_{1,p+2n}^{-1},
    \\
    & \ffbox{i}_p=Y_{i,p+i-1}Y_{i-1,p+i}^{-1} \qquad(2 \leq i \leq n),
    \\
    & \ffbox{\overline{i}}_p=Y_{i-1,p+2n-i}Y_{i,p+2n-i+1}^{-1}
    \qquad(2 \leq i \le n).
\end{aligned}
\end{equation*}}

\subsubsection{} First consider the case $\ell=1$. Let $M = Y_{1,0}Y_{0,1}^{-2}$. It follows from Corollary~\ref{corex} that $\mathcal M(M) \simeq \mathcal B(\varpi_\ell)$. Let $M'=\tilde{e}_0 (M)=Y_{0,-1}Y_{0,1}^{-1}$. The crystal graph of $\mathcal M(M)$ is given in Figure~\ref{fig:A22}.
We find $\tau_{2n} = z_{\ell}^{-1}$ and $\mathcal
M(M)/\tau_{2n}=\mathcal M_{I_0}(M)\sqcup \mathcal M_{I_0}(M')$.  
\begin{figure}[htbp]
\centering
\psset{xunit=.55mm,yunit=.55mm,runit=.55mm}
\psset{linewidth=0.3,dotsep=1,hatchwidth=0.3,hatchsep=1.5,shadowsize=1}
\psset{dotsize=0.7 2.5,dotscale=1 1,fillcolor=black}
\psset{arrowsize=1 2,arrowlength=1,arrowinset=0.25,tbarsize=0.7 5,bracketlength=0.15,rbracketlength=0.15}
\begin{pspicture}(0,0)(220,50)
\rput(105,40){$M'$}
\rput{0}(106,-500){\parametricplot[arrows=->]{78}{89.15}{ t cos 500.26 mul t sin 540.26 mul }}
\rput{0}(106,-500){\parametricplot[arrows=<-]{102.15}{91.5}{ t cos 500.26 mul t sin 540.26 mul }}
\rput(192,40){$\scriptscriptstyle 0[2n]$}
\rput(17,40){$\scriptscriptstyle 0$}
\rput(0,20){$\ffbox{1}_0$}
\rput(20,20){$\ffbox{2}_0$}
\psline{->}(5,20)(13,20)
\psline{->}(25,20)(33,20)
\psline{->}(85,20)(93,20)
\psline[linestyle=dashed,dash=1 1](40,20)(80,20)
\rput(100,20){$\ffbox{n}_0$}
\psline{->}(105,20)(113,20)
\rput(120,20){$\ffbox{\overline{n}}_0$}
\psline{->}(125,20)(133,20)
\psline{->}(185,20)(193,20)
\psline[linestyle=dashed,dash=1 1](140,20)(180,20)
\rput(200,20){$\ffbox{\overline{2}}_0$}
\psline{->}(205,20)(213,20)
\rput(220,20){$\ffbox{\overline{1}}_0$}
\rput(8,25){$\scriptscriptstyle 1$}
\rput(30,25){$\scriptscriptstyle 2$}
\rput(108,25){$\scriptscriptstyle n$}
\rput(129,25){$\scriptscriptstyle n-1$}
\rput(89,25){$\scriptscriptstyle n-1$}
\rput(190,25){$\scriptscriptstyle 2$}
\rput(210,25){$\scriptscriptstyle 1$}
\end{pspicture}
\caption{(Type $A_{2n}^{(2)}$) the crystal $\mathcal B(\varpi_1)$}
\label{fig:A22}
\end{figure}

\subsubsection{} Now we study $\mathcal{B}(\varpi_\ell)$ for $1\leq \ell\leq n$. Let
\(
  M_{0,0} = Y_{\ell,0} Y_{0,\ell}^{-2}
    = \ffbox{1}_{\ell-1}\ffbox{2}_{\ell-3} \cdots
      \ffbox{\ell}_{1-\ell}.
\)
It follows from Corollary~\ref{corex} that $\mathcal M(M_{0,0}) \simeq \mathcal B(\varpi_\ell)$. 

For $0\le j < \ell$, $0\le k < \ell$, let us define the monomial
$m_{T;j,k}$ associated with $T =
((i_1,\dots,i_{\ell-j-k}),(i_{j-2k+1},\cdots,i_{\ell-k}))\in
C_{\ell-k,\ell-j-k,n-\ell}$ by
\begin{enumerate} 
    \item $0\leq k\leq \ell-j-1$: 
                \begin{equation*}
                \begin{split}
                    m_{T;j,k}&=(Y_{0,\ell-2j}^{-1}Y_{0,\ell-2j-2k}) 
                    \prod_{a=1}^{\ell-j-k} \ffbox{i_a}_{-2j+\ell+1-2a-2k} 
                    \\
                    &\qquad\qquad\times
                \prod_{a=\ell-j-k+1}^{\ell-k} \ffbox{i_a}_{3\ell+1-2n-2j-2a-2k},
            \end{split}
            \end{equation*}

    \item $\ell-j\leq k\leq \ell-1$: 
            \begin{equation*}
                m_{T;j,k}=(Y_{0,\ell-2j}^{-1}Y_{0,-\ell}Y_{0,\ell-2n}^{-1}Y_{0,-2n+3\ell-2j-2k})
                \prod_{a=1}^{\ell-k} \ffbox{i_a}_{3\ell+1-2n-2j-2a-2k}. 
            \end{equation*}
\end{enumerate}
For $k = \ell$ we set $C_{0,-j,n-\ell} = \{ \emptyset\}$ and
define $m_{\emptyset;j,k}$ by the same formula as in (1),(2) where the
last product is understood as $1$.
We extend the definition of $m_{T;j,k}$ for all $j\in\mathbb Z$ so
that $m_{T;j+\ell,k} = \tau_{2n} m_{T;j,k}$.

We describe the action of $\tilde{e}_0$, $\tilde{f}_0$. We get that $\tilde{e}_0(m_{T;j,k})$ is equal to
\begin{equation*}
\begin{split} 
        \begin{cases}
                m_{(i_2,\cdots,i_{\ell-k});j,k+1}&\text{ if $i_1=1$ and $i_{\ell-k}\neq \overline{1}$,}
                                        \\ m_{(i_1,\cdots,i_{\ell-k},\overline{1});j+1,k-1}&\text{ if $i_1\neq 1$, $i_{\ell-k}\neq \overline{1}$ and $k>0$},
                                        \\0 &\text{ otherwise,}
                        \end{cases} 
\end{split}
\end{equation*}
and that $\tilde{f}_0(m_{T;j,k})$ is equal to
\begin{equation*}
\begin{split} 
        \begin{cases}
                   m_{(i_1,\cdots,i_{\ell-k-1});j-1,k + 1}&\text{ if $i_1 \neq 1$ and $i_{\ell-k} = \overline{1}$,}
                 \\ m_{(1,i_1,\cdots,i_{\ell-k});j,k-1}&\text{ if $i_1 \neq 1$, $i_{\ell-k} \neq \overline{1}$ and $k>0$,}
                 \\0 &\text{ otherwise.}
                          \end{cases} 
\end{split}
\end{equation*}
We have $\tau_{2n} = (z_{\ell})^{-\ell}$ and all
the monomials in $\mathcal M(M_{0,0})/\tau_{2n}$ are written as $m_{T;j,k}$.
The case $\ell = n$ is exceptional. We have $\tau_2 = z_{n}^{-1}$, so $\mathcal M(M_0)/\tau_2 \simeq \mathcal B(W(\varpi_n))$.
\begin{NB}
We use the same formula to define $m_{T;\ell,k} = \tau_{2n}(m_{T ; 0, k})$. For $0\le j\le \ell, 0\le k\le\ell$ we define $M_{j,k}=m_{(1,\cdots,\ell-k);j,k}$. It is clear that $z_\ell$ is given by
$\tau_{\ell - k,\ell-j-k-1,n-\ell}^{-1}\colon \mathcal M_{I_0}(M_{j,k})\to 
\mathcal M_{I_0}(M_{j+1,k})$ (for $0\le j <\ell$, $0\le k \le \ell$).
\end{NB}
For $\ell\neq n$, the crystal automorphism $z_\ell$ is given by
$\tau_{\ell - k,\ell-j-k-1,n-\ell}^{-1}$.
As an application, we have
\begin{equation*}
    \mathcal B(W(\varpi_\ell))\simeq
    \mathcal B_{I_0}(\varpi_\ell)
    \sqcup \mathcal B_{I_0}(\varpi_{\ell-1})\sqcup \cdots\sqcup
    \mathcal B_{I_0}(\varpi_1) \sqcup\mathcal B_{I_0}(0).
\end{equation*}
\begin{NB}
\begin{equation*}
    \mathcal B(W(\varpi_\ell))\simeq\mathcal M(M_{0,0})/z_\ell=\mathcal M_{I_0}(M_{0,0})\sqcup \mathcal M_{I_0}(M_{0,1})\sqcup        \cdots\sqcup \mathcal M_{I_0}(M_{0,\ell}).
\end{equation*}
\end{NB}


A conjectural description of the crystal of $\mathcal B(W(\varpi_\ell))$ was
proposed in \cite{oss}. As their description is given by relating the
crystal to an $A^{(1)}_{2n+1}$-crystal, it is not clear, at least to 
authors, whether their conjecture is true or not.

\begin{NB}
\subsubsection{} Finally we consider the case $\ell=n$. Let $M_{0} = Y_{n,0}Y_{0,n}^{-2}=\prod_{a=1}^{n} \ffbox{a}_{n+1-2a}$. It follows from Corollary~\ref{corex} that $\mathcal M(M_0) \simeq \mathcal B(\varpi_\ell)$. As $\tilde{e}_n\tilde{e}_{n-1}^2\tilde{e}_{n-2}^2\cdots\tilde{e}_0^2 M_0=\tau_{-2}(M_0)$, $\mathcal M(M_0)$ is preserved under $\tau_2$. For weight reason $z_\ell=\tau_{-2}$ (because $d_\ell=1$) and so $\mathcal M(M_0)/\tau_2\simeq \mathcal B(W(\varpi_\ell))$.

For $0\le k \le n$, let us define the monomial
$m_{T;k}$ associated with $T =
((i_1,\dots,i_{n-k})$  satisfying $i_a\in\mathbf B$ and $i_1\prec i_2 \prec \cdots \prec i_{n-k}$ by
\begin{equation*}m_{T;k}=(Y_{0,n}^{-1}Y_{0,n-2k})
\prod_{a=1}^{n-k} \ffbox{i_a}_{-2k+n+1-2a}.\end{equation*}
These and their $\tau_{2}$-images are the monomials appearing in 
$\mathcal M(M_{0})$. For $0\le k\le n$, we define $M_k=m_{(1,\cdots,n - k);k}$. As an application, we have
\begin{equation*}
    \mathcal M(M_0)/\tau_2=\mathcal M_{I_0}(M_0)\sqcup \mathcal M_{I_0}(M_1)\sqcup \cdots\sqcup\mathcal M_{I_0}(M_n).
\end{equation*}
\end{NB}

\subsection{Type $A_{2n}^{(2)\dagger}$ ($n\geq 1$)}


Let $\mathbf B = \{ 1,\dots,n,0,\overline{n},\dots, \overline{1}\}$. 
We give the ordering $\prec$ on the set $\mathbf B$
by
\begin{equation*}
  1 \prec 2 \prec \cdots \prec n \prec 0\prec \overline{n} \prec\cdots \prec \overline{2}\prec \overline{1}.
\end{equation*}

For $p\in \ZZ$, we define
{\allowdisplaybreaks
\begin{equation*}
\begin{aligned}[c]
    & \ffbox{i}_p = Y_{i-1,p+i}^{-1} Y_{i,p+{i-1}} \qquad(1 \le i \le n-1),
    \\
    & \ffbox{n}_p
     = Y_{n-1,p+n}^{-1} Y_{n,p+{n-1}}^2,
    \\
    & \ffbox{0}_p = Y_{n,p+{n}+1}^{-1} Y_{n,p+{n-1}},
    \\
    & \ffbox{\overline{n}}_p = Y_{n-1,p+n} Y_{n,p+{n}+1}^{-2},
    \\
    & \ffbox{\overline{i}}_p = Y_{i-1,p+{2n-i}} Y_{i,p+{2n+1-i}}^{-1}
    \qquad(1 \le i \le n-1).
  \end{aligned}
\end{equation*}}

\subsubsection{} 
\begin{figure}[htbp]
\centering
\psset{xunit=.55mm,yunit=.55mm,runit=.55mm}
\psset{linewidth=0.3,dotsep=1,hatchwidth=0.3,hatchsep=1.5,shadowsize=1}
\psset{dotsize=0.7 2.5,dotscale=1 1,fillcolor=black}
\psset{arrowsize=1 2,arrowlength=1,arrowinset=0.25,tbarsize=0.7 5,bracketlength=0.15,rbracketlength=0.15}
\begin{pspicture}(0,0)(220,30)
\rput(0,20){$\ffbox{1}_0$}
\rput(20,20){$\ffbox{2}_0$}
\psline{->}(5,20)(13,20)
\psline{->}(25,20)(33,20)
\psline{->}(75,20)(83,20)
\psline[linestyle=dashed,dash=1 1](40,20)(70,20)
\rput(90,20){$\ffbox{n}_0$}
\rput(110,20){$\ffbox{0}_0$}
\rput(130,20){$\ffbox{\overline{n}}_0$}
\psline{->}(135,20)(143,20)
\psline{->}(185,20)(193,20)
\psline[linestyle=dashed,dash=1 1](150,20)(180,20)
\rput(200,20){$\ffbox{\overline{2}}_0$}
\psline{->}(205,20)(213,20)
\rput(220,20){$\ffbox{\overline{1}}_0$}
\psline{->}(95,20)(103,20)
\psline{->}(115,20)(123,20)
\rput{0}(110.26,520){\parametricplot[arrows=<-]{-102.15}{-77.85}{ t cos 520.26 mul t sin 520.26 mul }}
\rput(8,25){$\scriptscriptstyle 1$}
\rput(30,25){$\scriptscriptstyle 2$}
\rput(98,25){$\scriptscriptstyle n$}
\rput(138,25){$\scriptscriptstyle n-1$}
\rput(118,25){$\scriptscriptstyle n$}
\rput(78,25){$\scriptscriptstyle n-1$}
\rput(190,25){$\scriptscriptstyle 2$}
\rput(210,25){$\scriptscriptstyle 1$}
\rput(110,3){$\scriptscriptstyle 0[2n]$}
\end{pspicture}
\caption{(Type $A_{2n}^{(2)\dagger}$) the crystal $\mathcal B(\varpi_1)$}
\label{fig:A2n2second}
\end{figure}
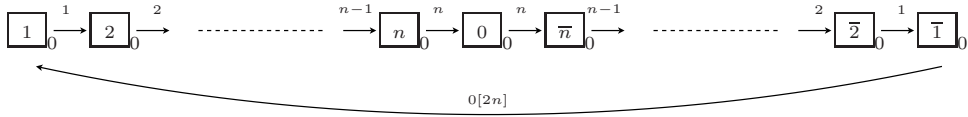

First consider the case $\ell=1$. Let $M = Y_{0,1}^{-1}Y_{1,0}$. It follows from Corollary~\ref{corex} that $\mathcal M(M) \simeq \mathcal B(\varpi_\ell)$.
The crystal graph of $\mathcal M(M)$ is given in Figure~\ref{fig:A2n2second}.
We find $\tau_{2n} = z_{\ell}^{-1}$ and $\mathcal
M(M)/\tau_{2n}=\mathcal M_{I_0}(M)$.

\subsubsection{} Now we study $\mathcal{B}(\varpi_\ell)$ for $1\leq \ell\leq n - 1$. Let
\(
  M_{0} = Y_{\ell,0} Y_{0,\ell}^{-1}
    = \ffbox{1}_{\ell-1}\ffbox{2}_{\ell-3} \cdots
      \ffbox{\ell}_{1-\ell}.
\)
It follows from Corollary~\ref{corex} that $\mathcal M(M_{0}) \simeq \mathcal B(\varpi_\ell)$.

For $0\le j < \ell$, let us define the monomial
$m_{T;j}$ associated with $T =
((i_1,\dots,i_{\ell - j}),(i_{\ell-j+1},\cdots,i_{\ell}))\in
B_{\ell,\ell-j,n-\ell}$ by
\begin{equation*}
m_{T;j} = \prod_{a=1}^{\ell-j} \ffbox{i_a}_{-2j+\ell+1-2a}\times \prod_{a=\ell-j+1}^\ell \ffbox{i_a}_{3\ell+1-2n-2j-2a}.
\end{equation*}
We extend the definition of $m_{T;j}$ for all $j\in\mathbb Z$ so
that $m_{T;j+\ell} = \tau_{2n} m_{T;j}$. 
We describe the action of $\tilde{e}_0$, $\tilde{f}_0$ by computation on monomials. We get that $\tilde{e}_0(m_{T;j})$ is equal to
\begin{equation*}
\begin{split} 
        \begin{cases}
        m_{(i_2,\cdots,i_{\ell},\overline{1});j+1}&\text{ if $i_1 = 1$ and $i_{\ell}\neq \overline{1}$,}
                \\0 &\text{ otherwise,}
  \end{cases} 
\end{split}
\end{equation*}
and that $\tilde{f}_0(m_{T;j})$ is equal to
\begin{equation*}
\begin{split} 
        \begin{cases}
        m_{(1,i_1,\cdots,i_{\ell-1});j-1}&\text{ if $i_1 \neq 1$ and $i_{\ell} = \overline{1}$,}
                \\0 &\text{ otherwise.}
  \end{cases} 
\end{split}
\end{equation*}
We have $\tau_{2n} = z_{\ell}^{-\ell}$ and all monomials in $\mathcal
M(M_{0})/\tau_{2n}$ are written as $m_{T;j}$.
The $P_{\mathrm{cl}}$-crystal automorphism $z_\ell$ is given by
$\tau_{\ell,\ell-j-1,n-\ell}^{-1}$.
As an application, we have
\begin{equation*}
    \mathcal B(W(\varpi_\ell))
    \simeq \mathcal B_{I_0}(\varpi_\ell).
\end{equation*}


A conjectural description of the crystal of $\mathcal B(W(\varpi_\ell))$ was
proposed in \cite{oss}. As their description is given by relating the
crystal to an $A^{(1)}_{2n+1}$-crystal, it is not clear, at least to 
authors, whether their conjecture is true or not.

\subsubsection{} Finally we consider the case $\ell=n$. Let
\(
  M = Y_{n,0}^2 Y_{0,\ell}^{-1}
    = \ffbox{1}_{n-1}\ffbox{2}_{n-3} \cdots
      \ffbox{n}_{1-n}.
\)
It follows from Corollary~\ref{corex} that $\mathcal M(M) \simeq
\mathcal B(\varpi_n)$.  Let us define the monomial $m_{T} =
\prod_{a=1}^{n} \ffbox{i_a}_{n+1-2a}$ associated with $T =
(i_1,\cdots,i_{n})$ satisfying (1) $i_a\in \mathbf B$ and $i_1\prec
i_2\prec\cdots \prec i_n$ but $0$ can be repeated, and (2) there is no
pair $a$, $b$ such that $i_a = k$, $i_b = \overline{k}$ and $b-a=n-k$.
The above exhausts all monomials in $\mathcal M_{I_0}(M)$ (see
\cite[Proposition~2.10]{kks}). We describe the action of $\tilde{e}_0$, $\tilde{f}_0$ on these monomials : we have
\begin{equation*}
\begin{split}
\tilde{e}_0(m_{T}) &=
        \begin{cases}
         \tau_{-2}(m_{(i_2,\cdots,i_n,\overline{1})}) &\text{ if $i_1 = 1$ and $i_n\neq\overline{1}$},
                \\  0&\text{ otherwise,}
  \end{cases}
\\
\tilde{f}_0(m_{T}) &=
        \begin{cases}
        \tau_2(m_{(1,i_1,\cdots,i_{n-1})}) &\text{ if $i_{n}=\overline{1}$ and $i_1\neq 1$,}
                \\  0&\text{ otherwise.}
  \end{cases} 
\end{split}
\end{equation*}
So the above exhausts all the monomials in $\mathcal
M(M)/\tau_{2}$. We have $\tau_2 = z_{n}^{-1}$, so $\mathcal
M(M)/\tau_2 \simeq \mathcal B(W(\varpi_n))$. As an application, we
have $\mathcal B(W(\varpi_n))\simeq \mathcal B_{I_0}(\varpi_n)$.
Note that $\varpi_n$ is identified with the {\it twice\/} of the
$n^{\mathrm{th}}$ fundamental weight of $\mathfrak g_{I_0}$.





\subsection{Type $A_{2n-1}^{(2)}$ ($n\geq 3$)}

Let $\mathbf B = \{ 1,\dots,n,\overline{n},\dots, \overline{1}\}$. 
We give the ordering $\prec$ on the set $\mathbf B$
by
\begin{equation*}
  1 \prec 2 \prec \cdots \prec n \prec \overline{n} \prec\cdots \prec \overline{2}\prec \overline{1}.
\end{equation*}

For $p\in \ZZ$, we define
{\allowdisplaybreaks\begin{equation*}
\begin{aligned}[c]
    & \ffbox{1}_p = Y_{0,p+2}^{-1} Y_{1,p}, \quad
     \ffbox{2}_p = Y_{0,p+2}^{-1} Y_{1,p+2}^{-1} Y_{2,p+1},
    \\
    & \ffbox{i}_p = Y_{i-1,p+i}^{-1} Y_{i,p+{i-1}} \qquad(3 \le i \le n),
    \\
    & \ffbox{\overline{i}}_p = Y_{i-1,p+{2n-i}} Y_{i,p+{2n+1-i}}^{-1}
    \qquad(3 \le i \le n),
    \\
    & \ffbox{\overline{2}}_p = Y_{0,p+2n-2} Y_{1,p+2n-2} Y_{2,p+2n-1}^{-1},
    \qquad
    \ffbox{\overline{1}}_p = Y_{0,p+2n-2} Y_{1,p+{2n}}^{-1}.
\end{aligned}
\end{equation*}}

\subsubsection{} First consider the case $\ell=1$. Let $M = Y_{1,0}Y_{0,2}^{-1}$. It follows from Corollary~\ref{corex} that $\mathcal M(M) \simeq \mathcal B(\varpi_1)$. The crystal graph of $\mathcal M(M)$ is given in Figure~\ref{fig:A2n12}.
We find $\tau_{2n-2} = z_{\ell}^{-1}$ and $\mathcal
M(M)/\tau_{2n-2}=\mathcal M_{I_0}(M)$.  
\begin{figure}[htbp]
\centering
\psset{xunit=.55mm,yunit=.55mm,runit=.55mm}
\psset{linewidth=0.3,dotsep=1,hatchwidth=0.3,hatchsep=1.5,shadowsize=1}
\psset{dotsize=0.7 2.5,dotscale=1 1,fillcolor=black}
\psset{arrowsize=1 2,arrowlength=1,arrowinset=0.25,tbarsize=0.7 5,bracketlength=0.15,rbracketlength=0.15}
\begin{pspicture}(0,0)(220,50)
\rput{0}(0,25){\qline(-0,0)(0,0)}\rput{0}(100,-480){\parametricplot[arrows=->]{78.85}{101.15}{ t cos 517.26 mul t sin 517.26 mul }}
\rput{0}(117.26,520){\parametricplot[arrows=<-]{-101.15}{-78.85}{ t cos 517.26 mul t sin 517.26 mul }}
\rput(0,20){$\ffbox{1}_0$}
\rput(20,20){$\ffbox{2}_0$}
\psline{->}(5,20)(13,20)
\psline{->}(25,20)(33,20)
\psline{->}(85,20)(93,20)
\psline[linestyle=dashed,dash=1 1](40,20)(80,20)
\rput(100,20){$\ffbox{n}_0$}
\psline{->}(105,20)(113,20)
\rput(120,20){$\ffbox{\overline{n}}_0$}
\psline{->}(125,20)(133,20)
\psline{->}(185,20)(193,20)
\psline[linestyle=dashed,dash=1 1](140,20)(180,20)
\rput(200,20){$\ffbox{\overline{2}}_0$}
\psline{->}(205,20)(213,20)
\rput(220,20){$\ffbox{\overline{1}}_0$}
\rput(8,25){$\scriptscriptstyle 1$}
\rput(30,25){$\scriptscriptstyle 2$}
\rput(108,25){$\scriptscriptstyle n$}
\rput(129,25){$\scriptscriptstyle n-1$}
\rput(89,25){$\scriptscriptstyle n-1$}
\rput(190,25){$\scriptscriptstyle 2$}
\rput(210,25){$\scriptscriptstyle 1$}
\rput(100,41){$\scriptscriptstyle 0[2n-2]$}
\rput(185,3){$\scriptscriptstyle 0[2n-2]$}
\end{pspicture}
\caption{(Type $A_{2n-1}^{(2)}$) the crystal $\mathcal B(\varpi_1)$}
\label{fig:A2n12}
\end{figure}

\subsubsection{} Now we study $\mathcal{B}(\varpi_\ell)$ for $2\leq \ell\leq n-1$. Let
\(
  M_{0,0} = Y_{\ell,0} Y_{0,\ell-1}^{-1} Y_{0,\ell+1}^{-1}
    = \ffbox{1}_{\ell-1}\ffbox{2}_{\ell-3} \cdots
      \ffbox{\ell}_{1-\ell}.
\)
As $\tilde{f}_2\tilde{f}_3\cdots\tilde{f}_\ell(M_{0,0})=Y_{0,\ell+1}^{-1}Y_{0,\ell-1}Y_{1,\ell}^{-1}Y_{\ell+1,1}$, we see as in Proposition~\ref{otherd} 
that $\mathcal M(M_{0,0}) \simeq \mathcal B(\varpi_\ell)$. 

For $0\le j < \ell$, $0\le k < \ell/2$, let us define the monomial
$m_{T;j,k}$ associated with $T =
((i_1,\dots,i_{j-2k}),(i_{j-2k+1},\cdots,i_{\ell-2k}))\in
C_{\ell-2k,j-2k,n-\ell-1}$ by
\begin{enumerate}
    \item $k < \lfloor j/2\rfloor$: 
                \begin{equation*}
                \begin{split}
                    m_{T;j,k}
                = &Y_{0,2n-\ell-4k+2j-1} Y_{0,2n-\ell+2j-1}^{-1} 
                \prod_{a=1}^{j-2k} \ffbox{i_a}_{2n-\ell-4k-2a+2j-1}
                \\ & \qquad\qquad\times
                \prod_{a=j-2k+1}^{\ell-2k} \ffbox{i_a}_{\ell-2(a-j+2k)+1},
                \end{split}
                \end{equation*}

    \item $j$ is odd and $k = (j-1)/2$:
                \begin{equation*}
                m_{T;j,(j-1)/2} = Y_{0,2n-\ell-3} Y_{0,2n-\ell+2j-1}^{-1} \ffbox{i_1}_{2n-\ell-1}
                \prod_{a=2}^{\ell-j+1} \ffbox{i_a}_{\ell-2a+3},
                \end{equation*}

    \item $j$ is even and $k \ge j/2$:
                \begin{equation*}
                m_{T;j,k} = Y_{0,\ell-4k+2j+1} Y_{0,\ell+1}^{-1}
            Y_{0,2n-\ell-1}Y_{0,2n-\ell+2j-1}^{-1}\prod_{a=1}^{\ell-2k} \ffbox{i_a}_{\ell-2a-4k+2j+1},
                \end{equation*}

    \item $j$ is odd and $k\ge (j+1)/2$:
                \begin{equation*}
                m_{T;j,k} = Y_{0,\ell-4k+2j+1} Y_{0,2n-\ell+2j-1}^{-1}
                Y_{1,\ell+1}^{-1} Y_{1,2n-\ell-1}
                \prod_{a=1}^{\ell-2k} \ffbox{i_a}_{\ell-2a-4k+2j+1}.
                \end{equation*}
\end{enumerate}
For $k = \ell/2$ we set $C_{0,j-\ell,n-\ell-1} = \{ \emptyset\}$ and
define $m_{\emptyset;j,k}$ by the same formula as in (3),(4) where the
last product is understood as $1$.
We extend the definition of $m_{T;j,k}$ for all $j\in\mathbb Z$ so
that $m_{T;j+\ell,k} = \tau_{2n-2} m_{T;j,k}$. 

We describe the action of $\tilde{e}_0$, $\tilde{f}_0$. We get that $\tilde{e}_0(m_{T;j,k})$ is equal to
\begin{equation*}
\begin{split} 
        \begin{cases}
        m_{(i_3,\cdots,i_{\ell-2k});j,k+1}&\text{ if $i_2=2$ and $i_{\ell-2k-1}\nsucceq \overline{2}$,}
                \\ m_{(i_1,\cdots,i_{\ell-2k},\overline{2},\overline{1});j-2,k-1}&\text{ if $i_2\npreceq 2$, $i_{\ell-2k}\nsucceq \overline{2}$ and $k > 0$},
                \\ m_{(i_2,\cdots,i_{\ell},\overline{3-i_1});j-1,0}&\text{ if $i_1\preceq 2$, $i_2\npreceq 2$, $i_{\ell}\nsucceq \overline{2}$ and $k=0$,}
                \\0 &\text{ otherwise,}
  \end{cases} 
\end{split}
\end{equation*}
and that $\tilde{f}_0(m_{T;j,k})$ is equal to
\begin{equation*}
\begin{split} 
        \begin{cases}
        m_{(1,2,i_1,\cdots,i_{\ell-2k});j,k - 1}&\text{ if $i_1 \npreceq 2$, $i_{\ell-2k-1}\nsucceq \overline{2}$ and $k > 0$},
                \\ m_{(i_1,\cdots,i_{\ell-2k-2});j+2,k+1}&\text{ if $i_{\ell-2k-1} = \overline{2}$ and $i_2\npreceq 2$},
                \\ m_{(3-\overline{i_{\ell}},i_1,\cdots,i_{\ell-1});j+1,0}&\text{ if $i_1 \npreceq 2$, $i_{\ell-1}\nsucceq \overline{2}$, $i_{\ell}\succeq \overline{2}$ and $k = 0$,}
                \\0 &\text{ otherwise.}
  \end{cases} 
\end{split}
\end{equation*}
We have $\tau_{2n-2} = (z_\ell)^{-\ell}$, and all
monomials in  $\mathcal
M(M_{0,0})/\tau_{2n-2}$ are written as $m_{T;j,k}$. The crystal automorphism $z_\ell$ is given by
$\tau_{\ell-2k,j-2k,n-\ell-1}^{-1}$.
\begin{NB}
These and their $\tau_{2n-2}$-images are the monomials appearing in 
$\mathcal M(M_{0,0})$. We use the same formula to define $m_{T;\ell,k} = \tau_{2n-2}(m_{T ; 0, k})$. For $0\le j\le \ell, 0\le k\le\lfloor \ell/2\rfloor$ we define $M_{j,k}=m_{(1,\cdots,\ell-2k);j,k}$. It is clear that $(z_\ell)^{-1}$ is given by
$\tau_{\ell - 2k,j-2k,n-\ell-1}\colon \mathcal M_{I_0}(M_{j,k})\to 
\mathcal M_{I_0}(M_{j+1,k})$ (for $0\le j <\ell$, $0\le k \le \lfloor \ell/2\rfloor$).
\end{NB}
As an application, we have
\begin{equation*}
\begin{split}
    \mathcal B(W(\varpi_\ell))
    \simeq 
   \mathcal B_{I_0}(\varpi_\ell) \sqcup
   \mathcal B_{I_0}(\varpi_{\ell-2}) \sqcup \cdots \sqcup
   \begin{cases}
   \mathcal B_{I_0}(\varpi_1) & \text{if $\ell$ is odd},
   \\
   \mathcal B_{I_0}(0) & \text{if $\ell$ is even}.
   \end{cases}
\end{split}
\end{equation*}
\begin{NB}
\begin{equation*}
    \simeq\mathcal M(M_{0,0})/z_\ell=\mathcal M_{I_0}(M_{0,0})\sqcup \mathcal M_{I_0}(M_{0,1})\sqcup        \cdots\sqcup \mathcal M_{I_0}(M_{0,\lfloor \ell/2\rfloor})
\end{equation*}
\end{NB}


A crystal base on $W(\varpi_\ell)$ was constructed in \cite{jmo}. A
key fact used there is that $W(\varpi_\ell)$ remains irreducible when
it is restricted to $\U_q(\overline{\Glie})$ for the finite dimensional
Lie algebra $\overline{\Glie}$ obtained by removing the vertex $n$. They
showed that the crystal base for the restriction is preserved also by
$\te_n$, $\tf_n$. By the uniqueness of the crystal base for
an irreducible $\U_q(\overline{\Glie})$-module we conclude that their
crystal base is isomorphic to the $\mathcal B(W(\varpi_\ell))$. However their
description of the Kashiwara operators was given in terms of
$\overline{\Glie}$, it is not obvious to compare our description to
theirs.


\subsubsection{}  Finally we consider the case $\ell=n$. Let $M_0=Y_{n,0}Y_{0,n-1}^{-2}=\prod_{a=1}^{n} \ffbox{a}_{n+1-2a}$. It follows from Corollary~\ref{corex} that $\mathcal M(M_0) \simeq \mathcal B(\varpi_\ell)$.
\begin{NB}
As $\tilde{e}_n\tilde{e}_{n-1}^2\tilde{e}_{n-2}^2\cdots\tilde{e}_1^2\tilde{e}_n\tilde{e}_{n-1}^2\tilde{e}_{n-2}^2\cdots\tilde{e}_2^2\tilde{e}_0^2 M_0=\tau_{-4}(M_0)$, $\mathcal M(M_0)$ is preserved under $\tau_4$. Moreover for weight reason $z_\ell=\tau_{-4}$ and so $\mathcal M(M_0)/\tau_4\simeq \mathcal{B}(W(\varpi_\ell))$. 
\end{NB}

For $0\le k \le \lfloor n/2\rfloor$, let us define the monomial
$m_{T;k}$ associated with $T =
(i_1,\dots,i_{n-2k})\in C_{n-2k,0,0}$ by
\begin{equation*}
   m_{T;k}= Y_{0,n-1}^{-1}Y_{0,n+1-4k} \times
   \prod_{a=1}^{n-2k} \ffbox{i_a}_{n+1-4k-2a},
\end{equation*}
where the case $n=2k$ is understood as before.

We describe the action of $\tilde{e}_0$, $\tilde{f}_0$. We get that $\tilde{e}_0(m_{T;k})$ is equal to
\begin{equation*}
\begin{split} 
        \begin{cases}
        m_{(i_3,\cdots,i_{n-2k});k+1}&\text{ if $i_2=2$ and $i_{n-2k-1}\nsucceq \overline{2}$,}
                \\ \tau_{-4}(m_{(i_1,\cdots,i_{n-2k},\overline{2},\overline{1});k-1})&\text{ if $i_2\npreceq 2$, $i_{n-2k}\nsucceq \overline{2}$ and $k > 0$},
                \\ \tau_{-2}(m_{(i_2,\cdots,i_n,\overline{3-i_1});0})&\text{ if $i_1\preceq 2$, $i_2\npreceq 2$, $i_n\nsucceq \overline{2}$ and $k=0$,}
                \\0 &\text{ otherwise,}
  \end{cases} 
\end{split}
\end{equation*}
and that $\tilde{f}_0(m_{T;k})$ is equal to
\begin{equation*}
\begin{split} 
        \begin{cases}
        m_{(1,2,i_1,\cdots,i_{n-2k});k - 1}&\text{ if $i_1 \npreceq 2$, $i_{n-2k-1}\nsucceq \overline{2}$ and $k > 0$},
                \\ \tau_4(m_{(i_1,\cdots,i_{n-2k-2});k+1})&\text{ if $i_{n-2k-1} = \overline{2}$ and $i_2\npreceq 2$},
                \\ \tau_2(m_{(3-\overline{i_{\ell}},i_1,\cdots,i_{\ell-1});0})&\text{ if $i_1 \npreceq 2$, $i_{n-1}\nsucceq \overline{2}$, $i_n\succeq \overline{2}$ and $k = 0$,}
                \\0 &\text{ otherwise.}
  \end{cases} 
\end{split}
\end{equation*}

We find that $z_\ell = \tau_{-2}$ and the monomials appearing in 
$\mathcal M(M_{0})/\tau_{2}$ are written as $m_{T;k}$. 
\begin{NB}
For $0\le k\le\lfloor \ell/2\rfloor$ we define $M_{k}=m_{(1,\cdots,\ell-2k);k}$.  
\end{NB}
As an application, we have
\begin{NB}
\begin{equation*}
    \mathcal B(W(\varpi_\ell))\simeq\mathcal M(M_{0})/\tau_2=\mathcal M_{I_0}(M_{0})\sqcup \mathcal M_{I_0}(M_{1})\sqcup        \cdots\sqcup \mathcal M_{I_0}(M_{\lfloor n/2\rfloor}).
\end{equation*}
\end{NB}
\begin{equation*}
    \mathcal B(W(\varpi_\ell))\simeq 
   \mathcal B_{I_0}(\varpi_\ell) \sqcup
   \mathcal B_{I_0}(\varpi_{\ell-2}) \sqcup \cdots \sqcup
   \begin{cases}
   \mathcal B_{I_0}(\varpi_1) & \text{if $\ell$ is odd},
   \\
   \mathcal B_{I_0}(0) & \text{if $\ell$ is even}.
   \end{cases}
\end{equation*}


\subsection{Type $D_{n+1}^{(2)}$ ($n\geq 2$)}

Let $\mathbf B = \{ 1,\dots,n,0,\overline{n},\dots, \overline{1}\}$. 
We give the ordering $\prec$ on the set $\mathbf B$
by
\begin{equation*}
  1 \prec 2 \prec \cdots \prec n \prec 0\prec \overline{n} \prec\cdots \prec \overline{2}\prec \overline{1}.
\end{equation*}

For $p\in \ZZ$, we define
{\allowdisplaybreaks
\begin{equation*}
\begin{aligned}[c]
& \ffbox{1}_p = Y_{1,p}Y_{0,p+1}^{-2},
\\
& \ffbox{i}_p = Y_{i,p+i-1}Y_{i-1,p+i}^{-1} \qquad(2 \le i \le n-1),
\\
& \ffbox{n}_p
 = Y_{n-1,p+n}^{-1}Y_{n,p+n-1}^2,
\\
& \ffbox{0}_p = Y_{n,p+n-1}Y_{n,p+n+1}^{-1},
\\
& \ffbox{\overline{n}}_p = Y_{n-1,p+n}Y_{n,p+n+1}^{-2},
\\
& \ffbox{\overline{i}}_p = Y_{i-1,p+2n-i}Y_{i,p+2n-i+1}^{-1}
\qquad(2\leq i\leq n-1),
\\
& \ffbox{\overline{1}}_p = Y_{0,p+2n-1}^2Y_{1,p+2n}^{-1}.
\end{aligned}
\end{equation*}}

\subsubsection{} First consider the case $\ell=1$. Let $M = Y_{1,0}Y_{0,1}^{-2}$. It follows from Corollary~\ref{corex} that $\mathcal M(M) \simeq \mathcal B(\varpi_\ell)$. Let $M'=\tilde{e}_0 (M)=Y_{0,-1}Y_{0,1}^{-1}$. The crystal graph of $\mathcal M(M)$ is given in Figure~\ref{fig:D2}.
We find $\tau_{2n} = z_{\ell}^{-1}$ and $\mathcal
M(M)/\tau_{2n}=\mathcal M_{I_0}(M)\sqcup \mathcal M_{I_0}(M')$.  
\begin{figure}[htbp]
\centering
\psset{xunit=.55mm,yunit=.55mm,runit=.55mm}
\psset{linewidth=0.3,dotsep=1,hatchwidth=0.3,hatchsep=1.5,shadowsize=1}
\psset{dotsize=0.7 2.5,dotscale=1 1,fillcolor=black}
\psset{arrowsize=1 2,arrowlength=1,arrowinset=0.25,tbarsize=0.7 5,bracketlength=0.15,rbracketlength=0.15}
\begin{pspicture}(0,0)(220,50)
\rput(105,40){$M'$}
\rput{0}(106,-500){\parametricplot[arrows=->]{78}{89.15}{ t cos 500.26 mul t sin 540.26 mul }}
\rput{0}(106,-500){\parametricplot[arrows=<-]{102.15}{91.5}{ t cos 500.26 mul t sin 540.26 mul }}
\rput(192,40){$\scriptscriptstyle 0[2n]$}
\rput(17,40){$\scriptscriptstyle 0$}
\rput(0,20){$\ffbox{1}_0$}
\rput(20,20){$\ffbox{2}_0$}
\psline{->}(5,20)(13,20)
\psline{->}(25,20)(33,20)
\psline{->}(75,20)(83,20)
\psline[linestyle=dashed,dash=1 1](40,20)(70,20)
\rput(90,20){$\ffbox{n}_p$}
\rput(110,20){$\ffbox{0}_p$}
\rput(130,20){$\ffbox{\overline{n}}_p$}
\psline{->}(135,20)(143,20)
\psline{->}(185,20)(193,20)
\psline[linestyle=dashed,dash=1 1](150,20)(180,20)
\rput(200,20){$\ffbox{\overline{2}}_p$}
\psline{->}(205,20)(213,20)
\rput(220,20){$\ffbox{\overline{1}}_p$}
\psline{->}(95,20)(103,20)
\psline{->}(115,20)(123,20)
\rput(8,25){$\scriptscriptstyle 1$}
\rput(30,25){$\scriptscriptstyle 2$}
\rput(98,25){$\scriptscriptstyle n$}
\rput(138,25){$\scriptscriptstyle n-1$}
\rput(118,25){$\scriptscriptstyle n$}
\rput(78,25){$\scriptscriptstyle n-1$}
\rput(190,25){$\scriptscriptstyle 2$}
\rput(210,25){$\scriptscriptstyle 1$}
\end{pspicture}
\caption{(Type $D_{n+1}^{(2)}$) the crystal $\mathcal B(\varpi_1)$}
\label{fig:D2}
\end{figure}

\subsubsection{} Now we study $\mathcal{B}(\varpi_\ell)$ for $1\leq \ell\leq n-1$. Let
\(
  M_{0,0} = Y_{\ell,0} Y_{0,\ell}^{-2}
    = \ffbox{1}_{\ell-1}\ffbox{2}_{\ell-3} \cdots
      \ffbox{\ell}_{1-\ell}.
\)
It follows from Corollary~\ref{corex} that $\mathcal M(M_{0,0}) \simeq \mathcal B(\varpi_\ell)$. 

For $0\le j < \ell$, $0\le k < \ell$, let us define the monomial
$m_{T;j,k}$ associated with $T =
((i_1,\dots,i_{\ell-j-k}),(i_{\ell-j-k+1},\cdots,i_{\ell-k}))\in
B_{\ell-k,\ell-j-k,n-\ell}$ by
\begin{enumerate} 
    \item $0\leq k\leq \ell-j-1$: 
                \begin{equation*}
                \begin{split}
                    m_{T;j,k}&=(Y_{0,\ell-2j}^{-1}Y_{0,\ell-2j-2k}) 
                    \prod_{a=1}^{\ell-j-k} \ffbox{i_a}_{-2j+\ell+1-2a-2k} 
                    \\
                    &\qquad\qquad\times
                \prod_{a=\ell-j-k+1}^{\ell-k} \ffbox{i_a}_{3\ell+1-2n-2j-2a-2k},
            \end{split}
            \end{equation*}

    \item $\ell-j\leq k\leq \ell-1$: 
            \begin{equation*}
                m_{T;j,k}=(Y_{0,\ell-2j}^{-1}Y_{0,-\ell}Y_{0,\ell-2n}^{-1}Y_{0,-2n+3\ell-2j-2k})
                \prod_{a=1}^{\ell-k} \ffbox{i_a}_{3\ell+1-2n-2j-2a-2k}. 
            \end{equation*}
\end{enumerate}
For $k = \ell$ we set $B_{0,-j,n-\ell} = \{ \emptyset\}$ and
define $m_{\emptyset;j,k}$ by the same formula as in (1),(2) where the
last product is understood as $1$.
We extend the definition of $m_{T;j,k}$ for all $j\in\mathbb Z$ so
that $m_{T;j+\ell,k} = \tau_{2n} m_{T;j,k}$. 

We describe the action of $\tilde{e}_0$, $\tilde{f}_0$. We get that $\tilde{e}_0(m_{T;j,k})$ is equal to
\begin{equation*}
\begin{split} 
        \begin{cases}
        m_{(i_2,\cdots,i_{\ell-k});j,k+1}&\text{ if $i_1=1$ and $i_{\ell-k}\neq \overline{1}$,}
                \\ m_{(i_1,\cdots,i_{\ell-k},\overline{1});j+1,k-1}&\text{ if $i_1\neq 1$, $i_{\ell-k}\neq \overline{1}$ and $k>0$},
                \\0 &\text{ otherwise,}
  \end{cases} 
\end{split}
\end{equation*}
and that $\tilde{f}_0(m_{T;j,k})$ is equal to
\begin{equation*}
\begin{split} 
        \begin{cases}
        m_{(i_1,\cdots,i_{\ell-k-1});j-1,k + 1}&\text{ if $i_1 \neq 1$ and $i_{\ell-k} = \overline{1}$,}
                \\ m_{(1,i_1,\cdots,i_{\ell-k});j,k-1}&\text{ if $i_1 \neq 1$, $i_{\ell-k} \neq \overline{1}$ and $k>0$,}
                \\0 &\text{ otherwise.}
  \end{cases} 
\end{split}
\end{equation*}
We have $\tau_{2n} = (z_\ell)^{-\ell}$ and all
monomials in $\mathcal M(M_{0,0})/\tau_{2n}$ are written as $m_{T;j,k}$. The crystal
automorphism $z_\ell$ is given by
$\tau_{\ell - k,\ell-j-k-1,n-\ell}^{-1}$.
\begin{NB}
We use the same formula to define $m_{T;\ell,k} = \tau_{2n}(m_{T ; 0, k})$. For $0\le j\le \ell, 0\le k\le\ell$ we define $M_{j,k}=m_{(1,\cdots,\ell-k);j,k}$. It is clear that $z_\ell$ is given by
$\tau_{\ell - k,\ell-j-k-1,n-\ell}^{-1}\colon \mathcal M_{I_0}(M_{j,k})\to 
\mathcal M_{I_0}(M_{j+1,k})$ (for $0\le j <\ell$, $0\le k \le \ell$).
\end{NB}
As an application, we have
\begin{equation*}
    \mathcal B(W(\varpi_\ell))\simeq
    \mathcal B_{I_0}(\varpi_\ell)\sqcup 
    \mathcal B_{I_0}(\varpi_{\ell-1})\sqcup
    \cdots\sqcup \mathcal B_{I_0}(\varpi_1)
    \sqcup \mathcal B_{I_0}(0).
\end{equation*}
\begin{NB}
\begin{equation*}
    \mathcal B(W(\varpi_\ell))\simeq\mathcal M(M_{0,0})/z_\ell=\mathcal M_{I_0}(M_{0,0})\sqcup \mathcal M_{I_0}(M_{0,1})\sqcup        \cdots\sqcup \mathcal M_{I_0}(M_{0,\ell}).
\end{equation*}
\end{NB}


A conjectural description of the crystal of $\mathcal B(W(\varpi_\ell))$ was
proposed in \cite{oss}. As their description is given by relating the
crystal to an $A^{(1)}_{2n+1}$-crystal, it is not clear, at least to 
authors, whether their conjecture is true or not.

\subsubsection{} Finally we consider the case $\ell=n$. Let $M =
Y_{n,0}Y_{0,n}^{-1}$. It follows from Corollary~\ref{corex} that
$\mathcal M(M) \simeq \mathcal B(\varpi_\ell)$.
\begin{NB}
 As $\tilde{e}_n\tilde{e}_{n-1}\tilde{e}_{n-2}\cdots\tilde{e}_1\tilde{e}_0 M = \tau_{-2}(M)$, $\mathcal M(M)$ is preserved under $\tau_2$. For weight reason $z_\ell=\tau_{-2}$ and so $\mathcal M(M)/\tau_2\simeq \mathcal B(W(\varpi_\ell))$. 
\end{NB}

Let
\begin{equation*}
\begin{split}
   \fhbox{i}_p & = 
   \begin{cases}
     Y_{i-1,p+{i-1}}^{-1} Y_{i,p+{i-2}} 
       & \text{if $1\le i\le n-1$},
   \\
     Y_{n-1,p+{n-1}}^{-1} & \text{if $i=n$},
   \\
     Y_{n,p+n} & \text{if $i=0$},
   \end{cases}
\\
   \fhbox{\overline{i}}_p &= 
   \begin{cases}
   Y_{0,p+2n} & \text{if $i=1$},
   \\
     1 & \text{if $2\le i\le n-1$},
   \\
   Y_{n,p+{n+2}}^{-2} 
     & \text{if $i=n$}.
   \end{cases}
\end{split}
\end{equation*}
Then the monomials appearing in $\mathcal M_{I_0}(M)$ are 
$m_T=\prod_{a=1}^{n+1} \fhbox{i_a}_{n+2-2a}$ associated with a tableau
$T = (i_1,\dots,i_{n+1})$ satisfying the conditions
\begin{enumerate}
    \item $i_a\in \mathbf B, \; i_1 \prec i_2 \prec \dots \prec i_{n+1}$,
    \item $i$ and $\overline{i}$ do not appear simultaneously.
\end{enumerate}

We describe the action of $\tilde{e}_0$, $\tilde{f}_0$ on these monomials. We have
\begin{equation*}
\begin{split}
\tilde{e}_0(m_{T}) &=
        \begin{cases}
         \tau_{-2}(m_{(i_2,\cdots,i_{n+1},\overline{1})}) &\text{ if $i_1 = 1$},
                \\  0&\text{ otherwise,}
  \end{cases}
\\
\tilde{f}_0(m_{T}) &=
        \begin{cases}
        \tau_2(m_{(1,i_1,\cdots,i_{n})}) &\text{ if $i_{n+1}=\overline{1}$,}
                \\  0&\text{ otherwise.}
  \end{cases} 
\end{split}
\end{equation*}
We have $\tau_2 = (z_\ell)^{-1}$ and the above monomials are those
appearing in $\mathcal M(M)/\tau_2$. As an application, we have
$\mathcal B(W(\varpi_n)) \simeq \mathcal B_{I_0}(\varpi_n)$. 


\section{Finite dimensional crystals -- exceptional types}

In this section we treat all exceptional cases (except some nodes of type $E_7^{(1)}$, $E_8^{(1)}$, and for one node of type $E_6^{(2)}$ where we do not get the decomposition in $I_0$-crystals at this moment). We
enumerate the nodes of the Dynkin diagram as explained in section \ref{numdyn}.

\subsection{Type $E_n^{(1)}$}

Recall $V_{I_0}(\lambda)$ denotes the irreducible
$\U_q(\Glie_{I_0})$-module with the highest weight $\lambda$. To save
the space, we write $i_p$ instead of $Y_{i,p}$ in some places.

\subsubsection{}

Let $\ell$ be a nonzero vertex with $a_\ell = 1$, i.e.,  $\ell = 1$ or
$5$ for $E_6^{(1)}$ and $\ell = 6$ for $E_7^{(1)}$. In these cases it
is known that the corresponding level $0$ fundamental representation
$W(\varpi_\ell)$ is restricted to the irreducible $\mathcal U_q(\mathfrak
g_{I_0})$-module $V_{I_0}(\varpi_\ell)$. Let us consider $\mathcal M(M)$
for $M = Y_{\ell,0} Y_{0,\theta_\ell}^{-1}$ where $\theta_\ell$ is the
distance of $0$ and $\ell$. By Corollary~\ref{corex} we have $\mathcal
M(M) \simeq \mathcal B(\varpi_\ell)$. Moreover an explicit calculation
shows that $Y_{\ell,p} Y_{0,\theta_\ell+p}^{-1} = \tau_p(M)$ appears
in $\mathcal M(M)$ where $p=6$ for $E_6^{(1)}$ and $p=8$ for
$E_7^{(1)}$. By the weight calculation 
we have $z_{\ell} = \tau_{-p}$. Hence  $\mathcal M(M)/\tau_p \simeq \mathcal
B(W(\varpi_\ell))$. We can check that all monomials are connected to
some $\tau_p^N(M)$ in the $I_0$-crystal. This recovers the above
mentioned result that $W(\varpi_\ell)$ is restricted to $V_{I_0}(\varpi_\ell)$.

\begin{NB}
Let us explain the $E_7^{(1)}$ case for an illustration. Let $m =
Y_{6,0} Y_{0,6}^{-1}$. Then a calculation shows $Y_{1,13}^{-1} Y_{6,8}
Y_{0,12}$ in $\mathcal B_{I_0}(m)$. Applying $\tilde f_0$, we get
$Y_{6,8} Y_{0,14}^{-1}$.
\end{NB}

Let us explain the $E_6^{(1)}$ case for an illustration. Let $M =
Y_{5,0} Y_{0,4}^{-1}$. Then a calculation shows that the following
27 monomials appear in $\mathcal B_{I_0}(M)$:

$5_0 0_4^{-1}$,
$4_1 5_2^{-1} 0_4^{-1}$,
$3_2 4_3^{-1} 0_4^{-1}$,
$6_3 2_3 3_4^{-1} 0_4^{-1}$,
$6_5^{-1} 2_3$,
$1_4 6_3 2_5^{-1} 0_4^{-1}$,
$1_4 6_5^{-1} 2_5^{-1} 3_4$,
$1_6^{-1} 6_3 0_4^{-1}$,
$1_6^{-1} 6_5^{-1} 3_4$,
$1_4 3_6^{-1} 4_5$,
$1_6^{-1} 2_5 3_6^{-1} 4_5$,
$1_4 4_7^{-1} 5_6$,
$2_7^{-1} 4_5$,
$1_6^{-1} 2_5 4_7^{-1} 5_6$,
$1_4 5_8^{-1}$,
$2_7^{-1} 3_6 4_7^{-1} 5_6$,
$1_6^{-1} 2_5 5_8^{-1}$,
$6_7 3_8^{-1} 5_6$,
$2_7^{-1} 3_6 5_8^{-1}$,
$6_9^{-1} 5_6 0_8$,
$6_7 3_8^{-1} 4_7 5_8^{-1}$,
$6_9^{-1} 4_7 5_8^{-1} 0_8$,
$6_7 4_9^{-1}$,
$6_9^{-1} 3_8 4_9^{-1} 0_8$,
$2_9 3_{10}^{-1} 0_8$,
$1_{10} 2_{11}^{-1} 0_8$,
$1_{12}^{-1} 0_8$.

Applying $\tilde f_0$ to $6_9^{-1} 5_6 0_8$, we get $5_6 0_{10}^{-1}
=\tau_6(M)$. It is also clear that all monomials are connected to
either $M$ or its $\tau_6$-images in the $I_0$-crystal.

\begin{rem}
  (1) For a level $0$ fundamental representation $W(\varpi_\ell)$ the
  corresponding quiver varieties are moduli spaces of vector bundles
  of rank $a_\ell$ on ALE spaces. In particular, they are moduli
  spaces of line bundles for the cases studied here. Then each
  component is a single point, and it is a geometric reason why
  $W(\varpi_\ell)$ is restricted to the irreducible representation of
  $\mathcal U_q(\mathfrak g_{I_0})$.

(2) This crystal has been studied in \cite{Magyar}.
\end{rem}

\subsubsection{}\label{ALE}

Let $\ell$ be the vertex adjacent to the vertex $0$, i.e., $\ell = 6$
for $E_6^{(1)}$, $1$ for $E_7^{(1)}$ and $E_8^{(1)}$. We have $a_\ell
= 2$. It is known that $W(\varpi_\ell)$ is restricted to the direct sum of
the adjoint representation $V_{I_0}(\varpi_\ell)$ and the trivial
representation $V_{I_0}(0)$ of $\mathcal U_q(\mathfrak g_{I_0})$. We can
check this, for example, by using the algorithm for the $t$--analog of
$q$--characters \cite{Naa}. All the coefficients of monomials are $1$
except one, whose coefficient is $1+t^2$. The exceptional monomial is
$Y_{3,5} Y_{3,7}^{-1}$ for $E_6^{(1)}$, $Y_{3,8} Y_{3,10}^{-1}$ for
$E_7^{(1)}$ and $Y_{5,14} Y_{5,16}^{-1}$ for $E_8^{(1)}$, if the
$l$--highest weight monomial is $Y_{\ell,0}$.

Let $M = Y_{\ell,0} Y_{0,1}^{-1} Y_{0,p}^{-1}$ where $p=5$ for
$E_6^{(1)}$, $7$ for $E_7^{(1)}$ and $11$ for $E_8^{(1)}$. 
We have 
\begin{gather*}
   E_6 : \quad \tilde f_3 \tilde f_6 M
   =Y_{4,2}Y_{2,2}Y_{3,3}^{-1}Y_{0,5}^{-1},
\\
   E_7 : \quad \tilde f_3 \tilde f_2 \tilde f_1 M =Y_{7,3}Y_{4,3}Y_{3,4}^{-1}Y_{0,7}^{-1},
\\
   E_8 : \quad \tilde f_5 \tilde f_4 \tilde f_3 \tilde f_2 \tilde f_1
   M
   =Y_{6,5}Y_{8,5}Y_{5,6}^{-1}Y_{0,11}^{-1}.
\end{gather*}
By the same argument as in the proof of Proposition~\ref{otherd}, we
see that $M$ is extremal. Therefore $\mathcal M(M) \simeq \mathcal
B(\varpi_\ell)$.
A direct calculation shows that the monomial corresponding to the lowest
weight vector in the adjoint representation is $m =
Y_{\ell,h^\vee}^{-1} Y_{0,h^\vee-p} Y_{0,h^\vee-1}$ where $h^\vee$ is
the dual Coxeter number, i.e., $h^\vee = 12$ for $E_6^{(1)}$, $18$ for
$E_7^{(1)}$ and $30$ for $E_8^{(1)}$. Applying $\tilde f_0$ to $m$,
we get $Y_{0,h^\vee-p} Y_{0,h^\vee+1}^{-1}$, which corresponds to the
trivial representation. 

We have 
\(
  \tilde f_0 \tilde e_{\ell} m 
  = Y_{0,h^\vee-p+2}^{-1} Y_{\ell, h^\vee-p+1} Y_{\ell,h^\vee-2}
  Y_{i,h^\vee-1}^{-1},
\)
where $i$ is the vertex adjacent to $\ell$ different from $0$. A direct
calculation shows that 
\[
  Y_{0,(h^\vee-p+3)/2}^{-1}
  Y_{\ell,(h^\vee-p+1)/2} Y_{\ell,(h^\vee+p-5)/2} Y_{i,(h^\vee+p-3)/2}^{-1}
  = \tau_{-(h^\vee-p+1)/2}(\tilde f_0 \tilde e_{\ell} m) 
\]
is in $\mathcal M_{I_0}(M)$. (Note that $(h^\vee-p+1)/2$ is $4$ for
$E_6^{(1)}$, $6$ for $E_7^{(1)}$ and $10$ for $E_8^{(1)}$.)  Therefore
$\tau_{(h^\vee-p+1)/2}(M)$ is contained in $\mathcal M(M)$. The weight
of $\tau_{(h^\vee-p+1)/2}(M)$ is equal to $\wt(M) - \delta$. Therefore
this is $z_\ell^{-1}(M)$ and we have $\tau_{(h^\vee-p+1)/2} =
z_\ell^{-1}$ and $\mathcal M(M)/\tau_{(h^\vee-p+1)/2} \simeq \mathcal
B(W(\varpi_\ell))$.

We can also check that 
\(
   \mathcal M(M)/ \tau_{(h^\vee-p+1)/2} \simeq
   \mathcal M_{I_0}(M)
   \sqcup\{
   Y_{0,h^\vee-p} Y_{0,h^\vee+1}^{-1}\}.
\)
Therefore we recover that $W(\varpi_\ell)$ is restricted to
$V_{I_0}(\varpi_\ell)\oplus V_{I_0}(0)$.

Let us give $E_6^{(1)}$ case for an illustration. The following
monomials appear in $\mathcal B_{I_0}(6_0 0_1^{-1} 0_5^{-1})$:

$6_0 0_1^{-1} 0_5^{-1}$,
$6_2^{-1} 3_1 0_5^{-1}$,
$2_2 3_3^{-1} 4_2 0_5^{-1}$,
$1_3 2_4^{-1} 4_2 0_5^{-1}$,
$2_2 4_4^{-1} 5_3 0_5^{-1}$,
$1_5^{-1} 4_2 0_5^{-1}$,
$1_3 2_4^{-1} 3_3 4_4^{-1} 5_3 0_5^{-1}$,
$2_2 5_5^{-1} 0_5^{-1}$,
$1_5^{-1} 3_3 4_4^{-1} 5_3 0_5^{-1}$,
$1_3 6_4 3_5^{-1} 5_3 0_5^{-1}$,
$1_3 2_4^{-1} 3_3 5_5^{-1} 0_5^{-1}$,
$1_5^{-1} 6_4 2_4 3_5^{-1} 5_3 0_5^{-1}$,
$1_5^{-1} 3_3 5_5^{-1} 0_5^{-1}$,
$1_3 6_6^{-1} 5_3$,
$1_3 6_4 3_5^{-1} 4_4 5_5^{-1} 0_5^{-1}$,
$1_5^{-1} 6_6^{-1} 2_4 5_3$,
$6_4 2_6^{-1} 5_3 0_5^{-1}$,
$1_5^{-1} 6_4 2_4 3_5^{-1} 4_4 5_5^{-1} 0_5^{-1}$,
$1_3 6_6^{-1} 4_4 5_5^{-1}$,
$1_3 6_4 4_6^{-1} 0_5^{-1}$,
$6_6^{-1} 2_6^{-1} 3_5 5_3$,
$1_5^{-1} 6_6^{-1} 2_4 4_4 5_5^{-1}$,
$6_4 2_6^{-1} 4_4 5_5^{-1} 0_5^{-1}$,
$1_5^{-1} 6_4 2_4 4_6^{-1} 0_5^{-1}$,
$1_3 6_6^{-1} 3_5 4_6^{-1}$,
$3_7^{-1} 4_6 5_3$,
$6_6^{-1} 2_6^{-1} 3_5 4_4 5_5^{-1}$,
$1_5^{-1} 6_6^{-1} 2_4 3_5 4_6^{-1}$,
$6_4 2_6^{-1} 3_5 4_6^{-1} 0_5^{-1}$,
$1_3 2_6 3_7^{-1}$,
$4_8^{-1} 5_3 5_7$,
$3_7^{-1} 4_4 4_6 5_5^{-1}$,
$6_6^{-1} 2_6^{-1} 3_5^2 4_6^{-1}$,
$1_5^{-1} 2_4 2_6 3_7^{-1}$,
$6_4 6_6 3_7^{-1} 0_5^{-1}$,
$1_3 1_7 2_8^{-1}$,
$5_3 5_9^{-1}$,
$4_4 4_8^{-1} 5_5^{-1} 5_7$,
$3_5 3_7^{-1}$,
$1_5^{-1} 1_7 2_4 2_8^{-1}$,
$6_4 6_8^{-1} 0_5^{-1}$,
$1_3 1_9^{-1}$,
$4_4 5_5^{-1} 5_9^{-1}$,
$3_5 4_6^{-1} 4_8^{-1} 5_7$,
$6_6 2_6 3_7^{-2} 4_6$,
$1_7 2_6^{-1} 2_8^{-1} 3_5$,
$6_6^{-1} 6_8^{-1} 3_5 0_7$,
$1_5^{-1} 1_9^{-1} 2_4$,
$3_5 4_6^{-1} 5_9^{-1}$,
$6_6 2_6 3_7^{-1} 4_8^{-1} 5_7$,
$6_8^{-1} 2_6 3_7^{-1} 4_6 0_7$,
$1_7 6_6 2_8^{-1} 3_7^{-1} 4_6$,
$1_9^{-1} 2_6^{-1} 3_5$,
$6_6 2_6 3_7^{-1} 5_9^{-1}$,
$6_8^{-1} 2_6 4_8^{-1} 5_7 0_7$,
$1_7 6_8^{-1} 2_8^{-1} 4_6 0_7$,
$1_7 6_6 2_8^{-1} 4_8^{-1} 5_7$,
$1_9^{-1} 6_6 3_7^{-1} 4_6$,
$6_8^{-1} 2_6 5_9^{-1} 0_7$,
$1_7 6_8^{-1} 2_8^{-1} 3_7 4_8^{-1} 5_7 0_7$,
$1_9^{-1} 6_8^{-1} 4_6 0_7$,
$1_7 6_6 2_8^{-1} 5_9^{-1}$,
$1_9^{-1} 6_6 4_8^{-1} 5_7$,
$1_7 6_8^{-1} 2_8^{-1} 3_7 5_9^{-1} 0_7$,
$1_9^{-1} 6_8^{-1} 3_7 4_8^{-1} 5_7 0_7$,
$1_7 3_9^{-1} 5_7 0_7$,
$1_9^{-1} 6_6 5_9^{-1}$,
$1_9^{-1} 6_8^{-1} 3_7 5_9^{-1} 0_7$,
$1_9^{-1} 2_8 3_9^{-1} 5_7 0_7$,
$1_7 3_9^{-1} 4_8 5_9^{-1} 0_7$,
$1_9^{-1} 2_8 3_9^{-1} 4_8 5_9^{-1} 0_7$,
$2_{10}^{-1} 5_7 0_7$,
$1_7 4_{10}^{-1} 0_7$,
$2_{10}^{-1} 4_8 5_9^{-1} 0_7$,
$1_9^{-1} 2_8 4_{10}^{-1} 0_7$,
$2_{10}^{-1} 3_9 4_{10}^{-1} 0_7$,
$6_{10} 3_{11}^{-1} 0_7$,
$6_{12}^{-1} 0_7 0_{11}$.

There is $6_4 6_6 3_7^{-1} 0_5^{-1}$ as claimed. We can also check
that all monomials are connected to either $M$, $0_7 0_{13}^{-1}$ or
their $\tau_4$-images in the $I_0$-crystal. 


\begin{rem}
(1) In this example, the corresponding quiver varieties are either a
single point or an ALE space of type $E_n$. The graded quiver
varieties, which are fixed point sets of a $\mathbb C^*$-action, are
single points or a complex projective line. The latter gives the
monomial with coefficient $1+t^2$.

(2) The crystal structure here is isomorphic to one studied recently
in \cite{bfkl}. As the crystal graph is connected, we conclude that
the crystal base constructed in \cite{bfkl} are isomorphic to
$\mathcal B(W(\varpi_\ell))$.
\end{rem}

\subsubsection{}

Let $\Glie = E_6^{(1)}$ and $\ell = 2$. The $t$--analog of
$q$--character of $W(\varpi_2)$ has 351 monomials among which the
following 27 monomials have coefficients $1+t^2$ and others have $1$:

$3_3 3_5^{-1} 5_3$,
$3_3 3_5^{-1} 4_4 5_5^{-1}$,
$3_3 4_6^{-1}$,
$6_4 2_4 3_5^{-1} 4_4 4_6^{-1}$,
$6_6^{-1} 2_4 4_4 4_6^{-1}$,
$1_5 6_4 2_6^{-1} 4_4 4_6^{-1}$,
$1_5 6_6^{-1} 2_6^{-1} 3_5 4_4 4_6^{-1}$,
$1_7^{-1} 6_4 4_4 4_6^{-1}$,
$1_7^{-1} 6_6^{-1} 3_5 4_4 4_6^{-1}$,
$1_5 3_7^{-1} 4_4$,
$1_7^{-1} 2_6 3_7^{-1} 4_4$,
$1_5 3_5 3_7^{-1} 4_6^{-1} 5_5$,
$2_8^{-1} 4_4$,
$1_7^{-1} 2_6 3_5 3_7^{-1} 4_6^{-1} 5_5$,
$1_5 3_5 3_7^{-1} 5_7^{-1}$,
$2_8^{-1} 3_5 4_6^{-1} 5_5$,
$1_7^{-1} 2_6 3_5 3_7^{-1} 5_7^{-1}$,
$2_8^{-1} 3_5 5_7^{-1}$,
$6_6 2_6 2_8^{-1} 3_7^{-1} 5_5$,
$6_6 2_6 2_8^{-1} 3_7^{-1} 4_6 5_7^{-1}$,
$6_8^{-1} 2_6 2_8^{-1} 5_5$,
$6_8^{-1} 2_6 2_8^{-1} 4_6 5_7^{-1}$,
$6_6 2_6 2_8^{-1} 4_8^{-1}$,
$6_8^{-1} 2_6 2_8^{-1} 3_7 4_8^{-1}$,
$2_6 3_9^{-1}$,
$1_7 2_8^{-1} 3_7 3_9^{-1}$,
$1_9^{-1} 3_7 3_9^{-1}$.

From this (or by other methods) we can see that $W(\varpi_2)$ is
restricted to $V_{I_0}(\varpi_2)\oplus V_{I_0}(\varpi_5)$.

Let us consider the monomial crystal $\mathcal M(M)$ with $M = 2_0
0_3^{-1} 0_5^{-1}$.
From $\tilde f_6\tilde f_3 \tilde f_2 M = 1_1 4_2 6_4^{-1} 0_5^{-1}$,
we see that $M$ is extremal by the argument in the proof of
Proposition~\ref{otherd}. Therefore $\mathcal M(M) \simeq \mathcal
B(\varpi_2)$.

There is a monomial 
\[
   m = 1_5 6_6^{-1} 6_8^{-1} 4_4 0_7
   = \tilde f_6 \tilde f_6 \tilde f_3 \tilde f_2
   \tilde f_3 \tilde f_2 \tilde f_1 \tilde f_4 \tilde f_5
   \tilde f_3 \tilde f_4 \tilde f_6 \tilde f_3 \tilde f_2 M
\]
in $\mathcal M_{I_0}(M)$. We have $\tilde e_2\tilde e_3 \tilde
e_6\tilde f_0 m = 1_3^{-1} 1_5 2_2 0_5^{-1} 0_9^{-1}$. By the weight
calculation, we find that this is $z_\ell^{-1}(M)$. Let us denote this
by $M_1$.

In $\mathcal M_{I_0}(M_1)$ we can find a monomial 
\[
   m' = 1_7 6_{10}^{-2} 4_8 0_9
   = \tilde f_6 \tilde f_6 \tilde f_3 \tilde f_2
   \tilde f_3 \tilde f_2 \tilde f_1 \tilde f_4 \tilde f_5
   \tilde f_3 \tilde f_4 \tilde f_6 \tilde f_3 \tilde f_2 M_1.
\]
We have $\tilde e_2\tilde e_3 \tilde e_6\tilde f_0 m' = 
2_6 0_9^{-1} 0_{11}^{-1} = \tau_6(M)$. This is equal to $z_\ell^{-2}(M)$.

We have
\[
   5_3 0_{11}^{-1}
   = \tilde e_5 \tilde e_4 \tilde e_3 \tilde e_6 \tilde e_0 \cdot \tau_6(M)
\]
in $\mathcal M(M)$. Write this $M_{0;1}$. Then $\mathcal
M_{I_0}(M_{0;1})$ consists of the following 27 monomials: 

$5_3 0_{11}^{-1}$,
$4_4 5_5^{-1} 0_{11}^{-1}$,
$3_5 4_6^{-1} 0_{11}^{-1}$,
$6_6 2_6 3_7^{-1} 0_{11}^{-1}$,
$6_8^{-1} 2_6 0_7 0_{11}^{-1}$,
$1_7 6_6 2_8^{-1} 0_{11}^{-1}$,
$1_7 6_8^{-1} 2_8^{-1} 3_7 0_7 0_{11}^{-1}$,
$1_9^{-1} 6_6 0_{11}^{-1}$,
$1_9^{-1} 6_8^{-1} 3_7 0_7 0_{11}^{-1}$,
$1_7 3_9^{-1} 4_8 0_7 0_{11}^{-1}$,
$1_9^{-1} 2_8 3_9^{-1} 4_8 0_7 0_{11}^{-1}$,
$1_7 4_{10}^{-1} 5_9 0_7 0_{11}^{-1}$,
$2_{10}^{-1} 4_8 0_7 0_{11}^{-1}$,
$1_9^{-1} 2_8 4_{10}^{-1} 5_9 0_7 0_{11}^{-1}$,
$1_7 5_{11}^{-1} 0_7 0_{11}^{-1}$,
$2_{10}^{-1} 3_9 4_{10}^{-1} 5_9 0_7 0_{11}^{-1}$,
$1_9^{-1} 2_8 5_{11}^{-1} 0_7 0_{11}^{-1}$,
$2_{10}^{-1} 3_9 5_{11}^{-1} 0_7 0_{11}^{-1}$,
$6_{10} 3_{11}^{-1} 5_9 0_7 0_{11}^{-1}$,
$6_{10} 3_{11}^{-1} 4_{10} 5_{11}^{-1} 0_7 0_{11}^{-1}$,
$6_{12}^{-1} 5_9 0_7 0_{11}^{-1}$,
$6_{12}^{-1} 4_{10} 5_{11}^{-1} 0_7$,
$6_{10} 4_{12}^{-1} 0_7 0_{11}^{-1}$,
$6_{12}^{-1} 3_{11} 4_{12}^{-1} 0_7$,
$2_{12} 3_{13}^{-1} 0_7$,
$1_{13} 2_{14}^{-1} 0_7$,
$1_{15}^{-1} 0_7$.

We have 
\(
  \tilde e_5 \tilde e_4 \tilde e_3 \tilde e_6 \tilde e_0 \cdot M_1
  = 1_3^{-1} 1_5 5_{-1} 0_9^{-1}.
\)
Set this $M_{1;1}$. Then $\mathcal M_{I_0}(M_{1;1})$ consists of

$1_3^{-1} 1_5 5_{-1} 0_9^{-1}$,
$1_3^{-1} 1_5 4_0 5_1^{-1} 0_9^{-1}$,
$1_3^{-1} 1_5 3_1 4_2^{-1} 0_9^{-1}$,
$1_3^{-1} 1_5 2_2 3_3^{-1} 6_2 0_9^{-1}$,
$1_3^{-1} 1_5 2_2 6_4^{-1} 0_9^{-1}$,
$1_5 2_4^{-1} 6_2 0_9^{-1}$,
$1_5 2_4^{-1} 3_3 6_4^{-1} 0_3 0_9^{-1}$,
$1_7^{-1} 2_4^{-1} 2_6 6_2 0_9^{-1}$,
$1_7^{-1} 2_4^{-1} 2_6 3_3 6_4^{-1} 0_3 0_9^{-1}$,
$1_5 3_5^{-1} 4_4 0_3 0_9^{-1}$,
$1_7^{-1} 2_6 3_5^{-1} 4_4 0_3 0_9^{-1}$,
$1_5 4_6^{-1} 5_5 0_3 0_9^{-1}$,
$2_8^{-1} 3_5^{-1} 3_7 4_4 0_3 0_9^{-1}$,
$1_7^{-1} 2_6 4_6^{-1} 5_5 0_3 0_9^{-1}$,
$1_5 5_7^{-1} 0_3 0_9^{-1}$,
$2_8^{-1} 3_7 4_6^{-1} 5_5 0_3 0_9^{-1}$,
$1_7^{-1} 2_6 5_7^{-1} 0_3 0_9^{-1}$,
$2_8^{-1} 3_7 5_7^{-1} 0_3 0_9^{-1}$,
$3_9^{-1} 4_6^{-1} 4_8 5_5 6_8 0_3 0_9^{-1}$,
$3_9^{-1} 4_8 5_7^{-1} 6_8 0_3 0_9^{-1}$,
$4_6^{-1} 4_8 5_5 6_{10}^{-1} 0_3$,
$4_8 5_7^{-1} 6_{10}^{-1} 0_3$,
$4_{10}^{-1} 5_7^{-1} 5_9 6_8 0_3 0_9^{-1}$,
$3_9 4_{10}^{-1} 5_7^{-1} 5_9 6_{10}^{-1} 0_3$,
$2_{10} 3_{11}^{-1} 5_7^{-1} 5_9 0_3$,
$1_{11} 2_{12}^{-1} 5_7^{-1} 5_9 0_3$,
$1_{13}^{-1} 5_7^{-1} 5_9 0_3$.

These have different weights, so there is only one way to make a
bijection to the above polynomials with coefficients $1+t^2$
preserving weights. It is the bijection given in order.

Also it should be possible to make the bijection between $\mathcal
M_{I_0}(M)$ and $\mathcal M_{I_0}(M_1)$ {\it explicit\/}, though we do
not do here, as both are 351 monomials.

Thus we have
\begin{equation*}
   \mathcal M(M)/\tau_6 \simeq 
   \mathcal M_{I_0}(M) \sqcup \mathcal M_{I_0}(M_1)
   \sqcup \mathcal M_{I_0}(M_{0;1}) \sqcup \mathcal M_{I_0}(M_{1;1}),
\end{equation*}
and we have a crystal isomorphism $\tau$ interchanging $\mathcal
M_{I_0}(M) \leftrightarrow \mathcal M_{I_0}(M_1)$ and $\mathcal
M_{I_0}(M_{0;1}) \leftrightarrow \mathcal M_{I_0}(M_{1;1})$. These
follow from the known results, but should be possible to check
directly from the above computation. 


\subsubsection{}\label{last}

Let $\Glie = E_6^{(1)}$ and $\ell = 3$. 
It is known that $W(\varpi_3)$ restricts to $V_{I_0}(\varpi_3)
\oplus V_{I_0}(\varpi_6)^{\oplus 2} \oplus V_{I_0}(\varpi_1+\varpi_5) \oplus
V_{I_0}(0)$ as a $\U_q(\Glie_{I_0})$-module.

\begin{NB}
The $t$--analog of $q$--character of $W(\varpi_3)$ has $2925$
monomials among which one ($6_5 6_7^{-1} 2_5 2_7^{-1} 4_5 4_7^{-1}$)
has the coefficient $1 + 3t^2 + 3t^4 + t^6$, $3$ has
$1 + 4t^2 + t^4$, $12$ has $1 + 3t^2 + t^4$, $44$ has
$1 + 2t^2 + t^4$, $16$ has $1 + t^2 + t^4$, $573$ has $1+t^2$.
\end{NB}

Let $M = 3_0 0_2^{-1} 0_4^{-1} 0_6^{-1}$. We have
\begin{equation*}
   m = \tilde f_6 \tilde f_3^2 \tilde f_6 \tilde f_4 \tilde f_2
   \tilde f_3 M = 1_2 2_3 3_4^{-1} 4_3 5_2 6_5^{-1} 0_6^{-1}.  
\end{equation*}
By the same argument as in the proof of Proposition~\ref{otherd}, we
see that $M$ is extremal. Therefore $\mathcal M(M) \simeq \mathcal
B(\varpi_3)$.

We have
\begin{equation*}
  m' = \tilde{e}_3\tilde{f}_6^2\tilde{f}_3^3\tilde{f}_4^2
  \tilde{f}_2^2\tilde{f}_5\tilde{f}_1 m
  = 1_4 3_4 5_4 6_5^{-1} 6_7^{-2} 0_6
\end{equation*}
in $\mathcal M_{I_0}(M)$. Then
\begin{equation*}
   \tilde e_3 \tilde e_4 \tilde e_2 \tilde e_3 \tilde e_6 \tilde e_6
   \tilde f_0 m' = 3_2 0_4^{-1} 0_6^{-1} 0_8^{-1}
   = \tau_2(M).
\end{equation*}
By the weight calculation, this is $z_\ell^{-1}(M)$, so we have
$z_{\ell}=\tau_{-2}$ and $\mathcal M(M)/\tau_2 \simeq \mathcal
B(W(\varpi_\ell))$. 

Let
\begin{equation*}
   M_1 = \tilde e_6 \tilde e_0 \tau_2(M)
   = 6_1 0_6^{-1} 0_8^{-1}.
\end{equation*}
Then $\mathcal M_{I_0}(M_1)$ is the crystal of the adjoint
representation of $\Glie_{I_0}$. By \ref{ALE} the lowest weight vector
is 
\(
   6_{13}^{-1} 0_8 0_{12} \times 0_2 0_8^{-1} = 
   6_{13}^{-1} 0_2 0_{12}.
\)
Applying $\tau_{-2}\tilde f_0$, we get $M_2 = 0_0 0_{12}^{-1}$. Applying
$\tilde f_0$ again, we get 
\(
   M_3 = 6_1 0_2^{-1} 0_{12}^{-1}.
\)
Looking at monomials in \ref{ALE}, we find $1_4 6_7^{-1} 5_4 0_6
0_{12}^{-1}$ in $\mathcal M_{I_0}(M_3)$. Applying $\tilde f_0$, we get
\(
   M_4 = 1_4 5_4 0_8^{-1} 0_{12}^{-1}.
\) 
This monomial generates the $I_0$-crystal of $V_{I_0}(\varpi_1+\varpi_5)$.

Thus 
\[
  \mathcal M(M)/\tau_2 = \mathcal M_{I_0}(M)\sqcup
  \mathcal M_{I_0}(M_1) \sqcup \mathcal M_{I_0}(M_2) \sqcup
  \mathcal M_{I_0}(M_3) \sqcup \mathcal M_{I_0}(M_4).
\]
This follows from $\operatorname{Res}W_0(\varpi_3) \simeq
V_{I_0}(\varpi_3) \oplus V_{I_0}(\varpi_6)^{\oplus 2} \oplus
V_{I_0}(\varpi_1+\varpi_5) \oplus V_{I_0}(0)$, but it is probably possible to
check directly from the above computation.


\begin{rem}
  The authors do not find the last two examples in the literature. One
  can probably check their {\it perfectness}, though we have not done yet.
\end{rem}

\subsection{Type $G_2^{(1)}$} 

\subsubsection{} First we consider $\ell=1$. Let $M=Y_{1,0}Y_{0,1}^{-1}Y_{0,3}^{-1}$. As $\tilde{f}_1(M)=Y_{1,2}^{-1}Y_{2,1}^3Y_{0,3}^{-1}$, we see as in Proposition~\ref{otherd} that $\mathcal M(M)\simeq \mathcal B(\varpi_\ell)$. 

As $\tilde{e}_1\tilde{e}_2^3\tilde{e}_1^2\tilde{e}_0\tilde{f}_1 M=\tau_{-2}(M)$, $\mathcal M(M)$ is preserved under $\tau_{-2}$. It has weight $\delta$, so $z_\ell=\tau_{-2}$ and hence $\mathcal M(m)/\tau_2\simeq \mathcal B(W(\varpi_\ell))$.

Let $M'=\tilde{e}_0 (M)=Y_{0,-1}Y_{0,3}^{-1}$. We have $\mathcal M_{I_0}(M')=\{M'\}$. The following $14$ monomials appear in $\mathcal{M}_{I_0}(M)$:

$1_00_1^{-1}0_3^{-1}$,
$2_1^31_2^{-1}0_3^{-1}$,
$2_1^22_3^{-1}0_3^{-1}$,
$2_12_3^{-2}1_20_3^{-1}$,
$2_3^{-3}1_2^20_3^{-1}$,
$2_12_31_4^{-1}$,
$1_21_4^{-1}$,
$2_12_5^{-1}$,
$1_4^{-2}2_3^30_3$,
$2_3^{-1}2_5^{-1}1_2$,
$2_3^22_5^{-1}1_4^{-1}0_3$,
$2_32_5^{-2}0_3$,
$2_5^{-3}0_31_4$,
$0_30_51_6^{-1}$.

By direct calculation, we find that these $14$ monomials and $M'$ are all monomials of $\mathcal M(M)/\tau_2$. As an application, we get
\begin{equation*}
    \mathcal B(W(\varpi_\ell)) \simeq
    \mathcal B_{I_0}(\varpi_\ell)\sqcup\mathcal B_{I_0}(0).
\end{equation*}


This crystal was described in \cite{Yamane,bfkl}. The crystal base
is isomorphic to ours by the same reason as in \ref{ALE}.

\subsubsection{} Now we consider the case $\ell=2$. Let $M=Y_{2,0}Y_{0,2}^{-1}$. It follows from Corollary~\ref{corex} that $\mathcal M(M)\simeq \mathcal B(\varpi_\ell)$. 
The following 7 monomials appear in $\mathcal{M}_{I_0}(M)$:

$M=2_00_2^{-1}$,
$m_2=1_12_2^{-1}0_2^{-1}$,
$m_3=1_3^{-1}2_2^2$,
$m_4=2_22_4^{-1}$,
$m_5=2_4^{-2}1_3$,
$m_6=1_5^{-1}2_40_4$,
$m_7=2_6^{-1}0_4$.

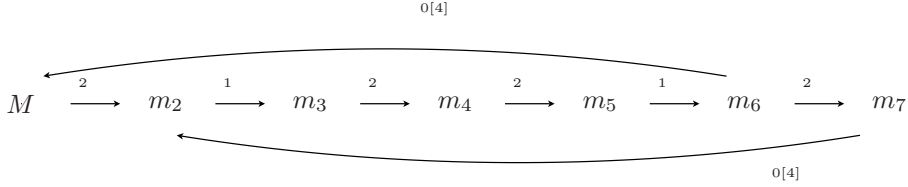
\begin{figure}[htbp]
\centering
\psset{xunit=.55mm,yunit=.55mm,runit=.55mm}
\psset{linewidth=0.3,dotsep=1,hatchwidth=0.3,hatchsep=1.5,shadowsize=1}
\psset{dotsize=0.7 2.5,dotscale=1 1,fillcolor=black}
\psset{arrowsize=1 2,arrowlength=1,arrowinset=0.25,tbarsize=0.7 5,bracketlength=0.15,rbracketlength=0.15}
\begin{pspicture}(0,0)(220,50)
\rput(0,20){$M$}
\psline{->}(12,20)(24,20)
\rput(35,20){$m_2$}
\psline{->}(47,20)(59,20)
\rput(70,20){$m_3$}
\psline{->}(82,20)(94,20)
\rput(105,20){$m_4$}
\psline{->}(117,20)(129,20)
\rput(140,20){$m_5$}
\psline{->}(152,20)(164,20)
\rput(175,20){$m_6$}
\psline{->}(187,20)(199,20)
\rput(210,20){$m_7$}
\rput{0}(0,20){\qline(-0,0)(0,0)}\rput{0}(88,-486){\parametricplot[arrows=->]{80.85}{99.15}{ t cos 519.26 mul t sin 519.26 mul }}
\rput{0}(120.26,525){\parametricplot[arrows=<-]{-99.15}{-80.85}{ t cos 519.26 mul t sin 519.26 mul }}
\rput(15,25){$\scriptscriptstyle 2$}
\rput(50,25){$\scriptscriptstyle 1$}
\rput(85,25){$\scriptscriptstyle 2$}
\rput(120,25){$\scriptscriptstyle 2$}
\rput(155,25){$\scriptscriptstyle 1$}
\rput(190,25){$\scriptscriptstyle 2$}
\rput(100,43){$\scriptscriptstyle 0[4]$}
\rput(185,3){$\scriptscriptstyle 0[4]$}
\end{pspicture}
\caption{(Type $G_2^{(1)}$) the crystal $\mathcal B(\varpi_2)$}
\label{fig:G2}
\end{figure}

The crystal graph of $\mathcal M(M)$ is given in Figure~\ref{fig:G2}.
We find $z_{\ell}=\tau_{-4}$ and $\mathcal M(M)/\tau_4=\mathcal M_{I_0}(M)$. 

The authors do not find a description of this crystal structure in the
literature (probably because it is not perfect), but one can easily
obtain it from the description of its $I_0$-crystal structure in
\cite{KM}.

\subsection{Type $F_4^{(1)}$} 

\subsubsection{} First let $\ell=1$ and $M=Y_{1,0}Y_{0,1}^{-1}Y_{0,5}^{-1}$. We have $\tilde{f}_2\tilde{f}_1 M=Y_{0,5}^{-1}Y_{2,3}^{-1}Y_{3,2}^2$ and so we see as in Proposition \ref{otherd} that $\mathcal M(M)\simeq \mathcal B(\varpi_\ell)$. 
As $\tilde{e}_1\tilde{e}_2\tilde{e}_3^2\tilde{e}_4^2\tilde{e}_2\tilde{e}_3\tilde{e}_3\tilde{e}_2\tilde{e}_1\tilde{e}_1\tilde{e}_0\tilde{f}_1 M=\tau_{-4}(M)$, $\mathcal M(M)$ is preserved under $\tau_4$, which has weight $\delta$. Therefore we have $z_\ell=\tau_{-4}$ and $\mathcal M(m)/\tau_4\simeq \mathcal B(W(\varpi_\ell))$.

Let $M'=\tilde{e}_0 (M)=Y_{0,-1}Y_{0,5}^{-1}$. We have $\mathcal M_{I_0}(M')=\{M'\}$. 

The following $52$ monomials appear in $\mathcal{M}_{I_0}(M)$:

$1_00_1^{-1}0_5^{-1}$, 
$1_2^{-1}2_10_5^{-1}$, 
$2_3^{-1}3_2^20_5^{-1}$, 
$3_23_4^{-1}4_30_5^{-1}$,
$2_33_4^{-2}4_3^20_5^{-1}$, $3_24_5^{-1}0_5^{-1}$,
$1_42_5^{-1}4_3^20_5^{-1}$, $2_33_4^{-1}4_34_5^{-1}0_5^{-1}$,
$1_6^{-1}4_3^2$, $1_42_5^{-1}4_34_5^{-1}3_40_5^{-1}$, $2_34_5^{-2}0_5^{-1}$,
$1_6^{-1}3_44_34_5^{-1}$, $1_42_5^{-1}4_5^{-2}3_4^20_5^{-1}$, $1_43_6^{-1}4_30_5^{-1}$,
$1_6^{-1}2_53_6^{-1}4_3$, $1_6^{-1}3_4^24_5^{-2}$, $1_43_43_6^{-1}4_5^{-1}0_5^{-1}$,
$2_7^{-1}3_64_3$, $1_6^{-1}2_53_6^{-1}3_44_5^{-1}$, $1_42_53_6^{-2}0_5^{-1}$,
$3_8^{-1}4_74_3$, $2_7^{-1}3_63_44_5^{-1}$, $1_6^{-1}2_5^23_6^{-2}$, $1_41_62_7^{-1}0_5^{-1}$,
$4_9^{-1}4_3$, $3_8^{-1}3_44_5^{-1}4_7$, $2_52_7^{-1}$, $1_41_8^{-1}0_70_5^{-1}$,
$3_44_5^{-1}4_9^{-1}$, $2_53_6^{-1}3_8^{-1}4_7$, $1_62_7^{-2}3_6^2$, $1_6^{-1}1_8^{-1}2_50_7$,
$2_53_6^{-1}4_9^{-1}$, $0_71_8^{-1}2_7^{-1}3_6^2$, $1_62_7^{-1}3_63_8^{-1}4_7$,
$1_62_7^{-1}3_64_9^{-1}$, $1_8^{-1}3_63_8^{-1}4_70_7$, $1_63_8^{-2}4_7^2$,
$1_8^{-1}3_64_9^{-1}0_7$, $1_63_8^{-1}4_74_9^{-1}$, $1_8^{-1}2_73_8^{-2}4_7^20_7$,
$1_8^{-1}2_73_8^{-1}4_74_9^{-1}0_7$, $1_64_9^{-2}$, $2_9^{-1}4_7^20_7$,
$2_9^{-1}3_84_74_9^{-1}0_7$, $1_8^{-1}2_74_9^{-2}0_7$,
$3_{10}^{-1}4_70_7$, $2_9^{-1}3_8^24_9^{-2}0_7$,
$3_83_{10}^{-1}4_9^{-1}0_7$,
$2_93_{10}^{-2}0_7$,
$1_{10}2_{11}^{-1}0_7$,
$1_{12}^{-1}0_70_{11}$.

These $52$ monomials, $M'$ and their $\tau_4$-images are all monomials of $\mathcal M(M)$. As an application we have
\begin{NB}
\begin{equation*}\mathcal M(M)/\tau_4=\mathcal M_{I_0}(M)\sqcup\mathcal M_{I_0}(M').\end{equation*}
\end{NB}
\begin{equation*}
   \mathcal B(W(\varpi_1)) \simeq 
   \mathcal B_{I_0}(\varpi_1)\sqcup\mathcal B_{I_0}(0).
\end{equation*}


The crystal base is isomorphic to one in \cite{bfkl} by the same
reason as in \ref{ALE}.

\subsubsection{} Let us consider $\ell=2$ and $M=Y_{2,0}Y_{0,2}^{-1}Y_{0,4}^{-1}Y_{0,6}^{-1}$. We have $\tilde{f}_1\tilde{f}_2^2\tilde{f}_3^2\tilde{f}_1\tilde{f}_2 M=Y_{4,2}^2Y_{3,3}^2Y_{2,4}^{-1}Y_{1,5}^{-1}Y_{0,6}^{-1}$ and so we see as in Proposition \ref{otherd} that $\mathcal M(M)\simeq \mathcal B(\varpi_\ell)$. 
As
$\tilde{e}_2\tilde{e}_3^2\tilde{e}_2\tilde{e}_1^2\tilde{e}_2\tilde{e}_4^2\tilde{e}_3^4\tilde{e}_2^3\tilde{e}_1^2
\tilde{e}_0\tilde{f}_1\tilde{f}_2^2\tilde{f}_3^2\tilde{f}_1\tilde{f}_2
M=\tau_{-2}(M)$, $\mathcal M(M)$ is preserved under $\tau_2$, which
has weight $\delta$. Therefore $\mathcal M(M)/\tau_2\simeq \mathcal B(W(\varpi_\ell))$.

Let $M_2=\tilde{e}_1\tilde{e}_0
M=Y_{0,4}^{-1}Y_{0,6}^{-1}Y_{1,-1}$. The following $52$ monomials
appear in $\mathcal{M}_{I_0}(M_2)$:

$1_{-1}0_4^{-1}0_6^{-1}$, 
$1_1^{-1}2_00_00_4^{-1}0_6^{-1}$, 
$2_2^{-1}3_1^20_00_4^{-1}0_6^{-1}$, 
$3_13_3^{-1}4_20_00_4^{-1}0_6^{-1}$,
$2_23_3^{-2}4_2^20_00_4^{-1}0_6^{-1}$, $3_14_4^{-1}0_00_4^{-1}0_6^{-1}$,
$1_32_4^{-1}4_2^20_00_4^{-1}0_6^{-1}$, $2_23_3^{-1}4_24_4^{-1}0_00_4^{-1}0_6^{-1}$,
$1_5^{-1}4_2^20_00_6^{-1}$, $1_32_4^{-1}4_24_4^{-1}3_30_5^{-1}0_00_4^{-1}0_6^{-1}$, $2_24_4^{-2}0_00_4^{-1}0_6^{-1}$,
$1_5^{-1}3_34_34_4^{-1}0_00_6^{-1}$, $1_32_4^{-1}4_4^{-2}3_3^20_5^{-1}0_00_4^{-1}0_6^{-1}$, $1_33_5^{-1}4_20_5^{-1}0_00_4^{-1}0_6^{-1}$,
$1_5^{-1}2_43_5^{-1}4_20_00_6^{-1}$, $1_5^{-1}3_3^24_4^{-2}0_00_6^{-1}$, $1_33_33_5^{-1}4_4^{-1}0_00_4^{-1}0_6^{-1}$,
$2_6^{-1}3_54_20_00_6^{-1}$, $1_5^{-1}2_43_5^{-1}3_34_4^{-1}0_00_6^{-1}$, $1_32_43_5^{-2}0_5^{-1}0_00_4^{-1}0_6^{-1}$,
$3_7^{-1}4_64_20_00_6^{-1}$, $2_6^{-1}3_53_34_4^{-1}0_00_6^{-1}$, $1_5^{-1}2_4^23_5^{-2}0_00_6^{-1}$, $1_31_52_6^{-1}0_5^{-1}0_00_4^{-1}0_6^{-1}$,
$4_8^{-1}4_20_00_6^{-1}$, $3_7^{-1}3_34_4^{-1}4_60_00_6^{-1}$, $2_42_6^{-1}0_00_6^{-1}$, $1_31_7^{-1}0_00_6^{-1}$,
$3_34_4^{-1}4_8^{-1}0_00_6^{-1}$, $2_43_5^{-1}3_7^{-1}4_60_00_6^{-1}$, $1_52_6^{-2}3_5^20_00_6^{-1}$, $1_5^{-1}1_7^{-1}2_40_00_6^{-1}$,
$2_43_5^{-1}4_8^{-1}0_00_6^{-1}$, $1_7^{-1}2_6^{-1}3_5^20_0$, $1_52_6^{-1}3_53_7^{-1}4_60_00_6^{-1}$,
$1_52_6{-1}3_54_8^{-1}0_00_6^{-1}$, $1_7^{-1}3_53_7^{-1}4_60_0$, $1_53_7^{-2}4_6^20_00_6^{-1}$,
$1_7^{-1}3_54_8^{-1}0_0$, $1_53_7^{-1}4_64_8^{-1}0_00_6^{-1}$, $1_7^{-1}2_63_7^{-2}4_6^20_0$,
$1_7^{-1}2_63_7^{-1}4_64_8^{-1}0_0$, $1_54_8^{-2}0_00_6^{-1}$, $2_8^{-1}4_6^20_0$,
$2_8^{-1}3_74_74_8^{-1}0_0$, $1_7^{-1}2_64_8^{-2}0_0$,
$3_{9}^{-1}4_60_0$, $2_8^{-1}3_7^24_8^{-2}0_0$,
$3_73_{9}^{-1}4_8^{-1}0_0$,
$2_83_{9}^{-2}0_0$,
$1_{9}2_{10}^{-1}0_0$,
$1_{11}^{-1}0_00_{10}$.

Let $M_3=\tilde{e}_1\tilde{e}_2\tilde{e}_3^2\tilde{e}_2\tilde{e}_4^2\tilde{e}_3^2\tilde{e}_2\tilde{e}_1\tilde{e}_0 M_2=Y_{1,-5}Y_{0,6}^{-1}Y_{0,-4}^{-1}$. The following $52$ monomials appear in $\mathcal{M}_{I_0}(M_3)$:

$1_{-5}0_6^{-1}0_{-4}^{-1}$, 
$1_{-3}^{-1}2_{-4}0_6^{-1}$, 
$2_{-2}^{-1}3_{-3}^20_6^{-1}$, 
$3_{-3}3_{-1}^{-1}4_{-2}0_6^{-1}$,
$2_{-2}3_{-1}^{-2}4_{-2}^20_6^{-1}$, $3_{-3}4_0^{-1}0_6^{-1}$,
$1_{-1}2_0^{-1}4_{-2}^20_6^{-1}$, $2_{-2}3_{-1}^{-1}4_{-2}4_0^{-1}0_6^{-1}$,
$1_1^{-1}4_{-2}^20_00_6^{-1}$, $1_{-1}2_0^{-1}4_{-2}4_0^{-1}3_{-1}0_6^{-1}$, $2_{-2}4_0^{-2}0_6^{-1}$,
$1_1^{-1}3_{-1}4_{-2}4_0^{-1}0_00_6^{-1}$, $1_{-1}2_0^{-1}4_0^{-2}3_{-1}^20_6^{-1}$, $1_{-1}3_1^{-1}4_{-2}0_6^{-1}$,
$1_1^{-1}2_03_1^{-1}4_{-2}0_00_6^{-1}$, $1_1^{-1}3_{-1}^24_0^{-2}0_00_6^{-1}$, $1_{-1}3_{-1}3_1^{-1}4_0^{-1}0_6^{-1}$,
$2_2^{-1}3_14_{-2}0_00_6^{-1}$, $1_1^{-1}2_03_1^{-1}3_{-1}4_0^{-1}0_00_6^{-1}$, $1_{-1}2_03_1^{-2}0_6^{-1}$,
$3_3^{-1}4_24_{-2}0_00_6^{-1}$, $2_2^{-1}3_13_{-1}4_0^{-1}0_00_6^{-1}$, $1_1^{-1}2_0^23_1^{-2}0_00_6^{-1}$, $1_{-1}1_12_2^{-1}0_6^{-1}$,
$4_4^{-1}4_{-2}0_00_6^{-1}$, $3_3^{-1}3_{-1}4_0^{-1}4_20_00_6^{-1}$, $2_02_2^{-1}0_00_6^{-1}$, $1_{-1}1_3^{-1}0_20_6^{-1}$,
$3_{-1}4_0^{-1}4_4^{-1}0_00_6^{-1}$, $2_03_1^{-1}3_3^{-1}4_20_00_6^{-1}$, $1_12_2^{-2}3_1^20_00_6^{-1}$, $1_1^{-1}1_3^{-1}2_00_20_00_6^{-1}$,
$2_03_1^{-1}4_4^{-1}0_00_6^{-1}$, $1_3^{-1}2_2^{-1}3_1^20_00_20_6^{-1}$, $1_12_2^{-1}3_13_2^{-1}4_20_00_6^{-1}$,
$1_12_2^{-1}3_14_4^{-1}0_00_6^{-1}$, $1_3^{-1}3_13_3^{-1}4_20_00_20_6^{-1}$, $1_13_3^{-2}4_2^20_00_6^{-1}$,
$1_3^{-1}3_14_4^{-1}0_20_00_6^{-1}$, $1_13_3^{-1}4_24_4^{-1}0_00_6^{-1}$, $1_3^{-1}2_23_3^{-2}4_2^20_20_00_6^{-1}$,
$1_3^{-1}2_23_3^{-1}4_24_4^{-1}0_20_00_6^{-1}$, $1_14_4^{-2}0_00_6^{-1}$, $2_4^{-1}4_2^20_20_00_6^{-1}$,
$2_4^{-1}3_34_24_4^{-1}0_20_00_6^{-1}$, $1_3^{-1}2_24_4^{-2}0_20_00_6^{-1}$,
$3_{5}^{-1}4_20_20_00_6^{-1}$, $2_4^{-1}3_3^24_4^{-2}0_20_00_6^{-1}$,
$3_33_{5}^{-1}4_4^{-1}0_20_00_6^{-1}$,
$2_43_{5}^{-2}0_20_00_6^{-1}$,
$1_{5}2_{6}^{-1}0_20_00_6^{-1}$,
$1_{7}^{-1}0_20_0$.

Let $M_4=\tilde{e}_0 M_3=Y_{0,-6}Y_{0,6}^{-1}$ and
$M_5=\tilde{e}_4^2\tilde{e}_3^2\tilde{e}_2^2\tilde{e}_1^2\tilde{e}_0\tilde{f}_1\tilde{f}_2
M=Y_{4,-2}^2Y_{0,2}^{-1}Y_{0,6}^{-1}$. We have
$\mathcal{M}_{I_0}(M_4)=\{M_4\}$. We do not give the list of monomials
of $\mathcal M_{I_0}(M)$ and $\mathcal M_{I_0}(M_5)$ (a total of
$1598$ monomials).

All monomials of $\mathcal M(M)/\tau_2$ are
connected to either $M$, $M_2$, $M_3$, $M_4$, $M_5$ in the
$I_0$-crystal (it is possible to check from the above computation; or
it also follows from $\text{Res}W(\varpi_2) = V_{I_0}(\varpi_2)\oplus
V_{I_0}(\varpi_1)^{\oplus 2}\oplus V_{I_0}(0)\oplus V_{I_0}(2\varpi_4)$
\begin{NB}
as it is clear that $M_2$ and $M_3$ are not in the same $\tau_2$-orbit
\end{NB}
). 
\begin{NB}
So we have
\begin{equation*}
    \mathcal M(M)/\tau_2=\mathcal M_{I_0}(M)\sqcup\mathcal M_{I_0}(M_2)\sqcup\mathcal M_{I_0}(M_3)\sqcup \mathcal
    M_{I_0}(M_4)\sqcup\mathcal M_{I_0}(M_5).
\end{equation*}
\end{NB}


\subsubsection{} Let us consider $\ell=3$ and $M=Y_{3,0}Y_{0,3}^{-1}Y_{0,5}^{-1}$. We have $\tilde{f}_1\tilde{f}_2\tilde{f}_3 M=Y_{4,1}Y_{3,2}Y_{1,4}^{-1}Y_{0,5}^{-1}$ and so we see as in Proposition \ref{otherd} that $\mathcal M(M)\simeq \mathcal B(\varpi_\ell)$. 
Let
$M_1=\tilde{e}_3\tilde{e}_2\tilde{e}_3\tilde{e}_4^2\tilde{e}_1\tilde{e}_2\tilde{e}_3^3\tilde{e}_2^2\tilde{e}_1^2\tilde{e}_0\tilde{f}_1\tilde{f}_2\tilde{f}_3
M=(Y_{4,-1}Y_{4,-3}^{-1})Y_{3,-4}Y_{0,3}^{-1}Y_{0,-1}^{-1}$. This has
weight $\wt M + \delta$ and hence $z_\ell(M)=M_1$.
As 
\begin{equation*}
\tilde{e}_3\tilde{e}_2\tilde{e}_3\tilde{e}_4^2\tilde{e}_1\tilde{e}_2\tilde{e}_3^3\tilde{e}_2^2\tilde{e}_1^2\tilde{e}_0\tilde{f}_1\tilde{f}_2\tilde{f}_3 M_1=Y_{3,-6}Y_{0,-1}^{-1}Y_{0,-3}^{-1}=\tau_{-6}(M),
\end{equation*}
$\mathcal M(m)$ is preserved under $\tau_6$ and we have
$(z_\ell)^{-2}=\tau_6$. 

Let us define the monomials $M_3=\tilde{e}_4\tilde{e}_3\tilde{e}_2\tilde{e}_1\tilde{e}_0 M=Y_{0,5}^{-1}Y_{4,-3}$ and  $M_4=\tilde{e}_4\tilde{e}_3\tilde{e}_2\tilde{e}_1\tilde{e}_0 M_1=Y_{4,-1}Y_{4,-3}^{-1}Y_{4,-7}Y_{0,3}^{-1}$. In particular, as $z_{\ell}$ is compatible with the operators $\tilde{e}_i$, it follows from $z_\ell(M)=M_1$ that $z_\ell(M_3)=M_4$.

The following $26$ monomials appear in $\mathcal{M}_{I_0}(M_3)$: 

$4_{-3}0_5^{-1}$,
$3_{-2}4_{-1}^{-1}0_5^{-1}$,
$2_{-1}3_0^{-1}0_5^{-1}$,
$1_02_1^{-1}3_00_5^{-1}$,
$1_03_2^{-1}4_10_5^{-1}$, $1_2^{-1}3_00_5^{-1}0_1$,
$1_04_3^{-1}0_5^{-1}$, $1_2^{-1}2_13_2^{-1}4_10_5^{-1}0_1$,
$1_2^{-1}2_14_3^{-1}0_5^{-1}0_1$, $2_3^{-1}3_24_10_5^{-1}0_1$,
$2_3^{-1}3_2^24_3^{-1}0_5^{-1}0_1$, $3_4^{-1}4_14_30_5^{-1}0_1$,
$3_23_4^{-1}0_5^{-1}0_1$, $4_14_5^{-1}0_5^{-1}0_1$,
$2_33_4^{-2}4_30_5^{-1}0_1$, $3_24_3^{-1}4_5^{-1}0_5^{-1}0_1$,
$1_42_5^{-1}4_30_5^{-1}0_1$, $2_33_4^{-1}4_5^{-1}0_5^{-1}0_1$,
$1_6^{-1}4_30_1$, $1_42_5^{-1}3_44_5^{-1}0_5^{-1}0_1$,
$1_6^{-1}3_44_8^{-1}0_1$, $1_43_6^{-1}0_5^{-1}0_1$,
$1_6^{-1}2_53_6^{-1}0_1$,
$2_{7}^{-1}3_60_1$,
$3_{8}^{-1}4_{7}0_1$,
$4_{9}^{-1}0_1$.

The following $26$ monomials appear in $\mathcal{M}_{I_0}(M_4)$:

$4_{-1}4_{-3}^{-1}4_{-7}0_3^{-1}$,
$3_{-6}4_{-5}^{-1}4_{-1}4_{-3}^{-1}0_3^{-1}$,
$2_{-5}3_{-4}^{-1}4_{-1}4_{-3}^{-1}0_3^{-1}$,
$1_{-4}2_{-3}^{-1}3_{-4}4_{-1}4_{-3}^{-1}0_3^{-1}$,
$1_{-4}3_{-2}^{-1}4_{-1}0_3^{-1}$, $1_{-2}^{-1}3_{-4}4_{-1}4_{-3}^{-1}0_{-3}0_3^{-1}$,
$1_{-4}3_{-2}^{-1}3_04_1^{-1}0_3^{-1}$, $1_{-2}^{-1}2_{-3}3_{-2}^{-1}4_{-1}0_{-3}0_3^{-1}$,
$1_{-2}^{-1}2_{-3}3_{-2}^{-1}3_04_1^{-1}0_{-3}0_3^{-1}$, $2_{-1}^{-1}3_{-2}4_{-1}0_{-3}0_3^{-1}$,
$2_{-1}^{-1}3_{-2}3_04_1^{-1}0_{-3}0_3^{-1}$, $3_0^{-1}4_{-1}^20_{-3}0_3^{-1}$,
$2_{-1}^{-1}2_13_{-2}3_2^{-1}0_{-3}0_3^{-1}$, $4_{-1}4_1^{-1}0_{-3}0_3^{-1}$,
$2_13_0^{-1}3_2^{-1}4_{-1}0_{-3}0_3^{-1}$, $3_04_1^{-2}0_{-3}0_3^{-1}$,
$1_22_3^{-1}3_0^{-1}3_24_{-1}0_{-3}0_3^{-1}$, $2_13_2^{-1}4_1^{-1}0_{-3}0_3^{-1}$,
$1_4^{-1}3_0^{-1}3_24_{-1}0_{-3}$, $1_22_3^{-1}3_24_1^{-1}0_{-3}0_3^{-1}$,
$1_4^{-1}3_24_1^{-1}0_{-3}$, $1_23_4^{-1}4_1^{-1}4_30_{-3}0_3^{-1}$,
$1_4^{-1}2_33_4^{-1}4_1^{-1}4_30_{-3}$,
$2_{5}^{-1}3_44_1^{-1}4_30_{-3}$,
$3_{6}^{-1}4_54_1^{-1}4_30_{-3}$,
$4_7^{-1}4_1^{-1}4_30_{-3}$.

The crystal isomorphism $z_{\ell}$ is given in order.

It should also be possible to make explicit the bijection between $\mathcal{M}_{I_0}(M)$ and $\mathcal{M}_{I_0}(M_1)$ (but to do not write it in the paper as there are $273$ monomials).

All monomials of $\mathcal M(M)/\tau_6$ are connected to either $M$,
 $M_1$, $M_3$, $M_4$ in the $I_0$-crystal (it is possible to check
 from the above computation; this follows also from
 $\text{Res}W_0(\varpi_3) = V_{I_0}(\varpi_3)\oplus
 V_{I_0}(\varpi_4)$).


\begin{NB}
As an application we have
\begin{equation*}
    \mathcal M(M)/\tau_6=\mathcal M_{I_0}(M)\sqcup\mathcal M_{I_0}(M_1)\sqcup\mathcal M_{I_0}(M_2)\sqcup\mathcal
    M_{I_0}(M_3),
\end{equation*}
and so 
\begin{equation*}
    \mathcal B(W(\varpi_\ell))\simeq \mathcal M_{I_0}(M)/z_\ell=\mathcal M_{I_0}(M)\sqcup\mathcal M_{I_0}(M_1).
\end{equation*}

 \end{NB}

\subsubsection{} Finally consider $\ell=4$ and
$M=Y_{4,0}Y_{0,4}^{-1}$. It follows from Corollary~\ref{corex} that
$\mathcal M(M)\simeq \mathcal B(\varpi_\ell)$. As
$\tilde{e}_4\tilde{e}_3\tilde{e}_2\tilde{e}_1\tilde{e}_3\tilde{e}_2\tilde{e}_4\tilde{e}_3^2\tilde{e}_2\tilde{e}_1\tilde{e}_0
M=\tau_{-6}(M)$, $\mathcal M(M)$ is preserved under $\tau_6$, which is
of weight $\delta$. So $z_\ell=\tau_{-6}$ and $\mathcal M(M)/\tau_6\simeq \mathcal B(W(\varpi_\ell))$.

The following $26$ monomials appear in $\mathcal{M}_{I_0}(M)$:

$4_00_4^{-1}$,
$3_14_2^{-1}0_4^{-1}$,
$2_23_3^{-1}0_4^{-1}$,
$1_32_4^{-1}3_30_4^{-1}$,
$1_33_5^{-1}4_40_4^{-1}$, $1_5^{-1}3_3$,
$1_34_6^{-1}0_4^{-1}$, $1_5^{-1}2_43_5^{-1}4_4$,
$1_5^{-1}2_44_6^{-1}$, $2_6^{-1}3_54_4$,
$2_6^{-1}3_5^24_6^{-1}$, $3_7^{-1}4_44_6$,
$3_53_7^{-1}$, $4_44_8^{-1}$,
$2_63_7^{-2}4_6$, $3_54_6^{-1}4_8^{-1}$,
$1_72_8^{-1}4_6$, $2_63_7^{-1}4_8^{-1}$,
$1_9^{-1}4_60_8$, $1_72_8^{-1}3_74_8^{-1}$,
$1_9^{-1}3_74_8^{-1}0_8$, $1_73_9^{-1}$,
$1_9^{-1}2_83_9^{-1}0_8$,
$2_{10}^{-1}3_90_8$,
$3_{11}^{-1}4_{10}0_8$,
$4_{12}^{-1}0_8$.

These are the monomials appearing in $\mathcal{M}(M)/\tau_6$. We thus
have $\mathcal B(W(\varpi_\ell)) \simeq \mathcal B_{I_0}(\varpi_\ell)$. 


\begin{rem}
The authors do not find the last three examples in the literature. One
can probably check whether they are {\it perfect} or not, though we
have not done yet.
\end{rem}

\begin{NB}
\subsection{Type $A_2^{(2)}$} Let $\ell=1$ and $m=Y_{1,0}Y_{0,1}^{-2}$. It follows from Corollary~\ref{corex} that $\mathcal B(m)\simeq \mathcal B(\varpi_\ell)$. 

As $\tilde{f}_0\tilde{f}_0\tilde{f}_1 m=\tau_2(m)$, $\mathcal M(m)$ is preserved under $\tau_2$. Moreover for weight reason $z_\ell=\tau_{-2}$ (the weight of $z_\ell$ is $d_\ell\delta=\delta=2\alpha_0+\alpha_1$) and so $\mathcal M(m)/\tau_2\simeq \mathcal B(W(\varpi_\ell))$.

Let $m'=\tilde{f}_1 m=Y_{1,2}^{-1}Y_{0,1}^2$ and $m''=\tilde{f}_0 m'=Y_{0,1}Y_{0,3}^{-1}$. We have $\mathcal M_{I_0}(m)=\{m, m'\}$, $\mathcal M_{I_0}(m'')=\{m''\}$. The crystal graph of $\mathcal M(m)$ is given in Figure~\ref{fig:A2}.
\begin{figure}[htbp]
\centering
\psset{xunit=1mm,yunit=1mm,runit=1mm}
\psset{linewidth=0.3,dotsep=1,hatchwidth=0.3,hatchsep=1.5,shadowsize=1}
\psset{dotsize=0.7 2.5,dotscale=1 1,fillcolor=black}
\psset{arrowsize=1 2,arrowlength=1,arrowinset=0.25,tbarsize=0.7 5,bracketlength=0.15,rbracketlength=0.15}
\begin{pspicture}(0,60)(160.78,85)
\rput(10,80){$m$}
\rput(60,80){$m'$}
\rput(110,80){$m''$}
\rput{90}(60,164.22){\parametricplot[arrows=<-]{150.26}{209.74}{ t cos 100.78 mul t sin 100.78 mul }}
\psline{->}(15,80)(55,80)
\psline{->}(65,80)(105,80)
\rput(30,82){$\scriptstyle 1$}
\rput(80,82){$\scriptstyle 0$}
\rput(60,65){$\scriptstyle 0[2]$}
\end{pspicture}
\caption{(Type $A_2^{(2)}$) the crystal $\mathcal B(\varpi_1)$}
\label{fig:A2}
\end{figure}
In particular all monomials of $\mathcal M(m)$ are connected to either $m$ or $m''$ or their $\tau_2$ images in the $I_0$-crystal, thus
\begin{equation*}\mathcal M(m)/\tau_2=\mathcal M_{I_0}(m)\sqcup\mathcal M_{I_0}(m'').\end{equation*}
\end{NB}

\subsection{Type $E_6^{(2)}$} 

\subsubsection{} First let $\ell=1$ and
$M=Y_{1,0}Y_{0,1}^{-1}Y_{0,5}^{-1}$. We have $\tilde{f}_2\tilde{f}_1
M=Y_{0,5}^{-1}Y_{2,3}^{-1}Y_{3,2}$ and so we see as in Proposition
\ref{otherd} that $\mathcal M(M)\simeq \mathcal B(\varpi_\ell)$. As
$\tilde{e}_1\tilde{e}_2\tilde{e}_3\tilde{e}_4\tilde{e}_2\tilde{e}_3\tilde{e}_2\tilde{e}_1\tilde{e}_1\tilde{e}_0\tilde{f}_1
M=\tau_{-4}(M)$, $\mathcal M(m)$ is preserved under $\tau_4$, which is
of weight $\delta$. Thus we have $z_\ell=\tau_{-4}$ and $\mathcal M(m)/\tau_4\simeq \mathcal B(W(\varpi_\ell))$.

Let $M'=\tilde{e}_0 M=Y_{0,-1}Y_{0,5}^{-1}$. We have $\mathcal M_{I_0}(M')=\{M'\}$. 

The following $26$ monomials appear in $\mathcal{M}_{I_0}(M)$:

$1_00_1^{-1}0_5^{-1}$,
$2_11_2^{-1}0_5^{-1}$,
$3_22_3^{-1}0_5^{-1}$,
$4_33_4^{-1}2_30_5^{-1}$,
$4_32_5^{-1}1_40_5^{-1}$, $4_5^{-1}2_30_5^{-1}$,
$4_31_6^{-1}$, $4_5^{-1}3_42_5^{-1}1_40_5^{-1}$,
$4_5^{-1}3_41_6^{-1}$, $3_6^{-1}2_51_40_5^{-1}$,
$3_6^{-1}2_5^21_6^{-1}$, $2_7^{-1}1_41_60_5^{-1}$,
$2_52_7^{-1}$, $1_41_8^{-1}0_5^{-1}0_7$,
$3_62_7^{-2}1_6$, $2_51_6^{-1}1_8^{-1}0_7$,
$4_73_8^{-1}1_6$, $3_62_7^{-1}1_8^{-1}0_7$,
$4_9^{-1}1_6$, $4_73_8^{-1}2_71_8^{-1}0_7$,
$4_9^{-1}2_71_8^{-1}0_7$, $4_72_9^{-1}0_7$,
$4_9^{-1}3_82_9^{-1}0_7$,
$3_{10}^{-1}2_90_7$,
$2_{11}^{-1}1_{10}0_7$,
$1_{12}^{-1}0_70_{11}$.

These $26$ monomials, $M'$ are all monomials of $\mathcal
M(M)/\tau_4$. As an application we have
\begin{equation*}
   \mathcal B(W(\varpi_1)) \simeq \mathcal B_{I_0}(\varpi_1) \sqcup 
   \mathcal B_{I_0}(0).
\end{equation*}
\begin{NB}
\begin{equation*}\mathcal M(M)/\tau_4=\mathcal M_{I_0}(M)\sqcup\mathcal M_{I_0}(M').\end{equation*}
\end{NB}


The crystal structure here is isomorphic to one studied recently
in \cite{bfkl}. As the crystal graph is connected, we conclude that
the crystal base constructed in \cite{bfkl} are isomorphic to
$\mathcal B(W(\varpi_\ell))$.

\subsubsection{} Now we consider $\ell=2$ and $M=Y_{2,0}Y_{0,2}^{-1}Y_{0,4}^{-1}Y_{0,6}^{-1}$. We have $\tilde{f}_1\tilde{f}_2^2\tilde{f}_3\tilde{f}_1\tilde{f}_2 M=Y_{4,2}Y_{3,3}Y_{2,4}^{-1}Y_{1,5}^{-1}Y_{0,6}^{-1}$ and so we see as in Proposition \ref{otherd} that $\mathcal M(M)\simeq \mathcal B(\varpi_\ell)$. 
As
$\tilde{e}_2\tilde{e}_3\tilde{e}_2\tilde{e}_1^2\tilde{e}_2\tilde{e}_4\tilde{e}_3^2\tilde{e}_2^3\tilde{e}_1^2\tilde{e}_0\tilde{f}_1\tilde{f}_2^2\tilde{f}_3\tilde{f}_1\tilde{f}_2
M=\tau_{-2}(M)$, $\mathcal M(M)$ is preserved under $\tau_2$, which is
of weight $\delta$. Therefore we have $z_\ell=\tau_{-2}$ and $\mathcal M(M)/\tau_2\simeq \mathcal B(W(\varpi_\ell))$.

Let $M_2=\tilde{e}_1\tilde{e}_0
M=Y_{1,-1}Y_{0,4}^{-1}Y_{0,6}^{-1}$. The following $26$ monomials
appear in $\mathcal{M}_{I_0}(M_2)$:

$1_{-1}0_4^{-1}0_6^{-1}$,
$2_01_1^{-1}0_00_4^{-1}0_6^{-1}$,
$3_12_2^{-1}0_00_4^{-1}0_6^{-1}$,
$4_23_3^{-1}2_20_00_4^{-1}0_6^{-1}$,
$4_22_4^{-1}1_30_00_4^{-1}0_6^{-1}$, $4_4^{-1}2_20_00_4^{-1}0_6^{-1}$,
$4_21_5^{-1}0_00_6^{-1}$, $4_4^{-1}3_32_4^{-1}1_30_00_4^{-1}0_6^{-1}$,
$4_4^{-1}3_31_5^{-1}0_00_6^{-1}$, $3_5^{-1}2_41_30_00_4^{-1}0_6^{-1}$,
$3_5^{-1}2_4^21_5^{-1}0_00_6^{-1}$, $2_6^{-1}1_31_50_00_4^{-1}0_6^{-1}$,
$2_42_6^{-1}0_00_6^{-1}$, $1_31_7^{-1}0_00_4^{-1}$,
$3_52_6^{-2}1_50_00_6^{-1}$, $2_41_5^{-1}1_7^{-1}0_0$,
$4_63_7^{-1}1_50_00_6^{-1}$, $3_52_6^{-1}1_7^{-1}0_0$,
$4_8^{-1}1_50_00_6^{-1}$, $4_63_7^{-1}2_61_7^{-1}0_0$,
$4_8^{-1}2_61_7^{-1}0_0$, $4_62_8^{-1}0_0$,
$4_8^{-1}3_72_8^{-1}0_0$,
$3_{9}^{-1}2_80_0$,
$2_{10}^{-1}1_{9}0_0$,
$1_{11}^{-1}0_00_{10}$.

Let
$M_3=\tilde{e}_1\tilde{e}_2\tilde{e}_3\tilde{e}_2\tilde{e}_4\tilde{e}_3\tilde{e}_2\tilde{e}_1\tilde{e}_0
M_2=Y_{1,-5}Y_{0,6}^{-1}Y_{0,-4}^{-1}$. The following $26$ monomials
appear in $\mathcal{M}_{I_0}(M_3)$:

$1_{-5}0_{-4}^{-1}0_6^{-1}$,
$2_{-4}1_{-3}^{-1}0_6^{-1}$,
$3_{-3}2_{-2}^{-1}0_6^{-1}$,
$4_{-2}3_{-1}^{-1}2_{-2}0_6^{-1}$,
$4_{-2}2_0^{-1}1_{-1}0_6^{-1}$, $4_0^{-1}2_{-2}0_6^{-1}$,
$4_{-2}1_1^{-1}0_00_6^{-1}$, $4_0^{-1}3_{-1}2_0^{-1}1_{-1}0_6^{-1}$,
$4_0^{-1}3_{-1}1_1^{-1}0_00_6^{-1}$, $3_1^{-1}2_01_{-1}0_6^{-1}$,
$3_1^{-1}2_0^21_1^{-1}0_00_6^{-1}$, $2_2^{-1}1_{-1}1_10_6^{-1}$,
$2_02_2^{-1}0_00_6^{-1}$, $1_{-1}1_3^{-1}0_6^{-1}0_2$,
$3_12_2^{-2}1_10_00_6^{-1}$, $2_01_1^{-1}1_3^{-1}0_6^{-1}0_00_2$,
$4_23_3^{-1}1_10_00_6^{-1}$, $3_12_2^{-1}1_3^{-1}0_00_6^{-1}0_2$,
$4_4^{-1}1_10_00_6^{-1}$, $4_23_3^{-1}2_21_3^{-1}0_00_6^{-1}0_2$,
$4_4^{-1}2_21_3^{-1}0_00_6^{-1}0_2$, $4_22_4^{-1}0_00_6^{-1}0_2$,
$4_4^{-1}3_32_4^{-1}0_00_6^{-1}0_2$,
$3_{5}^{-1}2_40_00_6^{-1}0_2$,
$2_{6}^{-1}1_{5}0_00_6^{-1}0_2$,
$1_{7}^{-1}0_00_2$.

Let $M_4=\tilde{e}_0 M_3=Y_{0,-6}Y_{0,6}^{-1}$. We have $\mathcal{M}_{I_0}(M_4)=\{M_4\}$.

Let
$M_5=\tilde{e}_4\tilde{e}_3\tilde{e}_2^2\tilde{e}_1^2\tilde{e}_0\tilde{f}_1\tilde{f}_2
M=Y_{4,-2}Y_{0,2}^{-1}Y_{0,6}^{-1}$. The following $52$ monomials
appear in $\mathcal{M}_{I_0}(M_5)$:

$4_{-2}0_2^{-1}0_6^{-1}$, 
$4_0^{-1}3_{-1}0_2^{-1}0_6^{-1}$, 
$3_1^{-1}2_0^20_2^{-1}0_6^{-1}$, 
$2_02_2^{-1}1_10_2^{-1}0_6^{-1}$,
$3_12_2^{-2}1_1^20_2^{-1}0_6^{-1}$, $2_01_3^{-1}0_6^{-1}$,
$4_23_3^{-1}1_1^20_2^{-1}0_6^{-1}$, $3_12_2^{-1}1_11_3^{-1}0_6^{-1}$,
$4_4^{-1}1_1^20_2^{-1}0_6^{-1}$, $4_23_4^{-1}1_11_3^{-1}2_20_6^{-1}$, $3_11_3^{-2}0_20_6^{-1}$,
$4_4^{-1}2_21_11_3^{-1}0_6^{-1}$, $4_23_3^{-1}1_3^{-2}2_2^20_20_6^{-1}$, $4_22_4^{-1}1_10_6^{-1}$,
$4_4^{-1}3_32_4^{-1}1_10_6^{-1}$, $4_4^{-1}2_2^21_3^{-2}0_20_6^{-1}$, $4_22_22_4^{-1}1_3^{-1}0_20_6^{-1}$,
$3_5^{-1}2_41_10_6^{-1}$, $4_4^{-1}3_32_4^{-1}2_21_3^{-1}0_20_6^{-1}$, $4_23_32_4^{-2}0_20_6^{-1}$,
$2_6^{-1}1_51_10_6^{-1}$, $3_5^{-1}2_42_21_3^{-1}0_20_6^{-1}$, $4_4^{-1}3_3^22_4^{-2}0_20_6^{-1}$, $4_24_43_5^{-1}0_20_6^{-1}$,
$1_7^{-1}1_1$, $2_6^{-1}2_21_3^{-1}1_50_20_6^{-1}$, $3_33_5^{-1}0_20_6^{-1}$, $4_24_6^{-1}0_20_6^{-1}$,
$2_21_3^{-1}1_7^{-1}0_2$, $3_32_4^{-1}2_6^{-1}1_50_20_6^{-1}$, $4_43_5^{-2}2_4^20_20_6^{-1}$, $4_4^{-1}4_6^{-1}3_30_20_6^{-1}$,
$3_32_4^{-1}1_7^{-1}0_2$, $4_6^{-1}3_5^{-1}2_4^20_20_6^{-1}$, $4_43_5^{-1}2_42_6^{-1}1_50_20_6^{-1}$,
$4_43_5^{-1}2_41_7^{-1}0_2$, $4_6^{-1}2_43_6^{-1}1_50_20_6^{-1}$, $4_42_6^{-2}1_5^20_20_6^{-1}$,
$4_6^{-1}2_41_7^{-1}0_2$, $4_42_6^{-1}1_51_7^{-1}0_2$, $4_6^{-1}3_52_6^{-2}1_5^20_20_6^{-1}$,
$4_6^{-1}3_52_6^{-1}1_51_7^{-1}0_2$, $4_41_7^{-2}0_20_6$, $3_7^{-1}1_5^20_20_6^{-1}$,
$3_7^{-1}2_61_51_7^{-1}0_2$, $4_6^{-1}3_51_7^{-2}0_20_6$,
$2_{8}^{-1}1_50_2$, $3_7^{-1}2_6^21_7^{-2}0_20_6$,
$2_62_{8}^{-1}1_7^{-1}0_20_6$,
$3_72_{8}^{-2}0_20_6$,
$4_{8}3_{9}^{-1}0_20_6$,
$4_{10}^{-1}0_20_6$.

We do not list the $273$ monomials of $\mathcal{M}_{I_0}(M)$, but we can check that all monomials of $\mathcal M(M)/\tau_2$ are connected to either $M$, $M_2$, $M_3$, $M_4$, $M_5$ in the $I_0$-crystal. As an application we have
\begin{equation*}
    \mathcal B(W(\varpi_2)) \simeq \mathcal B_{I_0}(\varpi_2)\sqcup
    \mathcal B_{I_0}(\varpi_1)\sqcup\mathcal B_{I_0}(\varpi_1)\sqcup
    \mathcal B_{I_0}(0)\sqcup\mathcal B_{I_0}(\varpi_4).
\end{equation*}


\begin{NB}
\begin{equation*}
    \mathcal M(M)/\tau_2=\mathcal M_{I_0}(M)\sqcup\mathcal M_{I_0}(M_2)\sqcup\mathcal M_{I_0}(M_3)\sqcup \mathcal   M_{I_0}(M_4)\sqcup\mathcal M_{I_0}(M_5).
\end{equation*}
\end{NB}

\subsubsection{} We consider $\ell=3$ and $M=Y_{3,0}Y_{0,3}^{-2}Y_{0,5}^{-2}$. We have $\tilde{f}_1^2\tilde{f}_2^2\tilde{f}_3 M=Y_{4,1}Y_{3,2}Y_{1,4}^{-2}Y_{0,5}^{-2}$ and so we see as in Proposition \ref{otherd} that $\mathcal M(M)\simeq \mathcal B(\varpi_\ell)$. 
As
$\tilde{e}_3\tilde{e}_2^2\tilde{e}_4\tilde{e}_3^2\tilde{e}_2^2\tilde{e}_1^4\tilde{e}_2^2\tilde{e}_3^2\tilde{e}_2^3\tilde{e}_4^3\tilde{e}_3^3\tilde{e}_2^4\tilde{e}_1^4\tilde{e}_0^4
M=\tau_{-6}(M)$, $\mathcal {M}(M)$ is preserved under $\tau_6$, which
is of weight $-4\delta = -2d_\ell \delta$. Therefore we have
$(z_\ell)^{-2}=\tau_6$. Let
$M_1=\tilde{e}_3\tilde{e}_2^2\tilde{e}_1^2\tilde{e}_3\tilde{e}_2^2\tilde{e}_4^2\tilde{e}_3^3\tilde{e}_2^4\tilde{e}_1^4\tilde{e}_0^2\tilde{f}_1^2\tilde{f}_2^2\tilde{f}_3
M=Y_{0,3}^{-2}Y_{0,-1}^{-2}Y_{3,-4}$. This has weight $\wt(M) +
d_\ell\delta$, hence we have $z_\ell(M)=M_1$.

We do not determine the $I_0$-crystal components of $\mathcal{M}(M)/\tau_6$ at this moment.

\subsubsection{} Finally let us consider $\ell=4$. Let $M=Y_{4,0}Y_{0,4}^{-2}$. It follows from Corollary~\ref{corex} that $\mathcal M(M)\simeq \mathcal B(\varpi_\ell)$. 
As
$\tilde{e}_4\tilde{e}_3\tilde{e}_2^2\tilde{e}_1^2\tilde{e}_3\tilde{e}_2^2\tilde{e}_4\tilde{e}_3^2\tilde{e}_2^2\tilde{e}_1^2\tilde{e}_0^2
M=\tau_{-6}(M)$, $\mathcal M(M)$ is preserved under $\tau_6$, which is
of weight $-2\delta = -d_\ell\delta$. Therefore we have $z_\ell=\tau_{-6}$ and $\mathcal M(M)/\tau_6\simeq \mathcal B(W(\varpi_\ell))$.

The following $52$ monomials appear in $\mathcal{M}_{I_0}(M)$:

$4_00_4^{-2}$, 
$2_2^{-1}3_10_4^{-2}$, 
$3_3^{-1}2_2^20_4^{-2}$, 
$2_22_4^{-1}1_30_4^{-2}$,
$3_32_4^{-2}1_3^20_4^{-2}$, $2_21_5^{-1}0_4^{-1}$,
$4_43_5^{-1}1_3^20_4^{-2}$, $3_32_4^{-1}1_31_5^{-1}0_4^{-1}$,
$4_6^{-1}1_3^20_4^{-2}$, $4_43_5^{-1}1_31_5^{-1}2_40_4^{-1}$, $3_31_5^{-2}$,
$4_6^{-1}2_41_31_5^{-1}0_4^{-1}$, $4_43_5^{-1}1_5^{-2}2_4^2$, $4_42_6^{-1}1_30_4^{-1}$,
$4_6^{-1}3_52_6^{-1}1_30_4^{-1}$, $4_6^{-1}2_4^21_5^{-2}$, $4_42_42_6^{-1}1_5^{-1}$,
$3_7^{-1}2_61_30_4^{-1}$, $4_6^{-1}3_52_6^{-1}2_41_5^{-1}$, $4_43_52_6^{-2}$,
$2_8^{-1}1_71_30_4^{-1}$, $3_7^{-1}2_62_41_5^{-1}$, $4_6^{-1}3_5^22_6^{-2}$, $4_44_63_7^{-1}$,
$1_9^{-1}1_30_80_4^{-1}$, $2_8^{-1}2_41_5^{-1}1_7$, $3_53_7^{-1}$, $4_44_8^{-1}$,
$2_41_5^{-1}1_9^{-1}0_8$, $3_52_6^{-1}2_8^{-1}1_7$, $4_63_7^{-2}2_6^2$, $4_6^{-1}4_8^{-1}3_5$,
$3_52_6^{-1}1_9^{-1}0_8$, $4_8^{-1}3_7^{-1}2_6^2$, $4_63_7^{-1}2_62_8^{-1}1_7$,
$4_63_7^{-1}2_61_9^{-1}0_8$, $4_8^{-1}2_62_8^{-1}1_7$, $4_62_8^{-2}1_7^2$,
$4_8^{-1}2_61_9^{-1}0_8$, $4_62_8^{-1}1_71_9^{-1}0_8$, $4_8^{-1}3_72_8^{-2}1_7^2$,
$4_8^{-1}3_72_8^{-1}1_71_9^{-1}0_8$, $4_61_9^{-2}0_8^2$, $3_9^{-1}1_7^2$,
$3_9^{-1}2_81_71_9^{-1}0_8$, $4_8^{-1}3_71_9^{-2}0_8^2$,
$2_{10}^{-1}1_70_8$, $3_9^{-1}2_8^21_9^{-2}0_8^2$,
$2_82_{10}^{-1}1_9^{-1}0_8^2$,
$3_92_{10}^{-2}0_8^2$,
$4_{10}3_{11}^{-1}0_8^2$,
$4_{12}^{-1}0_70_8^2$.

Let
$M'=\tilde{e}_1\tilde{e}_2\tilde{e}_3\tilde{e}_2\tilde{e}_1\tilde{e}_0
M=Y_{1,-3}Y_{0,4}^{-1}Y_{0,-2}^{-1}$. The following $26$ monomials
appear in $\mathcal{M}_{I_0}(M')$:

$1_{-3}0_{-2}^{-1}0_4^{-1}$,
$2_{-2}1_{-1}^{-1}0_4^{-1}$,
$3_{-1}2_0^{-1}0_4^{-1}$,
$4_03_1^{-1}2_00_4^{-1}$,
$4_02_2^{-1}1_10_4^{-1}$, $4_2^{-1}2_00_4^{-1}$,
$4_01_3^{-1}0_20_4^{-1}$, $4_2^{-1}3_12_2^{-1}1_10_4^{-1}$,
$4_2^{-1}3_11_3^{-1}0_20_4^{-1}$, $3_3^{-1}2_21_10_4^{-1}$,
$3_3^{-1}2_2^21_3^{-1}0_20_4^{-1}$, $2_4^{-1}1_11_30_4^{-1}$,
$2_22_4^{-1}0_20_4^{-1}$, $1_11_5^{-1}$,
$3_32_4^{-2}1_30_20_4^{-1}$, $2_21_3^{-1}1_5^{-1}0_2$,
$4_43_5^{-1}1_30_20_4^{-1}$, $3_32_4^{-1}1_5^{-1}0_2$,
$4_6^{-1}1_30_20_4^{-1}$, $4_43_5^{-1}2_41_5^{-1}0_2$,
$4_6^{-1}2_41_5^{-1}0_2$, $4_42_6^{-1}0_2$,
$4_6^{-1}3_52_6^{-1}0_2$,
$3_{7}^{-1}2_60_2$,
$2_{8}^{-1}1_{7}0_2$,
$1_{9}^{-1}0_20_8$.

Let $M''=\tilde{e}_0 M'=Y_{0,-4}Y_{0,4}^{-1}$. We have $\mathcal{M}_{I_0}(M'')=\{M''\}$.

The above exhausts all monomials of $\mathcal M(M)/\tau_6$. As an
application we have 
\begin{NB}
\begin{equation*}\mathcal M(M)/\tau_6=\mathcal M_{I_0}(M)\sqcup\mathcal M_{I_0}(M')\sqcup\mathcal M_{I_0}(M'').\end{equation*}
\end{NB}
\begin{equation*}
   \mathcal B(W(\varpi_4)) \simeq 
   \mathcal B_{I_0}(\varpi_4)\sqcup
   \mathcal B_{I_0}(\varpi_1)\sqcup
   \mathcal B_{I_0}(0).
\end{equation*}


\begin{rem}
The authors do not find the description of the examples $\ell = 2$,
$3$, $4$ in the literature.
\end{rem}

\subsection{Type $D_4^{(3)}$} 

\subsubsection{} First we consider $\ell=1$. Let $M = Y_{1,0}Y_{0,1}^{-1}Y_{0,3}^{-1}$. As $\tilde{f}_1 M=Y_{1,2}^{-2}Y_{0,3}^{-1}Y_{2,1}$ we see as in Proposition \ref{otherd} that $\mathcal M(M)\simeq \mathcal B(\varpi_\ell)$.
The following 7 monomials appear in $\mathcal{M}_{I_0}(M)$:

$M=1_00_1^{-1}0_3^{-1}$,
$m_2=1_2^{-1}0_3^{-1}2_1$,
$m_3=2_3^{-1}1_2^20_3^{-1}$,
$m_4=1_21_4^{-1}$,
$m_5=1_4^{-2}2_30_3$,
$m_6=2_5^{-1}1_40_3$,
$m_7=1_6^{-1}0_30_5$.

\begin{figure}[htbp]
\centering
\psset{xunit=.55mm,yunit=.55mm,runit=.55mm}
\psset{linewidth=0.3,dotsep=1,hatchwidth=0.3,hatchsep=1.5,shadowsize=1}
\psset{dotsize=0.7 2.5,dotscale=1 1,fillcolor=black}
\psset{arrowsize=1 2,arrowlength=1,arrowinset=0.25,tbarsize=0.7 5,bracketlength=0.15,rbracketlength=0.15}
\begin{pspicture}(0,0)(220,50)
\rput(0,20){$M$}
\psline{->}(12,20)(24,20)
\rput(35,20){$m_2$}
\psline{->}(47,20)(59,20)
\rput(70,20){$m_3$}
\psline{->}(82,20)(94,20)
\rput(105,20){$m_4$}
\psline{->}(117,20)(129,20)
\rput(140,20){$m_5$}
\psline{->}(152,20)(164,20)
\rput(175,20){$m_6$}
\psline{->}(187,20)(199,20)
\rput(210,20){$m_7$}
\rput(105,40){$M'$}
\rput{0}(0,20){\qline(-0,0)(0,0)}\rput{0}(88,-486){\parametricplot[arrows=->]{84.85}{95.15}{ t cos 515.26 mul t sin 515.26 mul }}
\rput{0}(120.26,528){\parametricplot[arrows=<-]{-95.15}{-84.85}{ t cos 515.26 mul t sin 515.26 mul }}
\rput{0}(106,-500){\parametricplot[arrows=->]{78}{89.15}{ t cos 500.26 mul t sin 540.26 mul }}
\rput{0}(106,-500){\parametricplot[arrows=<-]{102.15}{91.5}{ t cos 500.26 mul t sin 540.26 mul }}
\rput(15,25){$\scriptscriptstyle 1$}
\rput(50,25){$\scriptscriptstyle 2$}
\rput(85,25){$\scriptscriptstyle 1$}
\rput(120,25){$\scriptscriptstyle 1$}
\rput(155,25){$\scriptscriptstyle 2$}
\rput(190,25){$\scriptscriptstyle 1$}
\rput(120,30){$\scriptscriptstyle 0[2]$}
\rput(128,7){$\scriptscriptstyle 0[2]$}
\rput(192,40){$\scriptscriptstyle 0[2]$}
\rput(17,40){$\scriptscriptstyle 0[2]$}
\end{pspicture}
\caption{(Type $D_4^{(3)}$) the crystal $\mathcal B(\varpi_1)$}
\label{fig:D4}
\end{figure}
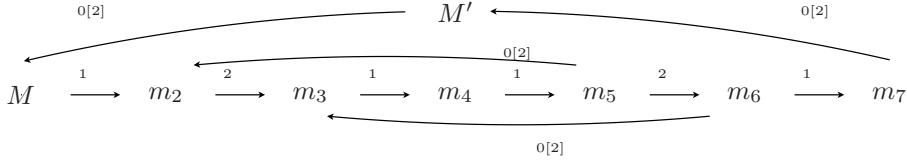

Let $M'=\tau_2(\tilde{e}_0 M)=Y_{0,1}Y_{0,5}^{-1}$. We have $\mathcal
M_{I_0}(M')=\{M'\}$. The crystal graph of $\mathcal M(M)$ is given in
Figure~\ref{fig:D4}.
We find that $z_l=\tau_{-2}$ and $\mathcal M(M)/\tau_2=\mathcal M_{I_0}(M)\sqcup\mathcal M_{I_0}(M')$.

This crystal was described in \cite{JM}.

\subsubsection{} Now we consider $\ell=2$. Let $M = Y_{2,0}Y_{0,2}^{-3}$. It follows from Corollary~\ref{corex} that $\mathcal M(M)\simeq \mathcal B(\varpi_\ell)$. 
As $\tilde{e}_2\tilde{e}_1^3\tilde{e}_2^2\tilde{e}_1^3\tilde{e}_0^3
M=\tau_{-4}(M)$, $\mathcal M(M)$ is preserved under $\tau_4$, which is
of weight $3\delta = d_2\delta$. Therefore  $z_\ell=\tau_{-4}$ and so $\mathcal M(m)/\tau_4\simeq \mathcal B(W(\varpi_\ell))$.

The following 14 monomials appear in $\mathcal{M}_{I_0}(M)$: 

$2_00_2^{-3}$,
$2_2^{-1}1_1^30_2^{-3}$,
$1_1^21_3^{-1}0_2^{-2}$,
$2_21_11_3^{-2}0_2^{-1}$,
$1_3^{-3}1_2^2$,
$2_4^{-1}1_11_30_2^{-1}$,
$2_22_4^{-1}$,
$1_11_5^{-1}0_2^{-1}0_4$,
$2_4^{-2}1_3^3$,
$1_3^{-1}1_5^{-1}2_20_4$,
$2_4^{-1}1_3^21_5^{-1}0_4$,
$1_31_5^{-2}0_4^2$,
$1_5^{-3}0_4^32_4$,
$2_6^{-1}0_4^3$.

Let $M_2=\tilde{e}_1\tilde{e}_0 M=Y_{0,2}^{-2}Y_{1,-1}$. The following
7 monomials appear in $\mathcal{M}_{I_0}(M_2)$:
 
$0_2^{-2}1_{-1}$,
$0_00_2^{-2}1_1^{-1}2_0$,
$0_00_2^{-2}1_1^22_2^{-1}$,
$0_00_2^{-1}1_11_3^{-1}$,
$1_3^{-2}2_20_0$,
$2_4^{-1}1_30_0$,
$1_5^{-1}0_00_4$.

Let $M_3=\tilde{e}_1\tilde{e}_2\tilde{e}_1\tilde{e}_0
M_2=Y_{1,-3}Y_{0,-2}^{-1}Y_{0,2}^{-1}$. The following 7 monomials
appear in $\mathcal{M}_{I_0}(M_3)$:

$0_{-2}^{-1}0_2^{-1}1_{-3}$,
$0_2^{-1}1_{-1}^{-1}2_{-2}$,
$0_2^{-1}1_{-1}^22_0^{-1}$,
$0_00_2^{-1}1_{-1}1_1^{-1}$,
$1_1^{-2}2_00_0^20_2^{-1}$,
$2_2^{-1}1_10_0^20_2^{-1}$,
$1_3^{-1}0_0^2$.

Let $M_4=\tilde{e}_0 M_3=Y_{0,-4}Y_{0,2}^{-1}$. We have $\mathcal M_{I_0}(M_4)=\{M_4\}$. 

By direct calculation we can see that all monomials of $\mathcal
M(M)/\tau_4$ are connected to either $M$ or $M_2$ or $M_3$ or $M_4$ in
the $I_0$-crystal. As an application we have
\begin{equation*}
   \mathcal B(W(\varpi_2)) \simeq
   \mathcal B_{I_0}(\varpi_2)\sqcup
   \mathcal B_{I_0}(\varpi_1)\sqcup\mathcal B_{I_0}(\varpi_1)
   \sqcup\mathcal B_{I_0}(0).
\end{equation*}
\begin{NB}
\begin{equation*}\mathcal M(M)/\tau_4 = \mathcal M_{I_0}(M)\sqcup\mathcal M_{I_0}(M_2)\sqcup\mathcal M_{I_0}(M_3)\sqcup\mathcal M_{I_0}(M_4).\end{equation*}
\end{NB}


The authors do not find the description of this example in the
literature.

\section{Discussions}

(1) As we saw in the simply-laced type examples (except the last one in \ref{last}) in this
paper, we can construct explicit bijections between monomial crystals
$\mathcal M(m)$ and the set $\mathcal C(m_0)$ of monomials in
$q$--characters counted with multiplicities. (Here $m_0$ is obtained
from $m$ by setting $Y_{0,*}$ as $1$.) Their origin is combinatorial
and we do not understand their representation theoretical meaning
yet. In the example in \ref{ALE}, the global crystal base element
corresponding to the exceptional monomial does not belong to a single
$l$--weight subspace.

Also we can check that the bijection is compatible with the crystal
structure in the following sense: Let $\mathcal M_{I_0}(m_0)$ be the
component of the monomial crystal for $\mathfrak g_{I_0}$ containing
$m_0$. Let $p\colon \mathcal M(m) \to \mathcal M_{I_0}(m_0)$ be the
composition of the above mentioned bijection and the map obtained by
forgetting multiplicities. Then $p$ is a morphism of the crystal (but
not strict). This is not true in general.

Counterexample: In the $q$--character of $W(\varpi_3)$ for $E_6$ we
have monomials $m_1=Y_{3,4}Y_{3,6}^{-1}Y_{4,3}Y_{2,5}^{-1}Y_{1,4}$
and $m_2=Y_{3,4}^2Y_{3,6}^{-1}Y_{4,5}^{-1}Y_{5,4}Y_{2,5}^{-1}Y_{1,4}$
with coefficients $1+ 2t^2 + t^4$ and $1 + t^2 + t^4$.
We have $\tilde{f}_4 m_1=m_2$ in the monomial crystal. If we had a
crystal morphism which preserves the weight, the 4 vectors corresponding
to $m_1$ would necessarily satisfy $\phi_4\geq 1$, and each of them would
be sent by $\tilde f_4$ to vectors corresponding to $m_2$. As there
are only 3 of them, we have a contradiction.

(2) In \cite{ns2} Naito-Sagaki proved that the crystal of
Lakshmibai-Seshadri paths of shape $\varpi_\ell$ is isomorphic to
$\mathcal B(\varpi_\ell)$. This result is better than Theorem~\ref{reallevel} in the sense that they determine all paths, not in a
recursive way as ours. Therefore it would be nice if we could give an
explicit map from the path crystal to the monomial crystal.

\end{document}